\documentclass[UTF8]{article}
	\usepackage{subfigure}
	\usepackage{float}
	\usepackage[graphicx]{realboxes}
	\usepackage[margin=1in]{geometry}  
	\usepackage{graphicx}
	\usepackage{amssymb}             
	\usepackage{amsmath}
	\usepackage{mathtools}
	\usepackage{slashed}            
	\usepackage{amsfonts}              
	\usepackage{amsthm}
	\usepackage{breqn}
	\usepackage{verbatim}            
	\usepackage{enumitem}
	\usepackage[colorlinks,allcolors=black]{hyperref}

	\usepackage[numbers]{natbib}
	\usepackage{mathtools}
	\usepackage{cases}
	\usepackage{fancyhdr}
	\usepackage{appendix}

	\theoremstyle{definition}
	\newtheorem{thm}{Theorem}[section]

	\newtheorem{defn}{Definition}[section]
	\newtheorem{lem}[thm]{Lemma}
	\newtheorem{prop}[thm]{Proposition}
	
	\newtheorem{cor}[thm]{Corollary}
	
	\newtheorem*{rmk}{Remark}
	
	\newtheorem*{pf}{Proof}
	
	\newcommand{\R}{\mathbb{R}}

	\newcommand{\T}{\mathbb{T}}

	\newcommand{\curl}{\text{curl }}
	\newcommand{\dive}{\text{div }}
	\newcommand{\Om}{\Omega}
	\newcommand{\q}{\quad}
	\newcommand{\p}{\partial}
	
	\newcommand{\DD}{\mathcal{D}}

	\newcommand{\nab}{\nabla}
	
	\newcommand{\lap}{\Delta}

	\newcommand{\di}{\text{div}\,}

	\newcommand{\PP}{\mathcal{P}}

    \newcommand{\nn}{\nabla^{\kappa}}
    \newcommand{\ak}{A^{\kappa}}
        \newcommand{\ve}{\varepsilon}
    
    \newcommand{\tphi}{\tilde{\phi}}
     \newcommand{\tPsi}{\tilde{\Psi}}
	\numberwithin{equation}{section}
	\usepackage{xcolor}

     \newcommand{\pb}{\bar\p}
     
\begin{document}
	\title{Local Well-Posedness of the Motion of Inviscid Liquid Crystals with a Free Surface Boundary}
	\date{}
	\author{Chenyun Luo\thanks{Department of Mathematics, The Chinese University of Hong Kong, Shatin, NT, Hong Kong. Email: \texttt{cluo@math.cuhk.edu.hk}. CL is partially supported by the Hong Kong RGC grants CUHK--14302922 and CUHK--14304424.}\,\,\,\,and Hang Yu\thanks{Department of Mathematics, The Chinese University of Hong Kong, Shatin, NT, Hong Kong. Email: \texttt{hyu@math.cuhk.edu.hk}.}}
	\maketitle
    \begin{abstract}
        In this article, we prove the local well-posedness of the free-boundary Lin--Liu equations describing the motion of inviscid nematic liquid crystals in the presence of surface tension in Lagrangian coordinates. It is well known that a priori energy estimates alone are insufficient for establishing local existence in free-boundary problems involving inviscid fluid equations, primarily due to the loss of symmetry in the linearized equations. The main challenge is to develop an effective approximate system of equations that is asymptotically consistent with the free-boundary Lin–Liu model expressed in the Lagrangian coordinates. This system must accurately capture the coupling between the fluid motion and the harmonic heat flows within the interior, as well as the regularity of the moving boundary. 
    \end{abstract}
\section{Introduction}
We study the motion of inviscid nematic liquid crystals occupying a moving domain $\mathcal{D}_t\in \mathbb{R}^3$ described by the simplified Ericksen--Leslie equations (also known as Lin--Liu equations): 
\begin{equation}\label{eu1}
\begin{cases}
        \p_{t}u+(u\cdot\nabla)u+\nabla P=-\nabla\cdot(\nabla d\odot\nabla d) \quad\text{in $\mathcal{D}$,}\\\nabla\cdot u=0\quad\text{in $\mathcal{D}$,}\\\p_{t}d+(u\cdot\nabla)d-\lap d=d|\nabla d|^2\quad\text{in $\mathcal{D}$,}
\end{cases}
    \end{equation}
    Here, $\mathcal{D}=\cup_{0\leq t\leq T}\{t\}\times\DD_{t}$,  $u:\DD\to\R^{3}$ is the velocity field, $d:\DD\to\mathbb{S}^2$ is the unit orientation vector fields of nematic liquid crystals, where $\nab d\cdot \nab d$ is a $2$-tensor with components $(\nabla d\odot\nabla d)_{ij}:=\nabla_{i}d\cdot\nabla_{j}d$, and $P:\DD\to \R$ is the pressure. 
    We equip \eqref{eu1} the initial conditions: 
    \begin{equation}
    \mathcal{D}_t|_{t=0} = \mathcal{D}_0, \quad 
    u(0,\cdot)=u_{0},\quad\text{and }\quad  
    d(0,\cdot)=d_{0}, 
\end{equation}
as well as the boundary conditions: 
\begin{equation}\label{BC}
    \begin{cases}
        (\p_{t}+v\cdot\nabla)|_{\p \mathcal{D}}\in \mathcal{T}(\p\mathcal{D}),\\P=\sigma H\quad\text{on $\p\mathcal{D}$},\\
        \nabla_{n}d=0 \quad\text{on $\p\mathcal{D}$},
    \end{cases}
\end{equation}
where $\p\DD=\cup_{0\leq t\leq T}\{t\}\times\p\DD_{t}$. 
In \eqref{BC}, $\mathcal{T}(\p \DD)$ is the tangent bundle of $\p \DD$, and $n$ is the outward-pointing unit normal of $\p \DD$. Also, we denote by $H$ the twice of mean curvature of the moving boundary, and by $\sigma>0$ the surface tension coefficient.

For the sake of simplicity and clean notations, we will only consider the case
\begin{align}\label{eu5}
    \DD_{0}=\T^{2}\times (0,1)
\end{align}
to avoid using local coordinate charts in a more general case. We will use $\p\DD^{+}(t)$ to represent the top moving boundary, and $\p\DD^{-}(t)=\{x_3=0\}$ to denote the fixed bottom. In other words, $\p\DD_t=\p\DD^{+}(t)\cup\p\DD^{-}(t)$. 
Additionally, we prescribe the natural slip boundary condition of the velocity on the fixed bottom, i.e., 
\begin{align}
u_{3}=0,\quad  \text{on } \p\DD^{-}(t).
\end{align}
\subsection{History and background}

Free boundary problems in fluid mechanics, nematic liquid crystals, and harmonic heat flow have been extensively studied over the past decades. Here, we mention only a few works relevant to our problem.

The system \eqref{eu1} was first introduced by Lin in \cite{lin1989nonlinear, lin1991nematic} as a simplified version of the Ericksen-Leslie model proposed by Ericksen \cite{ericksen1962hydrostatic} and Leslie \cite{leslie1968some}, which mathematically retains some basic nonlinearity structures. In \cite{lin1995nonparabolic, lin1996partial, lin2000existence}, Lin and Liu initiated the mathematical analysis of the system which established the global existence of the approximate solution in dimensions $2$ and $3$.  Since then, many research works have been established on its well-posedness, see \cite{hineman2013well, hong2011global, hong2012global, huang2014regularity, lei2014remarks, li2012global, lin2010uniqueness, lin2016global, wang2011well, wang2014global, xu2012global} and the references therein. We also refer to \cite{lin2014recent} for a survey of the mathematical developments in the theory of liquid crystals. Readers may also consult \cite{hong2014blow}, \cite{huang2012blow}, \cite{huang2016finite}, \cite{lai2022finite}, and \cite{lin2023nematic} for results concerning singularities and blow-up in related topics.

In addition, free-boundary problems in fluid mechanics have been extensively studied over the past few decades. For the incompressible Euler equation, earlier works mainly focused on irrotational fluids, also known as the water wave problem. In two dimensions, Nalimov \cite{nalimov1974cauchy} was the first to establish local well-posedness without surface tension for small initial data. See also Craig \cite{craig1985existence} and Yosihara \cite{yosihara1982gravity}. Under the Talor-sign condition, in dimensions $2$ and $3$, Wu \cite{wu1997well, wu1999well} developed the local existence for general initial data. The Taylor-sign condition is crucial in the case of $\sigma=0$. A well-known work by Ebin \cite{ebin1987equations} demonstrated the ill-posedness in the absence of surface tension. During the past two decades, water waves have been studied excessively, both in global and local solutions, for cases $\sigma = 0$ and $\sigma > 0$. See \cite{alazard2015global, berti2023birkhoff, berti2024hamiltonian, deng2017global, deng2025wave, germain2012global, germain2015global, ionescu2018global, ionescu2015global, ionescu2018global, wang2018global, wu2009almost, wu2011global}, etc. We refer readers to \cite{ionescu2018recent} for a survey on recent developments of the water waves.

On the other hand, for rotational initial data, Christodoulou and Lindblad \cite{christodoulou2000motion} established the a priori estimate in the case when $\sigma=0$ under the Taylor sign condition, followed by the well-posedness result by Lindblad \cite{lindblad2005well} using the Nash-Moser iteration. Later, Coutand and Shkoller \cite{coutand2007well} proved the local well-posedness in both $\sigma>0$ and $\sigma=0$ cases. The readers may refer to Zhang and Zhang \cite{zhang2008free}, Shatah and Zeng \cite{shatah2008geometry, shatah2008priori, shatah2011local}, Wang--Zhang--Zhao--Zheng \cite{wang2021local}, and Ifrim--Pineau--Tataru--Taylor \cite{ifrim2025sharp} for related topics. Additionally, we refer to \cite{beale1981initial, beale1984large, guo2013almost, guo2013decay, guo2013local, wang2020anisotropic} for the free-boundary problem in viscous fluids whose motion is described by the Navier--Stokes equations. 

In recent years, significant progress has been made in the study of free boundary problems coupled with other equations. It is worth noting the free boundary Magnetohydrodynamics (MHD) problem. Hao and Luo \cite{hao2014priori, hao2021well} initially established a priori estimates and linearized local well-posedness for ideal incompressible MHD. Gu and Wang \cite{gu2019construction} demonstrated local well-posedness for the case of vanishing surface tension. Luo and Zhang \cite{luo2020regularity} proved minimal regularity estimates of $H^{2.5+}$ for a small fluid domain. Gu, Luo, and Zhang \cite{gu2024local} were the first to establish local well-posedness for the non-zero surface tension case. Additionally, references are made to the works of Liu and Xin \cite{liu2025free} and Ifrim--Pineau--Tataru--Taylor \cite{ifrim2024sharp} for related topics. For compressible MHD, we recommend the works by Trakhinin and Wang \cite{trakhinin2018well, trakhinin2021well} for ideal MHD with and without surface tension, and Zhang \cite{zhang2019priori, zhang2023local} for resistive MHD and its corresponding inviscid limit.

Most studies on nematic liquid crystals primarily focus on cases involving viscosity and fixed boundary conditions, while research on inviscid liquid crystals is relatively limited. Additionally, the local well-posedness of free boundary problems in complex fluids is quite complex, especially when surface tension is present, which complicates the process of closing the energy estimates. In this article, we establish the local well-posedness and derive energy estimates for the free boundary problem of inviscid liquid crystals without any loss of regularity. This work serves as a crucial step toward further exploration of general nematic liquid crystals influenced by viscosity, both experimentally and theoretically.

\subsection{The Lagrangian description}
We reformulate the equations into Lagrangian coordinates to study a fixed domain question. Let $\Omega:=\T^2\times(0,1)$ be a fixed domain in $\R^3$, and $y=(y_1,y_2,y_3)\in \Omega$ be Lagrangian variables. Let $\eta(y,t):\Omega\to\DD_{t}$ be the solution of 
\begin{align*}
    \p_{t}\eta(y,t)=u(\eta(y,t),t)
\end{align*}
with initial condition
\begin{align*}
    \eta(y,0)=Id.
\end{align*}
We further set
\begin{align*}
\begin{aligned}
    v(y,t)&:=u(\eta(y,t),t),\\ q(y,t)&:=p(\eta(y,t),y), \\\phi(y,t)&:=d(\eta(y,t),t)\\A(y,t)&:=[D\eta]^{-1}
    \end{aligned}
\end{align*}
to denote the material velocity, pressure, orientation vector field of nematic liquid crystal and inverse of the deformation tensor in Lagrangian coordinates respectively. Since the fluid is incompressible, the Jacobian $J=\det[D\eta]=1$ in $\Omega\times[0,T]$. We introduce the moving top boundary and fixed bottom boundary respectively by $\Gamma=\{y_3=1\}$ and $\Gamma_{0}=\{y_3=0\}$, and set $N$ to be the normal outward direction of $\p\Omega=\Gamma\cup\Gamma_{0}$. 

On the fixed domain $\Omega$, the equations \eqref{eu1}-\eqref{eu5} can be reformulated as
\begin{equation}\label{la1}
    \begin{cases}
        \eta= e+\int_{0}^{t}v \quad\text{in $\Omega\times[0,T]$},\\
        \p_{t}v_{i}+A_{i}^{j}\p_{j}q=-A_{s}^{l}\p_{l}(A_{s}^{m}\p_{m}\phi_{n}A_{i}^{k}\p_{k}\phi_{n})\quad\text{in $\Omega\times[0,T]$},\\A_{i}^{j}\p_{j}v_{i}=0\quad\text{in $\Omega\times[0,T]$},\\\p_{t}\phi_{i}-A_{j}^{k}\p_{k}(A_{j}^{m}\p_{m} \phi_{i})=\phi_{i}(A_{s}^{l}\p_{l}\phi_{n}A_{s}^{m}\p_{m}\phi_{n})\quad\text{in $\Omega\times[0,T]$},\\ v^{3}=0\quad \textit{on $\Gamma_{0}\times[0,T]$},\\q \frac{A_{i}^{k}N_{k}}{|A_{j}^{l}N_{l}|}=-\sigma \lap_{h}(\eta) \quad\text{on $\Gamma\times[0,T]$},\\A_{i}^{l}N_{l}A_{i}^{k}\p_{k}\phi=0 \quad\text{on $\p\Omega\times[0,T]$},\\
        (\eta,v,\phi)=(e,u_{0},d_{0}) \quad\text{on $\Omega\times\{t=0\}$},\\
    \end{cases}
\end{equation}
where $e$ is the identity map, $N^{T}=(0,0,1)$ denotes the unit normal to $\Gamma$, $h_{\alpha\beta}=\p_{\alpha}\eta\cdot \p_{\beta}\eta$ represents the induced metric on $\Gamma$, $\lap_{h}$ serves as the corresponding boundary Laplacian operator and the Einstein summation convention is used. $v_{i}=v^{i}$ and $\phi_{i}=\phi^{i}$ are understood as the i-th competent of $u$ and $d$ in Euler setting. To clear the symbol, we a little abuse the notation by introducing the differential operator $\nabla_{i}=A_{i}^{j}\p_{j}$, $\nabla=(\nabla_1,\nabla_2,\nabla_3)$ and $\lap_{g}=\delta^{jl}A_{j}^{k}\p_{k}(A_{l}^{m}\p_{m})$. Since the orientation vector fields which we study in this paper only occur in pairs, the index of corresponding fields will be suppressed to simplify the expression if it will not cause any confusion. (i.e. Use $\nabla_{k}\phi^{i}\nabla_{l}\phi^{j}\delta_{ij}=\nabla_{k}\phi\nabla_{l}\phi$ to save notation.) Owing to the property of the incompressible flow, $J=1$. We obtain $\lap_{g}=\p_{i}(g^{ij}\p_{j})$. Simultaneously, we use the symbol div and curl in following sense:
\begin{align*}
\begin{aligned}
     \text{div }v=\delta^{ij}\p_{i}v_{j},\\
     \text{$\text{div}_{A}$ }v=\nabla_{A}\cdot v=\delta^{ik}A_{i}^{j}\p_{j}v_{k},\\
     (\text{curl }v)^i=\ve^{ijk}\p_{j}v_{k},\\
     (\text{$\text{curl}_{A} $ }v)^{i}=\nabla_{A}\times v=\ve^{ijk}A_{j}^{l}\p_{l}v_{k}
\end{aligned} 
\end{align*}
Under the above setting, we rewrite the equations \eqref{la1} into:
\begin{equation}\label{la2}
    \begin{cases}
        \eta= e+\int_{0}^{t}v \quad\text{in $\Omega\times[0,T]$},\\
        \p_{t}v+\nabla q=-\nabla\cdot(\nabla\phi\odot\nabla\phi)\quad\text{in $\Omega\times[0,T]$},\\ \nabla_{A}\cdot v=0\quad\text{in $\Omega\times[0,T]$},\\\p_{t}\phi-\lap_{g} \phi=\phi|\nabla\phi|^2\quad\text{in $\Omega\times[0,T]$},\\v^{3}=0 \quad \text{on $\Gamma_{0}\times[0,T]$},\\q\frac{A^{T}N}{|A^{T}N|}=-\sigma\lap_{h}(\eta) \quad\text{on $\Gamma\times[0,T]$},\\\nabla_{N}\phi=0 \quad\text{on $\Gamma\times[0,T]$},\\
        (\eta,v,\phi)=(e,u_{0},d_{0}) \quad\text{on $\Omega\times\{t=0\}$},\\
    \end{cases}
\end{equation}
where $h$ is the boundary metric induced by the metric $g$. Since $N=(0,0,1)$ on the moving boundary, we can simply write $A_{j}^{3}A_{j}^{i}\p_{i}\phi=g^{3i}\p_{i}\phi=0$ on $\Gamma$.

\subsection{Main theorem}
\thm [Main theorem]
Let $\sigma>0$ be fixed. Assuming that $u_{0}\in H^{5.5}(\Omega)$ with $\dive u_{0}=0$, and 
$$
\p_t^k d(0) \in H^{5.5-k}(\Omega), \quad k=0,1,2,3,\quad \p_{t}^{4}d(0)\in H^{1}(\Omega),
$$
 with $\nabla_{N}d_{0}=0$. Then there exists $T>0$, such that \eqref{la2} admits a unique solution $(v,q,\phi)$, satisfying  
\begin{align}\label{maintheorem1}
\begin{aligned}
     \sup_{t\in [0,T]}&(\|\p \mathcal{D}_{t}\|_{H^{6.5}}+\|v(t,\cdot)\|_{H^{5.5}(\Omega)}+\|\phi(t,\cdot)\|_{H^{5.5}(\Omega)}+\|q(t,\cdot)\|_{H^{4.5}(\Omega)}\\
     &+\|v_{tttt}(t,\cdot)\|_{L^{2}(\Omega)})+\|\p_{t}^{5}\phi\|_{L^{2}([0,T];L^{2}(\Omega))}\leq M_{0},
\end{aligned}   
\end{align}
where $M_{0}$ is a constant depending only on the
initial data. 

\rmk
The framework established by Coutand and Shkoller in \cite{coutand2007well} focuses on proving the local well-posedness of the free boundary problem coupled with the Euler equation. The existence result, dependent on initial data, is not optimal. Therefore, our objective does not aim for an optimal outcome with respect to the initial data $u_{0}$ and $d_{0}$ in this paper. For more refined results in the Lagrangian setting, readers may consult \cite{aydin2024construction} and \cite{disconzi2019priori}.

\subsection{Main difficulties and novelty of this article}
\paragraph*{The $L^2$-energy identity and a key cancellation scheme.}
We begin with the following energy identity, which is derived from testing the second equation of \eqref{la2} by $v$ and then integrating in $\Omega$:
\begin{align}\label{basigesti1}
\begin{aligned}    \frac{d}{dt}\frac{1}{2}\int_{\Omega}|v|^{2}+\sigma\int_{\Gamma}\sqrt{h}h^{ij}\p_{j}\eta^{k} \p_{i}v^{k}&=\underbrace{\int_{\Omega}p\nabla\cdot v}_{\dive v=0}+\int_{\Omega}\nabla_{i}\phi\nabla_{j}\phi\nabla_{i} v^{j}.
\end{aligned}    
\end{align}

Note that we can utilize the symmetry on the boundary:
\begin{align}
\begin{aligned}
      \sigma\int_{\Gamma}\frac{d}{dt}(\sqrt{h}h^{ij})\p_{i}\eta^{k}\p_{j}\eta^{k}&=\sigma\int_{\Gamma}\sqrt{h}(\frac{1}{2}h^{ij}h^{kl}-h^{lj}h^{ik})\p_{t}(h_{kl})h_{ij}\\&=\sigma\int_{\Gamma}\sqrt{h}h^{kl}\p_{t}(h_{kl})-\sqrt{h}h^{lk}\p_{t}(h_{kl})\\&=0.
\end{aligned} 
\end{align}
Here, $h_{ij}=\p_{i}\eta^{k}\p_{j}\eta^{k}$ is the induced boundary metric. We can rewrite the basis energy identity \eqref{basigesti1} as:
\begin{align}\label{basigesti2}
   \frac{d}{dt}\frac{1}{2}\int_{\Omega}|v|^{2}+\sigma\frac{1}{2}\frac{d}{dt}\int_{\Gamma}\underbrace{h^{ij}\sqrt{h}\p_{j}\eta^{k} \p_{i}\eta^{k}}_{2\sqrt{h}}&=\int_{\Omega}\nabla_{i}\phi\nabla_{j}\phi\nabla_{i} v^{j}. 
\end{align}

The calculation involves some geometric identities, see Section 2 for details. The second term appearing on the RHS of \eqref{basigesti2} needs to be cancelled. By testing $\p_{t} \phi$ on the fourth equation of \eqref{la2}
\begin{align}
    \int_{\Omega}|\p_{t}\phi|^{2}+\int_{\Omega}\nabla_{i} v^{k}\nabla_{k}\phi\nabla_{i}\phi+\frac{1}{2}\frac{d}{dt}\int_{\Omega}|\nabla\phi|^{2}=\int_{\Omega}\phi\p_{t}\phi|\nabla\phi|^{2}.
\end{align}
Owing to the identity $|\phi|=1$, we have $\phi\p_{t}\phi=0$. Thus, we obtain the $L^2$-energy identity of $\phi$: 

\begin{align}\label{basigesti3} 
    \frac{d}{dt}\frac{1}{2}\int_{\Omega}|\nabla\phi|^{2}+\int_{\Omega}|\p_{t}\phi|^{2}=-\int_{\Omega}\nabla_{i}v^{k}\nabla_{k}\phi\nabla_{i}\phi.
\end{align}
Further, there would be a cancellation on the RHS of \eqref{basigesti2} and \eqref{basigesti3}, which gives us
\begin{align}
    \frac{d}{dt}\frac{1}{2}\int_{\Omega}|v|^{2}+\sigma\frac{d}{dt}\int_{\Gamma}\sqrt{h}+ \frac{d}{dt}\frac{1}{2}\int_{\Omega}|\nabla\phi|^{2}+\int_{\Omega}|\p_{t}\phi|^{2}=0.
\end{align}

Additionally, in the higher-order energy, we will derive a similar energy estimate using cancellation. At the same time, we will treat the commutator terms as lower-order contributions. We denote $\p_{t}$, $\p_{1}$ and $\p_{2}$ by $\pb_{A}$. Taking $\pb_{A}^{4}$ on the second equation \eqref{la2}, the higher order energy estimate reads  
\begin{align}\label{basigesti5}
    \frac{1}{2}\frac{d}{dt}\int_{\Omega}|\pb_{A}^{4}v|^{2}+\sigma\int_{\Gamma}\pb_{A}^{4}(\sqrt{h}h^{ij}\p_{j}\eta^{k}) \p_{i}\pb_{A}^{4}v^{k}&=\int_{\Omega}\pb_{A}^{4}p\nabla\cdot \pb_{A}^{4}v+\int_{\Omega}\pb_{A}^{4}\nabla_{i}(\nabla_{i}\phi\nabla_{j}\phi) \pb_{A}^{4}v^{j}+\cdots,
\end{align}
where we used $\cdots$ to represent the lower-order terms. The terms $\sigma\int_{\Gamma}\pb_{A}^{4}(\sqrt{h}h^{ij}\p_{j}\eta^{k}) \p_{i}\pb_{A}^{4}v^{k}$ and $\int_{\Omega}\pb_{A}^{4}p\nabla\cdot \pb_{A}^{4}v$ have been estimated in \cite{coutand2007well} via a subtle cancellation scheme. See Section 6.  In this article, we derive a new cancellation relationship for the last term of \eqref{basigesti5}. The leading order term can be written as:
\begin{align}\label{basigesti6}
    \int_{\Omega}\pb_{A}^{4}\nabla_{i}(\nabla_{i}\phi\nabla_{j}\phi) \pb_{A}^{4}v^{j}=-\int_{\Omega}\pb_{A}^{4}(\nabla_{i}\phi)\nabla_{j}\phi \nabla_{i}\pb_{A}^{4}v^{j}-\int_{\Omega}\nabla_{i}\phi\pb_{A}^{4}(\nabla_{j}\phi) \nabla_{i}\pb_{A}^{4}v^{j}.
\end{align} 
The first term on the RHS of \eqref{basigesti6} can be controlled  directly through the parabolic estimate of $\lap_{A}\phi$. The cancellation scheme for the second term, $\int_{\Omega}\nabla_{i}\phi\pb_{A}^{4}(\nabla_{j}\phi) \nabla_{i}\pb_{A}^{4}v^{j}$, is a novel contribution of this article and serves as the key lemma in Section 4.

\paragraph*{No additional regularity for $\phi$ due to the moving boundary.}
The cancellation scheme is crucial because the regularity of the domain affects the parabolic estimate of \(\phi\). While the regularizing effect of surface tension contributes additional derivatives in the \(L^{2}\) energy estimate, it is still insufficient to fully close the energy estimate. This results in a loss of \(\frac{1}{2}\) derivatives. Specifically, in the Lagrangian framework, the term involving \(\lap_{A}\phi\) hinders the application of a standard parabolic estimate. As a result, the leading order of \(\phi\) is no longer two derivatives higher than \(\partial_{t}\phi\) and \(v\). Instead, the leading-order energy of \(\phi\) is limited by the regularity of \(\eta\).


Also, it is important to note that the regularizing effect does not yield additional derivatives. In other words, the regularity of $\eta$ ultimately aligns with that of $v$. For examples of this, refer to \cite{shatah2008geometry, shatah2008priori}. In our case, this will be established through an a posteriori estimate. However, the regularizing effect is understood in the context of $L^{2}$ estimates.

\paragraph*{The approximate system.}
We demonstrate the existence and uniqueness of solutions to the approximate system, along with the associated regularity. While our approach shares similarities with that of \cite{coutand2007well}, there are two key differences. First, we are required to construct a Picard iteration to establish the existence of a solution for the harmonic heat flow. Second, by mollifying the terms involving \(\phi\) and initial data \(u_{0}\), we introduce an approximate system that extends the framework in \cite{coutand2007well}, with the added feature of an \(L_t^\infty H^{12.5}\) estimate for \(v\).
\subsection{Organization of this article}
We devote Section 3 to constructing the approximate system and establishing the existence and uniqueness of its solutions. 
In Section 4, we establish a uniform estimate for $\phi$ with respect to $\kappa$. A key cancellation is derived in this section, which will be utilized in Section 6. In Sections 5 and 6, most of the calculations follow the framework of \cite{coutand2007well}. However, there are two notable differences. First, we provide a more precise estimate for $\|\p_{t}^{3} p\|_{1}$. Second, the key cancellation derived in Section 4 is applied to the $\p_{t}^{4}$-estimate in Section 6. Since the $\p_{t}^{4}$-estimate in \cite{coutand2007well} involves subtle dependencies on other energy norms, we provide detailed proofs for the reader's convenience and verification.

The following list of notations will be frequently used throughout this manuscript.
\paragraph*{List of notations}
\begin{enumerate}
    \item $\Omega=\T^{2}\times (0,1)$. $\Gamma=\T^{2}\times\{x_{3}=1\}$ is the free boundary and $\Gamma_{0}=\T^{2}\times \{x_{3}=0\}$.
    \item The Sobolev norms in $\Omega$: $\|\cdot\|_{H^{s}(\Omega)}=\|\cdot\|_{s}$. For $1\leq p\leq\infty$, $\|\cdot\|_{L^{p}(0,T;H^{s}(\Omega))}=\|\cdot\|_{L_{t}^{p}H^{s}(\Omega)}$.
    \item The Sobolev norms on $\p\Omega$: $\|\cdot\|_{H^{s}(\Gamma)}=|\cdot|_{s}$. For $1\leq p\leq\infty$, $\|\cdot\|_{L^{p}(0,T;H^{s}(\Omega))}=\|\cdot\|_{L_{t}^{p}H^{s}(\Gamma)}$.
    \item $\PP(\cdot)$: a generic non-decreasing continuous function.
    \item $\pb=\p_{1}$ or $\p_{2}$. $\bar{\lap}=\p_{1}^{2}+\p_{2}^{2}$.
\end{enumerate}

\section{Preliminary results}
\subsection{Sobolev space and the dual space}
In this paper, we will frequently use fractional Sobolev space. Since our domain is smooth, most common definitions on fractional Sobolev space coincide here. The reader may refer to \cite{adams2003sobolev}, \cite{leoni2023first}, \cite{lions2012non}, \cite{tartar2007introduction} and the references therein. When $s>0$, the space $H^{-s}$ denotes the dual space of $H^{s}$. We will use $(H_{0}^{1})^{*}$ to represent the dual space of $H_{0}^{1}$. We need the Kato-Ponce inequality (the fractional Leibnitz rule) to derive some inequalities later.
\lem[Kato-Ponce inequality]
Let $\Omega\subset \R^{3}$ be a smooth bounded domain. For $s\geq 0$, $2\leq p_{1},q_{2}<\infty$ satisfying the relationship $\frac{1}{2}=\frac{1}{p_{1}}+\frac{1}{p_{2}}=\frac{1}{q_{1}}+\frac{1}{q_{2}}$, $f\in W^{s,p_{1}}(\Omega)\cap L^{q_{1}}(\Omega)$ and $g\in L^{p_{2}}(\Omega)\cap W^{s,q_{2}}(\Omega)$, $fg\in H^{s}(\Omega)$ with the following estimate
\begin{align}
    \|fg\|_{H^{s}(\Omega)}\lesssim \|f\|_{W^{s,p_{1}}(\Omega)}\|g\|_{L^{p_{2}}(\Omega)}+\|f\|_{L^{q_{1}}(\Omega)}\|g\|_{W^{s,q_{2}}(\Omega)}.
\end{align}
In particular, owing to the Sobolev inequality, we have the following estimates
\begin{align}
    \|fg\|_{H^{0.5}(\Omega)}\lesssim\|f\|_{H^{0.5}(\Omega)}\|g\|_{H^{1.5+}(\Omega)},
\end{align}
and
\begin{align}
    \|fg\|_{H^{0.5}(\Omega)}\lesssim\|f\|_{H^{1}(\Omega)}\|g\|_{H^{1.5}(\Omega)}.
\end{align}
\pf
See \cite{grafakos2014kato}, \cite{kato1988commutator} and \cite{li2019kato} with appropriate modifications. (Sobolev extension with Kato-Ponce or Kato-Ponce with constrain) Actually, we only need the last two inequalities in our estimate, a rough version. Thus, one may also prove it via interpolation.

We denote the horizontal derivative by $\pb = \p_{1}$ or $\p_{2}$. In the context of the Lions-Magenes space, it is important to note that $\p$ is not, in general, a continuous operator from $H^{\frac{1}{2}}(\Omega)$ to $H^{-\frac{1}{2}}(\Omega)$. However, in our estimate, we will only encounter $\pb$ in cases where ``integrating 0.5 derivative by part" is required. To be precise, we need the following lemma
\lem
For any $F\in H^{\frac{1}{2}}(\Omega)$, we have
\begin{align}
    \|\pb F\|_{H^{-\frac{1}{2}}(\Omega)}\leq C\|F\|_{H^{\frac{1}{2}}(\Omega)},
\end{align}
for some generic constant $C$.
\pf 
See \cite{coutand2007well} and \cite{coutand2010simple}. 
\rmk
Instead of relying on this lemma, we could alternatively define $\pb = (1 - \bar{\lap})^{\frac{1}{2}}$, which provides a more convenient framework for handling a 0.5 derivative. However, this approach would require several additional and complex commutator estimates for the boundary analysis. For the sake of simplicity and convenience, we will continue to define $\pb = \p_{1}$ or $\p_{2}$.

\subsection{Horizontal convolution by layers}
While discussing functions on $\T^{2}$, we always view them as periodic functions on $\R^{2}$, with norms evaluated over a single period. Except for $\eta$, the functions $v$, $q$, and $\phi$ are periodic in $\R^{2}$. Although $\eta$ itself is not periodic, its derivative, $\p\eta$, is periodic. Throughout this paper, the norm of $\eta$ will always be taken over one period.
\defn[Horizontal mollifier]
For a fixed $\kappa_{0} > 0$ which is sufficiently small, we choose $\kappa \in (0, \kappa_{0})$. Let $f$ be an $L^{p}$ function ($1 \leq p \leq \infty$) defined on $\T^{2} \times (0,1)$, or a function such as $\eta \in L_{\text{loc}}^{p}(\R^{2})$
\begin{align}   \Lambda_{\kappa}f(x_{1},x_{2},x_{3})=\int_{\R}\int_{\R}\rho_{\kappa}(x_{1}-y_{1},x_{2}-y_{2})f(y_{1},y_{2},x_{3})dy_{1}dy_{2},
\end{align}
where $\rho_{\kappa}=c\frac{1}{\kappa^{2}}\rho(\frac{y_{1}}{\kappa},\frac{y_{2}}{\kappa})$ is a standard two dimensional smooth Friedrich mollifier (radial symmetric and $\|\rho_{\kappa}\|_{L^{1}(\R^{2})}=1$).

The following properties about the mollifier are well-known.
\lem
Suppose $f$ is a periodic function such that the norms in the following inequalities are well-defined.
\begin{align}
\begin{aligned}
\|\Lambda_{\kappa}f\|_{L^{2}(\T^{2})}&\lesssim \| f\|_{L^{2}(\T^{2})}\\\|\Lambda_{\kappa}f\|_{H^{s}(\T^{2})}&\lesssim \| f\|_{H^{s}(\T^{2})}\\\|\Lambda_{\kappa}f\|_{H^{s}(\T^{2})}&\lesssim_{\kappa} \| f\|_{L^{2}(\T^{2})}\\
  \|\Lambda_{\kappa}f-f\|_{L^{2}(\T^{2})}&\lesssim \kappa\| f\|_{H^{1}(\T^{2})}\\
  \|\Lambda_{\kappa}f-f\|_{L^{4}(\T^{2})}&\lesssim \kappa\| f\|_{H^{1.5}(\T^{2})}\\\|\Lambda_{\kappa}f-f\|_{L^{2}(\T^{2})}&\lesssim \sqrt{\kappa}\| f\|_{H^{0.5}(\T^{2})}\\\|\Lambda_{\kappa}f-f\|_{H^{0.5}(\T^{2})}&\lesssim \sqrt{\kappa}\| f\|_{H^{1}(\T^{2})}\\
    \|\Lambda_{\kappa}f-f\|_{L^{\infty}(\T^{2})}&\lesssim \sqrt{\kappa}\|f\|_{H^{1.5}(\T^{2})}\\
    \|\Lambda_{\kappa}f-f\|_{L^{\infty}(\T^{2})}&\lesssim \kappa\|f\|_{H^{2.5}(\T^{2})}
\end{aligned} 
\end{align}
\pf
See \cite{coutand2007well} and \cite{gu2024local}. We provide a sketchy proof for the last five inequalities. 
Using the definition, (taking $\tilde{f}$ as $f$ on $\R^{2}$)
\begin{align}
\begin{aligned}
    (\Lambda_{\kappa}\tilde{f}-\tilde{f})(x)&=\int_{\R^{2}}\rho_{\kappa}(y)(\tilde{f}(x-y)-\tilde{f}(x))dy\\&=\int_{\R^{2}}\rho(y)(\tilde{f}(x-\kappa y)-\tilde{f}(x))dy\\&=\int_{\R^{2}}\rho(y)\int_{0}^{1}\frac{d}{dt}\tilde{f}(x-t\kappa y)dtdy\\&=\int_{\R^{2}}-\kappa y\rho(y)\int_{0}^{1}D\tilde{f}(x-t\kappa y)dtdy
\end{aligned}   
\end{align}
Invoking Minkowski inequality and invariance of translation of $L^{p}$, we have
\begin{align}
\begin{aligned}
     \| (\Lambda_{\kappa}\tilde{f}-\tilde{f})(x)\|_{L^{2}(\T^{2})}&\lesssim\int_{0}^{1}\int_{\R^{2}}\kappa y\rho(y)\|D\tilde{f}(x-t\kappa y)\|_{L_{x}^{2}}dydt\\&\lesssim\kappa\|D\tilde{f}\|_{L^{2}(\T^{2})},
\end{aligned}   
\end{align}
\begin{align}
\begin{aligned}
     \| (\Lambda_{\kappa}\tilde{f}-\tilde{f})(x)\|_{L^{4}(\T^{2})}&\lesssim\int_{0}^{1}\int_{\R^{2}}\kappa y\rho(y)\|D\tilde{f}(x-t\kappa y)\|_{L_{x}^{4}}dydt\\&\lesssim\kappa\|D\tilde{f}\|_{L^{4}(\T^{2})},
\end{aligned}   
\end{align}
\begin{align}
\begin{aligned}
     \| (\Lambda_{\kappa}\tilde{f}-\tilde{f})(x)\|_{L^{\infty}(\T^{2})}&\lesssim\int_{0}^{1}\int_{\R^{2}}\kappa y\rho(y)\|D\tilde{f}(x-t\kappa y)\|_{L_{x}^{\infty}}dydt\\&\lesssim\kappa\|D\tilde{f}\|_{L^{\infty}(\T^{2})}.
\end{aligned}   
\end{align}
Therefore, the Sobolev embedding theorem implies that
\begin{align}
     \| (\Lambda_{\kappa}f-f)(x)\|_{L^{2}(\T^{2})}\lesssim\kappa\|f\|_{H^{1}(\T^{2})},
\end{align}
\begin{align}
     \| (\Lambda_{\kappa}f-f)(x)\|_{L^{4}(\T^{2})}\lesssim\kappa\|f\|_{H^{1.5}(\T^{2})},
\end{align}
\begin{align}
     \| (\Lambda_{\kappa}f-f)(x)\|_{L^{\infty}(\T^{2})}\lesssim\kappa\|f\|_{H^{2.5}(\T^{2})}.
\end{align}
Similarly, invoking the Morrey's estimate, (in our case $n=2$)
\begin{align}
\begin{aligned}
    |(\Lambda_{\kappa}\tilde{f}-\tilde{f})(x)|&\lesssim\int_{\R^{2}}\rho_{\kappa}(y)|\tilde{f}(x-y)-\tilde{f}(x)|dy\\&\lesssim\int_{\R^{2}}\rho_{\kappa}(y)\|D\tilde{f}\|_{L^{4}(\T^{2})}\|\frac{1}{r^{n-1}}\|_{L^{\frac{4}{3}}(B_{\kappa}(0))}dy\\&\lesssim\int_{\R^{2}}\rho_{\kappa}(y)\sqrt{\kappa}\|\tilde{f}\|_{H^{1.5}(\T^{2})}dy\\&\lesssim\sqrt{\kappa}\|\tilde{f}\|_{H^{1.5}(\T^{2})},
\end{aligned}    
\end{align}
Thus,
\begin{align}
    \|\Lambda_{\kappa}f-f\|_{L^{\infty}(\T^{2})}\lesssim\sqrt{\kappa}\|f\|_{H^{1.5}(\T^{2})}.
\end{align}
The inequality involving $H^{0.5}$ is an application of interpolation.
We verify this property via an equivalent definition of fractional Sobolev space (See \cite{adams2003sobolev},\cite{leoni2023first},\cite{lions2012non},\cite{tartar2007introduction})
\begin{align}
    \begin{aligned}
    (\Lambda_{\kappa}\tilde{f}-\tilde{f})(x)&=\int_{\R^{2}}\rho_{\kappa}(x-y)\frac{(\tilde{f}(y)-\tilde{f}(x))}{|y-x|^{\frac{3}{2}}}|y-x|^{\frac{3}{2}}dy\\&\lesssim(\int_{\R^{2}}|\rho_{\kappa}(x-y)|y-x|^{\frac{3}{2}}|^{2}dy)^{\frac{1}{2}}(\int_{\T^{2}}|\frac{(\tilde{f}(y)-\tilde{f}(x))}{|y-x|^{\frac{3}{2}}}|^{2}dy)^{\frac{1}{2}}\\&=(\int_{\R^{2}}\kappa\rho(x-y)|x-y|^{3}dy)^{\frac{1}{2}}(\int_{\T^{2}}\frac{|\tilde{f}(y)-\tilde{f}(x)|^{2}}{|y-x|^{3}}dy)^{\frac{1}{2}}
\end{aligned}   
\end{align}
Therefore,
\begin{align}
    \begin{aligned}
         \|\Lambda_{\kappa}\tilde{f}-\tilde{f}\|_{L^{2}}\lesssim \sqrt{\kappa}(\int_{\R^{2}}\int_{\R^{2}}\frac{|\tilde{f}(y)-\tilde{f}(x)|^{2}}{|y-x|^{3}}dydx)^{\frac{1}{2}}\lesssim\sqrt{\kappa}\|\tilde{f}\|_{H^{0.5}(\T^{2})}.
    \end{aligned}
\end{align}
Using a similar method, we can verify that
\begin{align}
      \|\Lambda_{\kappa}f-f\|_{H^{1}(\T^{2})}\lesssim \sqrt{\kappa}\|f\|_{H^{1.5}(\T^{2})}.
\end{align}
Hence, we can infer $\|\Lambda_{\kappa}f-f\|_{H^{0.5}(\T^{2})}\lesssim \sqrt{\kappa}\|f\|_{H^{1}(\T^{2})}$ from interpolation.
\rmk
Since we always calculate the norm in one period, a similar argument is also valid for $\eta$. 
\lem[Commutation type lemma]
Let $f\in H^{s}(\T^{2})$ for $s>1$ and $g\in L^{2}(\T^{2})$. Then there exists a constant C independent on $\kappa\in(0,\kappa_{0})$ s.t.
\begin{align}\label{commutationtypelemma}
    \|\Lambda_{\kappa}(fg)-f\Lambda_{\kappa}g\|_{L^{2}(\T^{2})}\leq C \kappa\|f\|_{H^{s+1}(\T^{2})}\|g\|_{L^{2}(\T^{2})}.
\end{align}
In particular, suppose $f\in H^{s}(\Omega)$ for $s>1$ and $g\in L^{2}(\Omega)$, the trace inequality implies that
\begin{align}
     \|\Lambda_{\kappa}(fg)-f\Lambda_{\kappa}g\|_{L^{2}(\Omega)}\leq C \kappa\|f\|_{H^{s+1.5}(\Omega)}\|g\|_{L^{2}(\Omega)}.
\end{align}
Further suppose $f\in W^{1,\infty}(\T^{2})$
\begin{align}\label{commutationtypelemma}
   \|\Lambda_{\kappa}(f\pb g)-f\Lambda_{\kappa}\pb g\|_{L^{2}(\T^{2})}\leq C \|f\|_{W^{1,\infty}(\T^{2})}\|g\|_{L^{2}(\T^{2})}.
\end{align}
\pf See \cite{coutand2007well} and \cite{coutand2010simple}.
\subsection{Elliptic estimates and trace estimates}
\lem[Quasi-linear elliptic estimate]
Let $s$ be an real number such that $1\leq s \leq 5.5$. Define an elliptic operator $L_{A}$ and corresponding Neumann boundary operator as
\begin{align}
    L_{A}&=\p_{i}(g^{ij}\p_{j}),\\
    B_{A}&=A_{i}^{j}A_{i}^{k}N_{j}\p_{k}
\end{align}
Let $f\in H^{s-2}(\Omega)$ and $g\in H^{s-1.5}(\Gamma)$. Suppose $u\in H^{s}(\Omega)$ is the solution of the following quasi-linear elliptic problem with Neumann condition.
\begin{align}
\begin{aligned}
     L_{A}u&=f\quad\textit{in $\Omega$},\\
    B_{A}u&=g\quad\textit{on $\p\Omega$}.
\end{aligned}  
\end{align}
Then tho following estimate holds:
\begin{align}\label{quasiellesti}
    \|u\|_{H^{s}(\Omega)}\lesssim_{s}\PP(\|\eta\|_{H^{5.5}(\Omega)})\|f\|_{H^{s-2}(\Omega)}+\|g\|_{H^{s-\frac{3}{2}}(\p\Omega)}+\|u\|_{L^{2}(\Omega)}.
\end{align}
\pf
A proof for the integer case can be found in \cite{ebenfeld2002L2}. The readers may derive the desired estimate through interpolation. Since our focus is on obtaining a solution over a short time interval, an alternative approach in our case is to consider:
\begin{align}
\begin{aligned}
     \lap u&= \p_{i}((\delta^{ij}-g^{ij})\p_{j}u)  +f\quad\textit{in $\Omega$},\\
     \p_{3}u&=(\delta^{3i}-g^{3i})\p_{i}u+g\quad \textit{on $\p\Omega$}.
\end{aligned} 
\end{align}
See \cite{disconzi2019priori} and the discussion therein.

We will utilize the following well-known curl-div decomposition.
\lem[Hodge elliptic decomposition]
Let $\Omega\subset\R^{3}$ be an $H^{r+1}$ domain with boundary $\Gamma$, $r>1.5$. If $F\in L^{2}(\Omega)$, $\curl F\in H^{s-1}(\Omega)$, $\di F\in H^{s-1}(\Omega)$, and $\pb F\cdot N|_{\Gamma}\in H^{s-\frac{3}{2}}(\Gamma)$ for $1\leq s\leq r$, then there exists a constant C depending only on $\Omega$ such that

\begin{align}\label{Hodgeelliptic}
    \|F\|_{H^{s}(\Omega)}\leq C(\|F\|_{L^{2}(\Omega}+\|\curl F\|_{H^{s-1}(\Omega)}+\|\di F\|_{H^{s-1}(\Omega)}+\|\pb F\cdot N\|_{H^{s-\frac{3}{2}}(\Gamma)}).
\end{align}
\pf See \cite{cheng2017solvability}.
\subsection{Geometric Identity}
\lem
Let $n$ be the unit outer normal to $\eta(\Gamma)$. We denote $\Pi=n\otimes n$ to be the normal projection operator. Let $\pb_{A}$ represent the $\pb$ or $\p_{t}$. We have the following geometric identities:
\begin{align}
    n&=\frac{A_{.}^{3}}{|A_{.}^{3}|},\\
    \p A_{i}^{j}&=-A_{i}^{l}\p_{l}\p \eta^{k}A_{k}^{j},\\
    \p_{t} A_{i}^{j}&=-A_{i}^{l}\p_{l}v^{k}A_{k}^{j},\\ \sqrt{h}&=|JA_{.}^{3}|,\\
    \Pi_{j}^{i}&=\delta^{i}_{j}-h^{\alpha\beta}\pb_{\alpha}\eta^{i}\pb_{\beta}\eta^{l}\delta_{jl},\\
\Pi_{i}^{j}&=\Pi_{k}^{j}\Pi_{i}^{k}  ,  \\\sqrt{h}\lap_{h}(\eta^{\alpha})&=\sqrt{h}h^{ij}\Pi_{\beta}^{\alpha}\pb_{i}\pb_{j}\eta^{\beta},\\ \pb_{A}(\sqrt{h}\lap_{h}\eta^{\alpha})&=\p_{i}(\sqrt{h}h^{ij}(\delta^{\alpha\lambda}-h^{lk}\pb_{l}\eta^{\lambda}\pb_{k}\eta^{\alpha})\pb_{j}\pb_{A}\eta^{\lambda}+\sqrt{h}(h^{ij}h^{kl}-h^{ik}h^{jl})\pb_{k}\eta^{\lambda}\pb_{A}\pb_{l}\eta^{\lambda}\pb_{j}\eta^{\alpha}),\label{derivativeofmeancur}\\\pb_{A}n_{j}&= -h^{kl}\pb_{k}\pb_{A}\eta^{i}n_{i}\pb_{l}\eta_{j}.
\end{align}
\pf
The proof of these formulas can be found in \cite{coutand2007well} and \cite{disconzi2019priori}.

\section{The smoothed $\kappa$-approximate system and functional framework}
Given $v\in L^{2}(0,T;H^{2}(\Omega))$, we introduce the approximating function $v^{\kappa}$ by
\begin{align}
\begin{aligned}
    \lap v^{\kappa}&= \lap(\rho_{\kappa}*E  v) \quad\textit{in $\Omega$},\\
    v^{\kappa}&= \Lambda_{\kappa}^{2}v \quad \textit{on $\p\Omega$},
\end{aligned}    
\end{align}
where $\rho_{\kappa}$ is standard mollifier on $\R^{3}$.
Then, we define
\begin{align}
    \eta^{\kappa}=e+\int_{0}^{T}v^{\kappa}
\end{align}
\lem
For some nonnegative integer $k$ and positive integer $l$ such that the following norms are well-defined, we have 
\begin{align}
\begin{aligned}     \|\p^{l}\p_{t}^{k}v^{\kappa}\|_{L^{2}(\Omega)}&\lesssim\|\p_{t}^{k}v\|_{H^{l}(\Omega)},\\
\|\p^{l}\p_{t}^{k}v^{\kappa}\|_{L^{2}(\Omega)}&\lesssim_{\kappa}\|\p_{t}^{k}v\|_{H^{1}(\Omega)},
\\\|\eta^{\kappa}\|_{H^{l}(\Omega)}&\lesssim \|\eta\|_{H^{l}(\Omega)}+\|e_{0}^{\kappa}\|_{H^{l}(\Omega)},\\\|\p_{t}^{k}(v^{\kappa}-v)\|_{H^{1}(\Omega)}&\lesssim \sqrt{\kappa}\|v\|_{H^{1.5}(\Omega)}.
\end{aligned}   
\end{align}
\pf
Note that if we only require $E$ to be a strong $m$-extension, it can be defined using the reflection method. In this case, $\p_{t}$ can be extended outside $\Omega$ and commutes with the extension operator. Thus, the first two inequalities is a direct consequence of elliptic estimate. Integrating on time, we have 
\begin{align}
    \begin{aligned}
        \lap\eta^{\kappa}&= \lap(\rho_{\kappa}*E\eta)-\lap e^{\kappa}\quad \textit{in $\Omega$},\\
        \eta^{\kappa}&=\Lambda_{\kappa}^{2}\eta \quad \textit{on $\p\Omega$}.
    \end{aligned}
\end{align}
Therefore, the elliptic estimate gives us the third estimate.
For the fourth estimate, we can compute that
\begin{align}
\begin{aligned}
     \|(\rho_{\kappa}*E  \p_{t}^{k}v)-\p_{t}^{k}v\|_{H^{1}(\Omega)}&\lesssim\|\rho_{\kappa}*E \p_{t}^{k} v-E\p_{t}^{k}v\|_{H^{1}}\lesssim\sqrt{\kappa}\|E\p_{t}^{k}v\|_{H^{1.5}}\lesssim\sqrt{\kappa}\|\p_{t}^{k}v\|_{H^{1.5}(\Omega)},\\
     \|\Lambda_{\kappa}^{2}\p_{t}^{k}v-\p_{t}^{k}v\|_{H^{\frac{1}{2}}(\p\Omega)}&\lesssim \sqrt{\kappa}\|\p_{t}^{k}v\|_{H^{1}(\p\Omega)}\lesssim\sqrt{\kappa}\|\p_{t}^{k}v\|_{H^{1.5}(\Omega)}.
\end{aligned}  
\end{align}
It immediately yields that
\begin{align}
    \|\p_{t}^{k}(v^{\kappa}-v)\|_{H^{1}(\Omega)}&\lesssim \sqrt{\kappa}\|v\|_{H^{1.5}(\Omega)}.
\end{align}

An important property is that $\eta^{\kappa}=\Lambda_{\kappa}^{2}\eta$ on $\p\Omega$, which preserves the doubly horizontally mollified structure on the boundary.
\defn
For a fixed small $\kappa_{0}$, for any $\kappa\in (0,\kappa_{0})$, we define the following notations:
\begin{align}
A^{\kappa}&=[D\eta^{\kappa}]^{-1},\\
\nabla^{\kappa}_{i}&=(\nabla_{A^{\kappa}})_{i}=(A^{\kappa})_{i}^{j}\p_{j}\\
\nabla^{\kappa}&=\nabla_{A^{\kappa}}=(\nabla^{\kappa}_{1},\nabla^{\kappa}_{2},\nabla^{\kappa}_{3})\\g_{\kappa}^{ij}&=(A^{\kappa})_{l}^{i}(A^{\kappa})_{l}^{j}\\
\lap_{g}^{\kappa}&=\p_{i}(g_{\kappa}^{ij}\p_{j})=\p_{i}((A^{\kappa})_{l}^{i}(A^{\kappa})_{l}^{j}\p_{j})\\\nabla_{N}^{\kappa}&=(\nabla_{A^{\kappa}})_{N}=\delta^{ik}(A^{\kappa})_{i}^{j}N_{j}(A^{\kappa})_{k}^{l}\p_{l}\\
J_{\kappa}&=\det \nabla \eta^{\kappa}.  
\end{align}

Now, 
we introduce the following $\kappa$-approximate system of \eqref{la2}:
\begin{align}\label{kappaapp}
    \begin{cases}
        \eta= e+\int_{0}^{t}v \quad\text{in $\Omega\times[0,T]$},\\
        \p_{t}v+\nabla^{\kappa} q=-\p_{t}\phi^{\kappa}\cdot \nabla^{\kappa}\phi^{\kappa}-\frac{1}{2}\nn|\nn\phi^{\kappa}|^{2}\quad\text{in $\Omega\times[0,T]$},\\ \nabla^{\kappa}\cdot v=0\quad\text{in $\Omega\times[0,T]$},\\\p_{t}\phi-\lap_{g}^{\kappa} \phi=\phi|\nabla^{\kappa}\phi|^2\quad\text{in $\Omega\times[0,T]$},\\v^{3}=0\quad \text{on $\Gamma_{0}\times[0,T]$},\\q =-\sigma(\frac{\sqrt{h}}{\sqrt{h_{\kappa}}} \lap_{h}(\eta) \cdot n_{\kappa})+\kappa\frac{1}{\sqrt{h_{\kappa}}}((1-\bar{\lap})(v\cdot n^{\kappa}))\quad\text{on $\Gamma\times[0,T]$},\\\nabla_{N}^{\kappa}\phi=0 \quad\text{on $\p\Omega\times[0,T]$},\\
        (\eta,v,\phi)=(e^{\kappa},u_{0},d_{0}) \quad\text{in $\Omega\times\{t=0\}$}.
    \end{cases}
\end{align}Later, we denote $(1-\bar{\lap})$ by $\bar{\lap}_{0}$.
$\phi^{\kappa}$ is the mollified $\phi$ by $\rho_{\kappa}*E \phi$ constrained on $\Omega$, where $E$ is the Sobolev extension operator.

Here, despite releasing the constrains on $\phi$ (i.e. we don't need to keep $|\phi^{\kappa}|=1$) in first velocity equation, we remain the form of $\p_{t}\phi^{\kappa}\cdot\nn\phi^{\kappa}$. i.e. 
\begin{align}\label{vkappat}
    \p_{t}v+\nabla^{\kappa} q=-\p_{t}\phi^{\kappa}\cdot \nabla^{\kappa}\phi^{\kappa}-\frac{1}{2}\nn|\nn\phi^{\kappa}|^{2}
\end{align}
\rmk
Later, we will utilize a cancellation structure to derive a uniform-in-$\kappa$ energy estimate. The term $\nn|\nn\phi^{\kappa}|^{2}$ will contribute a leading-order term. If we retain $\nn(\nn\phi^{\kappa} \odot \nn\phi^{\kappa})$, a similar proof for the cancellation can still be established. However, the curl estimate will become more lengthy and involve additional error analysis.

Recall that if we merely require $E$ to be a strong $m$ extension (m large enough), we can define $E\phi$ via the reflection method. Consequently, following a similar reason, $\p_{t}E\phi$ is well-defined provided that the regularity of $\phi$ is enough, and 
\begin{align}
\p_{t}E\phi=E\p_{t}\phi.
\end{align}
This ensures that $\p_{t}\phi^{\kappa}$ inherits the characteristics of standard mollification: for suitable positive integer $k,l$ such that the norms all well-defined
\begin{align}
    \begin{aligned}
        \|\p^{k}\p_{t}^{l}\phi^{\kappa}\|_{L^{2}(\Omega)}&\lesssim\|\p^{k}\p_{t}^{l}(\rho_{\kappa}*E\phi^{\kappa})\|_{L^{2}}=\|\rho_{\kappa}*\p^{k}\p_{t}^{l}E\phi\|_{L^{2}}\lesssim\|E\p_{t}^{l}\phi\|_{H^{k}}\lesssim\|\p_{t}^{l}\phi\|_{H^{k}(\Omega)}\\\|\p^{k}\p_{t}^{l}\phi^{\kappa}\|_{L^{2}(\Omega)}&\lesssim\|\p^{k}\p_{t}^{l}(\rho_{\kappa}*E\phi^{\kappa})\|_{L^{2}}=\|\p^{k}\rho_{\kappa}*E\p_{t}^{l}\phi\|_{L^{2}}\lesssim_{\kappa}\|E\p_{t}^{l}\phi\|_{L^{2}}\lesssim\|\p_{t}^{l}\phi\|_{L^{2}(\Omega)}\\\|\p_{t}^{l}\phi-\p_{t}^{l}\phi^{\kappa}\|_{L^{2}(\Omega)}&\lesssim\|\p_{t}^{l}E\phi-\p_{t}^{l}\rho_{\kappa}*E\phi\|_{L^{2}}\lesssim \kappa \|E\p_{t}^{l}\phi\|_{H^{1}}\lesssim \kappa \|\p_{t}^{l}\phi\|_{L^{2}(\Omega)}.        
    \end{aligned}
\end{align}
\prop\label{lwpappr}\label{lwpapp}
For each fixed $\kappa$, there exists a unique solution $(\eta,v,\phi)$ to \eqref{kappaapp} with $T_{\kappa}>0$ depending on $\|u_{0}\|_{H^{13.5}(\Omega)}$ and $\| d_{0}\|_{H^{9}(\Omega)}$, provided that $\nabla_{N}d_{0}=0$ on the boundary. Further, we have the following estimate
\begin{align}
    \sum_{k=0}^{5}\int_{0}^{T}\|\p_{t}^{k}\phi\|_{H^{10-2k}(\Omega)}^{2}+\sum_{k=0}^{4}\int_{0}^{T}\|\p_{t}^{k}v\|_{H^{13.5-3k}(\Omega)}^{2}+\sum_{k=0}^{3}\|\p_{t}^{k}p\|_{H^{14.5-3k}(\Omega)}^{2}+\sum_{k=0}^{4}\sup_t\|\p_{t}^{k}v\|_{12.5-3k}^{2}\lesssim C(\kappa,u_{0},d_{0}).
\end{align}
\rmk This $T_{\kappa}$ should be interpreted as a short time duration independent of other parameters that will emerge in Section 3. The energy estimate in Section 4 will demonstrate that a small $T$ can be attained, also independent of $\kappa$. This $T$ can be regarded as the consistent time interval throughout our paper, enabling us to employ a compactness argument at every step of the proof.
\pf
For $T>0$, we introduce a closed convex subset of $L^{2}(0,T;H^{13.5}(\Omega))$ as
\begin{align*}
    C_{T}=\{v\in L^{2}(0,T;H^{13.5}(\Omega))|\|v\|_{L^{2}_{t}H^{13.5}(\Omega)}\leq N_{0},  v_{0}=u_{0} \}
\end{align*}
We set the compatible conditions as follows:

Given $v_{0},d_{0}$, we require $\p_{t}d(0)=\lap d_{0}+d_{0}|\nabla d|^{2}$. $p(0)$ will be obtained via the elliptic equation:
\begin{align}
\begin{aligned}
    \lap p(0)&=-\nabla u_{0}:\nabla u_{0}-\nabla\cdot(\p_{t}d^{\kappa}(0)\nabla d_{0}^{\kappa})-\frac{1}{2}\lap|\nabla d_{0}^{\kappa}|^{2} \quad \textit{in $\Omega$},\\
    p(0)&=\kappa (1-\bar{\lap})(u_{0}^{3}) \quad \textit{on $\Gamma_{1}$},\\
    \nabla_{N}p(0)&=-\p_{t}d^{\kappa}(0)\nn_{N}d^{\kappa}-\frac{1}{2}\nn_{N}|\nn d^{\kappa}|^{2}\quad \textit{on $\Gamma_{0}$}.
\end{aligned}   
\end{align}
$\p_{t}v(0)=-\nabla p(0)-\p_{t}d_{0}^{\kappa}\nabla d_{0}^{\kappa}-\frac{1}{2}\nabla|\nabla d_{0}^{\kappa}|^{2}$. Inductively, after establishing the existence of $\phi$, $\p_{t}\phi$, $v$, $\p_{t}v$ and $p$, we define the initial values for their respective higher time derivatives. This inductive procedure is always valid, we don't need to worry about the existence of such a system.

For the sake of convenience, we always take $T\leq 1$. In the following section, $\bar{A}$, $\bar{A}^{\kappa}$ corresponds to $\bar{v}$. We further choose sufficiently small time interval such that $\sup_{t}(\|(\bar{A}^{\kappa})_{i}^{j}-\delta_{i}^{j}\|_{L^{\infty}(\Omega)}+\|(\bar{A}^{\kappa})_{i}^{j}-\delta_{i}^{j}\|_{L^{\infty}(\p\Omega)})\leq \ve_{1}$, where $\ve_{1}$ is a sufficient small number. $N_{0}$ will be determined later and will depend only on the initial data. If necessary, we may apply $T \leq 1$ to determine such an $N_{0}$. (i.e., $N_{0}$ is independent of the time interval) Similarly, $T$ will also be determined later and will depend solely on the initial data and $\kappa$ in Section 3. 
To make use of the estimate in \cite{coutand2007well}, we establish the following scheme:

1) Given $\bar{v}\in C_{T}$, we proved that there exists a unique solution $\phi$ for the equation:
\begin{align}\label{lwplinphi1}
\begin{cases}
     \p_{t}\phi-\lap_{\bar{A}^{\kappa}}\phi=\phi|\nabla_{\bar{A}^{\kappa}}\phi|^{2},\quad \textit{in $\Omega\times [0,T]$,}\\
     \bar{g}_{\kappa}^{3i}\p_{i}\phi=0,\quad\textit{on $\p\Om\times[0,T]$,}\\\phi=d_{0},\quad\textit{in $\Omega\times\{t=0\}$}
\end{cases}   
\end{align}

2)Provided with the same $\phi$, denote $|\phi|^{2}-1$ by $G$. We prove that there exists a unique solution of
\begin{align}\label{lwplinphi2}
\begin{cases}
     \p_{t}G-\lap_{\bar{A}^{\kappa}}G=G|\nabla_{\bar{A}^{\kappa}}\phi|^{2}, \quad\textit{in $\Omega\times[0,T]$},\\\bar{g}_{\kappa}^{3i}\p_{i}G=0,\quad\textit{on $\p\Omega\times[0,T]$},\\G=0,\quad\textit{in $\Omega\times\{t=0\}$}.
\end{cases}   
\end{align}
As a direct consequence, we can show $|\phi|=1$.

3)Plunge $\phi$ into the penalized problem
\begin{align}\label{penal1}
\begin{cases}
     \p_{t}v^{\ve}+\nabla_{\bar{A}^{\kappa}}q^{\ve}=-\p_{t}\phi^{\kappa}\cdot\nabla_{\bar{A}^{\kappa}}\phi^{\kappa}-\frac{1}{2}\nabla_{\bar{A}^{\kappa}}|\nabla_{\bar{A}^{\kappa}}\phi^{\kappa}|^{2},\quad\text{in $\Omega\times[0,T]$},\\\nabla_{\bar{A}^{\kappa}}\cdot v^{\ve}= -\ve q^{\ve},\quad \textit{in $\Omega\times[0,T]$},\\(v^{\ve})^{3}=0,\quad \text{on $\Gamma_{0}\times[0,T]$},\\q^{\ve} \bar{n}^{\kappa}_{i}=-\sigma(\frac{\sqrt{\bar{h}}}{\sqrt{\bar{h}_{\kappa}}} \lap_{\bar{h}}(\bar{\eta}) \cdot\bar{n}^{\kappa})\bar{n}^{\kappa}_{i}+\kappa\frac{1}{\sqrt{\bar{h}_{\kappa}}}(1-\bar{\lap})(v^{\ve}\cdot\bar{n}^{\kappa})\bar{n}^{\kappa}_{i},\quad\text{on $\Gamma\times[0,T]$},\\
    (\eta^{\ve},v^{\ve})=(Id,u_{0}),\quad\textit{in $\Omega\times\{t=0\}$}.
\end{cases} 
\end{align}
After using a standard Galerkin method to establish the existence of a weak solution for the penalized problem, we will perform a uniform $\ve$-estimate similar to the case in \cite{coutand2007well}.

4)After a compactness argument, we will show the existence of a weak solution to the linearized $\kappa$-problem and the corresponding regularity.
\begin{align}\label{penal2}
\begin{cases}
     \p_{t}v+\nabla_{\bar{A}^{\kappa}}q=-\p_{t}\phi^{\kappa}\cdot\nabla_{\bar{A}^{\kappa}}\phi^{\kappa}-\frac{1}{2}\nabla_{\bar{A}^{\kappa}}|\nabla_{\bar{A}^{\kappa}}\phi^{\kappa}|^{2},\quad\textit{in $\Omega\times[0,T]$},\\\nabla_{\bar{A}^{\kappa}}\cdot v= 0,\quad \textit{in $\Omega\times[0,T]$},\\v^{3}=0, \quad\text{on $\Gamma_{0}\times[0,T]$},\\q \bar{n}^{\kappa}_{i} =-\sigma(\frac{\sqrt{\bar{h}}}{\sqrt{\bar{h}_{\kappa}}} \lap_{\bar{h}}(\bar{\eta}) \cdot\bar{n})\bar{n}^{\kappa}_{i}+\kappa\frac{1}{\sqrt{\bar{h}_{\kappa}}}(1-\bar{\lap})(v\cdot\bar{n}^{\kappa})\bar{n}_{i}^{\kappa},\quad\text{on $\Gamma\times[0,T]$},\\
    (\eta,v)=(Id,u_{0}),\quad\textit{in $\Omega\times\{t=0\}$}.
\end{cases}
\end{align}

5) Finally, using a Tychonoff-type fixed point theorem, we will prove the existence of a unique solution for the $\kappa$-approximating problem \eqref{kappaapp}.
\subsection{Proof of \eqref{lwplinphi1}}

By a slight abuse of notation, in the following section, we will use $\phi$ to denote the solutions to different equations. The meaning of those variables is clear in the context.

\lem
There exists a unique solution $\phi$ to the following system:
\begin{align}\label{exlemforphi1}
\begin{cases}
     \p_{t}\phi-\lap_{\bar{A}^{\kappa}}\phi=F_{1},\quad \textit{in $\Omega\times [0,T]$,}\\
     \bar{g}_{\kappa}^{3i}\p_{i}\phi=F_{2},\quad\textit{on $\Gamma\times[0,T]$,}\\\phi=F_{0},\quad\textit{in $\Omega\times\{t=0\}$},
\end{cases}   
\end{align}
provided with $F_{0}\in H^{1}(\Omega)$, $F_{1}\in L^{2}(0,T;L^{2}(\Omega))$, and $F_{2}\in L^{2}(0,T;H^{\frac{1}{2}}(\p\Om))$. Further, we have the estimate:
\begin{align}
    \sup_{t}\|\phi\|_{H^{1}(\Omega)}+\|\phi\|_{L^{2}(0,T;H^{2}(\Omega))}+\|\p_{t}\phi\|_{L^{2}(0,T;L^{2}(\Omega))}\lesssim \|F_{1}\|_{L^{2}(0,T;L^{2}(\Omega))}+\|F_{2}\|_{L^{2}(0,T;H^{\frac{1}{2}}(\p\Om))}+\|F_{0}\|_{H^{1}(\Omega)}.
\end{align}
\pf
This lemma is widely known in standard parabolic theory, based on the Galerkin method. For detailed steps, Chapter III of \cite{ladyzhenskaia1968linear} can be consulted. (Actually, Theorem 5.1 with a slight modification).
\prop
There exists a small constant $T_{1}$ and a unique solution $\phi$ of equation \eqref{lwplinphi1} within the time interval $t\in(0, T_{1})$. Furthermore, we have the estimate
\begin{align}
     \sup_{t}(\|\phi\|_{H^{3}(\Omega)}+\|\p_{t}\phi\|_{L^{2}(\Omega)})+\|\phi\|_{L_{t}^{2}H^{4}(\Omega)}+\|\p_{t}\phi\|_{L_{t}^{2}H^{1}(\Omega)}+\|\p_{t}^{2}\phi\|_{L_{t}^{2}L^{2}(\Omega)}\lesssim M_{0},
\end{align}
where $M_{0}$ is a constant that only depends on $\|d_{0}\|_{H^{3}}$ and $N_{0}$.
\pf
Given the sufficiently high regularity of $\bar{v}$ in this proposition, we will omit the constant associated with $\bar{v}$. For example, we will write directly
\begin{align}    \|\nabla_{\bar{A}^{\kappa}}\phi\|_{H^{k}(\Omega)}\lesssim\|\bar{A}^{\kappa}\|_{H^{3}(\Omega)}\|\phi\|_{H^{k+1}(\Omega)}\lesssim\|\phi\|_{H^{k+1}(\Omega)},\quad\textit{k=0,1,2,3}.
\end{align}
Set
\begin{align}
    S_{1}=\{\phi|\|\phi\|_{L_{t}^{2}H^{4}_{x}(\Omega)}+\|\p_{t}\phi\|_{L_{t}^{2}H_{x}^{2}(\Omega)}+\sup_{t}(\|\phi\|_{H^{3}(\Omega)}+\|\p_{t}\phi\|_{H^{1}(\Omega)})\leq M_{0}, \phi|_{t=0}=\phi_{0},\p_{t}\phi|_{t=0}=\lap d_{0}+d_{0}|\nabla d_{0}|^{2}\}.
\end{align}
$M_{0}$ is to be determined later. Clearly, $S_{1}$ is a Banach space.
Now, given $\tilde{\phi}\in S_{1}$, we solve the linear system:
\begin{align}\label{exphifix1}
\begin{aligned}
     \p_{t}\phi-\lap_{\bar{A}^{\kappa}}\phi=\tilde{\phi}|\nabla_{\bar{A}^{\kappa}}\tilde{\phi}|^{2},\quad \textit{in $\Omega\times [0,T]$,}\\
     \bar{g}_{\kappa}^{3i}\p_{i}\phi=0,\quad\textit{on $\p\Omega\times[0,T]$,}\\\phi=d_{0},\quad\textit{in $\Omega\times\{t=0\}$}
\end{aligned}   
\end{align}
Utilizing the previous lemma, we can show the existence of $\phi$ as a direct consequence of the following estimate:
\begin{align}\label{lt2lx2phiexis}
 \begin{aligned}     \int_{0}^{t}\int_{\Omega}|\tilde{\phi}|^{2}|\nabla_{\bar{A}^{\kappa}}\tilde{\phi}|^{4}&\lesssim\int_{0}^{t}\|\tilde{\phi}\|_{L^{\infty}(\Omega)}^{2}\|\nabla_{\bar{A}^{\kappa}}\tilde{\phi}\|_{L^{4}(\Omega)}^{4}\\&\lesssim\int_{0}^{t}\|\tilde{\phi}\|_{H^{2}(\Omega)}^{2}\|\nabla_{\bar{A}^{\kappa}}\tilde{\phi}\|_{H^{1}(\Omega)}^{4}\\&\lesssim T \PP(M_{0}).
 \end{aligned}   
\end{align}
Thus,
\begin{align}
    \sup_{t}\|\phi\|_{H^{1}(\Omega)}+\|\phi\|_{L_{t}^{2}H^{2}(\Omega)}+\|\p_{t}\phi\|_{L^{2}_{t}L^{2}(\Omega)}\lesssim T^{\frac{1}{2}}\PP(M_{0})+\|d_{0}\|_{H^{1}(\Omega)}
\end{align}
Taking $\p_{t}$ on \eqref{exphifix1}, denoting $\p_{t} \phi$ by $\psi_{1}$, we can verify that $\psi_{1}$ is a weak solution to the following system: (Later, we will see the uniqueness of the solution. Thus, we use the same symbol to represent the variable and $\p_{t}\phi$)
\begin{align}\label{equationforphitildetpsi1}
    \begin{aligned}
         \p_{t}\psi_{1}-\lap_{\bar{A}^{\kappa}}\psi_{1}&=\p_{t}(\tilde{\phi}|\nabla_{\bar{A}^{\kappa}}\tilde{\phi}|^{2})+[\p_{t},\lap_{\bar{A}^{\kappa}}]\phi,\quad \textit{in $\Omega\times [0,T]$,}\\
     \bar{g}_{\kappa}^{3i}\p_{i}\psi_{1}&=-[\p_{t},\bar{g}_{\kappa}^{3i}]\p_{i}\phi,\quad\textit{on $\p\Omega\times[0,T]$,}\\\psi_{1}&= \lap d_{0}+d_{0}|\nabla d_{0}|^{2},\quad\textit{in $\Omega\times\{t=0\}$}
    \end{aligned}
\end{align}
And, by the Sobolev inequality, it indicates that:
\begin{align}\label{lt2lx2pbphi}
\begin{aligned}
     \|\p_{t}(\tphi|\nabla_{\bar{A}^{\kappa}}\tphi|^{2})\|_{L^{2}_{t}L_{x}^{2}}&\lesssim\|\p_{t}\tphi|\nabla_{\bar{A}^{\kappa}}\tphi|^{2}\|_{L^{2}_{t}L_{x}^{2}}+\|\tphi\nabla_{\bar{A}^{\kappa}}\tphi\p_{t}(\nabla_{\bar{A}^{\kappa}}\tphi)\|_{L^{2}_{t}L_{x}^{2}}\\&\lesssim (\int_{0}^{t}\|\p_{t}\tphi\|_{L^{6}}^{2}\|\nabla_{\bar{A}^{\kappa}}\tphi\|_{L^{6}}^{4})^{\frac{1}{2}}+(\int_{0}^{t}\|\tphi\|_{L^{\infty}}^{2}\|\nabla_{\bar{A}^{\kappa}}\tphi\|_{L^{\infty}}^{2}\|\p_{t}\nabla_{\bar{A}^{\kappa}}\tphi\|_{L^{2}}^{2})^{\frac{1}{2}}\\&\lesssim T^{\frac{1}{2}}\PP(M_{0}).
\end{aligned}   
\end{align}
Moreover,
\begin{align}\label{lt2lx2pbphi2}
\begin{aligned}
    \|[\p_{t},\lap_{\bar{\ak}}]\phi\|_{L_{t}^{2}L^{2}(\Omega)}=\|\p_{i}\p_{t}\bar{g}_{\kappa}^{ij}\p_{j}\phi+\p_{t}\bar{g}_{\kappa}^{ij}\p_{i}\p_{j}\phi\|_{L_{t}^{2}L^{2}(\Omega)}\lesssim \|\phi\|_{L_{t}^{2}H^{2}(\Omega)}\lesssim T^{\frac{1}{2}}\PP(M_{0})+\|d_{0}\|_{H^{3}(\Omega)}.
\end{aligned}   
\end{align}
For the boundary term $F_{2}=-\p_{t}\bar{g}_{\kappa}^{3i}\p_{i}\phi$, the trace inequality implies that
\begin{align}\label{lt2lx2pbphibdy}  \|\p_{t}\bar{g}_{\kappa}^{3i}\p_{i}\phi\|_{L_{t}^{2}H^{\frac{1}{2}}(\p\Omega)}\lesssim\|\p_{t}\bar{g}_{\kappa}^{3i}\|_{L_{t}^{\infty}H^{2}_{x}(\Omega)}\|\p\phi\|_{L_{t}^{2}H^{1}_{x}(\Omega)}\lesssim \|\phi\|_{L_{t}^{2}H^{2}(\Omega)}\lesssim T^{\frac{1}{2}}\PP(M_{0})+\|d_{0}\|_{H^{3}(\Omega)}.
\end{align}
Hence, combining Lemma 3.2,  \eqref{lt2lx2pbphi}, \eqref{lt2lx2pbphi2}, and \eqref{lt2lx2pbphibdy}, we have the existence and uniqueness of the solution and the following estimate:
\begin{align}\label{estimateforpsi1aspbphi}
    \sup_{t}\|\psi_{1}\|_{H^{1}(\Omega)}+\|\psi_{1}\|_{L_{t}^{2}H^{2}(\Omega)}+\|\p_{t}\psi_{1}\|_{L_{t}^{2}H^{1}(\Omega)}\leq C_{T}(T^{\frac{1}{2}}\PP(M_{0})+\|d_{0}\|_{H^{3}(\Omega)}).
\end{align}
Owing to our setting $T\leq 1$, $C_{T}$ can be bounded by a constant $C$ independent of time. Subsequently, we will consistently consider it as a constant without elaboration in analogous scenarios. We can directly check that $\p_{t}\phi$ is also the weak solution to \eqref{equationforphitildetpsi1}. Now, we can calculate that
\begin{align}    \|\tphi|\nabla_{\bar{A}^{\kappa}}\tphi|\|_{L^{2}_{t}H^{2}(\Omega)}\lesssim T^{\frac{1}{2}}\PP(\|\tphi\|_{L_{t}^{\infty}H^{3}(\Omega)})\lesssim T^{\frac{1}{2}}\PP(M_{0}),
\end{align}
the standard elliptic estimate gives us
\begin{align}\label{estimateforphilt4h4inifir}
    \|\phi\|_{L_{t}^{2}H^{4}(\Omega)}+\|\phi\|_{L_{t}^{\infty}H^{3}(\Omega)}\lesssim T^{\frac{1}{2}}\PP(M_{0})+\|d_{0}\|_{H^{3}(\Omega)}
\end{align}
Hence, combining Lemma 3.2, \eqref{lt2lx2phiexis}, \eqref{estimateforpsi1aspbphi}, and  \eqref{estimateforphilt4h4inifir} we have
\begin{align}
    \sup_{t}(\|\phi\|_{H^{3}(\Omega)}+\|\p_{t}\phi\|_{H^{1}(\Omega)})+\|\phi\|_{L_{t}^{2}H^{4}(\Omega)}+\|\p_{t}\phi\|_{L_{t}^{2}H^{2}(\Omega)}+\|\p_{t}^{2}\phi\|_{L_{t}^{2}L^{2}(\Omega)}\lesssim T^{\frac{1}{2}}\PP(M_{0})+\|d_{0}\|_{H^{3}(\Omega)}.
\end{align}

Now, selecting sufficiently large $M_{0}$ compared with $\|d_{0}\|_{H^{3}(\Omega)}$, and choosing sufficiently small $T_{1}$, the solution of \eqref{exphifix1} defines a map $M:S_{1}\to S_{1}$. To establish the existence theorem, it suffices to show that $M$ is a contraction. We still use $\phi^{(i)},\tilde{\phi}^{(i)}\in S_{1}$ , $i=1,2$ to express the relation $\phi^{(i)}=M(\tilde{\phi}^{(i)})$.
Let $\Psi=\phi^{(1)}-\phi^{(2)}$ and $\tilde{\Psi}=\tphi^{(1)}-\tphi^{(2)}$; then we obtain:
\begin{align}
     \begin{aligned}
         \p_{t}\Psi-\lap_{\bar{A}^{\kappa}}\Psi&=\tilde{\Psi}|\nabla_{\bar{A}^{\kappa}}\tilde{\phi}^{(1)}|^{2}+\tphi^{(2)}\nabla_{\bar{A}^{\kappa}}\tPsi\nabla_{\bar{A}^{\kappa}}(\tphi^{(1)}+\tphi^{(2)}),\quad \textit{in $\Omega\times [0,T]$,}\\
     \bar{g}_{\kappa}^{3i}\p_{i}\Psi&=0,\quad\textit{on $\p\Omega\times[0,T]$,}\\\Psi&=0,\quad\textit{in $\Omega\times\{t=0\}$}
    \end{aligned}
\end{align}
We will adopt a similar strategy to that of considering $\p_{t}\Psi$. From a direct calculation, it yields that
\begin{align}\label{estimateforPsiex1}
    \begin{aligned}        \|\tilde{\Psi}|\nabla_{\bar{A}^{\kappa}}\tilde{\phi}^{(1)}|^{2}\|_{L_{t}^{2}H_{x}^{2}}+\|\p_{t}(\tilde{\Psi}|\nabla_{\bar{A}^{\kappa}}\tilde{\phi}^{(1)}|^{2})\|_{L_{t}^{2}L_{x}^{2}}&\lesssim(\int_{0}^{t}\|\tPsi\|_{H^{2}}^{2}\|\nabla_{\bar{A}^{\kappa}}\tphi^{(1)}\|_{H^{2}}^{4})^{\frac{1}{2}}+(\int_{0}^{t}\|\p_{t}\tPsi\|_{L^{2}}^{2}\|\nabla_{\bar{A}^{\kappa}}\tphi^{(1)}\|_{H^{2}}^{4})^{\frac{1}{2}}\\&+(\int_{0}^{t}\|\tPsi\|_{H^{2}}^{2}\|\p_{t}\nabla_{\bar{A}^{\kappa}}\tphi^{(1)}\|_{L^{2}}^{2}\|\nabla_{\bar{A}^{\kappa}}\tphi^{(1)}\|_{H^{2}})^{\frac{1}{2}}\\&\lesssim T^{\frac{1}{2}}C(M_{0})(\|\p_{t}\tPsi\|_{L^{2}_{t}L^{2}}+\|\tPsi\|_{L_{t}^{\infty}H^{2}}),
    \end{aligned}
\end{align}
and a similar computation shows that
\begin{align}\label{estimateforPsiex2}
    \begin{aligned}     &\|\tphi^{(2)}\nabla_{\bar{A}^{\kappa}}\tPsi\nabla_{\bar{A}^{\kappa}}(\tphi^{(1)}+\tphi^{(2)})\|_{L_{t}^{2}H^{2}(\Omega)}+\|\p_{t}(\tphi^{(2)}\nabla_{\bar{A}^{\kappa}}\tPsi\nabla_{\bar{A}^{\kappa}}(\tphi^{(1)}+\tphi^{(2)}))\|_{L_{t}^{2}L^{2}(\Omega)}\\\lesssim &  (\int_{0}^{t}\|\tphi^{(2)}\|_{H^{2}}^{2}\|\nabla_{\bar{A}^{\kappa}}\tPsi\|_{H^{2}}^{2}\|\nabla_{\bar{A}^{\kappa}}(\tphi^{(1)}+\tphi^{(2)})\|_{H^{2}}^{2})^{\frac{1}{2}}+(\int_{0}^{t}\|\p_{t}\tphi^{(2)}\|_{L^{2}}^{2}\|\nabla_{\bar{A}^{\kappa}}\tPsi\|_{H^{2}}^{2}\|\nabla_{\bar{A}^{\kappa}}(\tphi^{(1)}+\tphi^{(2)})\|_{H^{2}}^{2})^{\frac{1}{2}}\\+&(\int_{0}^{t}\|\tphi^{(2)}\|_{H^{2}}^{2}\|\p_{t}\nabla_{\bar{A}^{\kappa}}\tPsi\|_{L^{2}}^{2}\|\nabla_{\bar{A}^{\kappa}}(\tphi^{(1)}+\tphi^{(2)})\|_{H^{2}}^{2})^{\frac{1}{2}}+(\int_{0}^{t}\|\tphi^{(2)}\|_{H^{2}}^{2}\|\nabla_{\bar{A}^{\kappa}}\tPsi\|_{H^{2}}^{2}\|\p_{t}\nabla_{\bar{A}^{\kappa}}(\tphi^{(1)}+\tphi^{(2)})\|_{L^{2}}^{2})^{\frac{1}{2}}\\\lesssim & T^{\frac{1}{2}}C(M_{0})(\|\tPsi\|_{L_{t}^{\infty}H^{3}}+\|\p_{t}\tPsi\|_{L_{t}^{\infty}H^{1}}).
    \end{aligned}
\end{align}
The estimate of commutator reads:
\begin{align} \label{estimateforPsiex3}   \|[\p_{t},\lap_{\bar{A}^{\kappa}}]\Psi\|_{L_{t}^{2}L^{2}}\lesssim \|\p_{i}\p_{t}\bar{g}_{\kappa}^{ij}\p_{j}\Psi+\p_{t}\bar{g}_{\kappa}^{ij}\p_{i}\p_{j}\Psi\|_{L_{t}^{2}L^{2}}\lesssim T^{\frac{1}{2}}\sup_{t}\|\Psi\|_{H^{2}}.
\end{align}
The boundary term can be controlled as
\begin{align}\label{estimateforPsiex4}
     \|\p_{t} \bar{g}_{\kappa}^{3i}\p_{i}\Psi\|_{L_{t}^{2}H^{\frac{1}{2}}(\p\Omega)}\lesssim T^{\frac{1}{2}}\sup_{t}\|\Psi\|_{H^{2}(\Omega)}.
\end{align}
Hence, combining \eqref{estimateforPsiex1}, \eqref{estimateforPsiex2}, \eqref{estimateforPsiex3}, and \eqref{estimateforPsiex4} 
, we have
\begin{align}
\begin{aligned}
     & \sup_{t}(\|\Psi\|_{H^{3}(\Omega)}+\|\p_{t}\Psi\|_{H^{1}(\Omega)})+\|\Psi\|_{L_{t}^{2}H^{4}(\Omega)}+\|\p_{t}\Psi\|_{L_{t}^{2}H^{2}(\Omega)}\\\lesssim & T^{\frac{1}{2}}C(M_{0})(  \sup_{t}\|\p_{t}\tPsi\|_{H^{1}(\Omega)}+\|\tPsi\|_{L_{t}^{\infty}H^{3}(\Omega)})+T^{\frac{1}{2}}\sup_{t}\|\Psi\|_{H^{2}}.
\end{aligned}  
\end{align}
Taking $T$ sufficiently small, we can ensure that $M$ is a contraction. A standard argument then yields the existence of a solution and the desired estimates (Choosing a constant slightly larger than $M_{0}$, which we still denote it by $M_{0}$). This completes the proof of the proposition.

\subsection{Proof of \eqref{lwplinphi2}}
Now, we verify that $|\phi|=1$. Test $G$ in \eqref{lwplinphi2},
\begin{align}
\begin{aligned}
      \int_{\Omega}G^{2}+\int_{0}^{t}\int_{\Omega}\bar{g}^{ij}\p_{i}G \p_{j}G&=\int_{0}^{t}\int_{\Omega}G^{2}|\nabla_{\bar{A}^{\kappa}}\phi|^{2}.
\end{aligned}  
\end{align}
Invoking the H\"{o}lder inequality, Young inequality and the Sobolev embedding theorem, it implies that
\begin{align}
    \begin{aligned}
        \sup_{t}\|G(t)\|_{L^{2}(\Omega)}^{2}+\int_{0}^{T}\| G\|_{H^{1}(\Omega)}^{2} &\lesssim\int_{0}^{T}\|G\|_{L^{3}(\Omega)}^{2}\|\nabla_{\bar{A}^{\kappa}}\phi\|_{L^{6}(\Omega)}^{2}\\&\lesssim\int_{0}^{T}\|G\|_{L^{2}(\Omega)}\|G\|_{H^{1}(\Omega)}\|\nabla_{\bar{A}^{\kappa}}\phi\|_{H^{1}(\Omega)}^{2}\\&\lesssim\int_{0}^{T}\|G\|_{L^{2}(\Omega)}\|G\|_{H^{1}(\Omega)}\PP(N_{0},M_{0})\\&\lesssim T^{\frac{1}{2}}\sup_{t\in [0,T]}\|G\|_{L^{2}(\Omega)}\|G\|_{L^{2}(0,T;H^{1}(\Omega))}\PP(N_{0},M_{0}).
    \end{aligned}
\end{align}
Choosing $T$ small enough, we can see that
\begin{align}
    G=0 \quad \textit{in $\Omega\times[0,T]$}.
\end{align}
Consequently, we have $|\phi|=1$.
\rmk
Alternatively, readers can refer to Theorem 3.2 in Chapter 3 of \cite{ladyzhenskaia1968linear}, where it is not difficult to determine the regularity of $\phi$ that meets the theorem's conditions.
\subsection{Proof of \eqref{penal1}}
Let $H$ denote $-\p_{t}\phi^{\kappa} \cdot \nabla_{\bar{A}^{\kappa}}\phi^{\kappa} - \frac{1}{2}\nabla_{\bar{A}^{\kappa}}|\nabla_{\bar{A}^{\kappa}}\phi^{\kappa}|^{2}$. The term related to $\phi^{\kappa}$ can be estimated using the properties of convolution, and $H$ can be regarded as smooth. The calculations in \cite{coutand2007well} are sufficient to justify the existence of the system. Here, we provide a sketchy proof, offering detailed explanations only where necessary modifications are required. Taking $(e_{l})_{l=1}^{\infty}$ as a basis of $H^{1.5}(\Omega)$ and an orthonormal basis of $L^{2}(\Omega)$, we define the approximating function 
\begin{align}
    w_{l}(t,x)=\sum_{k=1}^{l}y_{k}(t)e_{k}, \quad l\geq 2,
\end{align}
which satisfies the following ordinary differential equations:
\begin{equation}\label{testwlt}
\begin{cases}
     \int_{\Omega}\bar{J}_{\kappa} w_{lt}\psi -\int_{\Omega} \bar{J}_{\kappa} (\bar{A}^{\kappa})^{i}_{j}q_{l}\p_{i}\psi^{j}+\kappa\int_{\Gamma_{1}}w_{l}\cdot \bar{n}_{\kappa}\psi\cdot \bar{n}_{\kappa}+\sum_{i=1}^{2}\kappa\int_{\Gamma_{1}}\pb_{i} (w_{l}\cdot \bar{n}_{\kappa})\pb_{i}(\psi\cdot \bar{n}_{\kappa})\\=\int_{\Gamma_{1}}\sigma L_{\bar{h}}\bar{\eta}\cdot \bar{n}_{\kappa}\psi\cdot \bar{n}_{\kappa}+\int_{\Omega} \psi H \quad \textit{$\forall \psi\in span(e_{1},...,e_{l})$, in $[0,T]$},\\
     w_{l}(0)=(u_{0})_{l},\quad \textit{in $\Omega$,}
\end{cases}   
\end{equation}
where  $L_{\bar{h}}=\frac{\sqrt{\bar{h}}}{\sqrt{\bar{h}_{\kappa}}}$, $q_{l}=-\frac{1}{\ve}(\bar{A}_{\kappa})_{i}^{j}\p_{j}w_{l}^{i}$, and $(u_{0})_{l}$ denotes the $L^{2}(\Omega)$ projection of $u_{0}$ on $span(e_{1},...,e_{l})$. The Cauchy-Lipschitz theorem indicates the local well-posedness of $w_{l}$ on time interval $[0,T_{max}]$. Choosing $\psi=w_{l}$, we obtain the energy law for $\forall t \in (0,T_{max})$
\begin{align}
\begin{aligned}
     &\frac{1}{2}\|\bar{J}_{\kappa}^{\frac{1}{2}}w_{l}(t)\|_{L^{2}(\Omega)}^{2}+\kappa\int_{0}^{t}\|w_{l}\cdot \bar{n}_{\kappa}\|_{H^{1}(\Gamma_{1})}^{2}+\ve \int_{0}^{t}\|\bar{J}_{\kappa}^{\frac{1}{2}}q_{l}\|_{L^{2}(\Omega)}^{2}-\frac{1}{2}\int_{0}^{t}\int_{\Omega}\p_{t}\bar{J}_{\kappa}|w_{l}|^{2}\\=&\frac{1}{2}\|(u_{0})_{l}\|_{L^{2}(\Omega)}^{2}+\sigma \int_{0}^{t} \int_{\Gamma_{1}}L_{\bar{h}}\bar{\eta}\cdot \bar{n}_{\kappa}w_{l}\cdot \bar{n}_{\kappa}+\int_{0}^{t}\int_{\Omega}w_{l}H
\end{aligned}   
\end{align}

Followed by our assumption on $\bar{J}_{\kappa}$, we can choose a shorter time interval $t\in [0,T]$ if necessary such that
\begin{align}
    \begin{aligned}
        \sup_{t^{\prime}\in[0,t]}\|w_{l}\|_{L^{2}(\Omega)}^{2}+\kappa\int_{0}^{t}\|w_{l}\cdot \bar{n}_{\kappa}\|_{H^{1}(\Gamma_{1})}^{2}+\ve \int_{0}^{t}\|q_{l}\|_{L^{2}(\Omega)}^{2}\leq C(N_{0},M_{0})\quad \textit{$\forall t\in[0,T]$}.
    \end{aligned}
\end{align}
For every $\psi\in L^{2}(0,T;H^{\frac{3}{2}}(\Omega))$, we substitute it into equation \eqref{testwlt}, it is straightforward to verify that the integration is well-defined with an estimate
\begin{align}
    \|w_{lt}\|_{L^{2}_{t}H^{-\frac{3}{2}}(\Omega)}\leq C(N_{0},M_{0}).
\end{align}
Now, we can deduce that there exists a subsequence (still labeled as $l$)
\begin{align}
    \begin{aligned}
        w_{l}&\rightharpoonup^{*} w_{\ve}\quad\textit{in $L^{\infty}(0,T;L^{2}(\Omega))$},\\
        w_{l}&\rightharpoonup w_{\ve} \quad \text{in $L^{2}(0,T;L^{2}(\Omega))$},\\
         q_{l}&\rightharpoonup q_{\ve} \quad \text{in $L^{2}(0,T;L^{2}(\Omega))$},\\
          w_{l}\cdot \bar{n}_{\kappa}&\rightharpoonup s_{\ve} \quad \text{in $L^{2}(0,T;H^{1}(\Gamma_{1}))$},\\
          w_{lt}&\rightharpoonup r_{\ve}\quad \text{in $L^{2}(0,T;H^{-\frac{3}{2}}(\Omega))$}.
    \end{aligned}
\end{align}
Given $w_{\ve}\in L^{2}(0,T;L^{2}(\Omega))$, we regard $\dive_{\bar{A}^{\kappa}}w_{\ve}$ as an element of $L^{2}(0,T;(H_{0}^{1}(\Omega))^{*})$. By the definition of $q_{l}$, for any $\rho\in L_{t}^{2}H_{0}^{1}(\Omega)$, we also have
\begin{align}
    \begin{aligned}        \int_{0}^{T}\int_{\Omega}\dive_{\bar{A}^{\kappa}}w_{\ve} \rho&=-\int_{0}^{T}\int_{\Omega}\bar{J_{\kappa}}w_{\ve} \cdot\nabla_{\bar{A}^{\kappa}}(\bar{J}_{\kappa}^{-1}\rho)\\&=-\lim_{l\to\infty}\int_{0}^{T}\int_{\Omega}\bar{J_{\kappa}}w_{l} \cdot\nabla_{\bar{A}^{\kappa}}(\bar{J}_{\kappa}^{-1}\rho)\\&=\lim_{l\to\infty}\int_{0}^{T}\int_{\Omega}\dive_{\bar{A}^{\kappa}}w_{l} \rho\\&=-\lim_{l\to\infty}\int_{0}^{T}\int_{\Omega}\ve q_{\ve} \rho
    \end{aligned}
\end{align}
Thus, $\dive_{\bar{A}^{\kappa}}w_{\ve}=-\ve q_{\ve}$ in $L^{2}(0,T;(H_{0}^{1}(\Omega))^{*})$. Since $q_{\ve}\in L^{2}(0,T;L^{2}(\Omega))$, we have $\dive_{\bar{A}^{\kappa}}w_{\ve}\in L^{2}(0,T;L^{2}(\Omega))$ and
\begin{align}
    \dive_{\bar{A}^{\kappa}}(w_{l})\rightharpoonup  \dive_{\bar{A}^{\kappa}}(w_{\ve})\quad \text{in $L^{2}(0,T;L^{2}(\Omega))$}.
\end{align}
Similarly, for any $\rho^{\prime}\in L^{2}(0,T;L^{2}(\Gamma_{1}))$, we can find $\rho^{\prime}\in L^{2}(0,T;H^{1}(\Omega))$ (with a little abuse of notation, we choose the same symbol) such that $\rho^{\prime}$ is equal to its boundary trace. Now, integrating by part, we have
\begin{align}
    \begin{aligned}       \int_{0}^{T}\int_{\Omega}w_{\ve}\cdot \bar{n}_{\kappa}\rho^{\prime}&=\int_{0}^{T}\int_{\Omega}\dive_{\bar{A}^{\kappa}}w_{\ve} \frac{\rho^{\prime}}{|(\bar{A}^{\kappa})_{i}^{3}|}+\bar{J_{\kappa}}w_{\ve} \cdot\nabla_{\bar{A}^{\kappa}}(\frac{\rho^{\prime}}{|\bar{J}_{\kappa}(\bar{A}^{\kappa})_{i}^{3}|})\\&=\lim_{l\to\infty}\int_{0}^{T}\int_{\Omega}\dive_{\bar{A}^{\kappa}}w_{l} \frac{\rho^{\prime}}{|(\bar{A}^{\kappa})_{i}^{3}|}+\bar{J_{\kappa}}w_{l} \cdot\nabla_{\bar{A}^{\kappa}}(\frac{\rho^{\prime}}{|\bar{J}_{\kappa}(\bar{A}^{\kappa})_{i}^{3}|})\\&=\lim_{l\to\infty}\int_{0}^{T}\int_{\Omega}w_{l}\cdot \bar{n}_{\kappa}\rho^{\prime}\\&=\int_{0}^{T}\int_{\Omega}s_{\ve}\rho^{\prime}
    \end{aligned}
\end{align}
Consequently, as $s_{\ve}\in L^{2}(0,T;H^{1}(\Gamma_{1}))$, we have
\begin{align}
    w_{l}\cdot \bar{n}_{\kappa}\rightharpoonup w_{\ve}\cdot\bar{n}_{\kappa}  \quad \text{in $L^{2}(0,T;H^{1}(\Gamma_{1}))$ }.
\end{align}
We show $r_{\ve}=w_{\ve t}$ in a similar fashion:
\begin{align}
    \begin{aligned}        \int_{0}^{T}\int_{\Omega}r_{\ve}\rho_{1}(t)\rho_{2}(x)&=\lim_{l\to\infty}\int_{0}^{T}\int_{\Omega}w_{l t}\rho_{1}(t)\rho_{2}(x)\quad\textit{$\forall \rho_{1}\in C^{1}_{c}(0,T)$, $\forall \rho_{2}\in H^{\frac{3}{2}}(\Omega)$}\\&=-\lim_{l\to\infty}\int_{0}^{T}\int_{\Omega}w_{l}\p_{t}\rho_{1}(t)\rho_{2}(x)\\&=-\int_{0}^{T}\int_{\Omega}w_{\ve}\p_{t}\rho_{1}(t)\rho_{2}(x).
    \end{aligned}
\end{align}
Hence,
\begin{align}
     w_{lt}&\rightharpoonup r_{\ve}\quad \text{in $L^{2}(0,T;H^{-\frac{3}{2}}(\Omega))$}.
\end{align}
We summarize that
\begin{align}\label{estimateforwveinitiall2l2}
    \|w_{\ve}\|_{L^{2}(0,T;L^{2}(\Omega))}+\|\dive_{\bar{A}^{\kappa}}w_{\ve}\|_{L^{2}(0,T;L^{2}(\Omega))}+\|w_{\ve}\cdot \bar{n}_{\kappa}\|_{L^{2}(0,T;H^{1}(\Gamma_{1}))}\leq C(N_{0},M_{0}).
\end{align}
For any $\rho^{\prime\prime}\in L^{2}(0,T;H^{\frac{3}{2}}(\Omega))$, the weakly limit satisfies 
\begin{align}
\begin{aligned}
      &\int_{0}^{T}\int_{\Omega} \bar{J}_{\kappa}w_{\ve t}\rho^{\prime\prime} -  \int_{0}^{T}\int_{\Omega} \bar{J}_{\kappa} (\bar{A}^{\kappa})^{i}_{j}q_{\ve}\p_{i}(\rho^{\prime\prime})^{j}+\kappa  \int_{0}^{T}\int_{\Gamma_{1}}w_{\ve}\cdot \bar{n}_{\kappa}\rho^{\prime\prime}\cdot \bar{n}_{\kappa}+\sum_{i=1}^{2}\kappa  \int_{0}^{T}\int_{\Gamma_{1}}\pb_{i} (w_{\ve}\cdot \bar{n}_{\kappa})\pb_{i}(\rho^{\prime\prime}\cdot \bar{n}_{\kappa})\\= & \int_{0}^{T}\int_{\Gamma_{1}}\sigma L_{\bar{h}}\bar{\eta}\cdot \bar{n}_{\kappa}\rho^{\prime\prime}\cdot \bar{n}_{\kappa}+  \int_{0}^{T}\int_{\Omega} \rho^{\prime\prime} H.
\end{aligned} 
\end{align}
Thus, we can assert that 
\begin{align}\label{equationforwve}
    w_{\ve t}+\nabla_{\bar{A}^{\kappa}}q_{\ve}=H
\end{align}
in the sense of $L_{t}^{2}H^{-\frac{3}{2}}(\Omega)$. In other words, $(w_{\ve},q_{\ve})$ is the weak solution to the penalized problem. Since $\nabla_{\bar{A}^{\kappa}}q_{\ve} \in L^{2}(0,T;(H_{0}^{1}(\Omega))^{*})$, we can identify the equality \eqref{equationforwve} in the same space. Since $\curl_{\bar{A}^{\kappa}}w_{\ve}\in L^{2}(0,T;(H_{0}^{1}(\Omega))^{*})$, we have the relation
\begin{align} \label{equationforlowestcurlwve}   \curl_{\bar{A}^{\kappa}}w_{\ve}(t)=\curl u_{0}+\int_{0}^{t}\ve^{.ij}\p_{l}w_{\ve}^{j}\p_{t}(\bar{A}^{\kappa})_{i}^{l}+\int_{0}^{t}\curl_{\bar{A}^{\kappa}} H \quad\textit{in $L^{2}(0,T;(H_{0}^{1}(\Omega))^{*})$}.
\end{align}
The process of enhancing the regularity of $w_{\ve}$ will be segmented into the following steps: (1) enhancing interior regularity, (2) improving the regularity of $\curl_{\bar{A}^{\kappa}} w_{\ve}$ using a transport equation, (3) advancing boundary regularity through a difference quotient, and (4) combining steps (1)-(3) with estimates of $\dive_{\bar{A}^{\kappa}}w_{\ve}$ and Hodge's elliptic estimate to obtain $w_{\ve}\in L^{2}(0,T;H^{1}(\Omega))$.

Step 1: For any smooth subdomain $\Omega^{\prime}\subset\subset \Omega$, we can pick smooth function $\psi\in C_{c}^{\infty}(\Omega)$ such that $\psi =1$ in $\Omega^{\prime}$. Choosing $\epsilon^{\prime}$ small enough such that $\epsilon^{\prime}< dist(supp \psi, \Omega^{c})$, we denote standard Friederich's mollifier by $\rho_{\epsilon^{\prime}}$. For any $f\in W^{1.\infty}(\Omega)$, we investigate the equality in $H^{-1}(\Omega)$:
\begin{align}\label{estimforcommollwve}
\begin{aligned}
    \rho_{\epsilon^{\prime}}* (f\p(\psi w_{\ve}))- f\rho_{\epsilon^{\prime}}* (\p(\psi w_{\ve}))&=\int \p\rho_{\epsilon^{\prime}}(x-y)(f(x)-f(y))\psi w_{\ve}(y) dy\\&-\int  \rho_{\epsilon^{\prime}}(x-y)\p f(y) \psi w_{\ve}(y)dy.
\end{aligned}    
\end{align}
Since $w_{\ve t}\in L^{2}(0,T;(H_{0}^{1}(\Omega))^{*})$ and $w_{\ve }\in L^{2}(0,T;L^{2}(\Omega))$, we know \eqref{estimforcommollwve} is well-defined for a.e. $t\in [0,T]$.
The RHS of \eqref{estimforcommollwve} can be estimated in $L^{2}(\Omega)$ by mean value theorem and Young inequality:
\begin{align}\label{estimforcommollwve1}
\begin{aligned}
     & \|\int \p\rho_{\epsilon^{\prime}}(x-y)(f(x)-f(y))\psi w_{\ve}(y) dy-\int  \rho_{\epsilon^{\prime}}(x-y)\p f(y) \psi w_{\ve}(y)dy\|_{L^{2}(\Omega)}\\\lesssim &(\|\rho\|_{L^{1}}+\|\p\rho\|_{L^{1}})\|f\|_{W^{1,\infty}(\Omega)}\|w_{\ve}\|_{L^{2}(\Omega)}.
\end{aligned}  
\end{align}
From, \eqref{equationforlowestcurlwve}, we can write
\begin{align}
   \begin{aligned}
       \curl_{\bar{A}^{\kappa}}(\psi w_{\ve}(t))=\psi\curl u_{0}+\ve^{.ij}(\bar{A}^{\kappa})_{i}^{l}\p_{l}\psi w^{j}+\int_{0}^{t}\ve^{.ij}\psi\p_{l}w_{\ve}^{j}\p_{t}(\bar{A}^{\kappa})_{i}^{l}+\int_{0}^{t}\curl_{\bar{A}^{\kappa}} H \psi\quad\textit{in $L^{2}(0,T;(H_{0}^{1}(\Omega))^{*})$}
   \end{aligned} 
\end{align}
Invoking \eqref{estimforcommollwve1}, we can find that
\begin{align}    \curl_{\bar{A}^{\kappa}}\rho_{\epsilon^{\prime}}*(\psi w_{\ve}(t))=\psi\curl u_{0}+\ve^{.ij}(\bar{A}^{\kappa})_{i}^{l}\p_{l}\psi w^{j}+\int_{0}^{t}\ve^{.ij}\psi (\rho_{\epsilon^{\prime}}*\p_{l}w_{\ve}^{j})\p_{t}(\bar{A}^{\kappa})_{i}^{l}+\mathcal{R},
\end{align}
where $\|\mathcal{R}\|_{L_{t}^{2}L^{2}(\Omega)}\lesssim C(N_{0},M_{0})$.
The divergence type estimate can be derived in a similar strategy:
\begin{align}    \dive_{\bar{A}^{\kappa}}\rho_{\epsilon^{\prime}}*w_{\ve}+\mathcal{R},
\end{align}
where $\|\mathcal{R}\|_{L_{t}^{2}L^{2}(\Omega)}\lesssim C(N_{0},M_{0})$.
By virtue of \eqref{Hodgeelliptic}, we can show that
\begin{align}
    \|\rho_{\epsilon^{\prime}}*(\psi w_{\ve})\|_{H^{1}(\Omega)}\leq C(N_{0},M_{0})+C(N_{0},M_{0})\int_{0}^{t}\|\rho_{\epsilon^{\prime}}*(\psi w_{\ve})\|_{H^{1}(\Omega)} \quad \textit{for a.e. $t\in(0,T)$}.
\end{align}
Thus, for a $T$ small enough, we find the estimate independent on $\ve^{\prime}$ that
\begin{align}
     \|\rho_{\epsilon^{\prime}}*(\psi w_{\ve})\|_{L^{2}(0,T;H^{1}(\Omega))}\leq C(N_{0},M_{0}) .
\end{align}
It implies that $\psi w_{\ve} \in L^{2}(0,T;H^{1}(\Omega))$ and $w_{\ve}\in L^{2}(0,T;H^{1}(\Omega^{\prime}))$.

\rmk
In this instance, a modified type of the elliptic estimate \eqref{Hodgeelliptic} is employed. Readers can readily derive an appropriate version. Alternatively, we can utilize $(\bar{A}^{\kappa})_{i}^{j}-\delta_{i}^{j}=\int_{0}^{t}\p_{t}(\bar{A}^{\kappa})_{i}^{j}$ and introduce a small $T$ to arrive at a similar conclusion.

Let $\Lambda_{\epsilon^{\prime}}$ be the horizontal mollifier. We claim $\Lambda_{\epsilon^{\prime}}* w_{\ve}\in L^{2}(0,T;H^{1}(\Omega))$. By definition of the horizontal mollifier, we have
\begin{align}
    \|\pb\Lambda_{\epsilon^{\prime}}* w_{\ve}\|_{L^{2}(\T^{2}\times\{y\})}\lesssim_{\epsilon^{\prime}}\|w_{\ve}\|_{L^{2}(\T^{2}\times\{y\})}.
\end{align}
Integrating with respect to the variable $y$, we have
\begin{align}
      \|\pb\Lambda_{\epsilon^{\prime}}* w_{\ve}\|_{L^{2}(\Omega)}\lesssim_{\epsilon^{\prime}}\|w_{\ve}\|_{L^{2}(\Omega)}.
\end{align}
Now, similar to \cite{coutand2007well}, we use $\curl_{\bar{A}^{\kappa}}w_{\ve}$, $\dive_{\bar{A}^{\kappa}}w_{\ve}$ and $\pb w_{\ve}$ to express $\p_{3} w_{\ve}$:
\begin{align}   (\dive_{\bar{A}^{\kappa}}w_{\ve},\curl_{\bar{A}^{\kappa}}^{1}w_{\ve},\curl_{\bar{A}^{\kappa}}^{2}w_{\ve})^{T}=\sum_{i=1}^{3} M_{i}^{\kappa}\p_{i}w_{\ve},
\end{align}
where
\begin{align}
\begin{aligned}
      & M_{0}^{\kappa}= 
\begin{bmatrix}
(\bar{A}^{\kappa})_{1}^{1} & (\bar{A}^{\kappa})_{2}^{1} & (\bar{A}^{\kappa})_{3}^{1}\\
0 & -(\bar{A}^{\kappa})_{3}^{1} & (\bar{A}^{\kappa})_{2}^{1} \\
(\bar{A}^{\kappa})_{3}^{1} & 0 & -(\bar{A}^{\kappa})_{1}^{1}
\end{bmatrix}\quad
 M_{2}^{\kappa}= 
\begin{bmatrix}
(\bar{A}^{\kappa})_{1}^{2} & (\bar{A}^{\kappa})_{2}^{2} & (\bar{A}^{\kappa})_{3}^{2}\\
0 & -(\bar{A}^{\kappa})_{3}^{2} & (\bar{A}^{\kappa})_{2}^{2} \\
(\bar{A}^{\kappa})_{3}^{2} & 0 & -(\bar{A}^{\kappa})_{1}^{2}
\end{bmatrix}\\
& M_{3}^{\kappa}= 
\begin{bmatrix}
(\bar{A}^{\kappa})_{1}^{3} & (\bar{A}^{\kappa})_{2}^{3} & (\bar{A}^{\kappa})_{3}^{3}\\
0 & -(\bar{A}^{\kappa})_{3}^{3} & (\bar{A}^{\kappa})_{2}^{3} \\
(\bar{A}^{\kappa})_{3}^{3} & 0 & -(\bar{A}^{\kappa})_{1}^{3}.
\end{bmatrix}
\end{aligned}
\end{align}
A direct computation yields that
\begin{align}
    \det M_{3}^{\kappa}=(\bar{A}^{\kappa})_{3}^{3} \sum_{i=1}^{3}|(\bar{A}^{\kappa})_{i}^{3}|^{2}.
\end{align}
Thus, $M_{3}^{\kappa}$ is invertible. In brief, we express 
\begin{align}\label{expressionforp3}
    \p_{3}w_{\ve}=\tilde{M}\curl_{\bar{A}^{\kappa}} w_{\ve}+\tilde{N}\dive_{{\bar{A}^{\kappa}}}w_{\ve}+\sum_{i=1}^{2}\tilde{L}_{i}\p_{i}w_{\ve},
\end{align} where $\tilde{M}$, $\tilde{N}$, and $\tilde{L}_{i}$ denote matrices with entries consist of $\bar{A}^{\kappa}$.

 From \eqref{equationforwve}, it indicates that
\begin{align}
\begin{aligned}
   & \p_{t}(\curl_{\bar{A}^{\kappa}}w_{\ve})-\ve^{.ij}\p_{t}(\bar{A}^{\kappa})_{i}^{l}\p_{l}(w_{\ve})_{j}=\curl_{\bar{A}^{\kappa}}H,\\
     & \p_{t}(\curl_{\bar{A}^{\kappa}}w_{\ve})-\ve^{.ij}\p_{t}(\bar{A}^{\kappa})_{i}^{3}(\tilde{M}\curl_{\bar{A}^{\kappa}} w_{\ve}+\tilde{N}\dive_{{\bar{A}^{\kappa}}}w_{\ve}+\sum_{i=1}^{2}\tilde{L}_{i}\p_{i}w_{\ve})_{j}=\sum_{l=1}^{2}\ve^{.ij}\p_{t}(\bar{A}^{\kappa})_{i}^{l}\p_{l}(w_{\ve})_{j}+\curl_{\bar{A}^{\kappa}}H.
\end{aligned}   
\end{align}
In short, we can rewrite it as
\begin{align}
    \p_{t}(\curl_{\bar{A}^{\kappa}}w_{\ve})+\tilde{M}^{\prime}\curl_{\bar{A}^{\kappa}} w_{\ve}=\tilde{N}^{\prime}\dive_{{\bar{A}^{\kappa}}}w_{\ve}+\sum_{i=1}^{2}\tilde{L}_{i}^{\prime}\p_{i}w_{\ve}+\curl_{\bar{A}^{\kappa}}H,
\end{align}
where $\tilde{M}^{\prime}$, $\tilde{N}^{\prime}$, and $\tilde{L}_{i}^{\prime}$ are matrices in lower order. The transport equation can be solved as
\begin{align}
    \curl_{\bar{A}^{\kappa}}w_{\ve}= \tilde{P}^{-1}\curl u_{0}+\tilde{P}^{-1}\int_{0}^{t}\tilde{P}(\tilde{N}^{\prime}\dive_{{\bar{A}^{\kappa}}}w_{\ve}+\sum_{i=1}^{2}\tilde{L}_{i}^{\prime}\p_{i}w_{\ve}+\curl_{\bar{A}^{\kappa}}H),
\end{align}
with some matrix $\tilde{P}$. Plunge it back to \eqref{expressionforp3}:
\begin{align}\label{expressionforp3new}
    \p_{3}w_{\ve}=\tilde{M}\tilde{P}^{-1}\int_{0}^{t}\sum_{i=1}^{2}\tilde{P}\tilde{L}_{i}^{\prime}\p_{i}w_{\ve}+\sum_{i=1}^{2}\tilde{L}_{i}\p_{i}w_{\ve}+\mathcal{R},
\end{align}
where $\|\mathcal{R}\|_{L_{t}^{2}L^{2}(\Omega)}\lesssim C(N_{0},M_{0})$.
\rmk
In the absence of a direct estimate for the commutator of the horizontal mollifier and $\p_{3}$, we are required to represent $\p_{3}$ as $\pb$ and lower-order terms. The transport equation of ``vorticity" and the constraint of ``divergence-free" enable this strategy to be applicable for ``velocity". However, we mention that the relation \eqref{expressionforp3} is true for any function.

Now, we horizontally mollify \eqref{expressionforp3new}, it leads to the desired estimate
\begin{align}
    \| \p_{3}\Lambda_{\epsilon^{\prime}}* w_{\ve}\|_{L^{2}(0,T;L^{2}(\Omega))}\leq_{\epsilon^{\prime}} C(N_{0},M_{0}).
\end{align}
Hence, we arrive at
\begin{align}
      \| \Lambda_{\epsilon^{\prime}}* w_{\ve}\|_{L^{2}(0,T;H^{1}(\Omega))}\leq_{\epsilon^{\prime}} C(N_{0},M_{0}).
\end{align}
\rmk Up to this stage, we have already established the $L^{2}(0,T;H^{1}(\Omega))$ regularity of the mollified equation. However, it would be more convenient to use an elliptic estimate to form higher regularity in the following. Thus, we need to control the commutator of $\Lambda_{\ve}$ and $\dive_{\bar{A}^{\kappa}}$, $\curl_{\bar{A}^{\kappa}}$.

Thanks to the relation \eqref{expressionforp3new}, we have the estimate of commutator:
\begin{align}\label{goodestimateforcurldivandlamve}
\begin{aligned}
     & \|\Lambda_{\epsilon^{\prime}}*((\bar{A}^{\kappa})_{i}^{3}\p_{3}w_{\ve})-(\bar{A}^{\kappa})_{i}^{3}\Lambda_{\epsilon^{\prime}}*\p_{3}w_{\ve}\|_{L^{2}(0,T;L^{2}(\Omega))}\\\lesssim& \|[\Lambda_{\epsilon^{\prime}}*,(\bar{A}^{\kappa})_{i}^{3}](\tilde{M}\tilde{P}^{-1}\int_{0}^{t}\sum_{i=1}^{2}\tilde{P}\tilde{L}_{i}^{\prime}\p_{i}w_{\ve}+\sum_{i=1}^{2}\tilde{L}_{i}\p_{i}w_{\ve})\|_{L^{2}(0,T;L^{2}(\Omega))}+\mathcal{R}\\\leq & C(N_{0},M_{0}).
\end{aligned}  
\end{align}
It implies that $\curl_{\bar{A}^{\kappa}}$ and $\dive_{\bar{A}^{\kappa}}$ commute effectively with the horizontal mollifier $\Lambda_{\epsilon^{\prime}}$, given controlled by a constant $C(N_{0},M_{0})$. We express $\curl_{\bar{A}^{\kappa}}\Lambda_{\epsilon^{\prime}}*w_{\ve}$ and $\dive_{\bar{A}^{\kappa}}\Lambda_{\epsilon^{\prime}}*w_{\ve}$ into the combination of $\pb \Lambda_{\epsilon^{\prime}}*w_{\ve}$ and lower order term. If we can express $\p_{3}w_{\ve}$ into such forms, we can always achieve this. In our case, the domain is very simple. On $\{x_{3}=0\}$, the normal trace is always zero. On $\{x_{3}=1\}$, we apply normal trace estimate and commutation type lemma \eqref{commutationtypelemma} to obtain desired estimate. This scheme always holds for the higher regularity improvement. After this, we will use a Gr\"{o}nwall type inequality to obtain the uniform estimate.

Next, we apply horizontal mollification to the function. 
\begin{align}\label{curlakappalambdaepsprihorcolwve}
     \curl_{\bar{A}^{\kappa}}\Lambda_{\epsilon^{\prime}}*w_{\ve}= \tilde{P}^{-1}\int_{0}^{t}\tilde{P}(\tilde{N}^{\prime}\dive_{{\bar{A}^{\kappa}}}\Lambda_{\epsilon^{\prime}}*w_{\ve}+\sum_{i=1}^{2}\tilde{L}_{i}^{\prime}\p_{i}\Lambda_{\epsilon^{\prime}}*w_{\ve}+\Lambda_{\epsilon^{\prime}}*\curl_{\bar{A}^{\kappa}}H)+\mathcal{R},
\end{align}
where $\|\mathcal{R}\|_{L^{2}(0,T;L^{2}(\Omega))}\lesssim C(N_{0},M_{0})$. 
Now, we derive the estimate independent on $\epsilon^{\prime}$. Invoking \eqref{curlakappalambdaepsprihorcolwve}, we have
\begin{align}
    \|\curl_{\bar{A}^{\kappa}}\Lambda_{\epsilon^{\prime}}* w_{\ve}\|_{L^{2}(0,T;L^{2}(\Omega))}\leq C(N_{0},M_{0})+C(N_{0},M_{0})T\|\Lambda_{\epsilon^{\prime}}* w_{\ve}\|_{L^{2}(0,T;H^{1}(\Omega))}.
\end{align}
Since we always have the identity
\begin{align}
    \curl\Lambda_{\epsilon^{\prime}}* w_{\ve}=\curl_{\bar{A}^{\kappa}}\Lambda_{\epsilon^{\prime}}* w_{\ve}-\int_{0}^{t}\epsilon^{.ij}\p_{t}(\bar{A}^{\kappa})_{i}^{l}\p_{l}w_{\ve}^{j},
\end{align}
we can further write
\begin{align}\label{estimateforlamfwvecurl}
     \|\curl\Lambda_{\epsilon^{\prime}}* w_{\ve}\|_{L^{2}(0,T;L^{2}(\Omega))}\leq C(N_{0},M_{0})+C(N_{0},M_{0})T\|\Lambda_{\epsilon^{\prime}}* w_{\ve}\|_{L^{2}(0,T;H^{1}(\Omega))}.
\end{align}
On the other hand, in light of \eqref{goodestimateforcurldivandlamve}, \eqref{estimateforwveinitiall2l2} and time evolving equation, we have
\begin{align}\label{estimateforlamfwvediv}
    \|\dive\Lambda_{\epsilon^{\prime}}* w_{\ve}\|_{L^{2}(0,T;L^{2}(\Omega))}\leq C(N_{0},M_{0})+C(N_{0},M_{0})T\|\Lambda_{\epsilon^{\prime}}* w_{\ve}\|_{L^{2}(0,T;H^{1}(\Omega))}.
\end{align}
Notice that $\int_{0}^{t}\p_{t}\bar{n}^{\kappa}=\bar{n}^{\kappa}-N$, the commutator estimate and\eqref{estimateforwveinitiall2l2} give us
\begin{align}\label{estimateforlamfwvebdy}
    \|\pb \Lambda_{\ve^{\prime}}* w_{\ve}\cdot N\|_{L^{2}(0,T;H^{-\frac{1}{2}}(\Gamma_{1}))}\leq C(N_{0},M_{0})+C(N_{0},M_{0})T\|\Lambda_{\epsilon^{\prime}}* w_{\ve}\|_{L^{2}(0,T;H^{1}(\Omega))}.
\end{align}
\rmk
In our case, $N=(0,0,1)$; the condition facilitates a more concise proof. In the general case, one can apply the trace inequality to control the commutator term, see \cite{coutand2007well}.

Combining \eqref{Hodgeelliptic}, \eqref{estimateforlamfwvecurl}, \eqref{estimateforlamfwvediv}, and \eqref{estimateforlamfwvebdy}, we have
\begin{align}
     \| \Lambda_{\epsilon^{\prime}}* w_{\ve}\|_{L^{2}(0,T;H^{1}(\Omega))}\leq C(N_{0},M_{0})+C(N_{0},M_{0})T\|\Lambda_{\epsilon^{\prime}}* w_{\ve}\|_{L^{2}(0,T;H^{1}(\Omega))}.
\end{align}
Choosing small $T$ if necessary, we achieve that
\begin{align}
     \| \Lambda_{\epsilon^{\prime}}* w_{\ve}\|_{L^{2}(0,T;H^{1}(\Omega))}\leq C(N_{0},M_{0}).
\end{align}
Letting $\epsilon^{\prime}\to 0$, we obtain $ w_{\ve}\in L^{2}(0,T;H^{1}(\Omega))$ and
\begin{align}
   \|  w_{\ve}\|_{L^{2}(0,T;H^{1}(\Omega))}\leq C(N_{0},M_{0}).
\end{align}
\subsection{Proof of \eqref{penal2}}
For $p\in H^{-\frac{1}{2}}(\Omega)$, we define a linear operator $Q(t):H^{-\frac{1}{2}}(\Omega)\to H^{\frac{3}{2}}(\Omega)$ by Riesz representation theorem
\begin{align}
    \forall \varphi\in H^{\frac{3}{2}}(\Omega),\quad \langle p, \bar{J}^{\kappa}(\bar{A}^{\kappa})_{i}^{j}\p_{j}\varphi^{i}\rangle_{\frac{1}{2}}=(Q(t)p,\varphi)_{\frac{3}{2}}.
\end{align}
Clearly,
\begin{align}
    \|Q(t)p\|_{H^{\frac{3}{2}}(\Omega)}\lesssim \|p\|_{H^{-\frac{1}{2}}(\Omega)}.
\end{align}
\defn
\begin{align}
    \mathcal{V}_{\bar{v}}(t)=\{v\in L^{2}(\Omega)|\bar{J}^{\kappa}(\bar{A}^{\kappa})_{i}^{j}\p_{j}v^{i}(t)=0\}.
\end{align}
Due to $\mathcal{V}_{\bar{v}}(t)\cap H^{\frac{3}{2}}(\Omega)= R(Q(t))^{\perp}$, it implies that
\begin{align}\label{orthogonalrelationH32}
    H^{\frac{3}{2}}(\Omega)=R(Q(t))\oplus_{H^{\frac{3}{2}}(\Omega)} \mathcal{V}_{\bar{v}}(t)\cap H^{\frac{3}{2}}(\Omega).
\end{align}
\defn
\begin{align}
    X(t)=\{\psi\in H^{\frac{1}{2}}(\Omega)|\bar{J}^{\kappa}(\bar{A}^{\kappa})_{i}^{j}\p_{j}\psi^{i}\in H^{-\frac{1}{2}}(\Omega)\}.
\end{align}
We define a linear functional $Q^{\prime}(t):H^{\frac{1}{2}}(\Omega)\to X(t)$ via Riesz representation theorem, such that
\begin{align}
   \forall \varphi\in X(t),\quad \langle(\bar{A}^{\kappa})_{i}^{j}
   \p_{j}\varphi^{i},p\rangle_{H^{\frac{1}{2}} (\Omega)}=(Q^{\prime}(t)p,\varphi)_{X(t)}.
\end{align}
One can find that $R(Q(t))$ is closed in $X(t)$ and
\begin{align}
    \|p\|_{H^{\frac{1}{2}}(\Omega)}\lesssim \|Q(t)p\|_{X(t)},
\end{align}
where the norm of $X(t)$ is taken as a subspace of $H^{\frac{1}{2}}(\Omega)$. Owing to $\mathcal{V}_{\bar{v}}(t)\cap X(t)=R(Q(t)^{\prime})^{\perp}$, it yields that
\begin{align}
    X(t)=R(Q(t)^{\prime}))\oplus_{X(t)}\mathcal{V}_{\bar{v}}(t)\cap X(t).
\end{align}
Here, we will use two lemmas in \cite{coutand2007well}.
\lem[Lagrange Multiplier 1]
Let $\mathcal{L}(t)\in H^{-\frac{3}{2}}(\Omega)$ be such that $\mathcal{L}(t)\varphi=0$ for any $\varphi \in \mathcal{V}_{\bar{V}}(t)\cap H^{\frac{3}{2}}(\Omega)$. Then there exists a unique $q(t)\in H^{-\frac{1}{2}}(\Omega)$, which is termed the pressure function, satisfying
\begin{align}
    \forall \varphi\in H^{\frac{3}{2}}(\Omega), \quad \mathcal{L}(t)(\varphi)=\langle q(t), \bar{J}_{\kappa}(\bar{A}^{\kappa})_{i}^{j}\p_{j}\varphi^{i}\rangle_{H^{\frac{1}{2}}(\Omega)}.
\end{align}
Moreover, there is a $C>0$ (which does not depend on $t\in[0,T]$ and on the choice of $\bar{v}\in S_{0}$) such that
\begin{align}
    \|q(t)\|_{H^{-\frac{1}{2}}(\Omega)}\lesssim \|\mathcal{L}(t)\|_{H^{-\frac{3}{2}}(\Omega)}.
\end{align}
\lem[Lagrange Multiplier 2]
Let $\mathcal{L}(t)\in X(t)^{\prime}$ be such that $\mathcal{L}(t)\varphi=0$ for any $\varphi \in \mathcal{V}_{\bar{v}}(t)\cap H^{\frac{1}{2}}(\Omega)$. Then there exists a unique $q(t)\in H^{\frac{1}{2}}(\Omega)$, which is termed the pressure function, satisfying 
\begin{align}
    \forall \varphi\in X(t),\quad \mathcal{L}(t)(\varphi)=\langle\bar{J}_{\kappa}(\bar{A}^{\kappa})_{i}^{j}\p_{j}\varphi^{i},q(t)\rangle_{H^{\frac{1}{2}}(\Omega)}.
\end{align}
Moreover, there is a $C>0$ (which does not depend on $t\in[0,T]$ and on the choice of $\bar{v}\in S_{0}$) such that
\begin{align}
    \|q(t)\|_{H^{\frac{1}{2}}(\Omega)}\lesssim \|\mathcal{L}(t)\|_{X(t)^{\prime}}.
\end{align}
Recall that \eqref{estimateforwveinitiall2l2}
\begin{align}
    \|w_{\ve}\|_{L^{2}(0,T;L^{2}(\Omega))}+\frac{1}{\ve^{\frac{1}{2}}}\|\dive_{\bar{A}^{\kappa}}w_{\ve}\|_{L^{2}(0,T;L^{2}(\Omega))}+\|w_{\ve}\cdot \bar{n}_{\kappa}\|_{L^{2}(0,T;H^{1}(\Gamma_{1}))}\leq C(N_{0},M_{0}).
\end{align}
Passing $\ve\to 0$, we have
\begin{align}
    \begin{aligned}
         w_{\ve}&\rightharpoonup w \quad \text{in $L^{2}(0,T;L^{2}(\Omega))$},\\
           \dive_{\bar{A}^{\kappa}}(w_{\ve})&\rightharpoonup  \dive_{\bar{A}^{\kappa}}(w)\quad \text{in $L^{2}(0,T;L^{2}(\Omega))$},\\
           w_{\ve}\cdot \bar{n}_{\kappa}&\rightharpoonup w\cdot\bar{n}_{\kappa} \text{in $L^{2}(0,T;H^{1}(\Gamma_{1}))$}.
    \end{aligned}
\end{align}
The reason why $\dive_{\bar{A}^{\kappa}}(w)$ and $w\cdot\bar{n}_{\kappa}$ precisely represent the weak limit is akin to the explanation provided earlier. Further, by the definition of weak limit,
\begin{align}    \|\dive_{\bar{A}^{\kappa}}w\|_{L^{2}(0,T;L^{2}(\Omega))}=0.
\end{align}
We can employ a similar approach in step 2: 1) $\dive_{\bar{A}^{\kappa}}w=0$ provides us with $L_{t}^{2}H^{\frac{3}{2}}(\Omega)$ regularity of the divergence, 2) confirm the mollified function in $L^{2}_{t}H^{\frac{3}{2}}(\Omega)$, 3) find the uniform bound of the approximation in $L^{2}_{t}H^{\frac{3}{2}}(\Omega)$ via Gr\"{o}nwall-type inequality, 4) choose an appropriate small time and justify $w\in L^{2}_{t}H^{\frac{3}{2}}(\Omega)$, which leads us to:
\begin{align}
    \|w\|_{L^{2}(0,T;H^{\frac{3}{2}}(\Omega))}\leq C(N_{0},M_{0}).
\end{align}
For any $y\in L^{2}(0,T;H^{\frac{3}{2}}(\Omega))$ and $l=(\bar{A}^{\kappa})_{i}^{j}\p_{j}y^{i}$, we can find a solution to the following elliptic equation:
\begin{align}
    \begin{aligned}        (\bar{A}^{\kappa})_{i}^{j}\p_{j}((\bar{A}^{\kappa})_{i}^{k}\p_{k}\varphi)&=l\quad \textit{in $\Omega$,}\\
    \varphi&=0\quad\textit{on $\p\Omega$.}
    \end{aligned}
\end{align}
Let $e^{i}=(\bar{A}^{\kappa})_{i}^{k}\p_{k}\varphi$ and $s=y-e$. We check that $s,e\in L^{2}(0,T;H^{\frac{3}{2}}(\Omega))$ and by orthogonal relation \eqref{orthogonalrelationH32}:
\begin{align}
    \begin{aligned}
    (\bar{A}^{\kappa})_{i}^{j}\p_{j}s^{i}=0,\quad \textit{$s(t)\in \mathcal{V}_{\bar{v}}\cap H^{\frac{3}{2}}(\Omega)$ a.e. t},\\
       \|e\|_{L^{2}(0,T;H^{\frac{3}{2}}(\Omega))}+\|s\|_{L^{2}(0,T;H^{\frac{3}{2}}(\Omega)) }\lesssim \|y\|_{L^{2}(0,T;H^{\frac{3}{2}}(\Omega))}  .   
    \end{aligned}
\end{align}
Since $(\bar{A}^{\kappa})_{i}^{j}\p_{j}w^{i}=0$, we deduce that
\begin{align}
    \begin{aligned}    (\bar{A}^{\kappa})_{i}^{j}\p_{j}w_{t}^{i}&=-\p_{t}(\bar{A}^{\kappa})^{j}_{i}\p_{j}w^{i} \in L^{2}(0,T;H^{\frac{1}{2}}(\Omega))   \\
    \langle \bar{J}_{\kappa}w_{t},e\rangle_{H^{\frac{3}{2}}(\Omega)}&=\int_{\Omega}w_{t}\bar{J}_{\kappa}(\bar{A}^{\kappa})_{i}^{k}\p_{k}\varphi=-\int_{\Omega}\p_{k}w_{t}\bar{J}_{\kappa}(\bar{A}^{\kappa})_{i}^{k}\varphi\\&=(\p_{t}(\bar{A}^{\kappa})_{i}^{j}\p_{j}w^{i},\bar{J}_{\kappa}\varphi)_{L^{2}(\Omega)}
    \end{aligned}
\end{align}
$v$ also satisfies equation
\begin{align}
    \int_{0}^{T}\langle\bar{J}_{\kappa} w_{\ve t},s\rangle_{H^{\frac{3}{2}}(\Omega)}+\kappa\int_{0}^{T}(w_{\ve}\cdot \bar{n}_{\kappa},s\cdot\bar{n}_{\kappa})_{H^{1}(\Gamma_{1})}=\sigma\int_{0}^{T}(L_{\bar{h}}\bar{\eta}\cdot \bar{n}_{\kappa}, s\cdot \bar{n}_{\kappa})_{L^{2}(\Gamma_{1})}+\int_{0}^{T}(\bar{J}_{\kappa}H,s)_{L^{2}(\Omega)}.
\end{align}
Now, we can check that
\begin{align}
\begin{aligned}
     \lim_{\ve\to 0}\int_{0}^{T}\langle\bar{J}_{\kappa}w_{\ve t}, y\rangle_{H^{\frac{3}{2}}(\Omega)}&=\lim_{\ve\to 0}\int_{0}^{T}\langle\bar{J}_{\kappa}w_{\ve t}, e\rangle_{H^{\frac{3}{2}}(\Omega)}+\lim_{\ve\to 0}\int_{0}^{T}\langle\bar{J}_{\kappa}w_{\ve t}, s\rangle_{H^{\frac{3}{2}}(\Omega)}\\&=\int_{0}^{T}(\p_{t}(\bar{A}^{\kappa})_{i}^{j}\p_{j}w^{i},\bar{J}_{\kappa}\varphi)_{L^{2}(\Omega)}+\sigma\int_{0}^{T}(L_{\bar{h}}\bar{\eta}\cdot \bar{n}_{\kappa}, s\cdot \bar{n}_{\kappa})_{L^{2}(\Gamma_{1})}+\int_{0}^{T}(\bar{J}_{\kappa}H,s)_{L^{2}(\Omega)}\\&-\kappa\int_{0}^{T}(w\cdot \bar{n}_{\kappa},s\cdot\bar{n}_{\kappa})_{H^{1}(\Gamma_{1})}\\&\leq C(N_{0},M_{0})\|y\|_{L^{2}(0,T;H^{\frac{3}{2}}(\Omega))}.
\end{aligned}  
\end{align}
Hence, we deduce that for some small $\ve_{0}$ and $\ve<\ve_{0}$
\begin{align}
    \|w_{\ve t}\|_{L^{2}(0,T;H^{-\frac{3}{2}}(\Omega))}\leq C(N_{0},M_{0}).
\end{align}
Letting $\ve\to 0$, we conclude that
\begin{align}\label{weaklimitofwvet}
          w_{\ve t}\rightharpoonup w_{t}\quad \text{in $L^{2}(0,T;H^{-\frac{3}{2}}(\Omega))$}.
\end{align}
Likewise, we can confirm that $w_{t}$ represents the weak limit on the right hand side of \eqref{weaklimitofwvet}. Combining with $w\in L^{2}(0,T;H^{\frac{3}{2}}(\Omega))$, we can conclude that $w\in C(0,T;L^{2}(\Omega))$. Moreover, $w_{\ve}(0)=u_{0}$ gives us $w(0)=u_{0}$. Now, for all $\psi\in L^{2}(0,T;H^{\frac{3}{2}}(\Omega))$ such that $(\bar{A}^{\kappa})_{i}^{j}\p_{j}\psi^{i}=0$ a.e. $t$, we define distribution $\mathcal{L}(t)\in H^{-\frac{3}{2}}(\Omega)$ as
\begin{align}
\begin{aligned}
    \mathcal{L}(t)(\psi)&=\langle \bar{J}_{\kappa}w_{t},\psi\rangle_{H^{\frac{3}{2}}(\Omega)}+\kappa(w\cdot \bar{n}_{\kappa},\psi\cdot\bar{n}_{\kappa})_{H^{1}(\Gamma_{1})}\\-&\sigma(L_{\bar{h}}\bar{\eta}\cdot \bar{n}_{\kappa}, \psi\cdot \bar{n}_{\kappa})_{L^{2}(\Gamma_{1})}-(\bar{J}_{\kappa}H,\psi)_{L^{2}(\Omega)}\textit{a.e. $t$, $\forall \psi\in H^{\frac{3}{2}}(\Omega)$}.
\end{aligned}   
\end{align}
By Lagrange Multiplier 1, there exists $q\in H^{-\frac{1}{2}}(\Omega)$ such that
\begin{align}
\begin{aligned}
     \langle q(t),\bar{J}_{\kappa}(\bar{A}^{\kappa})_{i}^{j}\p_{j}\psi^{i}\rangle_{H^{\frac{1}{2}}(\Omega)}&=\langle \bar{J}_{\kappa}w_{t},\psi\rangle_{H^{\frac{3}{2}}(\Omega)}+\kappa(w\cdot \bar{n}_{\kappa},\psi\cdot\bar{n}_{\kappa})_{H^{1}(\Gamma_{1})}\\&+\sigma(L_{\bar{h}}\bar{\eta}\cdot \bar{n}_{\kappa}, \psi\cdot \bar{n}_{\kappa})_{L^{2}(\Gamma_{1})}+(\bar{J}_{\kappa}H,\psi)_{L^{2}(\Omega)}\quad\textit{$\forall \psi\in H^{\frac{3}{2}}(\Omega)$, a.e. $t$}.
\end{aligned}  
\end{align}
Now, for all $\tilde{\psi}\in L^{2}(0,T;H^{\frac{3}{2}}(\Omega))$, integrating with respect to time, it yields that
\begin{align}\label{equationidentityforvandvbar}
\begin{aligned}
     &\int_{0}^{T}\langle \bar{J}_{\kappa}w_{t},\tilde{\psi}\rangle_{H^{\frac{3}{2}}(\Omega)}  -\int_{0}^{T}\langle q(t),\bar{J}_{\kappa}(\bar{A}^{\kappa})_{i}^{j}\p_{j}\tilde{\psi}^{i}\rangle_{H^{\frac{1}{2}}(\Omega)}+\int_{0}^{T}\kappa(w\cdot \bar{n}_{\kappa},\tilde{\psi}\cdot\bar{n}_{\kappa})_{H^{1}(\Gamma_{1})}\\=&\int_{0}^{T}\sigma(L_{\bar{h}}\bar{\eta}\cdot \bar{n}_{\kappa}, \tilde{\psi}\cdot \bar{n}_{\kappa})_{L^{2}(\Gamma_{1})}+\int_{0}^{T}(\bar{J}_{\kappa}H,\tilde{\psi})_{L^{2}(\Omega)}
\end{aligned}
\end{align}
Hence, $w$ is the weak solution. Suppose $\tilde{w}$ is another weak solution in the space $L^{2}(0,T;H^{\frac{3}{2}}(\Omega))$ such that $\tilde{w}(0)=u_{0}$ and $\tilde{w}_{t}\in L^{2}(0,T;H^{\frac{3}{2}}(\Omega)^{\prime})$. Taking $\bar{w}=w-\tilde{w}$, we have
\begin{align}
    \int_{0}^{T}\langle \bar{J}_{\kappa}\bar{w}_{t},\tilde{\psi}\rangle_{H^{\frac{3}{2}}(\Omega)}+\int_{0}^{T}\kappa(\bar{w}\cdot \bar{n}_{\kappa},\tilde{\psi}\cdot\bar{n}_{\kappa})_{H^{1}(\Gamma_{1})}=0.
\end{align}
Taking $\tilde{\psi}=\bar{w}$, after a simple calculation, we have $\bar{w}=0$. As a consequence, the weak solution $w$ is unique. Later, we still denote the weak solution by $v$.
\subsection{Proof of \eqref{kappaapp}}
Now, we improve the regularity of $v$. The first step is to improve the regularity on the boundary. Then, we can utilize a similar strategy in the previous subsection to establish higher regularity. Similarly in \cite{coutand2007well}, we employ the difference quotient method. We use $D_{h}$ to represent the difference quotient operator in the direction $e_{1}$ or $e_{2}$. Invoking \eqref{equationidentityforvandvbar}, setting $\tilde{\psi}=D_{-h}D_{h}v$, it yields that
\begin{align}
       \begin{aligned}          &\underbrace{\int_{0}^{T}\langle \bar{J}_{\kappa}v_{t},D_{-h}D_{h}v\rangle_{H^{\frac{3}{2}}(\Omega)}}_{(1)}  \underbrace{-\int_{0}^{T}\langle q(t),\bar{J}_{\kappa}(\bar{A}^{\kappa})_{i}^{j}\p_{j}D_{-h}D_{h}v^{i}\rangle_{H^{\frac{1}{2}}(\Omega)}}_{(2)}\underbrace{+\int_{0}^{T}\kappa(v\cdot \bar{n}_{\kappa},D_{-h}D_{h}v\cdot\bar{n}_{\kappa})_{H^{1}(\Gamma_{1})}}_{(3)}\\=&\underbrace{\int_{0}^{T}\sigma(L_{\bar{h}}\bar{\eta}\cdot \bar{n}_{\kappa}, D_{-h}D_{h}v\cdot \bar{n}_{\kappa})_{L^{2}(\Gamma_{1})}}_{(4)}\underbrace{+\int_{0}^{T}(\bar{J}_{\kappa}H,D_{-h}D_{h}v)_{L^{2}(\Omega)}}_{(5)}
       \end{aligned}
\end{align}
We analyze them term by term: by Young's inequality,
\begin{align}\label{estimateforbra1}
    \begin{aligned}
        (1)&=\int_{0}^{T}\int_{\Omega}D_{h} (\bar{J}_{\kappa}v_{t}),D_{h}v\\&=\int_{0}^{T}\int_{\Omega}\bar{J}_{\kappa}(x+h)D_{h} v_{t}D_{h}v+\int_{0}^{T}\int_{\Omega}v_{t}D_{h}\bar{J}_{\kappa} D_{h}v\\&\geq\frac{1}{2}\|D_{h} v\|_{L^{2}(\Omega)}^{2}-C(\|u_{0}\|_{H^{1}(\Omega)})-C(N_{0})\|v\|_{L^{2}_{t}H^{1}(\Omega)}-\int_{0}^{T}\|v_{t}\|_{H^{-\frac{3}{2}}(\Omega)}\|D_{h}\bar{J}_{\kappa} D_{h}v\|_{H^{\frac{3}{2}}(\Omega)}\\&\geq\frac{1}{2}\|D_{h} v\|_{L^{2}(\Omega)}^{2}-C_{\delta}(N_{0})-\delta\int_{0}^{T}\|D_{h}v\|_{H^{\frac{3}{2}}(\Omega)}^{2}.
    \end{aligned}
\end{align}
Notice that
\begin{align}
\begin{aligned}
     \bar{J}_{\kappa}(\bar{A}^{\kappa})_{i}^{j}\p_{j}D_{-h}D_{h}v^{i}&= -D_{-h}D_{h}(\bar{J}_{\kappa}(\bar{A}^{\kappa})_{i}^{j})\p_{j}v^{i}- D_{h}(\bar{J}_{\kappa}(\bar{A}^{\kappa})_{i}^{j})D_{-h}\p_{j}v^{i}(x+h)\\&- D_{-h}(\bar{J}_{\kappa}(\bar{A}^{\kappa})_{i}^{j})D_{h}\p_{j}v^{i}(x-h).
\end{aligned}  
\end{align}
Thus, we have
\begin{align}\label{estimateforbra2}
\begin{aligned}
     |(2)|&\lesssim\int_{0}^{T}\|q\|_{H^{-\frac{1}{2}}}C(N_{0})(\|D_{h}\p v\|_{H^{\frac{1}{2}}(\Omega)}+\|\p v\|_{H^{\frac{1}{2}}} )\\&\lesssim C_{\delta}(N_{0},M_{0})+\delta\int_{0}^{T}\|D_{h}v\|_{H^{\frac{3}{2}}(\Omega)}^{2} .
\end{aligned} 
\end{align}
The estimate of $(3)$ reads:
\begin{align}
\begin{aligned}     (3)&=\sum_{i=1}^{2}\kappa\int_{0}^{T}\int_{\Omega}\pb_{i}(D_{h}v\cdot\bar{n}_{\kappa})\pb_{i}(D_{h}v\cdot\bar{n}_{\kappa})+\kappa\int_{0}^{T}\int_{\Omega}(D_{h}v\cdot\bar{n}_{\kappa})(D_{h}v\cdot\bar{n}_{\kappa})\\&+\kappa\int_{0}^{T}\int_{\Omega}v(x+h)\cdot D_{h}(\bar{n}_{\kappa}\bar{n}_{\kappa})\cdot D_{h}v-\sum_{i=1}^{2}\kappa\int_{0}^{T}\int_{\Omega}\pb_{i}(v\cdot\bar{n}_{\kappa})\pb_{i}(D_{h}v(x-h)\cdot D_{-h}\bar{n}_{\kappa})\\&-\sum_{i=1}^{2}\kappa\int_{0}^{T}\int_{\Omega}\pb_{i}(v(x-h)\cdot D_{h}\bar{n}_{\kappa})\pb_{i}(D_{h}v\cdot\bar{n}_{\kappa}).
\end{aligned}   
\end{align}
Applying Young's inequality again, we have
\begin{align}\label{estimateforbra3}
\begin{aligned}
     (3)&\geq \kappa\int_{0}^{T}\|D_{h}v\cdot \bar{n}_{\kappa}\|_{H^{1}(\Gamma_{1})}-C_{\delta}(N_{0},M_{0})-\delta\kappa\int_{0}^{T}\|D_{h}v\|_{H^{1}(\Gamma_{1})}^{2}\\&\geq\kappa\int_{0}^{T}\|D_{h}v\cdot \bar{n}_{\kappa}\|_{H^{1}(\Gamma_{1})}^{2}-C_{\delta}(N_{0},M_{0})-\delta\kappa\int_{0}^{T}\|D_{h}v\|_{H^{1.5}(\Omega)}^{2}
\end{aligned}   
\end{align}
$(4)$ and $(5)$ can be controlled directly by
\begin{align}\label{estimateforbra4}
\begin{aligned}
     |(4)|&\lesssim \sigma\int_{0}^{T}\|L_{\bar{h}}\bar{\eta}\cdot \bar{n}_{\kappa}\|_{L^{2}(\Gamma_{1})}\|D_{h}v\|_{H^{1}(\Gamma_{1})}\|\bar{n}_{\kappa}\|_{L^{\infty}(\Gamma_{1})}\\&\lesssim \delta\int_{0}^{T}\|D_{h}v\|_{H^{1.5}(\Omega)}^{2}+ C_{\delta}(N_{0},M_{0}),
\end{aligned}   
\end{align}
and
\begin{align}\label{estimateforbra5}
    \begin{aligned}
        |(5)|&\lesssim \int_{0}^{T}\|D_{h}(\bar{J}_{\kappa} H)\|_{L^{2}(\Omega)}\|D_{h}v\|_{L^{2}(\Omega)}\\&\lesssim  C(N_{0},M_{0})
    \end{aligned}
\end{align}
\rmk
The estimate of $(1)-(4)$ is identical to \cite{coutand2007well}. It can be verified that by introducing a higher-order term $D_{h_{1}}D_{h_{2}}...D_{h_{k}}$ to equation $(5)$ and evaluating its $L_{t}^{2}L^{2}{x}$ norm, it can still controlled by $C(N_{0},M_{0},\kappa)$. To see this, because of the horizontal mollification of $\bar{J}^{\kappa}$, $\bar{A}^{\kappa}$, and $\bar{n}^{\kappa}$, their respective estimates in the $\pb$ direction remain controllable. The terms involving $\phi^{\kappa}$ can be directly estimated using the properties of convolution.

Combining \eqref{estimateforbra1}, \eqref{estimateforbra2}, \eqref{estimateforbra3}, \eqref{estimateforbra4}, and \eqref{estimateforbra5}, we have
\begin{align}\label{estimateforDhvbdyini}
 \int_{0}^{T}\|D_{h}v\cdot \bar{n}_{\kappa}\|_{H^{1}(\Gamma_{1})}^{2}\lesssim\delta\int_{0}^{T}\|D_{h}v\|_{H^{1.5}(\Omega)}^{2}+ C_{\delta}(N_{0},M_{0}).
\end{align}
The $\dive_{\bar{A}^{\kappa}}$ and $\curl_{\bar{A}^{\kappa}}$ estimates read
\begin{align}\label{estimateforDhvdiveini}
\begin{aligned}
      \dive_{\bar{A}^{\kappa}} D_{h}v&=-D_{h}(\bar{A}^{\kappa})_{i}^{j}\p_{j}v^{i}(x+h)\\
      \| \dive_{\bar{A}^{\kappa}} D_{h}v\|_{L_{t}^{2}H^{0.5}(\Omega)}^{2}&\lesssim C(N_{0},M_{0}),
\end{aligned}  
\end{align}
and 
\begin{align}
    \begin{aligned}       \curl_{\bar{A}^{\kappa}}v(t)&=\curl u_{0}+\int_{0}^{t}\ve^{.ij}\p_{l}v^{j}\p_{t}(\bar{A}^{\kappa})_{i}^{l}-\int_{0}^{t}\curl_{\bar{A}^{\kappa}} (\p_{t}\phi^{\kappa}\nabla_{\bar{A}^{\kappa}}\phi^{\kappa}+\nabla_{\bar{A}^{\kappa}}|\nabla_{\bar{A}^{\kappa}}\phi^{\kappa}|^{2}) \\&=\curl u_{0}+\int_{0}^{t}\ve^{.ij}\p_{l}v^{j}\p_{t}(\bar{A}^{\kappa})_{i}^{l}-\int_{0}^{t}\ve^{.ij}(\bar{A}^{\kappa})_{i}^{k}\p_{t}\p_{k}\phi^{\kappa}\nabla_{\bar{A}^{\kappa}}\phi^{\kappa}.
    \end{aligned}
\end{align}
Moreover,
\begin{align}\label{estimateforDhvcurlini}
    \begin{aligned}
\curl_{\bar{A}^{\kappa}}D_{h}v(t)&=-\ve^{.ij}D_{h}(\bar{A}^{\kappa})_{i}^{k}\p_{k}v^{j}(x+h)+D_{h}\curl u_{0}\\&+\int_{0}^{t}\ve^{.ij}D_{h}(\p_{l}v^{j}\p_{t}(\bar{A}^{\kappa})_{i}^{l})-\int_{0}^{t}\ve^{.ij}D_{h}((\bar{A}^{\kappa})_{i}^{k}\p_{t}\p_{k}\phi^{\kappa}\nabla_{\bar{A}^{\kappa}}\phi^{\kappa})\\
\|\curl_{\bar{A}^{\kappa}}D_{h}v(t)\|_{L_{t}^{2}H^{0.5}}^{2}&\lesssim C_{\delta}(N_{0},M_{0})+\delta\int_{0}^{T}\|D_{h}v\|_{H^{1.5}(\Omega)}^{2}
    \end{aligned}
\end{align}
Hence, combining \eqref{estimateforDhvbdyini}, \eqref{estimateforDhvdiveini}, and \eqref{estimateforDhvcurlini}, we have
\begin{align}   \int_{0}^{T}\|D_{h}v\|_{H^{1.5}(\Omega)}^{2}\lesssim C_{\delta}(N_{0},M_{0})+\delta\int_{0}^{T}\|D_{h}v\|_{H^{1.5}(\Omega)}^{2}
\end{align}
Choosing $\delta$ small, we have
\begin{align}    \int_{0}^{T}\|D_{h}v\|_{H^{1.5}(\Omega)}^{2}\lesssim C_{\delta}(N_{0},M_{0}).
\end{align}
Consequently, the trace inequality shows that
\begin{align}  \int_{0}^{T}\|D_{h}v\|_{H^{1}(\Gamma_{1})}^{2}\lesssim C(N_{0},M_{0}).
\end{align}
Passing $h\to 0$, we get the desired inequality
\begin{align}
    \int_{0}^{T}\|v\|_{H^{2}(\Gamma_{1})}^{2}\lesssim C(N_{0},M_{0}).
\end{align}
Combining this with the argument in previous subsection, we have
\begin{align}
    \int_{0}^{T}\|v\|_{H^{2.5}(\Omega)}^{2}\lesssim C(N_{0},M_{0}).
\end{align}
Now, we estimate the pressure with the Lagrangian multiplier 2. For any $y\in X$, $\varphi$ is a solution of the elliptic problem
\begin{align}
    \begin{aligned}
        \bar{J}_{\kappa}(\bar{A}^{\kappa})_{i}^{j}\p_{j}(\bar{J}_{\kappa}(\bar{A}^{\kappa})_{i}^{k}\p_{k}\varphi)&=\bar{J}_{\kappa}(\bar{A}^{\kappa})_{i}^{k}\p_{k}y^{i}\quad \textit{in $H^{-\frac{1}{2}}(\Omega)$},\\
        \varphi&=0 \quad \textit{on $\p\Omega$}.
    \end{aligned}
\end{align}
By interpolation theorem, $\varphi$ is solvable and $\varphi\in H^{\frac{3}{2}}(\Omega)$. We still introduce $e^{i}=\bar{J}^{\kappa}(\bar{A}^{\kappa})_{i}^{k}\p_{k}\varphi$ and $s=y-e$, we have $e\in H^{\frac{1}{2}}(\Omega)$ and $s\in \mathcal{V}\cap H^{\frac{1}{2}}(\Omega)$. By orthogonal relationship, we have $\|e\|_{H^{\frac{1}{2}}(\Omega)}+\|s\|_{H^{\frac{1}{2}}(\Omega)}\leq \|y\|_{H^{\frac{1}{2}}(\Omega)}$. Similarly,
\begin{align}
    \langle\bar{J}_{\kappa}v_{t},e\rangle_{X}=(\p_{t}(\bar{J}_{\kappa}(\bar{A}^{\kappa})_{i}^{j})\p_{j}v^{i},\varphi)_{L^{2}(\Omega)}
\end{align}
and 
\begin{align}\label{equationforvtsX}
 \int_{0}^{T}  \langle\bar{J}_{\kappa}v_{t},s\rangle_{X} +\kappa\int_{0}^{T}\int_{\Gamma_{1}}(1-\lap)(v\cdot\bar{n}_{\kappa})s\cdot\bar{n}_{\kappa}=\sigma\int_{0}^{T}\int_{\Gamma_{1}}L_{\bar{h}}\bar{\eta}\cdot \bar{n}_{\kappa} s\cdot\bar{n}_{\kappa}+\int_{0}^{T}(\bar{J}_{\kappa}H,s)_{L^{2}(\Omega)}.
\end{align}
Notice that
\begin{align}
    \begin{aligned}
       \kappa\int_{0}^{T}\int_{\Gamma_{1}}(1-\bar{\lap})(v\cdot\bar{n}_{\kappa})s\cdot\bar{n}_{\kappa}&=   -\kappa\int_{0}^{T}\int_{\Omega}\bar{J}^{\kappa}(\bar{A}^{\kappa})_{i}^{j}\p_{j}(\frac{(1-\bar{\lap})(v\cdot\bar{n}_{\kappa})}{|\bar{J}^{\kappa}(\bar{A}^{\kappa})_{i}^{3}|})s^{i}\\&= -\kappa\int_{0}^{T}\int_{\Omega}\bar{J}^{\kappa}(\ak)_{i}^{j}\p_{j}(\frac{1}{|\bar{J}^{\kappa}(\bar{A}^{\kappa})_{i}^{3}|})(1-\bar{\lap})(v\cdot\bar{n}_{\kappa})s^{i}\\&- \kappa\int_{0}^{T}\int_{\Omega}\frac{\bar{J}^{\kappa}(\bar{A}^{\kappa})_{i}^{j}}{|\bar{J}^{\kappa}(\bar{A}^{\kappa})_{i}^{3}|}((1-\bar{\lap})\p_{j}(v\cdot\bar{n}_{\kappa}))s^{i}.
    \end{aligned}
\end{align}
Thus,
\begin{align}
\begin{aligned}
     | \kappa\int_{0}^{T}\int_{\Gamma_{1}}(1-\bar{\lap})(v\cdot\bar{n}_{\kappa})s\cdot\bar{n}_{\kappa}|&\lesssim\kappa\int_{0}^{T}\|\bar{J}^{\kappa}(\bar{A}^{\kappa})_{i}^{j}\p_{j}(\frac{1}{|\bar{J}^{\kappa}(\bar{A}^{\kappa})_{i}^{3}|})\|_{L^{\infty}(\Omega)}\|(1-\bar{\lap})(v\cdot\bar{n}_{\kappa})\|_{L^{2}(\Omega)}\|s^{i}\|_{L^{2}(\Omega)}\\&+\kappa\int_{0}^{T}\|((1-\bar{\lap})\p_{j}(v\cdot\bar{n}_{\kappa}))\|_{H^{-\frac{1}{2}}(\Omega)}\|\frac{\bar{J}^{\kappa}(\bar{A}^{\kappa})_{i}^{j}}{|\bar{J}^{\kappa}(\bar{A}^{\kappa})_{i}^{3}|}s^{i}\|_{H^{\frac{1}{2}}(\Omega)}\\&\lesssim C(N_{0},M_{0})+\kappa\int_{0}^{T}\|\p_{j}(v\cdot\bar{n}_{\kappa})\|_{H^{\frac{3}{2}}(\Omega)}^{2}\\&\lesssim C(N_{0},M_{0}).
\end{aligned}  
\end{align}
The right hand side of \eqref{equationforvtsX} can be directly controlled by $C(N_{0},M_{0})$. Hence, we can verify that $v_{t}\in L^{2}(0,T;X^{\prime})$ with
\begin{align}
    \int_{0}^{T}\|v_{t}\|_{X^{\prime}}^{2}\lesssim C(N_{0},M_{0}).
\end{align}
Invoking the Lagrangian multiplier 2, we can show that there exists a $q\in L^{2}(0,T;H^{\frac{1}{2}}(\Omega))$ with the estimate
\begin{align}
    \int_{0}^{T}\|q\|_{H^{\frac{1}{2}}(\Omega)}^{2}\lesssim C(N_{0},M_{0}),
\end{align}
such that $(v,q)$ is the weak solution to the approximating equation. Hence, by the uniqueness of the weak solution, we show that $\|q\|_{L^{2}_{t}H^{0.5}(\Omega)}\lesssim C(N_{0},M_{0})$. Now, we utilize the higher order difference quotient $D_{-h_{1}}D_{-h_{2}}D_{h_{2}}D_{h_{1}}v$ to test the equation. Using a similar method, we can show $v\in L^{2}(0,T;H^{\frac{7}{2}}(\Omega))$. The improvement of regularity of $q$ can be achieved by solving elliptic equation
\begin{align}
\begin{aligned}
     (\bar{A}^{\kappa})_{i}^{j}\p_{j}((\bar{A}^{\kappa})_{i}^{k}\p_{k}q)&= -\nabla_{\bar{A}^{\kappa}}\bar{v}^{\kappa}:\nabla_{\bar{A}^{\kappa}}v-\nabla_{\bar{A}^{\kappa}}\cdot(\phi^{\kappa}_{t}\nabla_{\bar{A}^{\kappa}} \phi^{\kappa})-\nabla_{\bar{A} ^{\kappa}}\cdot \nabla_{\bar{A} ^{\kappa}}|\nabla_{\bar{A} ^{\kappa}}\phi^{\kappa}|^{2}\quad\textit{in $\Omega$},\\
     q&=-\sigma(\frac{\sqrt{\bar{h}}}{\sqrt{\bar{h}_{\kappa}}} \lap_{\bar{h}}(\bar{\eta}) \cdot \bar{n}_{\kappa})+\kappa\frac{1}{\sqrt{\bar{h}_{\kappa}}}((1-\bar{\lap})(v\cdot \bar{n}^{\kappa})) \quad \textit{on $\Gamma_{1}$},\\
     (\bar{A}^{\kappa})_{i}^{3}(\bar{A}^{\kappa})_{j}^{i} \p_{i}q&=-\frac{1}{2}\nabla^{\kappa}_{N}|\nabla^{\kappa}\phi^{\kappa}|^{2} \quad \textit{on $\Gamma_{0}$}.
\end{aligned}
\end{align}
It's clear that
\begin{align}
    \int_{0}^{T}\|q\|_{H^{1.5}(\Omega)}^{2}\lesssim \PP(\|v\|_{H^{3.5}(\Omega)},\|\bar{v}^{\kappa}\|_{H^{4}(\Gamma_{1})},\|\nabla_{\bar{A}^{\kappa}}\phi^{\kappa}\|_{H^{2}(\Omega)})\lesssim C(N_{0},M_{0}),
\end{align}
and 
\begin{align}
    \int_{0}^{t}\|v_{t}\|_{H^{0.5}(\Omega)}^{2}\lesssim C(N_{0},M_{0}).
\end{align}
\rmk
Roughly speaking, the regularity of $q$ is approximately two derivatives less than that of $v$. Consequently, terms involving $\nabla_{\bar{A}^{\kappa}}\phi^{\kappa}$ are consistently more regular than $q$. The curl estimate primarily encompasses terms like $\nabla_{\bar{A}^{\kappa}}\phi_{t}^{\kappa}\nabla_{\bar{A}^{\kappa}}\phi^{\kappa}$. Estimates in the horizontal direction can be regarded as smooth. Hence, the framework established by Coutand \& Shkoller \cite{coutand2007well} is applicable in this context. The bootstrap argument stops when it encounters the regularity constraint $\lap_{\bar{h}}\bar{\eta} \cdot n^{\kappa}$.

After this bootstrap argument, we can obtain
\begin{align}
 \begin{aligned}   \|v\|_{L_{t}^{2}H^{13.5}}+\|q\|_{L_{t}^{2}H^{11.5}}+\|v_{t}\|_{L^{2}_{t}H^{10.5}}&\lesssim_{\kappa} C(N_{0},M_{0},\|u_{0}\|_{13.5}),
 \end{aligned}   
\end{align}
in particular,
\begin{align}
    \sup_{t}\|v\|_{12}\lesssim_{\kappa} C(N_{0},M_{0},\|u_{0}\|_{13.5}).
\end{align}
Now, we estimate the energy of the highest order term.
\begin{align}
    \begin{aligned}
        \curl v &=\curl u_{0}-\p_{l}v^{j}\int_{0}^{t}\ve^{.ij}\p_{t}(\bar{A}^{\kappa})_{i}^{j}+\int_{0}^{t}\ve^{.ij}\p_{l}v^{j}\p_{t}(\bar{A}^{\kappa})_{i}^{l}-\int_{0}^{t}\ve^{.ij}(\bar{A}^{\kappa})_{i}^{k}\p_{t}\p_{k}\phi^{\kappa}\nabla_{\bar{A}^{\kappa}}\phi^{\kappa}\\
        \dive v&=-\p_{j}v^{i}\int_{0}^{t}\p_{t}(\bar{A}^{\kappa})_{i}^{j}
    \end{aligned}
\end{align}
Thus, a similar estimate shows
\begin{align}\label{curldiveestimateforv125135}
\begin{aligned}
    \|\curl v\|_{12.5}&\lesssim \|u_{0}\|_{13.5}+\|v\|_{13.5}\int_{0}^{T}\|\p_{t}\bar{A}^{\kappa}\|_{L^{\infty}}+\|v\|_{L^{\infty}}\int_{0}^{T}\|\p_{t}\bar{A}^{\kappa}\|_{12.5}+\int_{0}^{T}\|\bar{A}^{\kappa}\p_{t}\p\phi^{\kappa}\nabla_{\bar{A}^{\kappa}}\phi^{\kappa}\|_{12.5}\\&\lesssim \|u_{0}\|_{13.5}+\|v\|_{13.5}T^{\frac{1}{2}}C(\kappa,N_{0},M_{0},\|u_{0}\|_{13.5}),\\\|\curl v\|_{11.5}&\lesssim\|v\|_{12.5}T^{\frac{1}{2}}C(\kappa,N_{0},M_{1},\|u_{0}\|_{12.5})+C(\kappa,N_{0},M_{0},\|u_{0}\|_{12.5}),\\
    \|\dive v\|_{12.5}&\lesssim\|v\|_{13.5} T^{\frac{1}{2}}C(\kappa,N_{0},M_{0},\|u_{0}\|_{13.5}),\\ \|\dive v\|_{11.5}&\lesssim\|v\|_{12.5} T^{\frac{1}{2}}C(\kappa,N_{0},M_{0},\|u_{0}\|_{12.5})+C(\kappa,N_{0},M_{0},\|u_{0}\|_{12.5}).
\end{aligned}   
\end{align}
We estimate $\sup_{t}\|v\|_{12.5}$ first. Taking $\pb^{10}$ on the velocity equation, testing $D_{-h}D_{h}\pb^{10}\p_{t}v$, we have
\begin{align}\label{dhdb10vtl2impltinftyv}
    \begin{aligned}
       & \int_{0}^{T}\|\bar{J}_{\kappa}D_{h}\pb^{10}v_{t}\|_{L^{2}(\Omega)}^{2}+\int_{0}^{T}\langle D_{h}\pb^{10}v\cdot  \bar{n}^{\kappa}D_{h}\pb^{10}v_{t}\cdot\bar{n}^{\kappa}\rangle_{H^{1}(\Gamma_{1})}\\&=\int_{0}^{T}\int_{\Omega}\pb^{10}q \bar{J}^{\kappa}\nabla_{\bar{A}^{\kappa}}\cdot D_{-h}D_{h}\pb^{10}v_{t}-\int_{0}^{T}\int_{\Omega}\pb\bar{J}^{\kappa}\pb^{9}v_{t}D_{-h}D_{h}\pb^{10}v_{t}-\int_{0}^{T}\int_{\Gamma_{1}}\pb^{9}(1-\lap)v\cdot\pb \bar{n}_{\kappa}D_{-h}D_{h}\pb^{10}v_{t} \\&+\int_{0}^{T}\int_{\Omega}\pb (\bar{J}^{\kappa}(\bar{A}^{\kappa})_{i}^{j})\p_{j}\pb^{9}q D_{-h}D_{h}\pb^{10}v_{t}+\int_{0}^{T}\int_{\Gamma_{1}}\sigma \pb^{10}(L_{\bar{h}}\bar{\eta}\cdot \bar{n}_{\kappa})D_{-h}D_{h}\pb^{10}v_{t}\cdot\bar{n}_{\kappa}+\mathcal{R}\\&=\int_{0}^{T}\int_{\Gamma_{1}}\sigma \pb^{10}(L_{\bar{h}}\bar{\eta}\cdot \bar{n}_{\kappa})D_{-h}D_{h}\pb^{10}v_{t}\cdot\bar{n}_{\kappa}+\mathcal{R},
    \end{aligned}
\end{align}
where $|\mathcal{R}|\lesssim C(\kappa,N_{0},M_{0},\|u_{0}\|_{12.5})$ is the lower order term. 
\rmk
We calculate the number of derivatives to control the commutators. For example, $2\times 10.5=21$ horizontal derivatives on $\p_{t}v$ has been controlled in $L_{t}^{2}L^{2}$. The commutators has at most $22-1-21$ derivatives on $\p_{t}v$. Thus, the we control the commutator term related to $\|D^{h}\pb^{10}\p_{t}v\|_{L_{t}^{2}L^{2}}$. Similarly, we can control the rest of commutators. As a result we would like to improve half derivative to obtain $\|\p_{t}v\|_{L^{2}H^{11}}$ and $\|v\|_{L_{t}^{\infty}H^{12}}$.

Moreover, we have
\begin{align}
\begin{aligned}
      &\int_{0}^{T}\int_{\Gamma_{1}}\sigma \pb^{10}(L_{\bar{h}}\bar{\eta}\cdot \bar{n}_{\kappa})D_{-h}D_{h}\pb^{10}v_{t}\cdot\bar{n}_{\kappa}\\=&\int_{\Gamma_{1}}\sigma \pb^{10}(L_{\bar{h}}\bar{\eta}\cdot \bar{n}_{\kappa})D_{-h}D_{h}\pb^{10}v\cdot\bar{n}_{\kappa}-\int_{0}^{T}\int_{\Gamma_{1}}\sigma \pb^{10}\p_{t}(L_{\bar{h}}\bar{\eta}\cdot \bar{n}_{\kappa})D_{-h}D_{h}\pb^{10}v\cdot\bar{n}_{\kappa}+\mathcal{R},
\end{aligned} 
\end{align}
where $|\mathcal{R}|\lesssim C(\kappa,N_{0},M_{0},\|u_{0}\|_{12.5})$ is the lower order term. Notice that
\begin{align}
   \sup_{t} \|L_{\bar{h}}\bar{\eta}\cdot \bar{n}_{\kappa}\|_{H^{11}(\Gamma_{1})}\lesssim\int_{0}^{t}\|\p_{t}(L_{\bar{h}}\bar{\eta}\cdot \bar{n}_{\kappa})\|_{H^{11}(\Gamma_{1})}\lesssim T^{\frac{1}{2}}C(\kappa,N_{0},M_{0},\|u_{0}\|_{12.5}).
\end{align}
\begin{align}
    |\int_{\Gamma_{1}}\sigma \pb^{10}(L_{\bar{h}}\bar{\eta}\cdot \bar{n}_{\kappa})D_{-h}D_{h}\pb^{10}v\cdot\bar{n}_{\kappa}|\lesssim \|D_{h}\pb^{10}v\cdot\bar{n}_{\kappa}\|_{L^{2}(\Gamma_{1})}T^{\frac{1}{2}}C(\kappa,N_{0},M_{0},\|u_{0}\|_{12.5}).
\end{align}
Consequently, we also have
\begin{align}\label{dhdb10vtl2impltinftyverror}
    |\int_{0}^{T}\int_{\Gamma_{1}}\sigma \pb^{10}(L_{\bar{h}}\bar{\eta}\cdot \bar{n}_{\kappa})D_{-h}D_{h}\pb^{10}v_{t}\cdot\bar{n}_{\kappa}|\lesssim C(\kappa,N_{0},M_{0},\|u_{0}\|_{12.5}).
\end{align}
Combining \eqref{dhdb10vtl2impltinftyv} and \eqref{dhdb10vtl2impltinftyverror}, we have
\begin{align}
   \int_{0}^{T}\|\bar{J}_{\kappa}D_{h}\pb^{10}v_{t}\|_{L^{2}(\Omega)}^{2}+\langle D_{h}\pb^{10}v\cdot  \bar{n}^{\kappa}D_{h}\pb^{10}v\cdot\bar{n}^{\kappa}\rangle_{H^{1}(\Gamma_{1})}\lesssim C(\kappa,N_{0},M_{0},\|u_{0}\|_{12.5}).
\end{align}
Letting $h\to 0$, we have
\begin{align}\label{estimateforpb11vt}
     \int_{0}^{T}\|\pb^{11}v_{t}\|_{L^{2}(\Omega)}^{2}+\sup_{t}\|\pb v\cdot  \bar{n}^{\kappa}\|_{H^{11}(\Gamma_{1})}\lesssim C(\kappa,N_{0},M_{0},\|u_{0}\|_{12.5}).
\end{align}
On the other hand,
\begin{align}
\begin{aligned}
      \curl v_{t}&=-\p_{l}v^{j}\int_{0}^{T}\ve^{.ij}(\p_{t}\bar{A}_{i}^{l})+\curl_{\bar{A}^{\kappa}}H,\\
      \dive v_{t}&=-\p_{i}v^{j}\int_{0}^{T}\p_{t}(\bar{A})_{j}^{i}.
\end{aligned}  
\end{align}
Therefore,
\begin{align}\label{estimateforcurldive11vt}
    \begin{aligned}
        \int_{0}^{T}\|\curl v_{t}\|_{H^{10}(\Omega)}^{2}+ \int_{0}^{T}\|\dive v_{t}\|_{H^{10}(\Omega)}^{2}&\lesssim C(\kappa,N_{0},M_{0},\|u_{0}\|_{12.5}).
    \end{aligned}
\end{align}
Combining \eqref{estimateforpb11vt}, \eqref{estimateforcurldive11vt} and normal trace inequality with elliptic estimate, we have
\begin{align}\label{estimateforvtH11}
    \int_{0}^{T}\|v_{t}\|_{H^{11}(\Omega)}^{2}\lesssim C(\kappa,N_{0},M_{0},\|u_{0}\|_{12.5}).
\end{align}
Combining \eqref{curldiveestimateforv125135}, \eqref{estimateforpb11vt}, and \eqref{estimateforvtH11}, we have
\begin{align}
    \sup_{t}\|v\|_{H^{12.5}(\Omega)}+\int_{0}^{T}\|v_{t}\|_{H^{11}(\Omega)}^{2}&\lesssim C(\kappa,N_{0},M_{0},\|u_{0}\|_{12.5}),
\end{align}
by choosing $T$ small enough.

Taking $\pb^{12}$ on the velocity equation, testing $\pb^{12}v$, the $L^{2}$ estimate suggests that 
\begin{align}   \label{estimateexpforcontrolv135lt2}
\begin{aligned}
    &\|\bar{J}^{\kappa}\pb^{12}v\|_{L^{2}}^{2}+\sigma\int_{0}^{T}\langle \pb^{12}v\cdot\bar{n}^{\kappa},\pb^{12}v\cdot\bar{n}^{\kappa}\rangle_{H^{1}(\Gamma_{1})}\\&=\|u_{0}\|_{12.5}+\int_{0}^{T}\int_{\Omega}\pb^{12}q \nabla_{\bar{A}^{\kappa}}\cdot\pb^{12}v+\int_{0}^{T}\int_{\Omega}\pb\bar{J}^{\kappa}\pb^{11}v_{t}\pb^{12}v+\int_{0}^{T}\int_{\Gamma_{1}}\pb^{11}(1-\lap)v\cdot \pb\bar{n}^{\kappa}\pb^{12}v\cdot \bar{n}^{\kappa}\\&+\int_{0}^{T}\int_{\Omega}\pb (\bar{A}^{\kappa})_{i}^{j}\p_{j}\pb^{11}q\pb^{12}v+\int_{0}^{T}\int_{\Gamma_{1}}\sigma \pb^{12}(L_{\bar{h}}\bar{\eta}\cdot \bar{n}_{\kappa})\pb^{12}v\cdot\bar{n}_{\kappa}+\mathcal{R},
\end{aligned}
\end{align}
 $|\mathcal{R}|\lesssim T^{\frac{1}{2}}C(\kappa,N_{0},M_{0},\|u_{0}\|_{13.5})$ is the lower order term. The $H^{-\frac{1}{2}}$-$H^{\frac{1}{2}}$ argument will lead us to
\begin{align}\label{estimateexpforcontrolv135lt21}
\begin{aligned}     &|\int_{0}^{T}\int_{\Omega}\pb\bar{J}^{\kappa}\pb^{11}v_{t}\pb^{12}v|+|+\int_{0}^{T}\int_{\Gamma_{1}}\pb^{11}(1-\lap)v\cdot \pb\bar{n}^{\kappa}\pb^{12}v\cdot \bar{n}^{\kappa}|+|\int_{0}^{T}\int_{\Omega}\pb (\bar{A}^{\kappa})_{i}^{j}\p_{j}\pb^{11}q\pb^{12}v|\\\lesssim & \int_{0}^{T}\|\pb\bar{J}^{\kappa}\|_{L^{\infty}(\Omega)}\|\pb^{11}v_{t}\|_{L^{2}(\Omega)}\|\pb^{12}v\|_{L^{2}(\Omega)}+\int_{0}^{T}\|\pb^{11}(1-\lap)v\cdot \pb\bar{n}^{\kappa}\|_{L^{2}(\Gamma_{1})}\|\pb^{12}v\cdot \bar{n}^{\kappa}\|_{L^{2}(\Gamma_{1})}\\+&\int_{0}^{T}\|\pb (\bar{A}^{\kappa})_{i}^{j}\|_{L^{\infty}(\Omega)}\|\p_{j}\pb^{11}q\|_{H^{-\frac{1}{2}}(\Omega)}\|\pb^{12}v\|_{H^{\frac{1}{2}}(\Omega)}\\\lesssim & T^{\frac{1}{2}}C(\kappa,N_{0},M_{0},\|u_{0}\|_{13.5}).
\end{aligned} 
\end{align}
Similarly, using $\dive_{\bar{A}^{\kappa}}v=0$ and $H^{-\frac{1}{2}}(\Omega)$-$H^{\frac{1}{2}}(\Omega)$ argument, we have
\begin{align}\label{estimateexpforcontrolv135lt22}
    |\int_{0}^{T}\int_{\Omega}\pb^{12}q \nabla_{\bar{A}^{\kappa}}\cdot\pb^{12}v|\lesssim T^{\frac{1}{2}}C(\kappa,N_{0},M_{0},\|u_{0}\|_{13.5}).
\end{align}
The last term can be estimated in a similar fashion
\begin{align}\label{estimateexpforcontrolv135lt23}
\begin{aligned}
     |\int_{0}^{T}\int_{\Gamma_{1}}\sigma \pb^{12}(L_{\bar{h}}\bar{\eta}\cdot \bar{n}_{\kappa})\pb^{12}v\cdot\bar{n}_{\kappa}|&\lesssim \sigma\int_{0}^{T}\|L_{\bar{h}}\bar{\eta}\cdot \bar{n}_{\kappa}\|_{H^{11}(\Gamma_{1})}\|\pb^{12}v\cdot\bar{n}_{\kappa}\|_{H^{1}(\Gamma_{1})}\\&\lesssim\sup_{t}\|L_{\bar{h}}\bar{\eta}\cdot \bar{n}_{\kappa}\|_{H^{11}(\Gamma_{1})}T^{\frac{1}{2}}C(\kappa,N_{0},M_{0},\|u_{0}\|_{13.5})\\&\lesssim T^{\frac{1}{2}}C(\kappa,N_{0},M_{0},\|u_{0}\|_{13.5})
\end{aligned}  
\end{align}
Hence, combining \eqref{curldiveestimateforv125135}, \eqref{estimateexpforcontrolv135lt2}, \eqref{estimateexpforcontrolv135lt21}, \eqref{estimateexpforcontrolv135lt22} and \eqref{estimateexpforcontrolv135lt23}, we have
\begin{align}
    \int_{0}^{T}\|v\|_{13.5}^{2}&\lesssim T^{\frac{1}{2}}C(\kappa,N_{0},M_{0},\|u_{0}\|_{13.5}).
\end{align}
Now, choosing $N_{0}$ large enough compared to $\|u_{0}\|_{18.5}$ and $T$ small enough such that $\|v\|_{L_{r}^{2}H^{13.5}}\leq N_{0}$. Thus, we define a functional $A:C_{T}\to C_{T}$. Similarly in \cite{coutand2007well}, it suffices to show $A$ is sequentially weakly lower semi-continuous. To this end, suppose $(\bar{v}^{n})_{n=0}^{\infty}$ is a sequence in $L^{2}(0,T;H^{15.5}(\Omega))$ such that
\begin{align}
    \bar{v}^{n} \rightharpoonup \bar{v} \quad \textit{in $L^{2}(0,T;H^{13.5}(\Omega))$}.
\end{align}
It's clear that $\bar{v}\in C_{T}$. By Rellich compactness theorem,
\begin{align}
     \bar{v}^{n} \to \bar{v} \quad \textit{in $L^{2}(0,T;H^{12.5}(\Omega))$}.
\end{align}
Now, since convolution operator, extension operator and Laplacian operator are linear, it's clear that 
\begin{align}
\begin{aligned}
     \|\bar{v}^{n}_{\kappa}-\bar{v}_{\kappa}\|_{L^{2}_{t}H^{12.5}(\Omega)}&\lesssim\|\rho_{\kappa}*E(\bar{v}^{n}-\bar{v})\|_{L_{t}^{2}H^{12.5}}+\|\Lambda_{\kappa}^{2}(\bar{v}^{n}-\bar{v})\|_{L_{t}^{2}H^{12}(\p\Omega)}\\&\lesssim \|\bar{v}^{n}-\bar{v}\|_{L_{t}^{2}H^{12.5}(\Omega)}\\&\to 0\quad \textit{when $n\to\infty$}.
\end{aligned}   
\end{align}
Hence,
\begin{align}
    \begin{aligned}
        \bar{v}^{n}_{\kappa}\to\bar{v}_{\kappa}\quad \textit{in $L^{2}(0,T;H^{12.5}(\Omega))$},\\
        \bar{\eta}^{n}_{\kappa}\to\bar{\eta}_{\kappa}\quad \textit{in $L^{2}(0,T;H^{12.5}(\Omega))$}.
    \end{aligned}
\end{align}
The readers may refer to the method of Theorem 4.5 in \cite{ladyzhenskaia1968linear} to show the convergence of $\phi^{n}$. Equivalently, we consider following system.
Let $(\bar{v}^{(1)},\phi^{(1)})$ and $(\bar{v}^{(2)},\phi^{(2)})$ be two pairs solutions to \eqref{lwplinphi1}. Still denoting $\phi^{(1)}-\phi^{(2)}$ by $\tPsi$, we have
\begin{align}
\begin{aligned}
     \p_{t}\tPsi-\p_{i}(\bar{g}^{(1)}_{ij}\p_{j}\tPsi)&=\p_{i}((\bar{g}_{ij}^{(1)}-\bar{g}_{ij}^{(2)})\p_{j}\phi^{(2)})+\tPsi|\nabla_{\bar{A}_{\kappa}^{(1)}}\phi^{(1)}|^{2}+\phi^{(2)}((\bar{A}_{\kappa}^{(1)})_{i}^{j}-(\bar{A}_{\kappa}^{(2)})_{i}^{j})\p_{j}\phi^{1}((\bar{A}_{\kappa}^{(1)})_{i}^{j}+(\bar{A}_{\kappa}^{(2)})_{i}^{j})\p_{j}\phi^{1}\\&+\phi^{(2)}\nabla_{\bar{A}_{\kappa}^{(2)}}\tPsi\nabla_{\bar{A}_{\kappa}^{(2)}}(\phi^{(1)}-\phi^{(2)})\quad \textit{in $\Omega$},\\
(g_{\kappa}^{(1)})^{3i}\p_{i}\tPsi&= ((g_{\kappa}^{(2)})^{3i}-(g_{\kappa}^{(1)})^{3i})\p_{i}\phi^{(2)}\quad \textit{on $\p\Omega$}.
\end{aligned}  
\end{align}
The standard parabolic estimate provides us the desired strong convergence. If we use $(\bar{v}^{n},\phi^{(n)})$, $(\bar{v},\phi)$ to represent corresponding solution pair, we have
\begin{align}
    \begin{aligned}
        \phi^{n}&\to \phi \quad \textit{in $L_{t}^{2}H^{4}(\Omega)$},\\  \phi^{n}&\to \phi \quad \textit{in $L_{t}^{\infty}H^{3}(\Omega)$},\\
        \p_{t}\phi^{n}&\to \p_{t}\phi \quad\textit{in $L_{t}^{2}H^{2}(\Omega)$},\\
        \p_{t}\phi^{n}&\to \p_{t}\phi \quad\textit{in $L_{t}^{\infty}H^{1}(\Omega)$}.
    \end{aligned}
\end{align}
Therefore, using the property of mollification, we have for some $s\geq 20$
\begin{align}
 \begin{aligned}   \phi^{\kappa(n)}&\to\phi^{\kappa} \quad \textit{in $L_{t}^{\infty}H^{s}(\Omega)$}\\ \p_{t}\phi^{\kappa(n)}&\to\p_{t}\phi^{\kappa} \quad \textit{in $L_{t}^{\infty}H^{s}(\Omega)$} 
 \end{aligned}   
\end{align}
Now, we write $v^{n}=A(\bar{v}^{n})$. Then Eberlein–Smulian theorem implies that
\begin{align}
    v^{n}\rightharpoonup v \quad \textit{in $L_{t}^{2}H^{13.5}(\Omega)$}.
\end{align}
Thus,
\begin{align}
      v^{n}\to v \quad \textit{in $L_{t}^{2}H^{12.5}(\Omega)$}.
\end{align}
We want to show $v=A(\bar{v})$. To this end, first we can verify that
\begin{align}    (\bar{A}^{\kappa})_{i}^{j}\p_{j}v^{i}=0,
\end{align}
which is a result from strong convergence. Moreover, suppose $p^{(n)}$ is solution to the following system:
\begin{align}
 \begin{aligned}  (\bar{A}^{\kappa(n)})_{i}^{j}\p_{j}(\bar{A}^{\kappa(n)})_{i}^{k}\p_{k}p^{(n)}  &= -(\bar{A}^{\kappa(n)})_{i}^{j}\p_{j}v^{(n)}_{k}(\bar{A}^{\kappa(n)})_{k}^{l}\p_{l}\bar{v}^{(n)}_{i}\\&- (\bar{A}^{\kappa(n)})_{i}^{j}\p_{j}(\phi_{t}^{(n)}\nabla_{\bar{A}^{\kappa(n)}}\phi^{(n)}+\frac{1}{2}\nabla_{\bar{A}^{\kappa(n)}}\cdot\nabla_{\bar{A}^{\kappa(n)}}|\nabla_{\bar{A}^{\kappa(n)}}\phi|^{2})\quad\textit{ in $\Omega$},
 \\p^{(n)}&=-\sigma\lap_{\bar{h}^{n}}\bar{\eta}^{(n)}\cdot \bar{n}_{\kappa}^{(n)}-\kappa\bar{\lap}_{0}(v^{(n)}\cdot \bar{n}_{\kappa}^{(n)})  \quad \textit{on $\Gamma_{1}$}.
 \end{aligned}   
\end{align}
Similarly, applying Eberlein–Smulian theorem and compactness theorem again, it yields that
\begin{align}
    \begin{aligned}
        p^{(n)}&\rightharpoonup p\quad \textit{in $L_{t}^{2}H^{11.5}(\Omega)$}\\ p^{(n)}&\to p\quad \textit{in $L_{t}^{2}H^{10.5}(\Omega)$}.
    \end{aligned}
\end{align}
The strong convergence immediately yields that
\begin{align}
 \begin{aligned}  (\bar{A}^{\kappa})_{i}^{j}\p_{j}(\bar{A}^{\kappa})_{i}^{k}\p_{k}p  &= -(\bar{A}^{\kappa})_{i}^{j}\p_{j}v_{k}(\bar{A}^{\kappa})_{k}^{l}\p_{l}\bar{v}_{i}- (\bar{A}^{\kappa})_{i}^{j}\p_{j}(\phi_{t}\nabla_{\bar{A}^{\kappa}}\phi+\frac{1}{2}\nabla_{\bar{A}^{\kappa}}\cdot\nabla_{\bar{A}^{\kappa}}|\nabla_{\bar{A}^{\kappa}}\phi|^{2})\textit{ in $\Omega$},
 \\p&=-\sigma\lap_{\bar{h}}\bar{\eta}\cdot \bar{n}_{\kappa}-\kappa\bar{\lap}_{0}(v\cdot \bar{n}_{\kappa})  \quad \textit{on $\Gamma_{1}$}.
 \end{aligned}   
\end{align}
Now, with use of the strong convergence of $\p_{t}v^{(n)}$ in $L_{t}^{2}H^{12.5}(\Omega)$ or by directly checking that $v$ is the weak solution to the \eqref{penal2}, we can deduce that $v=A(\bar{v})$. Hence, the Tychonoff fixed point theorem yields that there exists a solution to the $\kappa$-approximating problem.

Now, we consider the following steps: Given $\p_{t}v\in L^{2}(0,T;H^{11.5}(\Omega))$, we can establish $\phi\in W_{2}^{3,6}(Q_{T})$. Then, we can obtain the existence and regularity of $\p_{t}^{2}v\in L^{2}(0,T;H^{9.5}(\Omega))$. Finally, we can achieve that
\begin{align}
\begin{aligned}
      \p_{t}^{k}v&\in L_{t}^{2}H^{13.5-3k}(\Omega)\cap L_{t}^{\infty}H^{12.5-3k}(\Omega),\quad \textit{$k=0,1,2,3,4$},\\
       \p_{t}^{k}q&\in L_{t}^{2}H^{11.5-3k}(\Omega)\cap L_{t}^{\infty}H^{10.5-3k}(\Omega),\quad \textit{$k=0,1,2,3$},
       \\\phi&\in W_{2}^{5,10}(\Omega\times[0,T]). 
\end{aligned} 
\end{align}
\section{Energy estimate I: Estimate for $\phi$}
We will employ a similar strategy in \cite{coutand2007well}. And a dedicated cancellation regarding to $\phi$ will be applied to close the energy estimate. 
\defn
\begin{equation}
\begin{aligned}
    E^{\kappa}(t)&=\sup_{t}[\sum_{s=0}^{4}\|\p_{t}^{s}\eta\|_{H^{5.5-s}(\Omega)}^{2}+\|v_{tttt}\|_{L^{2}(\Omega)}^{2}+\sum_{s=0}^{3}\|\p_{t}^{s}\phi\|_{H^{4.5-s}(\Omega)}^{2}+\|\p_{t}^{4}\phi\|_{H^{1}(\Omega)}^{2}]\\&+\sum_{s=0}^{4}\int_{0}^{t}\|\p_{t}^{s}\phi\|_{H^{5.5-s}(\Omega)}^{2}+\int_{0}^{t}\|\p_{t}^{5}\phi\|_{L^{2}(\Omega)}^{2}+\|\sqrt{\kappa} \eta\|_{6.5}^{2}+\sum_{s=0}^{4}\int_{0}^{t}\|\sqrt{\kappa}\p_{t}^{s}v\|_{5.5-s}^{2}+1.
    \end{aligned}
\end{equation}

Since we only study the local well-posedness, we set $T_{\kappa}\leq 1$ first. 

\defn 
We define
\begin{align}    E_{\kappa}^{\prime}=\sup_{t}\sum_{s=0}^{3}\|\p_{t}^{s}\eta\|_{5.5-s}^{2}+\|\sqrt{\kappa}\eta\|_{6.5}^{2}+\sum_{k=0}^{3}\int_{0}^{t}\|\sqrt{\kappa}\p_{t}^{k}v\|_{5.5-k}^{2}
\end{align}
\rmk
In the logic sequence, we first finish the estimate of $E_{\kappa}^{\prime}$ and then $\|v_{tttt}\|_{L^{2}(\Omega)}$, $\|v_{ttt}\|_{1.5}$. However, the last two terms are more difficult to analyze. So, in our writing sequence, we deal with these two terms first.

In the following section, we will denote $\PP(E^{\kappa}(0),\|q(0)\|_{4.5}, \|q_{t}(0)\|_{3.5},\|q_{tt}(0)\|_{2.5},\|d(0)\|_{9}$ by $\PP_{0}$. We want to show that they can be bounded by $\PP(\|u(0)\|_{6},\|d(0)\|_{9},\sqrt{\kappa}\|u_{0}\|_{10.5})$. We consider the function satisfied by $q_{0}$.
\begin{align}
\begin{aligned}
     \lap q_{0}&=-\p_{j} v_{0}^{i}\p_{i}(u_{0}^{\kappa})^{j}-\p_{i}(\phi_{t}(0)\p_{i}\phi(0))-\lap(|\p\phi(0)|^{2}), \quad\textit{in $\Omega$},\\
    q(0)&= \kappa(1-\bar{\lap})v^{3}(0), \quad \textit{on $\Gamma_{1}$},\\
    \p_{N}q_{0}&=0 ,\quad\textit{on $\Gamma_{0}$}.
\end{aligned}  
\end{align}
We can see that
\begin{align}
\begin{aligned}
     \|q\|_{4.5}&\lesssim\|\p_{j} u_{0}^{i}\p_{i}(u_{0}^{\kappa})^{j}\|_{2.5}+\|\p_{i}(\phi_{t}(0)\p_{i}\phi(0))\|_{2.5}+\|\lap|\p\phi(0)|^{2}\|_{2.5}+\kappa\|(1-\bar{\lap})v^{3}(0)\|_{H^{4}(\Gamma_{1})}\\&\lesssim \PP(\|u(0)\|_{4.5},\|d(0)\|_{5.5})+\kappa\|v(0)\|_{6.5}.
\end{aligned}  
\end{align}
Similarly, we can derive estimates for $q_{t}$ and $q_{tt}$. Readers may refer to \cite{coutand2007well} for further discussions on removing extra regularity. However, since the initial data for $\phi$ is required to estimate $\p_{t}^{4}\phi$, it cannot be removed due to its parabolic structure. This requirement $\|u_{0}\|_{6}$ is higher than the velocity condition $\|u_{0}\|_{5.5}$ in \cite{coutand2007well}.

We now provide a list of functions which represent the lower-order terms in our analysis.
\lem
\begin{align}
     \sup_{t}\sum_{k=0}^{4}\|\p_{t}^{k}\ak\|_{4.5-k}\lesssim\PP(E^{\kappa}),
\end{align}
\begin{align}
    \sup_{t}(\sum_{k=0}^{3}\|\p_{t}^{k}\nn\phi\|_{H^{3.5-k}(\Omega)}+\|\p_{t}^{4}\nn\phi\|_{L^{2}(\Omega)})\lesssim\PP(E^{\kappa}).
\end{align}

\begin{align}
     \sup_{t}\sum_{k=0}^{4}\|\p_{t}^{k}g_{\kappa}^{ij}\|_{4.5-k}\lesssim\PP(E^{\kappa}),
\end{align}
\begin{align}
     \sup_{t}\sum_{k=0}^{3}\|\p_{t}^{k}\ak\|_{3.5-k}\lesssim\PP_{0}+T\PP(E^{\kappa}),
\end{align}
\begin{align}
  \sup_{t}\sum_{k=0}^{3}\|\p_{t}^{k}g_{\kappa}^{ij}\|_{3.5-k}\lesssim\PP_{0}+T\PP(E^{\kappa}),
\end{align}
\begin{align}
    \sup_{t}\sum_{k=0}^{2}(\|\p_{t}^{k}\nn v\|_{2.5-k}+\|\p_{t}^{k}\nn v^{\kappa}\|_{2.5-k})\lesssim\PP_{0}+T\PP(E^{\kappa}).
\end{align}
\begin{align}\label{linfinityestimateforphi}
    \sup_{t}\sum_{s=0}^{3}\|\p_{t}^{s}\phi\|_{H^{4.5-s}(\Omega)}^{2}\lesssim \PP_{0}+T\PP(E^{\kappa}).
\end{align}
\pf The estimate follows directly from the definition. We omit the details.

Similarly to the previous section, we further restrict $t$ to be sufficiently small so that $g_{\kappa}^{ij}$ remains close to $Id_{3\times 3}$ in our analysis.
\lem
Suppose there exists a constant $\tilde{N}_{0}$ independent on $\kappa$ such that $\sup_{t\in[0,T_{\kappa}]}E_{\kappa}\leq \tilde{N}_{0}$. Suppose $\ve>0$ is a small real number. Then there exists a time $T^{\prime}=\frac{\ve}{\PP(N_{0})}$(independent on $T_{\kappa}$) such that $\forall t\in[0. T^{\prime}]$, the following estimate is true:
\begin{align}\label{apprioriassumption}
\begin{aligned}
   \sup_{t\in[0,T^{\prime}]} \|(\ak)_{i}^{j}-\delta_{i}^{j}\|_{3.5}+\|(A)_{i}^{j}-\delta_{i}^{j}\|&\lesssim\ve,\\
    \sup_{t\in[0,T^{\prime}]} \|J_{\kappa}-1\|_{3.5}+\|J-1\|_{3.5}&\lesssim\ve,\\
    \sup_{t\in[0,T^{\prime}]} \|g_{\kappa}^{ij}-Id_{3\times3}\|_{3.5}+\|g^{ij}-Id_{3\times 3}\|_{3.5}&\lesssim\ve,\\
    \sup_{t\in[0,T^{\prime}]} \|h_{\kappa}^{ij}-Id_{2\times2}\|_{H^{3}(\p\Omega)}+\|h^{ij}-Id_{2\times 2}\|_{H^{3}(\p\Omega)}&\lesssim\ve. 
\end{aligned}  
\end{align}
\pf
We only take $A_{i}^{j}$ as an example. Notice that the constant in the norm is value of function taken at $t=0$. Therefore,
\begin{align}
    \|(A)_{i}^{j}-\delta_{i}^{j}\|_{3.5}\lesssim\int_{0}^{t}\|\p_{t}(A)_{i}^{j}\|_{3.5}\lesssim t\PP(E^{\kappa}).
\end{align}
If the time $t\leq\frac{\ve}{\PP(E^{\kappa})}$, the lemma follows immediately. The other term can be proved in a similar way.
\rmk
Note that if we have shown the existence of a $\kappa$-independent $t$ such that $E \lesssim \tilde{N}_{0}$, this a priori assumption can be immediately recovered. (i.e., once this is established, repeating the procedure will yield the same estimate. Therefore, we do not explicitly specify $\sup_{t}$ over a particular interval.)

\subsection{Estimate for $\phi$}
In the following section, we denote $\phi|\nn\phi|^2$ by $f$. Consider
\begin{align}\label{testphi}
    \p_{t}\phi-\p_{i}(g^{ij}_{\kappa}\p_{j}\phi)=f.
\end{align}
$f$ is always a lower order term:
\lem
\begin{align}\label{festi}
\sum_{k=0}^{4}\int_{0}^{t}\|\p_{t}^{k}f\|_{4.5-k}^{2}\lesssim T^{\frac{1}{2}}\PP(E^{\kappa}).
\end{align}
\pf
Our focus is solely on the term that includes $\p_{t}^{4}f$; the remaining terms are relatively easier to control.
\begin{align}
    \begin{aligned}        \p_{t}^{4}f&=\p_{t}^{4}\phi|\nn\phi|^{2}+4\p_{t}^{3}\phi(2\p_{t}\nn\phi\nn\phi)+6\p_{t}^{2}\phi(2\p_{t}^{2}\nn\phi\nn\phi+2\p_{t}\nn\phi\p_{t}\nn\phi)\\&+4\p_{t}\phi(2\p_{t}^{3}\nn\phi\nn\phi+6\p_{t}^{2}\nn\phi\p_{t}\nn\phi)+\phi(2\p_{t}^{4}\nn\phi\nn\phi+8\p_{t}^{3}\nn\phi\p_{t}\nn\phi+6\p_{t}^{2}\nn\phi\p_{t}^{2}\nn\phi),
    \end{aligned}
\end{align}
\begin{align}\label{estimatefor4tf}
    \begin{aligned}
        \|\p_{t}^{4}f\|_{H^{0.5}(\Omega)}&\lesssim\|\p_{t}^{4}\phi\|_{H^{0.5}(\Omega)}\|\nn\phi\|_{H^{2}(\Omega)}^{2}+\|\p_{t}^{3}\phi\|_{H^{0.5}(\Omega)}\|\p_{t}\nn\phi\|_{H^{2}(\Omega)}\|\nn\phi\|_{H^{2}(\Omega)}\\&+\|\p_{t}^{2}\phi\|_{H^{2}(\Omega)}\|\p_{t}^{2}\nn\phi\|_{H^{0.5}(\Omega)}\|\nn\phi\|_{H^{2}(\Omega)}\\&+\|\p_{t}^{2}\phi\|_{H^{2}(\Omega)}\|\p_{t}\nn\phi\|_{H^{0.5}(\Omega)}\|\p_{t}\nn\phi\|_{H^{2}(\Omega)}+\|\p_{t}\phi\|_{H^{2}(\Omega)}\|\p_{t}^{3}\nn\phi\|_{H^{0.5}(\Omega)}\|\nn\phi\|_{H^{2}(\Omega)}\\&+\|\p_{t}\phi\|_{H^{2}(\Omega)}\|\p_{t}^{2}\nn\phi\|_{H^{0.5}(\Omega)}\|\p_{t}\nn\phi\|_{H^{2}(\Omega)}+\|\phi\|_{H^{2}(\Omega)}\|\p_{t}^{4}\nn\phi\|_{H^{0.5}(\Omega)}\|\nn\phi\|_{H^{2}
        (\Omega)}\\&+\|\phi\|_{H^{2}(\Omega)}\|\p_{t}^{3}\nn\phi\|_{H^{0.5}(\Omega)}\|\p_{t}\nn\phi\|_{H^{2}(\Omega)}+\|\phi\|_{H^{2}(\Omega)}\|\p_{t}^{2}\nn\phi\|_{H^{1}(\Omega)}\|\p_{t}^{2}\nn\phi\|_{H^{1.5}(\Omega)}\\&\lesssim\PP(E^{\kappa})+\|\p_{t}^{4}\nn\phi\|_{H^{0.5}(\Omega)}(\PP_{0}+T\PP(E^{\kappa})),
    \end{aligned}
\end{align}
Thus,
\begin{align}
    \|\p_{t}^{4}f\|_{L_{t}^{2}H^{0.5}(\Omega)}\lesssim T^{\frac{1}{2}}\PP(E^{\kappa}).
\end{align}

We now begin our analysis with $\p_{t}^{4}\phi$.
\prop[Estimate for $\p_{t}^{4}\phi$]
\begin{align}\label{phiest1} \int_{0}^{t}\|\p_{t}^{4}\phi\|_{H^{1.5}(\Omega)}^{2}\lesssim\PP_{0}+T\PP(E^{\kappa}).
\end{align}
\begin{align}
    \|\sqrt{\kappa}\p_{t}^{4}\phi\|_{L^{2}_{t}H^{2}(\Omega)}\lesssim \PP_{0}+T\PP(E^{\kappa})+  \|\sqrt{\kappa}\p_{t}^{3} v\|_{L_{t}^{2}H^{2}(\Omega)}(\PP_{0}+T\PP(E^{\kappa})).
\end{align}
\pf
Taking $\p_{t}^{4}$ again on \eqref{testphi}, testing $\p_{t}^{4}\phi$, and integrating by part, it follows that
\begin{align}
    \begin{aligned}        \int_{\Omega}|\p_{t}^{4}\phi|^2+\int_{0}^{t}\int_{\Omega}\p_{t}^{4}(g^{ij}_{\kappa}\p_{j}\phi)\p_{t}^{4}\p_{i}\phi&=\|\p_{t}^{4}\phi(0)\|_{L^{2}(\Omega)}^{2}+\int_{0}^{t}\int_{\Omega}\p_{t}^{4}f\p_{t}^{4}\phi,
    \end{aligned}
\end{align}
where we use the boundary condition $\nabla_{N}\phi=0$. Now, we can express the second term as
\begin{align}\label{t4phiest1}
\begin{aligned}
    \int_{0}^{t}\int_{\Omega}\p_{t}^{4}(g^{ij}_{\kappa}\p_{j}\phi)\p_{t}^{4}\p_{i}\phi=\int_{0}^{t}\int_{\Omega}g^{ij}_{\kappa}\p_{t}^{4}\p_{j}\phi\p_{t}^{4}\p_{i}\phi +\int_{0}^{t}\int_{\Omega}[\p_{t}^{4},g_{\kappa}^{ij}]\p_{j}\phi\p_{t}^{4}\p_{i}\phi
\end{aligned} 
\end{align}
\begin{align}
\begin{aligned}
[\p_{t}^{4},g_{\kappa}^{ij}]\p_{j}\phi&=\p_{t}^{4}(g^{ij}_{\kappa})\p_{j}\phi+4\p_{t}^{3}(g_{\kappa}^{ij})\p_{t}\p_{j}
    \phi\\&+6\p_{t}^{2}(g_{\kappa}^{ij})\p_{t}^{2}\p_{j}\phi+4\p_{t}(g_{\kappa}^{ij})\p_{t}^{3}\p_{j}\phi
\end{aligned}    
\end{align}

With the help of assumption $g_{\kappa}^{ij}\geq \frac{1}{2}\delta^{ij}$, invoking Young's inequality, and choosing sufficiently small $\delta$, we can infer from \eqref{t4phiest1} that
\begin{align}   \int_{\Omega}|\p_{t}^{4}\phi|^2+\int_{0}^{t}\int_{\Omega}|\p_{t}^{4}\p\phi|^{2}&\lesssim\delta \int_{0}^{t}\int_{\Omega}|\p_{t}^{4}\p\phi|^{2}+\PP_{0}+C(\delta)(T\PP(E^{\kappa})).
\end{align}
After choosing small $\delta$, we can show that
\begin{align}\label{paraestphi1}
    \begin{aligned}
         \|\p_{t}^{4}\phi\|_{L^{2}(\Omega)}+\|\p_{t}^{4}\phi\|_{L^{2}(0,T;H^{1}(\Omega))}&\lesssim\|\p_{t}^{4}\phi(0)\|_{L^{2}(\Omega)}+\|[\p_{t}^{4},g_{\kappa}^{ij}]\p_{j}\phi\|_{L^{2}(0,T;L^{2}(\Omega))}+\|\p_{t}^{4}f\|_{L^{2}(0,T;L^{2}(\Omega))}
    \end{aligned}
\end{align}
Similarly, a standard parabolic estimate gives us
\begin{align}\label{paraestphi5}
 \|\p_{t}^{4}\phi\|_{L^{2}(0,t;H^{2}(\Omega))}\lesssim(1+\|g_{\kappa}^{ij}\|_{L^{\infty}_{t}H^{2}(\Omega)})(\|\p_{t}^{4}\phi(0)\|_{H^{1}(\Omega)}+\|[\p_{t}^{4},g_{\kappa}^{ij}]\p_{j}\phi\|_{L^{2}(0,t;H^{1}(\Omega))}+\|\p_{t}^{4}f\|_{L^{2}(0,t;L^{2}(\Omega))}).
\end{align}
We obtain the following estimate via interpolation between \eqref{paraestphi1} and \eqref{paraestphi5}:
\begin{align}\label{paraestiphi}
     \|\p_{t}^{4}\phi\|_{L^{2}(0,t;H^{1.5}(\Omega))}\lesssim(1+\|g_{\kappa}^{ij}\|_{L^{\infty}_{t}H^{2}(\Omega)})(\|\p_{t}^{4}\phi(0)\|_{H^{0.5}(\Omega)}+\|[\p_{t}^{4},g_{\kappa}^{ij}]\p_{j}\phi\|_{L^{2}(0,t;H^{0.5}(\Omega))}+\|\p_{t}^{4}f\|_{L^{2}(0,t;L^{2}(\Omega))})
\end{align}
Now, a direct computation shows that
\begin{align}
\begin{aligned}
    \|[\p_{t}^{4},g_{\kappa}^{ij}]\p_{j}\phi\|_{0.5}&\leq\|\p_{t}^{4}(g^{ij}_{\kappa})\|_{0.5}\|\p_{j}\phi\|_{2}+4\|\p_{t}^{3}(g_{\kappa}^{ij})\|_{0.5}\|\p_{t}\p_{j}
    \phi\|_{2}\\&+6\|\p_{t}^{2}(g_{\kappa}^{ij})\|_{1.5}\|\p_{t}^{2}\p_{j}\phi\|_{1}+4\|\p_{t}(g_{\kappa}^{ij})\|_{2}\|\p_{t}^{3}\p_{j}\phi\|_{0.5}\\&\lesssim\PP(E^{\kappa}), 
    \\
     \|[\p_{t}^{4},g_{\kappa}^{ij}]\p_{j}\phi\|_{1}&\leq\|\p_{t}^{4}(g^{ij}_{\kappa})\|_{1}\|\p_{j}\phi\|_{2}+4\|\p_{t}^{3}(g_{\kappa}^{ij})\|_{1}\|\p_{t}\p_{j}
    \phi\|_{2}\\&+6\|\p_{t}^{2}(g_{\kappa}^{ij})\|_{2}\|\p_{t}^{2}\p_{j}\phi\|_{1}+4\|\p_{t}(g_{\kappa}^{ij})\|_{2}\|\p_{t}^{3}\p_{j}\phi\|_{1}\\&\lesssim\PP(E^{\kappa})+(\PP_{0}+T\PP(E^{\kappa}))\|\p_{t}^{3}v\|_{H^{2}(\Omega)}, 
\end{aligned}
\end{align}
\begin{align}\label{paraestphi6}
     \|[\p_{t}^{4},g_{\kappa}^{ij}]\p_{j}\phi\|_{L^{2}(0,t;H^{0.5}(\Omega))}&\lesssim T^{\frac{1}{2}}\PP(E^{\kappa}),
\end{align}
Combining \eqref{festi}, \eqref{paraestiphi} and \eqref{paraestphi6}, we can show that
\begin{align} \int_{0}^{t}\|\p_{t}^{4}\phi\|_{H^{1.5}(\Omega)}^{2}\lesssim\PP_{0}+T\PP(E^{\kappa}).
\end{align}
On the other hand, we perform the $\kappa$-estimate for $\|\p_{t}^{4}\phi\|_{L^{2}_{t}H^{2}(\Omega)}$ as
\begin{align}
\begin{aligned}    \|\sqrt{\kappa}\p_{t}^{4}\phi\|_{L^{2}_{t}H^{2}(\Omega)}&\lesssim \PP_{0}+T\PP(E^{\kappa})+ \sqrt{\kappa} \|[\p_{t}^{4},g_{\kappa}^{ij}]\p_{j}\phi\|_{L^{2}(0,t;H^{1}(\Omega))}\\&\lesssim \PP_{0}+T\PP(E^{\kappa})+  \|\sqrt{\kappa}\p_{t}^{3} v\|_{L_{t}^{2}H^{2}(\Omega)}(\PP_{0}+T\PP(E^{\kappa})).
\end{aligned}
\end{align}

We now introduce the higher regularity of $\phi$, $\p_{t}\phi$, $\p_{t}^{2}\phi$, and $\p_{t}^{3}\phi$.
\lem[Quasi-linear elliptic estimate for $L_{t}^{2}$]Let $L_{A}$ be a divergence-form elliptic operator corresponding to a Neumann-boundary operator $B_{A}$ as
\begin{align}\label{quasiestitl2}
    L_{A}=\p_{i}(g^{ij}_{\kappa}\p_{j})\\
    B_{A}=(\ak)_{i}^{3}(\ak)_{i}^{j}\p_{j}
\end{align}
Then, we can consider the elliptic problem ($k=0,1,2,3$) as 
\begin{align}
    L_{A}(\p_{t}^{k}\phi)&=\p_{t}^{k}(\p_{t}\phi-f)+[L_{A},\p_{t}^{k}]\phi\\
    B_{A}(\p_{t}^{k}\phi)&=-[\p_{t}^{k},g_{\kappa}^{3i}]\p_{i}\phi.
\end{align}
The commutator terms read
\begin{align}
    [L_{A},\p_{t}]\phi&=-\p_{i}\p_{t}(g_{\kappa}^{ij})\p_{j}\phi-\p_{t}g_{\kappa}^{ij}\p_{i}\p_{j}\phi,\\
  - [\p_{t},g_{\kappa}^{3i}]\p_{i}\phi&=-\p_{t}g_{\kappa}^{3i}\p_{i}\phi
\end{align}
\begin{align}
\begin{aligned}
    [L_{A},\p_{t}^{2}]\phi&=-\p_{i}\p_{t}^{2}(g_{\kappa}^{ij})\p_{j}\phi-2\p_{i}\p_{t}(g_{\kappa}^{ij})\p_{t}\p_{j}\phi-\p_{t}^{2}g_{\kappa}^{ij}\p_{i}\p_{j}\phi-2\p_{t}g_{\kappa}^{ij}\p_{t}\p_{i}\p_{j}\phi\\
     - [\p_{t}^{2},g_{\kappa}^{3i}]\p_{i}\phi&=-\p_{t}^{2}g_{\kappa}^{3i}\p_{i}\phi-2\p_{t}g_{\kappa}^{3i}\p_{t}\p_{i}\phi
\end{aligned}    
\end{align}
\begin{align}
\begin{aligned}
     [L_{A},\p_{t}^{3}]\phi&=-\p_{i}\p_{t}^{3}g_{\kappa}^{ij}\p_{j}\phi-\p_{t}^{3}g_{\kappa}^{ij}\p_{i}\p_{j}\phi-3\p_{i}\p_{t}^{2}g_{\kappa}^{ij}\p_{t}\p_{j}\phi\\&-3\p_{t}^{2}g_{\kappa}^{ij}\p_{t}\p_{i}\p_{j}\phi-3\p_{i}\p_{t}g_{\kappa}^{ij}\p_{t}^{2}\p_{j}\phi-3\p_{t}g_{\kappa}^{ij}\p_{t}^{2}\p_{i}\p_{j}\phi\\\\
     - [\p_{t}^{3},g_{\kappa}^{3i}]\p_{i}\phi&=-\p_{t}^{3}g_{\kappa}^{3i}\p_{i}\phi-3\p_{t}^{2}g_{\kappa}^{3i}\p_{t}\p_{i}\phi-3\p_{t}g_{\kappa}^{3i}\p_{t}^{2}\p_{i}\phi
\end{aligned} 
\end{align}
The following elliptic estimate holds:
\begin{align}\label{l2estimateforphi}  
\begin{aligned}
\sum_{k=0}^{3}\|\p_{t}^{k}\phi\|_{5.5-k}&\lesssim\PP(E^{\kappa}),\\
     \sum_{k=0}^{3}\int_{0}^{t}\|\p_{t}^{k}\phi\|_{5.5-k}^{2}&\lesssim \PP_{0}+T\PP(E^{\kappa}).
\end{aligned}
\end{align}
\pf
From a direct computation, using Sobolev trace lemma, we can show that
\begin{align}
\begin{aligned}
     \|[L_{A},\p_{t}]\phi\|_{2.5}&\lesssim\|\p_{i}\p_{t}(g_{\kappa}^{ij})\|_{2.5}\|\p_{j}\phi\|_{2.5}+\|\p_{t}g_{\kappa}^{ij}\|_{2.5}\|\p_{i}\p_{j}\phi\|_{2.5}\lesssim\PP(E^{\kappa}),\\
     \|
  - [\p_{t},g_{\kappa}^{3i}]\p_{i}\phi\|_{H^{3}(\p\Omega)}&\lesssim\|\p_{t}g_{\kappa}^{3i}\p_{i}\phi\|_{3.5}\lesssim\PP(E^{\kappa})
\end{aligned}    
\end{align}
\begin{align}
    \begin{aligned}
         \|[L_{A},\p_{t}^{2}]\phi\|_{1.5}&=\|\p_{i}\p_{t}^{2}(g_{\kappa}^{ij})\|_{1.5}\|\|\p_{j}\phi\|_{2}+\|\p_{i}\p_{t}(g_{\kappa}^{ij})\|_{1.5}\|\p_{t}\p_{j}\phi\|_{2}+\|\p_{t}^{2}g_{\kappa}^{ij}\|_{1.5}\|\p_{i}\p_{j}\phi\|_{2}\\&+\|\p_{t}g_{\kappa}^{ij}\|_{2}\|\p_{t}\p_{i}\p_{j}\phi\|_{1.5}\\&\lesssim\PP(E^{\kappa}), \\
         \| - [\p_{t}^{2},g_{\kappa}^{3i}]\p_{i}\phi\|_{H^{2}(\p\Omega)}&\lesssim\|-\p_{t}^{2}g_{\kappa}^{3i}\p_{i}\phi-2\p_{t}g_{\kappa}^{3i}\p_{t}\p_{i}\phi\|_{2.5}\lesssim\PP(E^{\kappa}).
    \end{aligned}
\end{align}
\begin{align}\label{comestil2phit3}
    \begin{aligned}
         \|[L_{A},\p_{t}^{3}]\phi\|_{0.5}&\lesssim\|\p_{i}\p_{t}^{3}g_{\kappa}^{ij}\|_{0.5}\|\p_{j}\phi\|_{2}+\|\p_{t}^{3}g_{\kappa}^{ij}\|_{0.5}\|\p_{i}\p_{j}\phi\|_{2}+\|\p_{i}\p_{t}^{2}g_{\kappa}^{ij}\|_{0.5}\|\p_{t}\p_{j}\phi\|_{2}\\&\|\p_{t}^{2}g_{\kappa}^{ij}\|_{2}\|\p_{t}\p_{i}\p_{j}\phi\|_{0.5}+\|\p_{i}\p_{t}g_{\kappa}^{ij}\|_{2}\|\p_{t}^{2}\p_{j}\phi\|_{0.5}+\|\p_{t}g_{\kappa}^{ij}\|_{2}\|\p_{t}^{2}\p_{i}\p_{j}\phi\|_{0.5}\\&\lesssim\PP(E^{\kappa}).\\
          \|- [\p_{t}^{3},g_{\kappa}^{3i}]\p_{i}\phi\|_{H^{1}(\p\Omega)}&\lesssim\|-\p_{t}^{3}g_{\kappa}^{3i}\p_{i}\phi-3\p_{t}^{2}g_{\kappa}^{3i}\p_{t}\p_{i}\phi-3\p_{t}g_{\kappa}^{3i}\p_{t}^{2}\p_{i}\phi\|_{1.5}\\&\lesssim\|\p_{t}^{3}g_{\kappa}^{3i}\|_{1.5}\|\p_{i}\phi\|_{2}+\|\p_{t}^{2}g_{\kappa}^{3i}\|_{1.5}\|\p_{t}\p_{i}\phi\|_{2}+\|\p_{t}g_{\kappa}^{3i}\|_{2}\|\p_{t}^{2}\p_{i}\phi\|_{1.5}\\&\lesssim\PP(E^{\kappa}).
    \end{aligned}
\end{align}
We commence the estimate at $k=3$. Invoking the Lemma 2.1, it yields that
\begin{align}
\begin{aligned}
     \|\p_{t}^{3}\phi\|_{H^{2.5}(\Omega)}&\lesssim\PP(E^{\kappa})(\|\p_{t}^{4}\phi\|_{{H^{0.5}}(\Omega)}+\|\p_{t}^{3}(\phi|\nn\phi|^{2})\|_{H^{0.5}(\Omega)}+\|[\p_{t}^{3},L_{A}]\phi\|_{H^{0.5}(\Omega)})\\&+\|- [\p_{t}^{3},g_{\kappa}^{3i}]\p_{i}\phi\|_{H^{1}(\p\Omega)}+\|\p_{t}^{3}\phi\|_{L^{2}(\Omega)}.
\end{aligned}   
\end{align}
Owing to \eqref{festi} and \eqref{comestil2phit3}, we have
\begin{align}\label{H2esriphit3}
     \|\p_{t}^{3}\phi\|_{H^{2.5}(\Omega)}\lesssim\PP(E^{\kappa}).
\end{align}

\begin{align}
    \int_{0}^{t}\|\p_{t}^{3}\phi\|_{H^{2.5}(\Omega)}^{2}\lesssim \PP_{0}+T\PP(E^{\kappa}).
\end{align}
Inductively, applying an analogous argument, we deduce that
\begin{align*}
\begin{aligned}
\sum_{k=0}^{3}\|\p_{t}^{k}\phi\|_{5.5-k}&\lesssim\PP(E^{\kappa}),\\
     \sum_{k=0}^{3}\int_{0}^{t}\|\p_{t}^{k}\phi\|_{5.5-k}^{2}&\lesssim \PP_{0}+T\PP(E^{\kappa}).
\end{aligned}  
\end{align*}

To close the energy estimate, we also need two estimates of $\phi$ involving $\kappa$.
\lem[Estimate of $\|\p_{t}^{3}\phi\|_{L_{t}^{\infty}H^{2}}$ using $E_{\kappa}^{\prime}$]
\begin{align}   \sup_{t}\|\p_{t}^{3}\phi\|_{2}\lesssim\PP(E_{\kappa}^{\prime})+\PP_{0}+T\PP(E^{\kappa})
\end{align}
\pf
Similar to \eqref{comestil2phit3}, we have
\begin{align}
    \begin{aligned}
         \|[L_{A},\p_{t}^{3}]\phi\|_{L^{2}(\Omega)}&\lesssim\|\p_{i}\p_{t}^{3}g_{\kappa}^{ij}\|_{L^{2}(\Omega)}\|\p_{j}\phi\|_{2}+\|\p_{t}^{3}g_{\kappa}^{ij}\|_{L^{2}(\Omega)}\|\p_{i}\p_{j}\phi\|_{2}+\|\p_{i}\p_{t}^{2}g_{\kappa}^{ij}\|_{L^{2}(\Omega)}\|\p_{t}\p_{j}\phi\|_{2}\\&+\|\p_{t}^{2}g_{\kappa}^{ij}\|_{2}\|\p_{t}\p_{i}\p_{j}\phi\|_{L^{2}(\Omega)}+\|\p_{i}\p_{t}g_{\kappa}^{ij}\|_{2}\|\p_{t}^{2}\p_{j}\phi\|_{L^{2}(\Omega)}+\|\p_{t}g_{\kappa}^{ij}\|_{2}\|\p_{t}^{2}\p_{i}\p_{j}\phi\|_{L^{2}(\Omega)}\\&\lesssim\PP(E_{\kappa}^{\prime})+\PP_{0}+T\PP(E^{\kappa}).\\
          \|- [\p_{t}^{3},g_{\kappa}^{3i}]\p_{i}\phi\|_{H^{0.5}(\p\Omega)}&\lesssim\|-\p_{t}^{3}g_{\kappa}^{3i}\p_{i}\phi-3\p_{t}^{2}g_{\kappa}^{3i}\p_{t}\p_{i}\phi-3\p_{t}g_{\kappa}^{3i}\p_{t}^{2}\p_{i}\phi\|_{1}\\&\lesssim\|\p_{t}^{3}g_{\kappa}^{3i}\|_{1}\|\p_{i}\phi\|_{2}+\|\p_{t}^{2}g_{\kappa}^{3i}\|_{1}\|\p_{t}\p_{i}\phi\|_{2}+\|\p_{t}g_{\kappa}^{3i}\|_{2}\|\p_{t}^{2}\p_{i}\phi\|_{1}\\&\lesssim\PP(E_{\kappa}^{\prime})+\PP_{0}+T\PP(E^{\kappa}).
    \end{aligned}
\end{align}
The norm involving $\phi$ has been controlled, and the norm consist of $g_{\kappa}$ can be bounded by $E_{\kappa}^{\prime}$.
Now, since
\begin{align}
\begin{aligned}
      \|\p_{t}^{3}f\|_{L^{2}(\Omega)}&\lesssim\PP_{0}+\int_{0}^{t}\|\p_{t}^{4}f\|_{L^{2}(\Omega)}\lesssim\PP_{0}+T\PP(E^{\kappa})
\end{aligned} 
\end{align}
, and $\|\p_{t}^{4}\phi\|_{L^{2}(\Omega)}$ has been controlled, we can write
\begin{align}
\begin{aligned}
     \|\p_{t}^{3}\phi\|_{H^{2}(\Omega)}&\lesssim\PP(E_{\kappa}^{\prime})(\|\p_{t}^{4}\phi\|_{L^{2}(\Omega)}+\|\p_{t}^{3}(\phi|\nn\phi|^{2})\|_{L^{2}(\Omega)}+\|[\p_{t}^{3},L_{A}]\phi\|_{L^{2}(\Omega)}\\&+\|- [\p_{t}^{3},g_{\kappa}^{3i}]\p_{i}\phi\|_{H^{0.5}(\p\Omega)}+\|\p_{t}^{3}\phi\|_{L^{2}(\Omega)}\\&\lesssim\PP(E_{\kappa}^{\prime})+\PP_{0}+T\PP(E^{\kappa}).
\end{aligned}   
\end{align}

\lem[$\kappa$ estimate for $\|\phi\|_{L^{2}_{t}H^{6.5}(\Omega)}$]
\begin{align}   \int_{0}^{t}\|\sqrt{\kappa}\phi\|_{6.5}^{2}\lesssim \PP_{0}+T\PP(E^{\kappa}).
\end{align}
\pf
We consider following elliptic system:
\begin{align}
\begin{aligned}
    \lap\phi&=\p_{t}\phi+\p_{i}((\delta^{ij}-g_{\kappa}^{ij})\p_{j}\phi)-f \quad \textit{in $\Omega$},\\
     \p_{3}\phi&=(\delta_{j}^{3}-(A^{\kappa})_{i}^{3}(A^{\kappa})_{i}^{j})\p_{j}\phi \quad \textit{on $\p\Omega$}.
\end{aligned}  
\end{align}
The Sobolev inequality shows that
\begin{align}
\begin{aligned}
    \|\p_{i}((\delta^{ij}-g_{\kappa}^{ij})\p_{j}\phi)\|_{4.5}&\lesssim\|\delta^{ij}-g_{\kappa}^{ij}\|_{5.5}\|\p\phi\|_{L^{\infty}(\Omega)}+\|\delta^{ij}-g_{\kappa}^{ij}\|_{L^{\infty}(\Omega)}\|\p\phi\|_{5.5}\\&\lesssim \|\eta\|_{6.5}(\PP_{0}+T\PP(E^{\kappa}))+\delta \|\phi\|_{6.5}+\PP_{0}+T\PP(E^{\kappa}).
\end{aligned}    
\end{align}
Similarly, the trace inequality implies that
\begin{align}
    \begin{aligned}
       \| (\delta_{j}^{3}-(A^{\kappa})_{i}^{3}(A^{\kappa})_{i}^{j})\p_{j}\phi\|_{H^{5}(\p\Omega)}&\lesssim \|\delta_{j}^{3}-(A^{\kappa})_{i}^{3}(A^{\kappa})_{i}^{j}\|_{H^{5}(\p\Omega)}\|\p\phi\|_{L^{\infty}(\p\Omega)}\\&+\|\delta_{j}^{3}-(A^{\kappa})_{i}^{3}(A^{\kappa})_{i}^{j}\|_{L^{\infty}(\p\Omega)}\|\p\phi\|_{H^{5}(\p\Omega)}\\&\lesssim \|\eta\|_{6.5}(\PP_{0}+T\PP(E^{\kappa}))+\delta \|\phi\|_{6.5}+\PP_{0}+T\PP(E^{\kappa}).
    \end{aligned}
\end{align}
Therefore, the elliptic estimate shows that
\begin{align}
    \begin{aligned}        
    \|\phi\|_{6.5}&\lesssim \|\p_{t}\phi\|_{4.5}+ \|\p_{i}((\delta^{ij}-g_{\kappa}^{ij})\p_{j}\phi)\|_{4.5}+\|f\|_{4.5}+  \| (\delta_{j}^{3}-(A^{\kappa})_{i}^{3}(A^{\kappa})_{i}^{j})\p_{j}\phi\|_{H^{5}(\p\Omega)}\\
    \|\sqrt{\kappa}\phi\|_{L_{t}^{2}H^{6.5}}&\lesssim T^{\frac{1}{2}}\|\sqrt{\kappa}\eta\|_{6.5}(\PP_{0}+T\PP(E^{\kappa}))+\PP_{0}+T\PP(E^{\kappa})\\&\lesssim\PP_{0}+T\PP(E^{\kappa}).
    \end{aligned}
\end{align}

We provide the following corollary for $\int_{0}^{t}\|\sqrt{\kappa}\p_{t}^{4}\nn\phi\|_{H^{1}(\Omega)}^{2}$ which will be used later.
\cor
\begin{align}   \label{estimateforsqrtkappapt4nablaphih1}
\int_{0}^{t}\|\sqrt{\kappa}\p_{t}^{4}\nn\phi\|_{H^{1}(\Omega)}^{2}\lesssim \PP_{0}+T\PP(E^{\kappa})+(\PP_{0}+T\PP(E^{\kappa}))\|\sqrt{\kappa}\p_{t}^{3}v\|_{L_{t}^{2}H^{2}(\Omega)}^{2}.
\end{align}
\pf
We compute it directly as
\begin{align}
\begin{aligned}
     \|\sqrt{\kappa}\p_{t}^{4}\nn\phi\|_{H^{1}(\Omega)}^{2}&\lesssim \|\sqrt{\kappa}\p_{t}^{4}\ak\|_{H^{1}(\Omega)}\|\p\phi\|_{H^{2}(\Omega)}+\|\ak\|_{H^{2}(\Omega)}\|\sqrt{\kappa}\p_{t}^{4}\p\phi\|_{H^{1}(\Omega)}+\mathcal{R}\\&\lesssim \|\sqrt{\kappa}\p_{t}^{3}v\|_{H^{2}(\Omega)}(\PP_{0}+T\PP(E^{\kappa}))+\mathcal{R},
\end{aligned}
\end{align}
where $\int_{0}^{T}|\mathcal{R}|^{2}\lesssim \PP_{0}+T\PP(E^{\kappa})$ is the lower order term. Integrating on time, we immediately obtain the desired estimate
\begin{align}
    \int_{0}^{t}\|\sqrt{\kappa}\p_{t}^{4}\nn\phi\|_{H^{1}(\Omega)}^{2}\lesssim \PP_{0}+T\PP(E^{\kappa})+(\PP_{0}+T\PP(E^{\kappa}))\|\sqrt{\kappa}\p_{t}^{3}v\|_{L_{t}^{2}H^{2}(\Omega)}^{2}.
\end{align}

Next, we introduce the key lemma in our article.
\prop[key cancellation]
\begin{align}\label{keycancellation1}
  \int_{0}^{t}\int_{\Omega}J_{\kappa}\nn_{l}\p_{t}^{4}v_{\kappa}^{\mu}\nn_{l}\phi\p_{t}^{4}\nn_{\mu}\phi= -\int_{0}^{t}\int_{\Omega}J_{\kappa}|\p_{t}^{5}\phi|^{2}-\frac{1}{2}\int_{\Omega}J_{\kappa}g_{\kappa}^{ij}\p_{t}^{4}\p_{j}\phi\p_{t}^{4}\p_{i}\phi+\mathcal{R},
\end{align}
where $|\mathcal{R}|\lesssim\delta\|\p_{t}^{4}\phi\|_{1}^{2}+\delta\int_{0}^{t}\|\p_{t}^{5}\phi\|_{L^{2}(\Omega)}^{2}+\delta\|\p_{t}^{3}v\|_{1.5}^{2}+C(\delta)(\PP_{0}+T\PP(E^{\kappa}))$.
\pf
Taking $\p_{t}^{4}$ on \eqref{testphi}, testing $J_{\kappa}\p_{t}^{5}\phi$, and integrating by part in space, it implies that
\begin{align}\label{t5estphi1}   \underbrace{\int_{0}^{t}\int_{\Omega}J_{\kappa}|\p_{t}^{5}\phi|^2}_{K_{1}}\underbrace{+\int_{0}^{t}\int_{\Omega}J_{\kappa}\p_{t}^{4}(g^{ij}_{\kappa}\p_{j}\phi)\p_{t}^{5}\p_{i}\phi}_{K_{2}}&=\PP_{0}+\underbrace{\int_{0}^{t}\int_{\Omega}J_{\kappa}\p_{t}^{4} f\p_{t}^{5}\phi}_{K_{3}}\underbrace{-\int_{0}^{t}\int_{\Omega}\p_{i}J_{\kappa}\p_{t}^{4}(g^{ij}_{\kappa}\p_{j}\phi)\p_{t}^{5}\phi}_{K_{4}},
\end{align}
where we apply the Neumann boundary condition on $\phi$. After a direct calculation, the second term $K_{2}$ on the LHS reads
\begin{align}\label{t5estphi2}
 K_{2}=\underbrace{\int_{0}^{t}\int_{\Omega}J_{\kappa}(\ak)_{l}^{i}(\ak)_{l}^{j}\p_{j}\p_{t}^{4}\phi\p_{t}^{5}\p_{i}\phi}_{K_{21}}\underbrace{+\int_{0}^{t}\int_{\Omega}J_{\kappa}\p_{t}^{4}((\ak)_{l}^{i}(\ak)_{l}^{j})\p_{j}\phi\p_{t}^{5}\p_{i}\phi}_{K_{22}}+\mathcal{R},   
\end{align}
$\mathcal{R}$ can be written in the form of
\begin{align}
    \int_{0}^{t}\int_{\Omega}J_{\kappa}\p_{t}^{s_{1}}(g_{\kappa}^{ij})\p_{t}^{s_{2}}(\p_{j}\phi)\p_{t}^{5}(\p_{i}\phi),\quad  \textit{$s_{1}$, $s_{2}\geq 1$}, \quad s_{1}+s_{2}=4.
\end{align}
Integrating by part in time, 
\begin{align}
\begin{aligned}
   \mathcal{R}&\sim\sum_{s_{1}+s_{2}=2}   \int_{\Omega}J_{\kappa}\p_{t}^{s_{1}}\p_{t}(g_{\kappa}^{ij})\p_{t}^{s_{2}}\p_{t}(\p_{j}\phi)\p_{t}^{4}(\p_{i}\phi)|_{0}^{t}-   \int_{0}^{t}\int_{\Omega}J_{\kappa}\p_{t}(\p_{t}^{s_{1}}\p_{t}(g_{\kappa}^{ij})\p_{t}^{s_{2}}\p_{t}(\p_{j}\phi))\p_{t}^{4}(\p_{i}\phi)\\&-\sum_{s_{1}+s_{2}=2}   \int_{0}^{t}\int_{\Omega}\p_{t}J_{\kappa}\p_{t}^{s_{1}}\p_{t}(g_{\kappa}^{ij})\p_{t}^{s_{2}}\p_{t}(\p_{j}\phi)\p_{t}^{4}(\p_{i}\phi),
\end{aligned}     
\end{align}
Recall that using our notation, $\p\phi$, $g_{\kappa}$ have at least order 4 (correspondingly, have been controlled in order 3). We directly calculate that
\begin{align}\label{estimateforK2expanresidue}
    \begin{aligned}
        |\mathcal{R}|&\lesssim\PP_{0}+\sum_{s_{1}+s_{2}=2}   \|J_{\kappa}\|_{L^{\infty}(\Omega)}\|\p_{t}^{s_{1}}\p_{t}(g_{\kappa}^{ij})\p_{t}^{s_{2}}\p_{t}(\p_{j}\phi)\|_{L^{2}(\Omega)}\|\p_{t}^{4}\p_{i}\phi\|_{L^{2}(\Omega)}\\&+ \sum_{s_{1}+s_{2}=2}   \int_{0}^{t}\|J_{\kappa}\|_{L^{\infty}(\Omega)}\|\p_{t}(\p_{t}^{s_{1}}\p_{t}(g_{\kappa}^{ij})\p_{t}^{s_{2}}\p_{t}(\p_{j}\phi))\|_{L^{2}(\Omega)}\|\p_{t}^{4}\p_{i}\phi\|_{L^{2}(\Omega)}\\&+\sum_{s_{1}+s_{2}=2}   \int_{0}^{t}\|\pb J_{\kappa}\|_{L^{\infty}(\Omega)}\|\p_{t}(\p_{t}^{s_{1}}\p_{t}(g_{\kappa}^{ij})\p_{t}^{s_{2}}\p_{t}(\p_{j}\phi))\|_{L^{2}(\Omega)}\|\p_{t}^{4}(\p_{i}\phi)\|_{L^{2}(\Omega)}\\&\lesssim\delta\|\p_{t}^{4}\phi\|_{1}^{2}+C(\delta)(\PP_{0}+T\PP(E^{\kappa})).
    \end{aligned}
\end{align}
Similarly, we can write the major part in $K_{21}$ as
\begin{align}\label{expansionforK21}
\begin{aligned}
      K_{21}&=\int_{0}^{t}\int_{\Omega}\frac{1}{2}J_{\kappa}(\ak)_{l}^{i}(\ak)_{l}^{j}\p_{t}(\p_{j}\p_{t}^{4}\phi\p_{t}^{4}\p_{i}\phi)\\&=\int_{\Omega}\frac{1}{2}J_{\kappa}g_{\kappa}^{ij}\p_{j}\p_{t}^{4}\phi\p_{t}^{4}\p_{i}\phi\underbrace{-\int_{0}^{t}\int_{\Omega}\frac{1}{2}\p_{t}(J_{\kappa}(\ak)_{l}^{i}(\ak)_{l}^{j})\p_{j}\p_{t}^{4}\phi\p_{t}^{4}\p_{i}\phi}_{K_{21,1}}.
\end{aligned}  
\end{align}

One can compute that
\begin{align}\label{estimateforK211}
    |K_{21,1}|&\lesssim\int_{0}^{t}\frac{1}{2}\|\p_{t}(J_{\kappa}(\ak)_{l}^{i}(\ak)_{l}^{j})\|_{L^{\infty}(\Omega)}\|\p_{j}\p_{t}^{4}\phi\|_{L^{2}(\Omega)}\|\p_{t}^{4}\p_{i}\phi\|_{L^{2}(\Omega)}\\&\lesssim\PP_{0}+T\PP(E^{\kappa}).
\end{align}
Thus, combining \eqref{t5estphi2}, \eqref{estimateforK2expanresidue}, \eqref{expansionforK21}, \eqref{estimateforK211}, we can express $K_{2}$ as
\begin{align}\label{expansionforK2}
    K_{2}=\int_{\Omega}\frac{1}{2}J_{\kappa}g_{\kappa}^{ij}\p_{t}^{4}\p_{j}\phi\p_{t}^{4}\p_{i}\phi+\underbrace{\int_{0}^{t}\int_{\Omega}J_{\kappa}\p_{t}^{4}((\ak)_{l}^{i}(\ak)_{l}^{j})\p_{j}\phi\p_{t}^{5}\p_{i}\phi}_{K_{22}}+\mathcal{R},
\end{align}
where $|\mathcal{R}|\lesssim\delta\|\p_{t}^{4}\phi\|_{1}^{2}+C(\delta)(\PP_{0}+T\PP(E^{\kappa}))$.
Now, we estimate $K_{3}$ in \eqref{t5estphi1}, applying \eqref{estimatefor4tf}
\begin{align}\label{estimateforK3}
\begin{aligned} 
|K_{3}|&\lesssim\int_{0}^{t}\|J_{\kappa}\|_{L^{\infty}(\Omega)}\|\p_{t}^{4} f\|_{L^{2}(\Omega)}\|\p_{t}^{5}\phi\|_{L^{2}(\Omega)}\\&\lesssim\delta\int_{0}^{t}\|\p_{t}^{5}\phi\|_{L^{2}(\Omega)}^{2}+C(\delta)(\PP_{0}+T\PP(E^{\kappa})).
\end{aligned}    
\end{align}
The estimate for $K_{4}$ is similar.
Using our order argument,
\begin{align}\label{estimateforK4}
    \begin{aligned}        |K_{4}|&\lesssim\sum_{k=0}^{4}C_{4}^{k}\int_{0}^{t}\|\p_{i}J_{\kappa}\|_{L^{\infty}(\Omega)}\|\p_{t}^{k}g^{ij}_{\kappa}\p_{t}^{4-k}\p_{j}\phi\|_{L^{2}(\Omega)}\|\p_{t}^{5}\phi\|_{L^{2}(\Omega)}\\&\lesssim\delta\int_{0}^{t}\|\p_{t}^{5}\phi\|_{L^{2}(\Omega)}^{2}+C(\delta)(\PP_{0}+T\PP(E^{\kappa})).
    \end{aligned}
\end{align}
Combining \eqref{t5estphi1}, \eqref{expansionforK2}, \eqref{estimateforK3}, and \eqref{estimateforK4}, we can summarize that
\begin{align}\label{expansionforcancellation1}
    \underbrace{\int_{0}^{t}\int_{\Omega}J_{\kappa}\p_{t}^{4}((\ak)_{l}^{i}(\ak)_{l}^{j})\p_{j}\phi\p_{t}^{5}\p_{i}\phi}_{K_{22}}=-\int_{0}^{t}\int_{\Omega}J_{\kappa}|\p_{t}^{5}\phi|^2-\int_{\Omega}\frac{1}{2}J_{\kappa}g_{\kappa}^{ij}\p_{t}^{4}\p_{j}\phi\p_{t}^{4}\p_{i}\phi+\mathcal{R},
\end{align}
where $|\mathcal{R}|\lesssim\delta\|\p_{t}^{4}\phi\|_{1}^{2}+\delta\int_{0}^{t}\|\p_{t}^{5}\phi\|_{L^{2}(\Omega)}^{2}+C(\delta)(\PP_{0}+T\PP(E^{\kappa}))$. 
Now, we calculate the major part of $K_{22}$.
\begin{align}\label{expansionforK221123}
    \begin{aligned}
        K_{22}&=\underbrace{\int_{0}^{t}\int_{\Omega}J_{\kappa}\p_{t}^{4}(\ak)_{l}^{i}(\ak)_{l}^{j}\p_{j}\phi\p_{t}^{5}\p_{i}\phi}_{K_{22,1}}\underbrace{+\int_{0}^{t}\int_{\Omega}J_{\kappa}(\ak)_{l}^{i}\p_{t}^{4}(\ak)_{l}^{j}\p_{j}\phi\p_{t}^{5}\p_{i}\phi}_{K_{22,2}}\\&\underbrace{+\sum_{s=1}^{3}C_{4}^{s}\int_{0}^{t}\int_{\Omega}J_{\kappa}\p_{t}^{s}(\ak)_{l}^{i}\p_{t}^{4-s}(\ak)_{l}^{j}\p_{j}\phi\p_{t}^{5}\p_{i}\phi}_{K_{22,3}}
    \end{aligned}
\end{align}
We estimate $K_{22,3}$ first. Integrating by part in time, 
\begin{align}
    \begin{aligned}
        K_{22,3}&=\sum_{s=0}^{2}C_{4}^{s+1}\int_{\Omega}J_{\kappa}\p_{j}\phi\p_{t}^{s}\p_{t}(\ak)_{l}^{i}\p_{t}^{2-s}\p_{t}(\ak)_{l}^{j}\p_{t}^{4}\p_{i}\phi|_{0}^{t}\\&-\sum_{s=0}^{2}C_{4}^{s+1}\int_{0}^{t}\int_{\Omega}\p_{t}(J_{\kappa}\p_{j}\phi)\p_{t}^{s}
        \p_{t}(\ak)_{l}^{i}\p_{t}^{2-s}\p_{t}(\ak)_{l}^{j}\p_{t}^{4}\p_{i}\phi\\&-\sum_{s=0}^{2}C_{4}^{s+1}\int_{0}^{t}\int_{\Omega}J_{\kappa}\p_{j}\phi\p_{t}(\p_{t}^{s}\p_{t}(\ak)_{l}^{i}\p_{t}^{2-s}\p_{t}(\ak)_{l}^{j})\p_{t}^{4}\p_{i}\phi.
    \end{aligned}
\end{align}
As a consequence,
\begin{align}\label{estimateforK223}
    \begin{aligned}
        |K_{22,3}|&\lesssim\PP_{0}+\sum_{s=0}^{2}C_{4}^{s+1}\|J_{\kappa}\p_{j}\phi\|_{L^{\infty}(\Omega)}\|\p_{t}^{s}
        \p_{t}(\ak)_{l}^{i}\p_{t}^{2-s}\p_{t}(\ak)_{l}^{j}\|_{L^{2}(\Omega)}\|\p_{t}^{4}\p_{i}\phi\|_{L^{2}(\Omega)}\\&+\sum_{s=0}^{2}C_{4}^{s+1}\int_{0}^{t}\|\p_{t}(J_{\kappa}\p_{j}\phi)\|_{L^{\infty}(\Omega)}\|\p_{t}^{s}
        \p_{t}(\ak)_{l}^{i}\p_{t}^{2-s}\p_{t}(\ak)_{l}^{j}\|_{L^{2}(\Omega)}\|\p_{t}^{4}\p_{i}\phi\|_{L^{2}(\Omega)}\\&+\sum_{s=0}^{2}C_{4}^{s+1}\int_{0}^{t}\|J_{\kappa}\p_{j}\phi\|_{L^{\infty}(\Omega)}\|\p_{t}(\p_{t}^{s}\p_{t}(\ak)_{l}^{i}\p_{t}^{2-s}\p_{t}(\ak)_{l}^{j})\|_{L^{2}(\Omega)}\|\p_{t}^{4}\p_{i}\phi\|_{L^{2}(\Omega)}\\&\lesssim\delta\|\p_{t}^{4}\phi\|_{H^{1}(\Omega)}^{2}+C(\delta)(\PP_{0}+T\PP(E^{\kappa})).
    \end{aligned}
\end{align}
From a direct calculation, we calculate the leading order term of $K_{22,1}$ and $K_{22,2}$
\begin{align}\label{expansionforK221}
    \begin{aligned}
        K_{22,1}&=\underbrace{-\int_{0}^{t}\int_{\Omega}J_{\kappa}(\ak)_{l}^{k}\p_{k}\p_{t}^{3}v_{\kappa}^{\mu}(\ak)_{\mu}^{i}(\ak)_{l}^{j}\p_{j}\phi\p_{t}^{5}\p_{i}\phi}_{K_{22,11}}\\&\underbrace{-\sum_{s=1}^{3}C_{3}^{s}\int_{0}^{t}\int_{\Omega}J_{\kappa}\p_{t}^{s}((\ak)_{\mu}^{i}(\ak)_{l}^{k})\p_{t}^{3-s}\p_{k}v_{\kappa}^{\mu}(\ak)_{l}^{j}\p_{j}\phi\p_{t}^{5}\p_{i}\phi}_{K_{22,12}}
    \end{aligned}
\end{align}
Integrating by part in time, 
\begin{align}
    \begin{aligned}
        K_{22,12}&=-\sum_{s=0}^{2}C_{3}^{s+1}\int_{\Omega}J_{\kappa}(\ak)_{l}^{j}\p_{j}\phi\p_{t}^{s}\p_{t}((\ak)_{\mu}^{i}(\ak)_{l}^{k})\p_{t}^{2-s}\p_{k}v_{\kappa}^{\mu}\p_{t}^{4}\p_{i}\phi|_{0}^{t}\\&+\sum_{s=0}^{2}C_{3}^{s+1}\int_{0}^{t}\int_{\Omega}\p_{t}(J_{\kappa}(\ak)_{l}^{j}\p_{j}\phi)\p_{t}^{s}\p_{t}((\ak)_{\mu}^{i}(\ak)_{l}^{k})\p_{t}^{2-s}\p_{k}v_{\kappa}^{\mu}\p_{t}^{4}\p_{i}\phi\\&+\sum_{s=0}^{2}C_{3}^{s+1}\int_{0}^{t}\int_{\Omega}J_{\kappa}(\ak)_{l}^{j}\p_{j}\phi\p_{t}(\p_{t}^{s}\p_{t}((\ak)_{\mu}^{i}(\ak)_{l}^{k})\p_{t}^{2-s}\p_{k}v_{\kappa}^{\mu})\p_{t}^{4}\p_{i}\phi.
    \end{aligned}
\end{align}
Thus,
\begin{align}\label{estimateforK2212}
    \begin{aligned}        |K_{22,12}|&\lesssim\PP_{0}+\sum_{s=0}^{2}C_{3}^{s+1}\|J_{\kappa}(\ak)_{l}^{j}\p_{j}\phi\|_{L^{\infty}(\Omega)}\|\p_{t}^{s}\p_{t}((\ak)_{\mu}^{i}(\ak)_{l}^{k})\p_{t}^{2-s}\p_{k}v_{\kappa}^{\mu}\|_{L^{2}(\Omega)}\|\p_{t}^{4}\p_{i}\phi\|_{L^{2}(\Omega)}\\&\sum_{s=0}^{2}C_{3}^{s+1}\int_{0}^{t}\|\p_{t}(J_{\kappa}(\ak)_{l}^{j}\p_{j}\phi)\|_{L^{\infty}(\Omega)}\|\p_{t}^{s}\p_{t}((\ak)_{\mu}^{i}(\ak)_{l}^{k})\p_{t}^{2-s}\p_{k}v_{\kappa}^{\mu}\|_{L^{2}(\Omega)}\|\p_{t}^{4}\p_{i}\phi\|_{L^{2}(\Omega)}\\&\sum_{s=0}^{2}C_{3}^{s+1}\int_{0}^{t}\|J_{\kappa}(\ak)_{l}^{j}\p_{j}\phi\|_{L^{\infty}(\Omega)}\|\p_{t}(\p_{t}^{s}\p_{t}((\ak)_{\mu}^{i}(\ak)_{l}^{k})\p_{t}^{2-s}\p_{k}v_{\kappa}^{\mu})\|_{L^{2}(\Omega)}\|\p_{t}^{4}\p_{i}\phi\|_{L^{2}(\Omega)}\\&\lesssim\delta\|\p_{t}^{4}\phi\|_{1}^{2}+C(\delta)(\PP_{0}+T\PP(E^{\kappa})).   \end{aligned}
\end{align}
Similarly, 
\begin{align}\label{expansionforK222}
    \begin{aligned}
        K_{22,2}&=\underbrace{-\int_{0}^{t}\int_{\Omega}J_{\kappa}(\ak)_{l}^{i}(\ak)_{l}^{k}\p_{t}^{3}\p_{k} v_{\kappa}^{\mu}(\ak)_{\mu}^{j}\p_{j}\phi\p_{t}^{5}\p_{i}\phi}_{K_{22,21}}\\&\underbrace{-\sum_{s=1}^{3}C_{3}^{s}\int_{0}^{t}\int_{\Omega}J_{\kappa}(\ak)_{l}^{i}\p_{t}^{s}((\ak)_{l}^{k}(\ak)_{\mu}^{j})\p_{t}^{3-s}\p_{k}v_{\kappa}^{\mu}\p_{j}\phi\p_{t}^{5}\p_{i}\phi}_{K_{22,22}}.
    \end{aligned}
\end{align}
Integrating by part in time,
\begin{align}
    \begin{aligned}
        K_{22,22}&=-\sum_{s=0}^{2}C_{3}^{s+1}\int_{\Omega}J_{\kappa}(\ak)_{l}^{i}\p_{j}\phi\p_{t}^{s}\p_{t}((\ak)_{l}^{k}(\ak)_{\mu}^{j})\p_{t}^{2-s}\p_{k}v_{\kappa}^{\mu}\p_{t}^{4}\p_{i}\phi|_{0}^{t}\\&+\sum_{s=0}^{2}C_{3}^{s+1}\int_{0}^{t}\int_{\Omega}\p_{t}(J_{\kappa}(\ak)_{l}^{i}\p_{j}\phi)\p_{t}^{s}\p_{t}((\ak)_{l}^{k}(\ak)_{\mu}^{j})\p_{t}^{2-s}\p_{k}v_{\kappa}^{\mu}\p_{t}^{4}\p_{i}\phi\\&+\sum_{s=0}^{2}C_{3}^{s+1}\int_{0}^{t}\int_{\Omega}J_{\kappa}(\ak)_{l}^{i}\p_{j}\phi\p_{t}(\p_{t}^{s}\p_{t}((\ak)_{l}^{k}(\ak)_{\mu}^{j})\p_{t}^{2-s}\p_{k}v_{\kappa}^{\mu})\p_{t}^{4}\p_{i}\phi.
    \end{aligned}
\end{align}
Consequently,
\begin{align}\label{estimateforK2222}
    \begin{aligned}
        |K_{22,22}|&\lesssim\PP_{0}+\sum_{s=0}^{2}C_{3}^{s+1}\|J_{\kappa}(\ak)_{l}^{i}\p_{j}\phi\|_{L^{\infty}(\Omega)}\|\p_{t}^{s}\p_{t}((\ak)_{l}^{k}(\ak)_{\mu}^{j})\p_{t}^{2-s}\p_{k}v_{\kappa}^{\mu}\|_{L^{2}(\Omega)}\|\p_{t}^{4}\p_{i}\phi\|_{L^{2}(\Omega)}\\&+\sum_{s=0}^{2}C_{3}^{s+1}\int_{0}^{t}\|\p_{t}(J_{\kappa}(\ak)_{l}^{i}\p_{j}\phi)\|_{L^{\infty}(\Omega)}\|\p_{t}^{s}\p_{t}((\ak)_{l}^{k}(\ak)_{\mu}^{j})\p_{t}^{2-s}\p_{k}v_{\kappa}^{\mu}\|_{L^{2}(\Omega)}\|\p_{t}^{4}\p_{i}\phi\|_{L^{2}(\Omega)}\\&+\sum_{s=0}^{2}C_{3}^{s+1}\int_{0}^{t}\|J_{\kappa}(\ak)_{l}^{i}\p_{j}\phi\|_{L^{\infty}(\Omega)}\|\p_{t}(\p_{t}^{s}\p_{t}((\ak)_{l}^{k}(\ak)_{\mu}^{j})\p_{t}^{2-s}\p_{k}v_{\kappa}^{\mu})\|_{L^{2}(\Omega)}\|\p_{t}^{4}\p_{i}\phi\|_{L^{2}(\Omega)}\\&\lesssim\delta\|\p_{t}^{4}\phi\|_{1}^{2}+C(\delta)(\PP_{0}+T\PP(E^{\kappa})).
    \end{aligned}
\end{align}
Combining \eqref{expansionforcancellation1}, \eqref{estimateforK223}, \eqref{expansionforK221}, \eqref{estimateforK2212}, \eqref{expansionforK222}, \eqref{estimateforK2222}, we can conclude that
\begin{align}\label{expansionforcancellation2}
\begin{aligned}
    &\underbrace{-\int_{0}^{t}\int_{\Omega}J_{\kappa}\nn_{l}\p_{t}^{3}v_{\kappa}^{\mu}\nn_{l}\phi\nn_{\mu}\p_{t}^{5}\phi}_{K_{22,11}}\underbrace{-\int_{0}^{t}\int_{\Omega}J_{\kappa}\nn_{l}\p_{t}^{3} v_{\kappa}^{\mu}\nn_{\mu}\phi\nn_{l}\p_{t}^{5}\phi}_{K_{22,21}}\\&=-\int_{0}^{t}\int_{\Omega}J_{\kappa}|\p_{t}^{5}\phi|^2-\int_{\Omega}\frac{1}{2}J_{\kappa}g_{\kappa}^{ij}\p_{t}^{4}\p_{j}\phi\p_{t}^{4}\p_{i}\phi+\mathcal{R},
\end{aligned}   
\end{align}

where $|\mathcal{R}|\lesssim\delta\|\p_{t}^{4}\phi\|_{1}^{2}+\delta\int_{0}^{t}\|\p_{t}^{5}\phi\|_{L^{2}(\Omega)}^{2}+C(\delta)(\PP_{0}+T\PP(E^{\kappa}))$. Now, taking $\p_{t}^{4}$ on \eqref{testphi} again, testing $J_{\kappa}\p_{t}^{4}v_{\kappa}^{l}\nn_{l}\phi$, and integrating by part, we can see that
\begin{align} \label{t5estphi4}   \underbrace{\int_{0}^{t}\int_{\Omega}J_{\kappa}\p_{t}^{5}\phi\p_{t}^{4}v_{\kappa}^{l}\nn_{l}\phi}_{J_{1}}+\underbrace{\int_{0}^{t}\int_{\Omega}\p_{t}^{4}(g_{\kappa}^{ij}\p_{j}\phi)\p_{i}(J_{\kappa}\p_{t}^{4}v_{\kappa}^{l}\nn_{l}\phi)}_{J_{2}}=\PP_{0}\underbrace{+\int_{0}^{t}\int_{\Omega}\p_{t}^{4} f\p_{t}^{4}v_{\kappa}^{l}\nn_{l}\phi}_{J_{3}}
\end{align}
We can estimate $J_{1}$ and $J_{3}$ in a direct way:
\begin{align}\label{estimateforJ1}
    \begin{aligned}
        |J_{1}|&\lesssim\int_{0}^{t}\|J_{\kappa}\|_{L^{\infty(\Omega)}}\|\p_{t}^{5}\phi\|_{L^{2}(\Omega)}\|\p_{t}^{4}v_{\kappa}^{l}\|_{L^{2}(\Omega)}\|\nn_{l}\phi\|_{L^{\infty}(\Omega)}\\&\lesssim\delta\int_{0}^{t}\|\p_{t}^{5}\phi\|_{L^{2}(\Omega)}^{2}+C(\delta)(\PP_{0}+T\PP(E^{\kappa})),
    \end{aligned}
\end{align}
and
\begin{align}\label{estimateforJ3}
    \begin{aligned}
        |J_{3}|&\lesssim+\int_{0}^{t}\|\p_{t}^{4}f\|_{L^{2}(\Omega)}\|\p_{t}^{4}v_{\kappa}^{l}\|_{L^{2}(\Omega)}\|\nn_{l}\phi\|_{L^{\infty}(\Omega)}\\&\lesssim\PP_{0}+T\PP(E^{\kappa}).
    \end{aligned}
\end{align}
Now, the expansion of $J_{2}$ can be written as
\begin{align}\label{expansionforJ2}
\begin{aligned}  J_{2}&=\underbrace{\int_{0}^{t}\int_{\Omega}\p_{t}^{4}(g_{\kappa}^{ij}\p_{j}\phi)\p_{i}(J_{\kappa}\nn_{l}\phi)\p_{t}^{4}v_{\kappa}^{l}}_{J_{21}}\\&\underbrace{+\int_{0}^{t}\int_{\Omega}\p_{t}^{4}(\ak)_{k}^{i}(\ak)_{k}^{j}\p_{j}\phi J_{\kappa}\nn_{l}\phi\p_{i}\p_{t}^{4} v_{\kappa}^{l}}_{J_{22}}\\&\underbrace{+\int_{0}^{t}\int_{\Omega}\p_{t}^{4}((\ak)_{k}^{j}\p_{j}\phi)J_{\kappa}\nn_{l}\phi (\ak)_{k}^{i}\p_{i}\p_{t}^{4}v_{\kappa}^{l}}_{J_{23}}\\ &\underbrace{+\sum_{s=0}^{2}C_{4}^{s+1}\int_{0}^{t}\int_{\Omega}\p_{t}^{s}\p_{t}(\ak)_{k}^{i}\p_{t}^{2-s}\p_{t}((\ak)_{k}^{j}\p_{j}\phi) J_{\kappa}\nn_{l}\phi\p_{i}\p_{t}^{4} v_{\kappa}^{l}}_{J_{24}},
\end{aligned}  
\end{align}
Using the order of $g_{\kappa}^{ij}$ and $\p\phi$ again,
\begin{align}\label{estimateforJ21}
    \begin{aligned}
       |J_{21}|&\lesssim \int_{0}^{t}\|\p_{t}^{4}(g_{\kappa}^{ij}\p_{j}\phi)\|_{L^{2}(\Omega)}\|\p_{i}(J_{\kappa}\nn_{l}\phi)\|_{L^{\infty}(\Omega)}\|\p_{t}^{4}v_{\kappa}^{l}\|_{L^{2}(\Omega)}\\&\lesssim\PP_{0}+T\PP(E^{\kappa}).
    \end{aligned}
\end{align}
For $J_{24}$, integrating by part in time,
\begin{align}
    \begin{aligned}
        J_{24}&=\sum_{s=0}^{2}C_{4}^{s+1}\int_{\Omega}\p_{t}^{s}\p_{t}(\ak)_{k}^{i}\p_{t}^{2-s}\p_{t}((\ak)_{k}^{j}\p_{j}\phi) J_{\kappa}\nn_{l}\phi\p_{i}\p_{t}^{3} v_{\kappa}^{l}|_{0}^{t}\\&-\sum_{s=0}^{2}C_{4}^{s+1}\int_{0}^{t}\int_{\Omega}\p_{t}(J_{\kappa}\nn_{l}\phi)\p_{t}^{s}\p_{t}(\ak)_{k}^{i}\p_{t}^{2-s}\p_{t}((\ak)_{k}^{j}\p_{j}\phi) \p_{i}\p_{t}^{3} v_{\kappa}^{l}\\&-\sum_{s=0}^{2}C_{4}^{s+1}\int_{0}^{t}\int_{\Omega}J_{\kappa}\nn_{l}\phi\p_{t}(\p_{t}^{s}\p_{t}(\ak)_{k}^{i}\p_{t}^{2-s}\p_{t}((\ak)_{k}^{j}\p_{j}\phi)) \p_{i}\p_{t}^{3} v_{\kappa}^{l}.
    \end{aligned}
\end{align}
Thus,
\begin{align}\label{estimateforJ24}
    \begin{aligned}
        |J_{24}|&\lesssim\PP_{0}+\sum_{s=0}^{2}C_{4}^{s+1}\|\p_{t}^{s}\p_{t}(\ak)_{k}^{i}\p_{t}^{2-s}\p_{t}((\ak)_{k}^{j}\p_{j}\phi)\|_{L^{2}(\Omega)} \|J_{\kappa}\nn_{l}\phi\|_{L^{\infty}(\Omega)}\|\p_{i}\p_{t}^{3} v_{\kappa}^{l}\|_{L^{2}(\Omega)}\\&+\sum_{s=0}^{2}C_{4}^{s+1}\int_{0}^{t}\|\p_{t}^{s}\p_{t}(\ak)_{k}^{i}\p_{t}^{2-s}\p_{t}((\ak)_{k}^{j}\p_{j}\phi)\|_{L^{2}(\Omega)}\| \p_{t}(J_{\kappa}\nn_{l}\phi)\|_{L^{\infty}(\Omega)}\|\p_{i}\p_{t}^{3} v_{\kappa}^{l}\|_{L^{2}(\Omega)}\\&+\sum_{s=0}^{2}C_{4}^{s+1}\int_{0}^{t}\|\p_{t}(\p_{t}^{s}\p_{t}(\ak)_{k}^{i}\p_{t}^{2-s}\p_{t}((\ak)_{k}^{j}\p_{j}\phi))\|_{L^{2}(\Omega)}\| J_{\kappa}\nn_{l}\phi\|_{L^{\infty}(\Omega)}\|\p_{i}\p_{t}^{3}v_{\kappa}^{l}\|_{L^{2}(\Omega)}\\&\lesssim \delta\|\p_{t}^{3}v\|_{1}^{2}+C(\delta)(\PP_{0}+T\PP(E^{\kappa})).
    \end{aligned}
\end{align}
Now, we calculate the leading order of $J_{22}$:
\begin{align}\label{expansionforJ22}
    \begin{aligned}
        J_{22}&=\underbrace{-\int_{0}^{t}\int_{\Omega}(\ak)_{k}^{\alpha}\p_{\alpha}\p_{t}^{3}v_{\kappa}^{\beta}(\ak)_{\beta}^{i}(\ak)_{k}^{j}\p_{j}\phi J_{\kappa}\nn_{l}\phi\p_{i}\p_{t}^{4}v_{\kappa}^{l}}_{J_{22,1}}\\&\underbrace{-\sum_{s=1}^{3}C_{3}^{s}\int_{0}^{t}\int_{\Omega}\p_{t}^{s}((\ak)_{k}^{\alpha}(\ak)_{\beta}^{i})\p_{t}^{3-s}\p_{\alpha} v_{\kappa}^{\beta}(\ak)_{k}^{j}\p_{j}\phi J_{\kappa}\nn_{l}\phi\p_{i}\p_{t}^{4} v_{\kappa}^{l}}_{J_{22,2}},
    \end{aligned}
\end{align}
Also, estimate $J_{22,2}$ via integrating by part in time,
\begin{align}
    \begin{aligned}
       J_{22,2}&=-\sum_{s=0}^{2}C_{3}^{s+1}\int_{\Omega}(\ak)_{k}^{j}\p_{j}\phi J_{\kappa}\nn_{l}\phi\p_{t}^{s}\p_{t}((\ak)_{k}^{\alpha}(\ak)_{\beta}^{i})\p_{t}^{2-s}\p_{\alpha} v_{\kappa}^{\beta}\p_{i}\p_{t}^{3} v_{\kappa}^{l}|_{0}^{t}\\&+\sum_{s=0}^{2}C_{3}^{s+1}\int_{0}^{t}\int_{\Omega}\p_{t}((\ak)_{k}^{j}\p_{j}\phi J_{\kappa}\nn_{l}\phi)\p_{t}^{s}\p_{t}((\ak)_{k}^{\alpha}(\ak)_{\beta}^{i})\p_{t}^{2-s}\p_{\alpha} v_{\kappa}^{\beta}\p_{i}\p_{t}^{3} v_{\kappa}^{l}\\&+\sum_{s=0}^{2}C_{3}^{s+1}\int_{0}^{t}\int_{\Omega}(\ak)_{k}^{j}\p_{j}\phi J_{\kappa}\nn_{l}\phi\p_{t}(\p_{t}^{s}\p_{t}((\ak)_{k}^{\alpha}(\ak)_{\beta}^{i})\p_{t}^{2-s}\p_{\alpha} v_{\kappa}^{\beta})\p_{i}\p_{t}^{3} v_{\kappa}^{l}.
    \end{aligned}
\end{align}
Hence,
\begin{align}\label{estimateforJ222}
    \begin{aligned}
        |J_{22,2}|&\lesssim \sum_{s=0}^{2}C_{3}^{s+1}\|\p_{t}^{s}\p_{t}((\ak)_{k}^{\alpha}(\ak)_{\beta}^{i})\p_{t}^{2-s}\p_{\alpha} v_{\kappa}^{\beta}\|_{L^{2}(\Omega)}\|(\ak)_{k}^{j}\p_{j}\phi J_{\kappa}\nn_{l}\phi\|_{L^{\infty}(\Omega)}\|\p_{i}\p_{t}^{3} v_{\kappa}^{l}\|_{L^{2}(\Omega)}\\&+\sum_{s=0}^{2}C_{3}^{s+1}\int_{0}^{t}\|\p_{t}^{s}\p_{t}((\ak)_{k}^{\alpha}(\ak)_{\beta}^{i})\p_{t}^{2-s}\p_{\alpha} v_{\kappa}^{\beta}\|_{L^{2}(\Omega)}\|\p_{t}((\ak)_{k}^{j}\p_{j}\phi J_{\kappa}\nn_{l}\phi)\|_{L^{\infty}(\Omega)}\|\p_{i}\p_{t}^{3} v_{\kappa}^{l}\|_{L^{2}(\Omega)}\\&+\sum_{s=0}^{2}C_{3}^{s+1}\int_{0}^{t}\|\p_{t}(\p_{t}^{s}\p_{t}((\ak)_{k}^{\alpha}(\ak)_{\beta}^{i})\p_{t}^{2-s}\p_{\alpha} v_{\kappa}^{\beta})\|_{L^{2}(\Omega)}\|(\ak)_{k}^{j}\p_{j}\phi J_{\kappa}\nn_{l}\phi\|_{L^{\infty}(\Omega)}\|\p_{i}\p_{t}^{3} v_{\kappa}^{l}\|_{L^{2}(\Omega)}\\&\lesssim\delta\|\p_{t}^{3}v\|_{1}^{2}+C(\delta)(\PP_{0}+T\PP(E^{\kappa})).
    \end{aligned}
\end{align}

Now, combining \eqref{estimateforJ1}, \eqref{estimateforJ3}, \eqref{expansionforJ2},\eqref{estimateforJ21}, \eqref{estimateforJ24}, \eqref{expansionforJ22}, \eqref{estimateforJ222}. We can reformulate \eqref{t5estphi4} as

\begin{align}\label{t5estphi5}
    \begin{aligned}
      &\underbrace{-\int_{0}^{t}\int_{\Omega}\p_{t}^{4}((\ak)_{k}^{j}\p_{j}\phi)J_{\kappa}\nn_{l}\phi (\ak)_{k}^{i}\p_{i}\p_{t}^{4}v_{\kappa}^{l}}_{-J_{23}}\underbrace{+\int_{0}^{t}\int_{\Omega}(\ak)_{k}^{\alpha}\p_{\alpha}\p_{t}^{3}v_{\kappa}^{\beta}(\ak)_{\beta}^{i}(\ak)_{k}^{j}\p_{j}\phi J_{\kappa}\nn_{l}\phi\p_{i}\p_{t}^{4}v_{\kappa}^{l}}_{-J_{22,1}}=\mathcal{R},
    \end{aligned}
\end{align}
where $|\mathcal{R}|\lesssim\delta\|\p_{t}^{4}\phi\|_{1}^{2}+\delta\int_{0}^{t}\|\p_{t}^{5}\phi\|_{L^{2}(\Omega)}^{2}+\delta\|\p_{t}^{3}v\|_{1}^{2}+C(\delta)(\PP_{0}+T\PP(E^{\kappa}))$.
Now, adding \eqref{expansionforcancellation2} on \eqref{t5estphi5}, it yields that
\begin{align}\label{expansionforcancellation3}
    \begin{aligned}       &\underbrace{-\int_{0}^{t}\int_{\Omega}J_{\kappa}\nn_{l}\p_{t}^{3}v_{\kappa}^{\mu}\nn_{l}\phi\nn_{\mu}\p_{t}^{5}\phi}_{K_{22,11}}\underbrace{-\int_{0}^{t}\int_{\Omega}J_{\kappa}\nn_{l}\p_{t}^{3} v_{\kappa}^{\mu}\nn_{\mu}\phi\nn_{l}\p_{t}^{5}\phi}_{K_{22,21}}\\&-\int_{0}^{t}\int_{\Omega}\p_{t}^{4}((\ak)_{k}^{j}\p_{j}\phi)J_{\kappa}\nn_{l}\phi (\ak)_{k}^{i}\p_{i}\p_{t}^{4}v_{\kappa}^{l}+\int_{0}^{t}\int_{\Omega}(\ak)_{k}^{\alpha}\p_{\alpha}\p_{t}^{3}v_{\kappa}^{\beta}(\ak)_{\beta}^{i}(\ak)_{k}^{j}\p_{j}\phi J_{\kappa}\nn_{l}\phi\p_{i}\p_{t}^{4}v_{\kappa}^{l} \\= &-\int_{0}^{t}\int_{\Omega}J_{\kappa}|\p_{t}^{5}\phi|^{2}-\frac{1}{2}\int_{\Omega}J_{\kappa}g_{\kappa}^{ij}\p_{t}^{4}\p_{j}\phi\p_{t}^{4}\p_{i}\phi+\mathcal{R},
    \end{aligned}
\end{align}
where $|\mathcal{R}|\lesssim\delta\|\p_{t}^{4}\phi\|_{1}^{2}+\delta\int_{0}^{t}\|\p_{t}^{5}\phi\|_{L^{2}(\Omega)}^{2}+\delta\|\p_{t}^{3}v\|_{1}^{2}+C(\delta)(\PP_{0}+T\PP(E^{\kappa}))$.
Notice that the LHS of \eqref{expansionforcancellation3} can be written as
\begin{align}
    \begin{aligned}
      \underbrace{-\int_{0}^{t}\int_{\Omega}J_{\kappa}\nn_{l}\p_{t}^{3}v_{\kappa}^{\mu}\nn_{l}\phi\nn_{\mu}\p_{t}^{5}\phi+\int_{0}^{t}\int_{\Omega}\nn_{k}\p_{t}^{3}v_{\kappa}^{\beta}\nn_{k}\phi J_{\kappa}\nn_{l}\phi\nn_{\beta}\p_{t}^{4}v_{\kappa}^{l} }_{L_{1}}\\\underbrace{-\int_{0}^{t}\int_{\Omega}J_{\kappa}\nn_{l}\p_{t}^{3} v_{\kappa}^{\mu}\nn_{\mu}\phi\nn_{l}\p_{t}^{5}\phi-\int_{0}^{t}\int_{\Omega}\p_{t}^{4}(\nn_{k}\phi)J_{\kappa}\nn_{l}\phi \nn_{k}\p_{t}^{4}v_{\kappa}^{l}}_{L_{2}}
    \end{aligned}
\end{align}
Notice that
\begin{align}\label{expansionforL1}
    \begin{aligned}
        L_{1}&=\underbrace{-\int_{0}^{t}\int_{\Omega}J_{\kappa}\nn_{l}\p_{t}^{3}v_{\kappa}^{\mu}\nn_{l}\phi\p_{t}^{5}\nn_{\mu}\phi}_{L_{11}}\underbrace{+\sum_{s=1}^{4}C_{5}^{s}\int_{0}^{t}\int_{\Omega}J_{\kappa}\nn_{l}\p_{t}^{3}v_{\kappa}^{\mu}\nn_{l}\phi\p_{t}^{s}(\ak)_{\mu}^{i}\p_{t}^{5-s}\p_{i}\phi}_{L_{12}}\\&\underbrace{-\sum_{s=1}^{4}C_{4}^{s}\int_{0}^{t}\int_{\Omega}J_{\kappa}\nn_{l}\p_{t}^{3}v_{\kappa}^{\mu}\nn_{l}\phi\p_{t}^{s}((\ak)_{\mu}^{j}(\ak)_{k}^{i})\p_{t}^{4-s}\p_{j}v_{\kappa}^{k}\p_{i}\phi}_{L_{13}}.
    \end{aligned}
\end{align}
We can compute that
\begin{align}\label{estimateforL13}
    \begin{aligned}
        |L_{13}|&\lesssim\sum_{s=0}^{3}C_{4}^{s+1}\int_{0}^{t}\|J_{\kappa}\nn_{l}\phi\p_{i}\phi\|_{L^{\infty}(\Omega)}\|\nn_{l}\p_{t}^{3}v_{\kappa}^{\mu}\|_{L^{2}(\Omega)}\|\p_{t}^{s}\p_{t}((\ak)_{\mu}^{j}(\ak)_{k}^{i})\p_{t}^{3-s}\p_{j}v_{\kappa}^{k}\|_{L^{2}(\Omega)}\\&\lesssim\PP_{0}+T\PP(E^{\kappa}).
    \end{aligned}
\end{align}
Considering $L_{12}$
we can control it as
\begin{align}\label{estimateforL12}
    \begin{aligned}
        |L_{12}|&\lesssim\sum_{s=0}^{3}C_{5}^{s+1}\int_{0}^{t}\|J_{\kappa}\nn_{l}\phi\|_{L^{\infty}(\Omega)}\|\nn_{l}\p_{t}^{3}v_{\kappa}^{\mu}\|_{L^{2}(\Omega)}\|\p_{t}^{s}\p_{t}(\ak)_{\mu}^{i}\p_{t}^{3-s}\p_{t}\p_{i}\phi\|_{L^{2}(\Omega)}\\&\lesssim\PP_{0}+T\PP(E^{\kappa}).
    \end{aligned}
\end{align}
The remainder of $L_{1}$, $L_{11}$ contributes the key cancellation term in our article as
\begin{align}\label{expansionforL11}
    \begin{aligned}
        L_{11}&=\underbrace{-\int_{\Omega}J_{\kappa}\nn_{l}\p_{t}^{3}v_{\kappa}^{\mu}\nn_{l}\phi\p_{t}^{4}\nn_{\mu}\phi|_{0}^{t}}_{L_{11,1}}\underbrace{+\int_{0}^{t}\int_{\Omega}J_{\kappa}\nn_{l}\p_{t}^{4}v_{\kappa}^{\mu}\nn_{l}\phi\p_{t}^{4}\nn_{\mu}\phi}_{L_{11,2}}\\&\underbrace{+\int_{0}^{t}\int_{\Omega}\p_{t}(J_{\kappa}(\ak)_{l}^{i}\nn_{l}\phi)\p_{i}\p_{t}^{3}v_{\kappa}^{\mu}\p_{t}^{4}\nn_{\mu}\phi}_{L_{11,3}}.
    \end{aligned}
\end{align}
We have
\begin{align}\label{estimateforL113}
    \begin{aligned}
       |L_{11,3}|&\lesssim \int_{0}^{t}\|\p_{t}(J_{\kappa}(\ak)_{l}^{i}\nn_{l}\phi)\|_{L^{\infty}(\Omega)}\|\p_{i}\p_{t}^{3}v_{\kappa}^{\mu}\|_{L^{2}(\Omega)}\|\p_{t}^{4}\nn_{\mu}\phi\|_{L^{2}(\Omega)}\\&\lesssim\PP_{0}+T\PP(E^{\kappa}),
    \end{aligned}
\end{align}
and 
\begin{align}\label{expansionforL111}
\begin{aligned}
     L_{11,1}&=\underbrace{-\int_{\Omega}J_{\kappa}\nn_{l}\p_{t}^{3}v_{\kappa}^{\mu}\nn_{l}\phi\nn_{\mu}\p_{t}^{4}\phi|_{0}^{t}}_{L_{11,11}}\underbrace{-\sum_{s=1}^{3} C_{4}^{s}\int_{\Omega}J_{\kappa}\nn_{l}\p_{t}^{3}v_{\kappa}^{\mu}\nn_{l}\phi\p_{t}^{s}(\ak)_{\mu}^{i}\p_{t}^{4-s}\p_{i}\phi|_{0}^{t}}_{L_{11,12}}\\&\underbrace{-\int_{\Omega}J_{\kappa}\nn_{l}\p_{t}^{3}v_{\kappa}^{\mu}\nn_{l}\phi\p_{t}^{4}(\ak)_{\mu}^{i}\p_{i}\phi|_{0}^{t}}_{L_{11,13}}.
\end{aligned}  
\end{align}
From direct calculation,
\begin{align}
    \begin{aligned}        L_{11,13}&=\int_{\Omega}J_{\kappa}\nn_{l}\p_{t}^{3}v_{\kappa}^{\mu}\nn_{l}\phi(\ak)_{\mu}^{j}(\ak)_{k}^{i}\p_{t}^{3}\p_{j}v_{\kappa}^{k}\p_{i}\phi|_{0}^{t}\\&+\sum_{s=1}^{3}C_{3}^{s}\int_{\Omega}J_{\kappa}\nn_{l}\p_{t}^{3}v_{\kappa}^{\mu}\nn_{l}\phi\p_{t}^{s}((\ak)_{\mu}^{j}(\ak)_{k}^{i})\p_{t}^{3-s}\p_{j}v_{\kappa}^{k}\p_{i}\phi|_{0}^{t}.
    \end{aligned}
\end{align}
A standard $L^{2}$ estimate shows that
\begin{align}
    \begin{aligned}
        |L_{11,13}|&\lesssim\PP_{0}+\|J_{\kappa}(\ak)_{l}^{\nu}\nn_{l}\phi(\ak)_{\mu}^{j}(\ak)_{k}^{i}\p_{i}\phi\|_{L^{\infty}(\Omega)}\|\p_{\nu}\p_{t}^{3}v_{\kappa}^{\mu}\|_{L^{2}(\Omega)}\|\p_{t}^{3}\p_{j}v_{\kappa}^{k}\|_{L^{2}(\Omega)}\\&+\sum_{s=0}^{2}C_{3}^{s+1}\|J_{\kappa}(\ak)_{l}^{\nu}\nn_{l}\phi\p_{i}\phi\|_{L^{\infty}(\Omega)}\|\p_{\nu}\p_{t}^{3}v_{\kappa}^{\mu}\|_{L^{2}(\Omega)}\|\p_{t}^{s}\p_{t}((\ak)_{\mu}^{j}(\ak)_{k}^{i})\p_{t}^{2-s}\p_{j}v_{\kappa}^{k}\|_{L^{2}(\Omega)}\\&\lesssim\PP_{0}+(\PP_{0}+T\PP(E^{\kappa}))\|\p_{t}^{3}v\|_{1}^{2}+\|\p_{t}^{3}v\|_{1}(\PP_{0}+T\PP(E^{\kappa})).
    \end{aligned}
\end{align}
Invoking Sobolev interpolation theorem,
\begin{align}
    \begin{aligned}
        \|\p_{t}^{3}v\|_{1}\lesssim\|\p_{t}^{3}v\|_{L^{2}(\Omega)}^{\frac{1}{3}}\|\p_{t}^{3}v\|_{1.5}^{\frac{2}{3}},
    \end{aligned}
\end{align}
and Young's inequality, we have
\begin{align}\label{estimateforL1113}
    \begin{aligned}
        |L_{11,13}|&\lesssim\PP_{0}+(\PP_{0}+T\PP(E^{\kappa}))\|\p_{t}^{3}v\|_{L^{2}(\Omega)}^{\frac{2}{3}}\|\p_{t}^{3}v\|_{1.5}^{\frac{4}{3}}+\delta\|\p_{t}^{3}v\|_{1}^{2}+C(\delta)(\PP_{0}+T\PP(E^{\kappa}))\\&\lesssim\PP_{0}+C(\delta)(\PP_{0}+T\PP(E^{\kappa}))\|\p_{t}^{3}v\|_{L^{2}(\Omega)}^{2}+\delta\|\p_{t}^{3}v\|_{1.5}^{2}+\delta\|\p_{t}^{3}v\|_{1}^{2}+C(\delta)(\PP_{0}+T\PP(E^{\kappa}))\\&\lesssim\delta\|\p_{t}^{3}v\|_{1.5}^{2}+C(\delta)(\PP_{0}+T\PP(E^{\kappa})).
    \end{aligned}
\end{align}
Similarly,
\begin{align}\label{estimateforL1112}
    \begin{aligned}
        |L_{11,12}|&\lesssim\PP_{0}+\sum_{s=0}^{2} C_{4}^{s+1}\|J_{\kappa}(\ak)_{l}^{\nu}\nn_{l}\phi\|_{L^{\infty}(\Omega)}\|\p_{\nu}\p_{t}^{3}v_{\kappa}^{\mu}\|_{L^{2}(\Omega)}\|\p_{t}^{s}\p_{t}(\ak)_{\mu}^{i}\p_{t}^{3-s}\p_{i}\phi\|_{L^{2}(\Omega)}\\&\lesssim\delta\|\p_{t}^{3}v\|_{1}^{2}+C(\delta)(\PP_{0}+T\PP(E^{\kappa})),
    \end{aligned}
\end{align}

\begin{align}\label{estimateforL1111}
    \begin{aligned}
        |L_{11,11}|&\lesssim\PP_{0}+\int_{\Omega}\|J_{\kappa}(\ak)_{l}^{\nu}\nn_{l}\phi(\ak)_{\mu}^{j}\|_{L^{\infty}(\Omega)}\|\p_{\nu}\p_{t}^{3}v_{\kappa}^{\mu}\|_{L^{2}(\Omega)}\|\p_{j}\p_{t}^{4}\phi\|_{L^{2}(\Omega)}\\&\lesssim\PP_{0}+\delta\|\p_{t}^{4}\phi\|_{H^{1}(\Omega)}+C(\delta)\|\p_{t}^{3}v\|_{1}^{2}(\PP_{0}+T\PP(E^{\kappa}))\\&\lesssim\delta\|\p_{t}^{4}\phi\|_{H^{1}(\Omega)}+\delta\|\p_{t}^{3}v\|_{1.5}^{2}+C(\delta)(\PP_{0}+T\PP(E^{\kappa})).
    \end{aligned}
\end{align}
Now, combining \eqref{expansionforL111}, \eqref{estimateforL1113}, \eqref{estimateforL1112} and \eqref{estimateforL1111}, we have
\begin{align}\label{estimateforL111}
    |L_{11,1}|\lesssim\delta\|\p_{t}^{4}\phi\|_{H^{1}(\Omega)}+\delta\|\p_{t}^{3}v\|_{1.5}^{2}+C(\delta)(\PP_{0}+T\PP(E^{\kappa})).
\end{align}
Combining \eqref{expansionforL11}, \eqref{estimateforL113} and \eqref{estimateforL111}, we have
\begin{align}  \label{estimateforL11} L_{11}=\underbrace{\int_{0}^{t}\int_{\Omega}J_{\kappa}\nn_{l}\p_{t}^{4}v_{\kappa}^{\mu}\nn_{l}\phi\p_{t}^{4}\nn_{\mu}\phi}_{L_{11,2}}+\mathcal{R}
\end{align}
where $|\mathcal{R}|\lesssim\delta\|\p_{t}^{4}\phi\|_{H^{1}(\Omega)}+\delta\|\p_{t}^{3}v\|_{1.5}^{2}+C(\delta)(\PP_{0}+T\PP(E^{\kappa}))$.
Combining \eqref{expansionforL1}, \eqref{estimateforL13}, \eqref{estimateforL12} and \eqref{estimateforL11}, we have
\begin{align}\label{estimateforL1}
    \begin{aligned}  L_{1}=\int_{0}^{t}\int_{\Omega}J_{\kappa}\nn_{l}\p_{t}^{4}v_{\kappa}^{\mu}\nn_{l}\phi\p_{t}^{4}\nn_{\mu}\phi+\mathcal{R},
    \end{aligned}
\end{align}
where $|\mathcal{R}|\lesssim\delta\|\p_{t}^{4}\phi\|_{H^{1}(\Omega)}+\delta\|\p_{t}^{3}v\|_{1.5}^{2}+C(\delta)(\PP_{0}+T\PP(E^{\kappa}))$.
The $L_{2}$ can be controlled directly:
\begin{align}\label{expansionforL2}
    \begin{aligned}
        L_{2}&=-\int_{0}^{t}\int_{\Omega}J_{\kappa}\nn_{l}\p_{t}^{3} v_{\kappa}^{\mu}\nn_{\mu}\phi\nn_{l}\p_{t}^{5}\phi-\int_{\Omega}\p_{t}^{4}(\nn_{k}\phi)J_{\kappa}\nn_{l}\phi \nn_{k}\p_{t}^{3}v_{\kappa}^{l}|_{0}^{t}\\&+\int_{0}^{t}\int_{\Omega}J_{\kappa}\p_{t}^{5}(\nn_{k}\phi) \nn_{k}\p_{t}^{3} v_{\kappa}^{l}\nn_{l}\phi+\int_{0}^{t}\int_{\Omega}\p_{t}(J_{\kappa}(\ak)_{k}^{i}\nn_{l}\phi)\p_{t}^{4}(\nn_{k}\phi) \p_{i}\p_{t}^{3} v_{\kappa}^{l}\\&=-\int_{0}^{t}\int_{\Omega}J_{\kappa}\nn_{k}\p_{t}^{3} v_{\kappa}^{l}\nn_{l}\phi\nn_{k}\p_{t}^{5}\phi+\int_{0}^{t}\int_{\Omega}J_{\kappa}\nn_{k}\p_{t}^{5}\phi \nn_{k}\p_{t}^{3} v_{\kappa}^{l}\nn_{l}\phi\\&\underbrace{+\sum_{s=1}^{5}C_{5}^{s}\int_{0}^{t}\int_{\Omega}J_{\kappa}\p_{t}^{s}(\ak)_{k}^{i}\p_{t}^{5-s}\p_{i}\phi \nn_{k}\p_{t}^{3} v_{\kappa}^{l}\nn_{l}\phi}_{L_{21}}\\&\underbrace{\int_{0}^{t}\int_{\Omega}\p_{t}(J_{\kappa}(\ak)_{k}^{i}\nn_{l}\phi)\p_{t}^{4}(\nn_{k}\phi) \p_{i}\p_{t}^{3} v_{\kappa}^{l}}_{L_{22}}\underbrace{-\int_{\Omega}\p_{t}^{4}(\nn_{k}\phi)J_{\kappa}\nn_{l}\phi \nn_{k}\p_{t}^{3}v_{\kappa}^{l}|_{0}^{t}}_{L_{23}}
    \end{aligned}
\end{align}
It suffices to control $L_{21}$-$L_{23}$ term by term.
\begin{align}\label{expansionforL21}
    \begin{aligned}
       L_{21}&=\underbrace{\int_{0}^{t}\int_{\Omega}J_{\kappa}\p_{t}^{5}(\ak)_{k}^{i}\p_{i}\phi \nn_{k}\p_{t}^{3} v_{\kappa}^{l}\nn_{l}\phi}_{L_{21,1}}\\& \underbrace{+\sum_{s=0}^{3}C_{5}^{s+1}\int_{0}^{t}\int_{\Omega}J_{\kappa}\p_{t}^{s}\p_{t}(\ak)_{k}^{i}\p_{t}^{4-s}\p_{i}\phi \nn_{k}\p_{t}^{3} v_{\kappa}^{l}\nn_{l}\phi}_{L_{21,2}},
    \end{aligned}
\end{align}
From a direct calculation 
\begin{align}\label{estimateforL212}
    \begin{aligned}        |L_{21,2}|&\lesssim\sum_{s=0}^{3}C_{5}^{s+1}\int_{0}^{t}\|J_{\kappa}\nn_{l}\phi\|_{L^{\infty}(\Omega)}\|\p_{t}^{s}\p_{t}(\ak)_{k}^{i}\p_{t}^{3-s}\p_{t}\p_{i}\phi\|_{L^{2}(\Omega)}\| \nn_{k}\p_{t}^{3} v_{\kappa}^{l}\|_{L^{2}(\Omega)}\\&\lesssim\PP_{0}+T\PP(E^{\kappa}).
    \end{aligned}
\end{align}
We expand $L_{21,1}$ as
\begin{align}\label{expansionforL211}
    \begin{aligned}
        L_{21,1}&=\underbrace{-\int_{0}^{t}\int_{\Omega}J_{\kappa}(\ak)_{k}^{\mu}\p_{t}^{4}\p_{\mu}v_{\kappa}^{j}(\ak)_{j}^{i}\p_{i}\phi \nn_{k}\p_{t}^{3} v_{\kappa}^{l}\nn_{l}\phi}_{L_{21,11}}\\&\underbrace{-\sum_{s=1}^{4}C_{4}^{s}\int_{0}^{t}\int_{\Omega}J_{\kappa}\p_{t}^{s}((\ak)_{k}^{\mu}(\ak)_{j}^{i})\p_{t}^{4-s}\p_{\mu}v_{\kappa}^{j}\p_{i}\phi \nn_{k}\p_{t}^{3} v_{\kappa}^{l}\nn_{l}\phi}_{L_{21,12}}.
    \end{aligned}
\end{align}
From a direct computation,
\begin{align}\label{estimateforL2112}
    \begin{aligned}
        |L_{21,12}|&\lesssim\sum_{s=0}^{3}C_{4}^{s+1}\int_{0}^{t}\|J_{\kappa}\p_{i}\phi\nn_{l}\phi\|_{L^{\infty}(\Omega)}\|\p_{t}^{s}\p_{t}((\ak)_{k}^{\mu}(\ak)_{j}^{i})\p_{t}^{3-s}\p_{\mu}v_{\kappa}^{j}\|_{L^{2}(\Omega)} \|\nn_{k}\p_{t}^{3} v_{\kappa}^{l}\|_{L^{2}(\Omega)}\\&\lesssim\PP_{0}+T\PP(E^{\kappa}).
    \end{aligned}
\end{align}
For $L_{21,11}$, integrating by part in time,
we have
\begin{align}
    \begin{aligned}
        L_{21,11}&=-\frac{1}{2}\int_{0}^{t}\int_{\Omega}J_{\kappa}(\ak)_{k}^{\mu}(\ak)_{k}^{\nu}\p_{t}(\p_{t}^{3}\p_{\mu}v_{\kappa}^{j} \p_{\nu}\p_{t}^{3} v_{\kappa}^{l})\nn_{j}\phi\nn_{l}\phi\\&=-\frac{1}{2}\int_{\Omega}J_{\kappa}(\ak)_{k}^{\mu}(\ak)_{k}^{\nu}\p_{t}^{3}\p_{\mu}v_{\kappa}^{j} \p_{\nu}\p_{t}^{3} v_{\kappa}^{l}\nn_{j}\phi\nn_{l}\phi|_{0}^{t}+\frac{1}{2}\int_{0}^{t}\p_{t}(J_{\kappa}(\ak)_{k}^{\mu}(\ak)_{k}^{\nu}\nn_{j}\phi\nn_{l}\phi)\p_{t}^{3}\p_{\mu}v_{\kappa}^{j} \p_{\nu}\p_{t}^{3} v_{\kappa}^{l}
    \end{aligned}
\end{align}
We can deduce that
\begin{align}\label{estimateforL2111}
    \begin{aligned}
       |L_{21,11}|&\lesssim \PP_{0}+ \frac{1}{2}\|J_{\kappa}(\ak)_{k}^{\mu}(\ak)_{k}^{\nu}\nn_{j}\phi\nn_{l}\phi\|_{L^{\infty}(\Omega)}\|\p_{t}^{3}\p_{\mu}v_{\kappa}^{j}\|_{L^{2}(\Omega)}\| \p_{\nu}\p_{t}^{3} v_{\kappa}^{l}\|_{L^{2}(\Omega)}\\&+\frac{1}{2}\int_{0}^{t}\|\p_{t}(J_{\kappa}(\ak)_{k}^{\mu}(\ak)_{k}^{\nu}\nn_{j}\phi\nn_{l}\phi)\|_{L^{\infty}(\Omega)}\|\p_{t}^{3}\p_{\mu}v_{\kappa}^{j}\|_{L^{2}(\Omega)}\| \p_{\nu}\p_{t}^{3} v_{\kappa}^{l}\|_{L^{2}(\Omega)}\\&\lesssim(\PP_{0}+T\PP(E^{\kappa}))\times\|\p_{t}^{3}v\|_{H^{1}(\Omega)}^{2}+\PP_{0}+T\PP(E^{\kappa})\\&\lesssim\delta\|\p_{t}^{3}v\|_{1.5}^{2}+C(\delta)(\PP_{0}+T\PP(E^{\kappa})).
    \end{aligned}
\end{align}

Combining \eqref{expansionforL211}, \eqref{estimateforL2112} and \eqref{estimateforL2111}, we have
\begin{align}\label{estimateforL211}
    |L_{21,1}|\lesssim\delta\|\p_{t}^{3}v\|_{1.5}^{2}+C(\delta)(\PP_{0}+T\PP(E^{\kappa})).
\end{align}
Combining \eqref{expansionforL21}, \eqref{estimateforL212}, \eqref{estimateforL211}, we have
\begin{align}\label{estimateforL21}
    |L_{21}|\lesssim\delta\|\p_{t}^{3}v\|_{1.5}^{2}+C(\delta)(\PP_{0}+T\PP(E^{\kappa})).
\end{align}
Now, we estimate $L_{22}$,
\begin{align}\label{estimateforL22}
    \begin{aligned}
        |L_{22}|&\lesssim\int_{0}^{t}\|\p_{t}(J_{\kappa}(\ak)_{k}^{i}\nn_{l}\phi)\|_{L^{\infty}(\Omega)}\|\p_{t}^{4}(\nn_{k}\phi)\|_{L^{2}(\Omega)} \|\p_{i}\p_{t}^{3} v_{\kappa}^{l}\|_{L^{2}(\Omega)}\\&\lesssim\PP_{0}+T\PP(E^{\kappa}).
    \end{aligned}
\end{align}

and we use a similar strategy to express $L_{23}$ as
\begin{align}
    \begin{aligned}
        L_{23}&=-\int_{\Omega}(\ak)_{k}^{\mu}\p_{\mu}\p_{t}^{4}\phi J_{\kappa}\nn_{l}\phi (\ak)_{k}^{\nu}\p_{\nu}\p_{t}^{3}v_{\kappa}^{l}|_{0}^{t}-\int_{\Omega}\p_{t}^{4}(\ak)_{k}^{\mu}\p_{\mu}\phi J_{\kappa}\nn_{l}\phi \nn_{k}\p_{t}^{3}v_{\kappa}^{l}|_{0}^{t}\\&-\sum_{s=1}^{3}C_{4}^{s}\int_{\Omega}\p_{t}^{s}(\ak)_{k}^{\mu}\p_{\mu}\p_{t}^{4-s}\phi J_{\kappa}\nn_{l}\phi (\ak)_{k}^{\nu}\p_{\nu}\p_{t}^{3}v_{\kappa}^{l}|_{0}^{t}\\&=-\int_{\Omega}(\ak)_{k}^{\mu}\p_{\mu}\p_{t}^{4}\phi J_{\kappa}\nn_{l}\phi (\ak)_{k}^{\nu}\p_{\nu}\p_{t}^{3}v_{\kappa}^{l}|_{0}^{t}+\int_{\Omega}((\ak)_{k}^{i}\p_{t}^{3}\p_{i}v_{\kappa}^{j}(\ak)_{j}^{\mu})\p_{\mu}\phi J_{\kappa}\nn_{l}\phi \nn_{k}\p_{t}^{3}v_{\kappa}^{l}|_{0}^{t}\\&+\sum_{s=1}^{3}C_{3}^{s}\int_{\Omega}\p_{t}^{s}((\ak)_{k}^{i}(\ak)_{j}^{\mu}))\p_{t}^{3-s}\p_{i}v_{\kappa}^{j}\p_{\mu}\phi J_{\kappa}\nn_{l}\phi \nn_{k}\p_{t}^{3}v_{\kappa}^{l}|_{0}^{t}\\&-\sum_{s=1}^{3}C_{4}^{s}\int_{\Omega}\p_{t}^{s}(\ak)_{k}^{\mu}\p_{\mu}\p_{t}^{4-s}\phi J_{\kappa}\nn_{l}\phi (\ak)_{k}^{\nu}\p_{\nu}\p_{t}^{3}v_{\kappa}^{l}|_{0}^{t}.
    \end{aligned}
\end{align}
Thus,
\begin{align}\label{estimateforL23}
    \begin{aligned}
        |L_{23}|&\lesssim\PP_{0}+\|(\ak)_{k}^{\mu}J_{\kappa}\nn_{l}\phi (\ak)_{k}^{\nu}\|_{L^{\infty}(\Omega)}\|\p_{\mu}\p_{t}^{4}\phi\|_{L^{2}(\Omega)} \|\p_{\nu}\p_{t}^{3}v_{\kappa}^{l}\|_{L^{2}(\Omega)}\\&+\|(\ak)_{k}^{i}(\ak)_{j}^{\mu}\p_{\mu}\phi J_{\kappa}\nn_{l}\phi (\ak)_{k}^{\nu}\|_{L^{\infty}(\Omega)}\|\p_{t}^{3}\p_{i}v_{\kappa}^{j}\|_{L^{2}(\Omega)}\|\p_{\nu}\p_{t}^{3}v_{\kappa}^{l}\|_{L^{2}(\Omega)}\\&+\sum_{s=0}^{2}C_{3}^{s+1}\|\p_{t}^{s}\p_{t}((\ak)_{k}^{i}(\ak)_{j}^{\mu}))\p_{t}^{2-s}\p_{i}v_{\kappa}^{j}\|_{L^{2}(\Omega)}\|\p_{\mu}\phi J_{\kappa}\nn_{l}\phi (\ak)_{k}^{\nu}\|_{L^{\infty}(\Omega)}\|\p_{\nu}\p_{t}^{3}v_{\kappa}^{l}\|_{L^{2}(\Omega)}\\&+\sum_{s=0}^{2}C_{4}^{s+1}\|\p_{t}^{s}\p_{t}(\ak)_{k}^{\mu}\p_{\mu}\p_{t}^{3-s}\phi\|_{L^{2}(\Omega)}\| J_{\kappa}\nn_{l}\phi (\ak)_{k}^{\nu}\|_{L^{\infty}(\Omega)}\|\p_{\nu}\p_{t}^{3}v_{\kappa}^{l}\|_{L^{2}(\Omega)}\\&\lesssim\delta\|\p_{t}^{4}\phi\|_{H^{1}(\Omega)}+\delta\|\p_{t}^{3}v\|_{1.5}^{2}+C(\delta)(\PP_{0}+T\PP(E^{\kappa})).
    \end{aligned}
\end{align}
Combining \eqref{expansionforL2}, \eqref{estimateforL21}, \eqref{estimateforL22} and \eqref{estimateforL23}, we have
\begin{align}\label{estimateforL2}
    |L_{2}|\lesssim\delta\|\p_{t}^{4}\phi\|_{H^{1}(\Omega)}+\delta\|\p_{t}^{3}v\|_{1.5}^{2}+C(\delta)(\PP_{0}+T\PP(E^{\kappa})).
\end{align}
Combining \eqref{expansionforcancellation3}, \eqref{estimateforL1}, \eqref{estimateforL2}, we have
\begin{align}
  \int_{0}^{t}\int_{\Omega}J_{\kappa}\nn_{l}\p_{t}^{4}v_{\kappa}^{\mu}\nn_{l}\phi\p_{t}^{4}\nn_{\mu}\phi= -\int_{0}^{t}\int_{\Omega}J_{\kappa}|\p_{t}^{5}\phi|^{2}-\frac{1}{2}\int_{\Omega}J_{\kappa}g_{\kappa}^{ij}\p_{t}^{4}\p_{j}\phi\p_{t}^{4}\p_{i}\phi+\mathcal{R},
\end{align}
where $|\mathcal{R}|\lesssim\delta\|\p_{t}^{4}\phi\|_{1}^{2}+\delta\int_{0}^{t}\|\p_{t}^{5}\phi\|_{L^{2}(\Omega)}^{2}+\delta\|\p_{t}^{3}v\|_{1.5}^{2}+C(\delta)(\PP_{0}+T\PP(E^{\kappa}))$.
\rmk This proposition will provide a key cancellation which is crucial in the full-time derivative estimate for $v$. For the estimate involving the spatial derivative for $\phi$, we can see that it is self-contained. And an $H^{-\frac{1}{2}}$-$H^{\frac{1}{2}}$ duality may be employed to derive the rest of the estimate for $v$. Indeed, if we adjust our framework to start our estimate from $\|v_{ttt}\|_{H^{1}(\Omega)}$ instead of $\|v_{tttt}\|_{L^{2}}$, this cancellation will be much easier to obtain.

\section{Energy estimate II: Estimate for $\curl$ and $\dive$}
\subsection{Curl estimate for $v$, $v_{t}$, $v_{tt}$, $v_{ttt}$, 
$\eta$}

\prop
\begin{align}\label{curlest}
    \sup_{t}(\sum_{k=0}^{4}\|\curl \p_{t}^{k}\eta\|_{4.5-k}^{2} +\|\sqrt{\kappa}\curl \eta\|_{5.5}^{2})+\sum_{k=0}^{4}\int_{0}^{t}\|\sqrt{\kappa}\curl \p_{t}^{k}v\|_{4.5-k}^{2}\lesssim \PP_{0}+T\PP(E^{\kappa}).
\end{align}
\pf
The major difference arises from the necessity of estimate $\phi$ in contrast to the methodology used in \cite{coutand2007well}. Our focus will be only on the modified segment of the proof.
Taking $\curl_{\eta^{\kappa}}$ on the second equation of \eqref{vkappat} shows that 
\begin{align}\label{vt1}
    (\curl_{\eta^{\kappa}}v_{t})^{i}=-\epsilon^{ijk}\nabla_{j}^{\kappa}\p_{t}\phi^{\kappa}\cdot \nabla_{k}^{\kappa}\phi^{\kappa},
\end{align}
Changing the sequence of $\p_{t}$ and $\curl_{\eta^{\kappa}}$, we can arrive at
\begin{align}
    \p_{t}(\curl_{\eta^{\kappa}}v)^{i}+\epsilon^{inm}\nn_{n}(v^{\kappa})^{s}\nn_{s}v_{m}=-\epsilon^{ijk}\nabla_{j}^{\kappa}\p_{t}\phi^{\kappa}\cdot \nabla_{k}^{\kappa}\phi^{\kappa},
\end{align}
After integrating on time, we obtain
\begin{align}\label{v1}
   \curl_{\eta^{\kappa}}v=\curl u_{0}-\int_{0}^{t}\epsilon^{.nm}\nabla_{n}^{\kappa}v^{\kappa}\cdot\nabla^{\kappa}v_{m}-\int_{0}^{t}\epsilon^{.jk}\nabla_{j}^{\kappa}\p_{t}\phi^{\kappa}\cdot \nabla_{k}^{\kappa}\phi^{\kappa}.
\end{align}
From \eqref{v1}, it is clear that
\begin{align}\label{curlv}
   \curl v+\epsilon^{.sw}\p_{i}v_{w}\int_{0}^{t}\p_{t}(\ak)_{s}^{i}=\curl u_{0}-\int_{0}^{t}\epsilon^{.nm}\nabla_{n}^{\kappa}v^{\kappa}\cdot\nabla^{\kappa}v_{m}-\int_{0}^{t}\epsilon^{.jk}\nabla_{j}^{\kappa}\p_{t}\phi^{\kappa}\cdot \nabla_{k}^{\kappa}\phi^{\kappa}.
\end{align}
Given that the precise expression of the equation is not important, we can write
\begin{align}
\begin{aligned}
      \curl v\sim \curl u_{0}+\p v\int_{0}^{t}\p_{t}\ak+\int_{0}^{t}\nn v^{\kappa}\nn v+\int_{0}^{t}\nn\p_{t}\phi^{\kappa}\nn\phi^{\kappa}
\end{aligned}
\end{align}
Integrating on time again, it follows that
\begin{align}
\begin{aligned}
     \curl \eta&\sim t \curl u_{0}+\int_{0}^{t}\p v\int_{0}^{t^\prime}\p_{t}\ak+\int_{0}^{t}\int_{0}^{t^\prime}\nn v^{\kappa}\nn v+\int_{0}^{t}\int_{0}^{t^\prime}\nn\p_{t}\phi^{\kappa}\nn\phi^{\kappa}
\end{aligned} 
\end{align}
Taking spatial derivative $D$, and integrating on time by part, we can write
\begin{align}
   \begin{aligned}
    D \curl \eta&\sim t D\curl u_{0}+\int_{0}^{t}\p Dv\int_{0}^{t^\prime}\p_{t}\ak+\int_{0}^{t}\p v\int_{0}^{t^\prime}\p_{t}D\ak+\int_{0}^{t}\int_{0}^{t^\prime}D\nn v^{\kappa}\nn v+\int_{0}^{t}\int_{0}^{t^\prime}\nn v^{\kappa}D\nn v\\&+\int_{0}^{t}\int_{0}^{t^{\prime}}D(\nabla^{\kappa}\p_{t}\phi^{\kappa}\nabla^{\kappa}\phi^{\kappa})\\&\sim \underbrace{t D\curl u_{0}+\p D\eta\int_{0}^{t}\p_{t}\ak+\int_{0}^{t}\p D\eta\p_{t}\ak+\int_{0}^{t}\p vD\ak+\int_{0}^{t}\int_{0}^{t^\prime}\ak D\ak D v^{\kappa}D v}_{I}\\&\underbrace{+\int_{0}^{t}\int_{0}^{t^\prime}\ak\ak D^{2}v^{\kappa}D v+\int_{0}^{t}\int_{0}^{t^\prime}\ak\ak Dv^{\kappa}D^{2} v}_{II}\\&\underbrace{+\int_{0}^{t}\int_{0}^{t^{\prime}}D(\nabla^{\kappa}\p_{t}\phi^{\kappa}\nabla^{\kappa}\phi^{\kappa})}_{III}
\end{aligned}  
\end{align}
Notice that
\begin{align}\label{approcurlphi}
\begin{aligned}
\int_{0}^{t}\int_{0}^{t^{\prime}}D(\nabla^{\kappa}\p_{t}\phi^{\kappa}\nabla^{\kappa}\phi^{\kappa})&\sim \int_{0}^{t}\int_{0}^{t^{\prime}}D(\p_{t}\nn\phi^{\kappa}\nn\phi^{\kappa})+\int_{0}^{t}\int_{0}^{t^{\prime}}D(A^{\kappa}Dv^{\kappa}\nn\phi^{\kappa}\nn\phi^{\kappa})\\&\sim\int_{0}^{t}\int_{0}^{t^{\prime}}D\p_{t}\nn\phi^{\kappa}\nn\phi^{\kappa}+\int_{0}^{t}\int_{0}^{t^{\prime}}\p_{t}\nn\phi^{\kappa} D\nn\phi^{\kappa}\\&+\int_{0}^{t}\int_{0}^{t^{\prime}}A^{\kappa}D^{2}v^{\kappa}\nn\phi^{\kappa}\nn\phi^{\kappa}+\int_{0}^{t}\int_{0}^{t^{\prime}}A^{\kappa}Dv^{\kappa}D\nn\phi^{\kappa}\nn\phi^{\kappa}\\&+\int_{0}^{t}\int_{0}^{t^{\prime}}DA^{\kappa}Dv^{\kappa}\nn\phi^{\kappa}\nn\phi^{\kappa}\\&\sim\underbrace{\int_{0}^{t}D\nn\phi^{\kappa}\nn\phi^{\kappa}|_{0}^{t^{\prime}}+\int_{0}^{t}\int_{0}^{t^{\prime}}\p_{t}\nn\phi^{\kappa} D\nn\phi^{\kappa}}_{III_{1}}\\&\underbrace{+\int_{0}^{t}\int_{0}^{t^{\prime}}DA^{\kappa}Dv^{\kappa}\nn\phi^{\kappa}\nn\phi^{\kappa}+\int_{0}^{t}\int_{0}^{t^{\prime}}A^{\kappa}Dv^{\kappa}D\nn\phi^{\kappa}\nn\phi^{\kappa}}_{III_{2}}\\&\underbrace{+\int_{0}^{t}A^{\kappa}D^{2}\eta^{\kappa}\nn\phi^{\kappa}\nn\phi^{\kappa}+\int_{0}^{t}\int_{0}^{t^{\prime}}\p_{t}A^{\kappa}D^{2}\eta^{\kappa}\nn\phi^{\kappa}\nn\phi^{\kappa}+\int_{0}^{t}\int_{0}^{t^{\prime}}A^{\kappa}D^{2}\eta^{\kappa}\p_{t}\nn\phi^{\kappa}\nn\phi^{\kappa}}_{III_{3}}
\end{aligned}     
\end{align}
We can control the term $III_{1}$, $III_{2}$ and $III_{3}$ is a direct way. Therefore,
\begin{align}
\begin{aligned}
    \|III\|_{3.5}&\lesssim t\PP_{0}+(\int_{0}^{t}\|D\nn\phi^{\kappa}\|_{3.5}^{2})^{\frac{1}{2}}(\int_{0}^{t}\|\nn\phi^{\kappa}\|_{3.5}^{2})^{\frac{1}{2}}+\int_{0}^{t}(\int_{0}^{t^{\prime}}\|\p_{t}\nn\phi^{\kappa}\| _{3.5}^{2}+\int_{0}^{t^{\prime}}\|D\nn\phi^{\kappa}\|_{3.5}^{2})\\&+\int_{0}^{t}\int_{0}^{t^{\prime}}\PP(E^{\kappa})+\int_{0}^{t}(\int_{0}^{t^{\prime}}\|D\nn\phi^{\kappa}\|_{3.5}^{2}+\int_{0}^{t^{\prime}}\PP(E^{\kappa}))\\&+\int_{0}^{t}\PP(E^{\kappa})+\int_{0}^{t}\int_{0}^{t^{\prime}}\PP(E^{\kappa})+\int_{0}^{t}(\int_{0}^{t^{\prime}}\|\p_{t}\nn\phi^{\kappa}\|_{3.5}^{2}+\int_{0}^{t^{\prime}}\PP(E^{\kappa}))\\&\lesssim T\PP(E^{\kappa}).
\end{aligned}    
\end{align}
The control of the terms corresponding to $\eta$ (I and II) is the same as \cite{coutand2007well}, which reads 
\begin{align}
\begin{aligned}
     \|I\|_{3.5}&\lesssim t \|D\curl u_{0}\|_{3.5}+\|\p D\eta\|_{3.5}\int_{0}^{t}\|\p_{t}\ak\|_{2.5}+\int_{0}^{t}\|\p D\eta\|_{3.5}\|\p_{t}\ak\|_{3.5}\\&+\int_{0}^{t}\|\p v\|_{3.5}\|D\ak\|_{3.5}+\int_{0}^{t}\int_{0}^{t^\prime}\|\ak\|_{3.5}\| D\ak\|_{3.5}\| D v^{\kappa}\|_{3.5}\|D v\|_{3.5}\\&\lesssim T\PP(E^{\kappa}),
\end{aligned}  
\end{align}
and 
\begin{align}
\begin{aligned}
    D^{2}(II)&\sim\underbrace{ \int_{0}^{t}\int_{0}^{t^\prime}\ak\ak D^{4}v^{\kappa}D v+\int_{0}^{t}\int_{0}^{t^\prime}\ak\ak Dv^{\kappa}D^{4} v}_{II_{2,1}}+ \underbrace{\int_{0}^{t}\int_{0}^{t^\prime} D^{3}v^{\kappa} \ak\ak D^{2} v+\int_{0}^{t}\int_{0}^{t^\prime}D^{3} v\ak\ak D^{2}v^{\kappa}}_{II_{2,2}}\\&\underbrace{+\int_{0}^{t}\int_{0}^{t^\prime}D^{3}v^{\kappa}D(\ak\ak) D v+\int_{0}^{t}\int_{0}^{t^\prime}Dv^{\kappa}D(\ak\ak) D^{3} v}_{II_{2,3}}\\&\underbrace{+\int_{0}^{t}\int_{0}^{t^\prime}D^{2}v^{\kappa}D^{2}(\ak\ak) D v+\int_{0}^{t}\int_{0}^{t^\prime}D^{2}(\ak\ak) Dv^{\kappa}D^{2} v+\int_{0}^{t}\int_{0}^{t^{\prime}}D^{2}v^{\kappa}D(\ak\ak) D^{2} v}_{II_{2,4}},
\end{aligned}    
\end{align}
 and integrating by part on time;
 \begin{align}
     II_{2,1}=\int_{0}^{t}\ak\ak D^{4}\eta^{\kappa}D v+\int_{0}^{t}\ak\ak Dv^{\kappa}D^{4} \eta+\int_{0}^{t}\int_{0}^{t^\prime} D^{4}\eta^{\kappa}\p_{t}(\ak\ak D v)+\int_{0}^{t}\int_{0}^{t^\prime}D^{4} \eta\p_{t}(\ak\ak Dv^{\kappa})
 \end{align}
From a direct computation,
\begin{align}
\begin{aligned}
    \|II_{2,1}\|_{1.5}&\lesssim\int_{0}^{t}\|D^{4}\eta^{\kappa}\|_{1.5}\|\ak\ak D v\|_{2}+\int_{0}^{t}\|\ak\ak Dv^{\kappa}\|_{2}\|D^{4} \eta\|_{1.5}\\&+\int_{0}^{t}\int_{0}^{t^\prime} \|D^{4}\eta^{\kappa}\|_{1.5}\|\p_{t}(\ak\ak D v)\|_{1.5}+\int_{0}^{t}\int_{0}^{t^\prime}\|D^{4} \eta\|_{1.5}\|\p_{t}(\ak\ak Dv^{\kappa})\|_{2}\\&\lesssim T\PP(E^{\kappa}).
\end{aligned}  
\end{align}
The rest of terms of $II$ are lower order, $\|II_{2,2}+II_{2,3}+II_{2,4}\|_{1.5}\lesssim T\PP(E^{\kappa})$.
Immediately, we obtain
\begin{align}
    \|D\curl \eta\|_{3.5}\lesssim T\PP(E^{\kappa}).
\end{align}
We can conclude that
\begin{align}
    \|\curl \eta\|_{4.5}^{2}\lesssim \PP_{0}+T\PP(E^{\kappa}).
\end{align}
Now we focus on the estimate of $v$ and the time derivatives on $v$. 
\begin{align}
\begin{aligned}
      \curl v\sim \curl u_{0}+\p v\int_{0}^{t}\p_{t}\ak+\int_{0}^{t}\nn v^{\kappa}\nn v+\int_{0}^{t}\nn\p_{t}\phi^{\kappa}\nn\phi^{\kappa},
\end{aligned}
\end{align}

\begin{align}
    \curl v_{t}\sim\nn v^{\kappa}\nn v+\p_{t} \p\phi^{\kappa}\ak\nn\phi^{\kappa}+\p v_{t}\int_{0}^{t}\p_{t}\ak,
\end{align}
\begin{align}
    \curl v_{tt}\sim\nn v_{t}\nn v^{\kappa}+\p_{t}(\nn v^{\kappa}\nn v)+\p_{t}(\p_{t} \p\phi^{\kappa}\ak\nn\phi^{\kappa})+\p v_{tt}\int_{0}^{t}\p_{t}\ak,
\end{align}
\begin{align}
\begin{aligned}
     \curl v_{ttt}&\sim\nn v_{tt}\nn v^{\kappa}+\p_{t}(\nn v_{t}\nn v^{\kappa})+\p_{t}^{2}(\nn v^{\kappa}\nn v)+\p_{t}^{2}(\p_{t} \p\phi^{\kappa}\ak\nn\phi^{\kappa})+\p v_{ttt}\int_{0}^{t}\p_{t}\ak.
\end{aligned}   
\end{align}
\begin{align}
\begin{aligned}
     \curl v_{tttt}&\sim\p v_{ttt}\p_{t}\ak+\p_{t}(\nn v_{tt}\nn v^{\kappa})+\p_{t}^{2}(\nn v_{t}\nn v^{\kappa})+\p_{t}^{3}(\nn v^{\kappa}\nn v)\\&+\p_{t}^{3}(\p_{t} \p\phi^{\kappa}\ak\nn\phi^{\kappa})+\p v_{tttt}\int_{0}^{t}\p_{t}\ak.
\end{aligned}   
\end{align}
Roughly speaking, the time-spatial mixed norms of $\phi$ are of order 4.5 in $L_{t}^{\infty}$ and of order 5.5 in $L_{t}^2$. Consequently, the term involving $\phi^{\kappa}$ is controlled by the established estimate of $\phi$. The term related to $v$ is handled in the same manner as in \cite{coutand2007well}. Therefore, we present the estimate below without further elaboration.
\begin{align}
    \begin{aligned}
         \|\curl v\|_{3.5}&\lesssim \|\curl u_{0}\|_{3.5}+\|\p v\|_{3.5}\int_{0}^{t}\|\p_{t}\ak\|_{3.5}+\int_{0}^{t}\|\nn v^{\kappa}\nn v\|_{3.5}+\int_{0}^{t}\|\p\p_{t}\phi^{\kappa}\|_{3.5}\|\ak\nn\phi^{\kappa}\|_{3.5}\\&\lesssim \|\curl u_{0}\|_{3.5}+\|\p v\|_{3.5}\int_{0}^{t}\|\p_{t}\ak\|_{3.5}+\int_{0}^{t}\|\nn v^{\kappa}\nn v\|_{3.5}+\int_{0}^{t}\|\p\p_{t}\phi^{\kappa}\|_{3.5}^{2}+\int_{0}^{t}\|\ak\nn\phi^{\kappa}\|_{3.5}^{2}\\&\lesssim\PP_{0}+T\PP(E^{\kappa}),
    \end{aligned}
\end{align}
\begin{align}
    \|\curl v_{t}\|_{2.5}&\lesssim\|\nn v^{\kappa}\nn v\|_{2.5}+\|\p_{t} \p\phi^{\kappa}\|_{2.5}\|\ak\nn\phi^{\kappa}\|_{2.5}+\|\p v_{t}\|_{2.5}\int_{0}^{t}\|\p_{t}\ak\|_{2.5}\\&\lesssim\PP_{0}+T\PP(E^{\kappa}),
\end{align}
\begin{align}
\begin{aligned}
    \|\curl v_{tt}\|_{1.5}&\lesssim\| \p_{t}\p v\|_{1.5}\|\ak\nn v^{\kappa}\|_{2}+\|\p_{t}(\nn v^{\kappa}\nn v)\|_{1.5}+\|\p_{t}(\p_{t} \p\phi^{\kappa}\ak\nn\phi^{\kappa})\|_{1.5}+\|\p v_{tt}\|_{1.5}\int_{0}^{t}\|\p_{t}\ak\|_{2}\\&\lesssim \PP_{0}+T\PP(E^{\kappa}),
\end{aligned}    
\end{align}
\begin{align}
\begin{aligned}
     \|\curl v_{ttt}\|_{0.5}&\lesssim\|\p_{t}^{2} \p v\|_{0.5}\|\ak\nn v^{\kappa}\|_{2}+\|\p_{t}(\p_{t}\p v\ak\nn v^{\kappa})\|_{0.5}+\|\p_{t}^{2}(\nn v^{\kappa}\nn v)\|_{0.5}\\&+\|\p_{t}^{2}(\p_{t} \p\phi^{\kappa}\ak\nn\phi^{\kappa})\|_{0.5}+\|\p v_{ttt}\|_{0.5}\int_{0}^{t}\|\p_{t}\ak\|_{2}\\&\lesssim\PP_{0}+T\PP(E^{\kappa}).
\end{aligned}   
\end{align}
We now turn to controlling the terms involving $\kappa$. It is evident that the terms associated with $\phi$ are better controlled for a similar reason. To be precise,
\begin{align}
    \begin{aligned}        \|\sqrt{\kappa}D\curl\eta\|_{4.5}\lesssim\|\sqrt{\kappa}I\|_{4.5}+\|\sqrt{\kappa}II\|_{4.5}+\|\sqrt{\kappa}III\|_{4.5},
    \end{aligned}
\end{align}
As
\begin{align}
\begin{aligned}
     D^{2}(III_{1})&\sim \int_{0}^{t}D^{3}\nn\phi^{\kappa}\nn\phi^{\kappa}+\int_{0}^{t}D^{2}\nn\phi^{\kappa} D\nn\phi^{\kappa}+\int_{0}^{t}D^{2}(D\nn\phi^{\kappa}\nn\phi^{\kappa}(0))\\&+\int_{0}^{t}\int_{0}^{t^{\prime}}D^{3}\nn\phi^{\kappa}\p_{t}\nn\phi^{\kappa}+D^{2}\nn\phi^{\kappa}\p_{t}D\nn\phi^{\kappa}+\int_{0}^{t}D^{2}\nn\phi^{\kappa} D\nn\phi^{\kappa}|_{0}^{t^{\prime}}
\end{aligned}   
\end{align}
\begin{align}
    \begin{aligned}
       \|\sqrt{\kappa}D^{2}(III_{1})\|_{2.5}  &\lesssim T\PP_{0}+\int_{0}^{t} \|\sqrt{\kappa}D^{3}\nn\phi^{\kappa}\|_{2.5} \|\nn\phi^{\kappa}\|_{2.5}+\int_{0}^{t} \|\sqrt{\kappa}D^{2}\nn\phi^{\kappa} D\nn\phi^{\kappa}\|_{2.5}\\&+\int_{0}^{t}\int_{0}^{t^{\prime}} \|\sqrt{\kappa}D^{3}\nn\phi^{\kappa}\|_{2.5} \|\p_{t}\nn\phi^{\kappa}\|_{2.5}+ \|\sqrt{\kappa} D^{2}\nn\phi^{\kappa}\p_{t}D\nn\phi^{\kappa}\|_{2.5}\\&\lesssim T^{\frac{1}{2}}\PP(E^{\kappa}).
    \end{aligned}
\end{align}
\begin{align}
    \begin{aligned}       DIII_{2}&\sim\int_{0}^{t}\int_{0}^{t^{\prime}}D^{2}A^{\kappa}Dv^{\kappa}\nn\phi^{\kappa}\nn\phi^{\kappa}+\int_{0}^{t}\int_{0}^{t^{\prime}}DA^{\kappa}D^{2}v^{\kappa}\nn\phi^{\kappa}\nn\phi^{\kappa}+\int_{0}^{t}\int_{0}^{t^{\prime}}DA^{\kappa}Dv^{\kappa}D\nn\phi^{\kappa}\nn\phi^{\kappa}\\&+\int_{0}^{t}\int_{0}^{t^{\prime}}A^{\kappa}D^{2}v^{\kappa}D\nn\phi^{\kappa}\nn\phi^{\kappa}+\int_{0}^{t}\int_{0}^{t^{\prime}}A^{\kappa}Dv^{\kappa}D^{2}\nn\phi^{\kappa}\nn\phi^{\kappa}+\int_{0}^{t}\int_{0}^{t^{\prime}}A^{\kappa}Dv^{\kappa}D\nn\phi^{\kappa} D\nn\phi^{\kappa},
    \end{aligned}
\end{align}
we have
\begin{align}
    \begin{aligned}         \|\sqrt{\kappa}DIII_{2}\|_{3.5}&\lesssim\int_{0}^{t}\int_{0}^{t^{\prime}} \|\sqrt{\kappa}D^{2}A^{\kappa}\|_{3.5}  \|Dv^{\kappa}\nn\phi^{\kappa}\nn\phi^{\kappa}\|_{3.5} +\int_{0}^{t}\int_{0}^{t^{\prime}} \|\sqrt{\kappa}D^{2}v^{\kappa}\|_{3.5}  \|DA^{\kappa}\nn\phi^{\kappa}\nn\phi^{\kappa}\|_{3.5} \\&+\int_{0}^{t}\int_{0}^{t^{\prime}} \|\sqrt{\kappa}DA^{\kappa}Dv^{\kappa}D\nn\phi^{\kappa}\nn\phi^{\kappa}\|_{3.5} +\int_{0}^{t}\int_{0}^{t^{\prime}} \|\sqrt{\kappa}D^{2}v^{\kappa}\|_{3.5} \| A^{\kappa} D\nn\phi^{\kappa}\nn\phi^{\kappa}\|_{3.5}\\&+\int_{0}^{t}\int_{0}^{t^{\prime}} \|\sqrt{\kappa}D^{2}\nn\phi^{\kappa}\|_{3.5}  \|A^{\kappa}Dv^{\kappa}\nn\phi^{\kappa}\|_{3.5} +\int_{0}^{t}\int_{0}^{t^{\prime}} \|\sqrt{\kappa}A^{\kappa}Dv^{\kappa}D\nn\phi^{\kappa} D\nn\phi^{\kappa}\|_{3.5} \\&\lesssim T\PP(E^{\kappa}),
    \end{aligned}
\end{align}
as an application of Cauchy inequality. 
\begin{align}
    \begin{aligned}        DIII_{3}&\sim\int_{0}^{t}A^{\kappa}D^{3}\eta^{\kappa}\nn\phi^{\kappa}\nn\phi^{\kappa}+\int_{0}^{t}DA^{\kappa}D^{2}\eta^{\kappa}\nn\phi^{\kappa}\nn\phi^{\kappa}+\int_{0}^{t}A^{\kappa}D^{2}\eta^{\kappa}D\nn\phi^{\kappa}\nn\phi^{\kappa}\\&+\int_{0}^{t}\int_{0}^{t^{\prime}}\p_{t}A^{\kappa}D^{3}\eta^{\kappa}\nn\phi^{\kappa}\nn\phi^{\kappa}+\int_{0}^{t}\int_{0}^{t^{\prime}} D^{2}v^{\kappa}\ak\ak D^{2}\eta^{\kappa}\nn\phi^{\kappa}\nn\phi^{\kappa}\\&+\int_{0}^{t}\int_{0}^{t^{\prime}}D v^{\kappa} D\ak\ak D^{2}\eta^{\kappa}\nn\phi^{\kappa}\nn\phi^{\kappa}+\int_{0}^{t}\int_{0}^{t^{\prime}}\p_{t}A^{\kappa}D^{2}\eta^{\kappa}D\nn\phi^{\kappa}\nn\phi^{\kappa}\\&+\int_{0}^{t}\int_{0}^{t^{\prime}}A^{\kappa}D^{3}\eta^{\kappa}\p_{t}\nn\phi^{\kappa}\nn\phi^{\kappa}+\int_{0}^{t}\int_{0}^{t^{\prime}}DA^{\kappa}D^{2}\eta^{\kappa}\p_{t}\nn\phi^{\kappa}\nn\phi^{\kappa}\\&+\int_{0}^{t}\int_{0}^{t^{\prime}}A^{\kappa}D^{2}\eta^{\kappa}\p_{t}\nn\phi^{\kappa} D\nn\phi^{\kappa}+\int_{0}^{t}A^{\kappa}D^{2}\eta^{\kappa}D\nn\phi^{\kappa}\nn\phi^{\kappa}|_{0}^{t^{\prime}}+\int_{0}^{t}\int_{0}^{t^{\prime}}A^{\kappa}D^{2}v^{\kappa}D\nn\phi^{\kappa}\nn\phi^{\kappa}\\&+\int_{0}^{t}\int_{0}^{t^{\prime}}D^{2}\eta^{\kappa}D\nn\phi^{\kappa}\p_{t}(A^{\kappa}\nn\phi^{\kappa}),
    \end{aligned}
\end{align}
\begin{align}
    \begin{aligned}        \|\sqrt{\kappa}DIII_{3}\|_{3.5}&\lesssim \int_{0}^{t}\|\sqrt{\kappa}D^{3}\eta^{\kappa}\|_{3.5}\|A^{\kappa}\nn\phi^{\kappa}\nn\phi^{\kappa}\|_{3.5}+\int_{0}^{t}\|\sqrt{\kappa}DA^{\kappa}D^{2}\eta^{\kappa}\nn\phi^{\kappa}\nn\phi^{\kappa}\|_{3.5}\\&+\int_{0}^{t}\|\sqrt{\kappa}A^{\kappa}D^{2}\eta^{\kappa}D\nn\phi^{\kappa}\nn\phi^{\kappa}\|_{3.5}+\int_{0}^{t}\int_{0}^{t^{\prime}}\|\sqrt{\kappa}D^{3}\eta^{\kappa}\|_{3.5}\|\sqrt{\kappa}\p_{t}A^{\kappa}\nn\phi^{\kappa}\nn\phi^{\kappa}\|_{3.5}\\&+\int_{0}^{t}\int_{0}^{t^{\prime}} \|\sqrt{\kappa}D^{2}v^{\kappa}\|_{3.5}\|\ak\ak D^{2}\eta^{\kappa}\nn\phi^{\kappa}\nn\phi^{\kappa}\|_{3.5}\\&+\int_{0}^{t}\int_{0}^{t^{\prime}}\|\sqrt{\kappa}D v^{\kappa} D\ak\ak D^{2}\eta^{\kappa}\nn\phi^{\kappa}\nn\phi^{\kappa}\|_{3.5}+\int_{0}^{t}\int_{0}^{t^{\prime}}\|\sqrt{\kappa}\p_{t}A^{\kappa}D^{2}\eta^{\kappa}D\nn\phi^{\kappa}\nn\phi^{\kappa}\|_{3.5}\\&+\int_{0}^{t}\int_{0}^{t^{\prime}}\|\sqrt{\kappa}D^{3}\eta^{\kappa}\|_{3.5}\|A^{\kappa}\p_{t}\nn\phi^{\kappa}\nn\phi^{\kappa}\|_{3.5}+\int_{0}^{t}\int_{0}^{t^{\prime}}\|\sqrt{\kappa}DA^{\kappa}D^{2}\eta^{\kappa}\p_{t}\nn\phi^{\kappa}\nn\phi^{\kappa}\|_{3.5}\\&+\int_{0}^{t}\int_{0}^{t^{\prime}}\|\sqrt{\kappa}A^{\kappa}D^{2}\eta^{\kappa}\p_{t}\nn\phi^{\kappa} D\nn\phi^{\kappa}\|_{3.5}+\int_{0}^{t}\int_{0}^{t^{\prime}}\|\sqrt{\kappa}D^{2}v^{\kappa}\|_{3.5}\|A^{\kappa}D\nn\phi^{\kappa}\nn\phi^{\kappa}\|_{3.5}\\&+\int_{0}^{t}\int_{0}^{t^{\prime}}\|\sqrt{\kappa}D^{2}\eta^{\kappa}D\nn\phi^{\kappa}\p_{t}(A^{\kappa}\nn\phi^{\kappa})\|_{3.5}+T\PP_{0}\\&\lesssim T^{\frac{1}{2}}\PP(E^{\kappa}),
    \end{aligned}
\end{align}
invoking the Cauchy's inequality. The estimate of $\|\sqrt{\kappa}I\|_{4.5}$ and $\|\sqrt{\kappa}II\|_{4.5}$ is the same as \cite{coutand2007well}. We omit the details.
Hence, we can write that $\|\sqrt{\kappa}\curl \eta\|_{5.5}\lesssim T^{\frac{1}{2}}\PP(E^{\kappa})$.

Due to a similar reason, we can also control the terms related to $\int_{0}^{t}\|\sqrt{\kappa}\curl \p_{t}^{k}v\|_{4.5-k}^{2}$, $k=0,1,2,3,4$. 
\begin{align}
\begin{aligned}
     D\curl v&\sim D\curl u_{0}+D^{2}v\int_{0}^{t}\p_{t}\ak+Dv\int_{0}^{t}(D^{2}v^{\kappa}\ak\ak+Dv^{\kappa}D\ak\ak)+D^{2}\phi^{\kappa}\ak\ak D\phi^{\kappa}|_{0}^{t}\\&+\int_{0}^{t}D^{2}\phi^{\kappa}\p_{t}(\ak\ak) D\phi^{\kappa}+\int_{0}^{t}D^{2}\phi^{\kappa}\ak\ak \p_{t}D\phi^{\kappa}+\int_{0}^{t}D\p_{t}\phi^{\kappa} D(\ak\ak) D\phi^{\kappa},
\end{aligned}   
\end{align}
\begin{align}
    \begin{aligned}
       \|\sqrt{\kappa}D\curl v\|_{3.5}&\lesssim \PP_{0}+\|\sqrt{\kappa}D^{2}v\|_{3.5}\int_{0}^{t}\|\p_{t}\ak\|_{3.5}\\&+\|D v\|_{3.5}\int_{0}^{t}(\|\sqrt{\kappa}D^{2} v^{\kappa}\|_{3.5}\|\ak\ak\|_{3.5}+\sqrt{\kappa}\|D v^{\kappa}D\ak\ak\|_{3.5})\\&+\|\sqrt{\kappa}D^{2}\phi^{\kappa}\|_{3.5}(\|\ak\ak D\phi^{\kappa}\|_{3.5})+\sqrt{\kappa}\int_{0}^{t}\|D^{2}\phi^{\kappa}\|_{3.5}\|\p_{t}(\ak\ak) D\phi^{\kappa}\|_{3.5}\\&+\sqrt{\kappa}\int_{0}^{t}\|D^{2}\phi^{\kappa}\|_{3.5}\|\ak\ak \|_{3.5}\|\p_{t}D\phi^{\kappa}\|_{3.5}+\sqrt{\kappa}\int_{0}^{t}\|D\p_{t}\phi^{\kappa}\|_{3.5}\| D(\ak\ak) D\phi^{\kappa}\|_{3.5}\\&\lesssim\|\sqrt{\kappa}D^{2}v\|_{3.5}T\PP(E^{\kappa})+\PP_{0}+T\PP(E^{\kappa}).
    \end{aligned}
\end{align}
We obtain
\begin{align}
    \int_{0}^{t}\|\sqrt{\kappa}\curl v\|_{4.5}^{2}\lesssim \PP_{0}+T\PP(E^{\kappa}).
\end{align}
Similarly,
\begin{align}
\begin{aligned}
     \|\sqrt{\kappa}\curl v_{t}\|_{3.5}&\lesssim\sqrt{\kappa}\|\ak\ak Dv^{\kappa}D v\|_{3.5}+\sqrt{\kappa}\|\p_{t}\p\phi^{\kappa}\|_{3.5}\|\ak\ak D\phi^{\kappa}\|_{3.5}+\|\sqrt{\kappa}Dv_{t}\|_{3.5}\int_{0}^{t}\|\p_{t}\ak\|_{3.5}\\&\lesssim\sqrt{\kappa}\PP(E^{\kappa})+\|\sqrt{\kappa}D v_{t}\|_{3.5}(T\PP(E^{\kappa})),
\end{aligned} 
\end{align}
\begin{align}
    \begin{aligned}
       \int_{0}^{t} \|\sqrt{\kappa}\curl v_{t}\|_{3.5} ^{2}\lesssim \PP_{0}+T\PP(E^{\kappa}),
    \end{aligned}
\end{align}

\begin{align}
    \begin{aligned}
        \|\sqrt{\kappa}\curl v_{tt}\|_{2.5}&\lesssim \sqrt{\kappa}\|D v_{t}\ak\ak D v^{\kappa}\|_{2.5}+\sqrt{\kappa}\|\p_{t}(\nn v^{\kappa}\nn v)\|_{2.5}+\sqrt{\kappa}\|\p_{t}^{2}D\phi^{\kappa}\|_{2.5}\|\ak\ak D\phi^{\kappa}\|_{2.5}\\&+\sqrt{\kappa}\|\p_{t}D\phi^{\kappa}\p_{t}(\ak\ak D\phi^{\kappa})\|_{2.5}+\|\sqrt{\kappa}D v_{tt}\|_{2.5}\int_{0}^{t}\|\p_{t}\ak\|_{2.5}\\&\lesssim \sqrt{\kappa}(\PP_{0}+\PP(E^{\kappa}))+\|\sqrt{\kappa}Dv_{tt}\|_{2.5}T\PP(E^{\kappa}),
    \end{aligned}
\end{align}
\begin{align}
    \begin{aligned}
         \int_{0}^{t}\|\sqrt{\kappa}\curl v_{tt}\|_{2.5}^{2}&\lesssim \PP_{0}+T\PP(E^{\kappa}).
    \end{aligned}
\end{align},
\begin{align}
    \begin{aligned}
        \|\sqrt{\kappa}\curl v_{ttt}\|_{1.5}&\lesssim \sqrt{\kappa}\|\p_{t}^{2}Dv\ak\ak Dv^{\kappa}\|_{1.5}+\sqrt{\kappa}\|\p_{t}(Dv_{t}\ak\ak D v^{\kappa})\|_{1.5}\\&+\sqrt{\kappa}\|\p_{t}^{2}(\nn v^{\kappa}\nn v)\|_{0.5}+\sqrt{\kappa}\|\p_{t}^{3}\p\phi^{\kappa}\|_{1.5}\|\ak\ak D\phi^{\kappa}\|_{2}\\&+\sqrt{\kappa}\sum_{i+j=1}\|\p_{t}^{i}\p_{t}\p\phi^{\kappa}\p_{t}^{j}\p_{t}(\ak\ak D\phi^{\kappa})\|_{1.5}+\sqrt{\kappa}\|D v_{ttt}\|_{1.5}\int_{0}^{t}\|\p_{t}\ak\|_{2}\\&\lesssim \sqrt{\kappa}(\PP_{0}+\PP(E^{\kappa}))+\sqrt{\kappa}\|D v_{ttt}\|_{1.5}T\PP(E^{\kappa}),
    \end{aligned}
\end{align}
\begin{align}
        \int_{0}^{t}  \|\sqrt{\kappa}\curl v_{ttt}\|_{1.5}^{2}\lesssim \PP_{0}+T\PP(E^{\kappa}),
\end{align}
\begin{align}
\begin{aligned}
      \|\sqrt{\kappa}\curl v_{tttt}\|_{0.5}&\lesssim \sqrt{\kappa}(\|D v_{ttt}\|_{0.5}\|\p_{t}\ak\|_{2}+\|\p_{t}(\p_{t}Dv_{t}\ak\ak Dv^{\kappa})\|_{0.5})\\&+\sqrt{\kappa}(\|\p_{t}^{2}(D v_{t}\ak\ak D v^{\kappa})\|_{0.5}+\|\p_{t}^{3}(\nn v^{\kappa}\nn v)\|_{0.5})\\&+\sqrt{\kappa}\|\p_{t}^{4}D\phi^{\kappa}\|_{0.5}\|\ak\ak D\phi^{\kappa}\|_{2}+\sqrt{\kappa}\sum_{i+j=2}\|\p_{t}^{i}\p_{t}D\phi^{\kappa} \p_{t}^{j}\p_{t}(\ak\ak D\phi^{\kappa})\|_{0.5}\\&+\|\sqrt{\kappa}D v_{tttt}\|_{0.5}\int_{0}^{t}\|\p_{t}\ak\|_{2}\\&\lesssim \sqrt{\kappa}(\PP_{0}+\PP(E^{\kappa}))+\sqrt{\kappa}\|\p_{t}^{4}D\phi^{\kappa}\|_{0.5}(\PP_{0}+T\PP(E^{\kappa}))+\|\sqrt{\kappa}v_{tttt}\|_{1.5}(T\PP(E^{\kappa})),
\end{aligned}
\end{align}
\begin{align}
    \begin{aligned}
       \int_{0}^{t} \|\sqrt{\kappa}\curl v_{tttt}\|_{0.5}^{2}\lesssim\PP_{0}+T\PP(E^{\kappa}). 
    \end{aligned}
\end{align}

This completes the curl estimate \eqref{curlest}.
\subsection{Div estimate}
\prop
\begin{align}\label{divest}
    \sup_{t}(\sum_{k=0}^{4}\|\dive \p_{t}^{k}\eta\|_{4.5-k}^{2} +\|\sqrt{\kappa}\dive \eta\|_{5.5}^{2})+\sum_{k=0}^{4}\int_{0}^{t}\|\sqrt{\kappa}\dive \p_{t}^{k}v\|_{4.5-k}^{2}\lesssim \delta E^{\kappa}+C(\delta)(\PP_{0}+T^{\frac{1}{2}}\PP(E^{\kappa})).
\end{align}
\pf
Since the estimate is totally based on the incompressible condition, an identical argument in \cite{coutand2007well} is enough. i.e. No term is related $\phi$. (Actually, it is exactly the same estimate of terms without $\phi$ as above). To see this,
\begin{align}
\begin{aligned}
    \p_{t}( (\ak)_{i}^{j}D\p_{j}\eta^{i})&=\p_{t}(\ak)_{i}^{j}D\p_{j}\eta^{i}+(\ak)_{i}^{j}D\p_{j}v^{i}\\&=\p_{t}(\ak)_{i}^{j}D\p_{j}\eta^{i}-D(\ak)_{i}^{j}\p_{j}v^{i}.
\end{aligned}   
\end{align}
Thus,
\begin{align}
  D \di \eta =D\p_{j}\eta^{i}\int_{0}^{t}\p_{t}(\ak)_{i}^{j}+\int_{0}^{t} \p_{t}(\ak)_{i}^{j}D\p_{j}\eta^{i}-\int_{0}^{t}D(\ak)_{i}^{j}\p_{j}v^{i}.
\end{align}
Now, other estimate is similar.
\section{Energy estimate III: Estimate for $\pb_{A}$}

We first establish the regularity of the pressure. For $\|q\|_{4.5}$, $\|q_{t}\|_{3.5}$, and $\|q_{tt}\|_{2.5}$, we apply the elliptic estimate provided in \eqref{quasiellesti} under the Neumann boundary condition. Regarding the data $\|q\|_{L^{2}}$, $\|q_{t}\|_{L^{2}}$, and $\|q_{tt}\|_{L^{2}}$, the bounds obtained from a rough elliptic estimate under the Dirichlet boundary condition for $\|q\|_{H^1}$, $\|q_{t}\|_{H^1}$, and $\|q{tt}\|_{H^{1}}$ are sufficient.

However, when estimating $\|q_{ttt}\|_1$, we require the additional assistance of $\bar{q} = \frac{1}{|\Omega|} \int q$ to provide a sharper bound. If we were to rely only on the elliptic estimate with Dirichlet boundary conditions for $\|q_{ttt}\|_{1}$ to control $\|q_{ttt}\|_{L^{2}}$, a regularity loss would occur at the boundary. Specifically, $\|\bar{\Delta} \p_t^3 \eta\|_{H^{\frac{1}{2}}(\partial\Omega)} \sim \|v_{tt}\|_{H^{2.5}(\partial\Omega)}$, while our energy functional only contains $\|\p_{t}^{2} v\|_{H^{2.5}(\Omega)}$. This discrepancy highlights the need for a more refined approach. An appropriate form of lemma was proposed by Ignotova \& Kukavica\cite{ignatova2016local} stated as following:
\lem[Low regularity elliptic estimates]
Let $\Omega$ be a smooth bounded domain. $\|b^{ij}\|_{\infty}\leq M$ satisfies strong elliptic condition. i.e. $b^{ij}\xi_{i}\xi_{j}\geq c|\xi|^{2}$, $\forall \xi\in\R^{n}$ for some constant $c>0$. Let $q$ be an $H^{1}$ weak solution of the
\begin{align}\label{lemmaforq}
    \p_{i}(b^{ij}\p_{j}q)&=\di \pi\quad\textit{in $\Omega$,}\\
    b_{mk}\p_{k}qN^{m}&=g\quad\textit{on $\p \Omega$,}
\end{align}
where $\pi$, $\di \pi\in L^{2}(\Omega)$ and $g\in H^{-\frac{1}{2}}(\Omega)$ with the compatibility condition
\begin{align}
    \int_{\p\Omega}(\pi\cdot N-g)=0,
\end{align}
then we have
\begin{align}
    \|q-\bar{q}\|_{H^{1}(\Omega)}\lesssim \|\pi\|_{L^{2}(\Omega)}+\|g-\pi\cdot N\|_{H^{-\frac{1}{2}}(\p\Omega)}.
\end{align}
\pf See \cite{ignatova2016local}. Alternatively, we prove the lemma as follows:
First, we use the definition of weak solution.
\begin{align}    \int_{\Omega}b^{ij}\p_{j}q\p_{i}\phi=\int_{\p\Omega}g\phi -\int_{\Omega}\di \pi \phi,\quad \text{$\forall\phi\in H^{1}(\Omega)$.}
\end{align}
Since $\pi$,$\di \pi\in L^{2}(\Omega)$, the normal trace of $\pi$ is well defined. We can write
\begin{align}    \int_{\Omega}b^{ij}\p_{j}q\p_{i}\phi=\int_{\p\Omega}g\phi +\int_{\Omega} \pi\cdot\nabla \phi-\int_{\p\Omega}\pi\cdot N\phi.
\end{align}
Setting $\phi=q$ and utilizing Sobolev trace lemma, we find
\begin{align}
\begin{aligned}
    \|\nabla q\|_{L^{2}(\Omega)}^{2}&\lesssim\|g-\pi\cdot N\|_{H^{-\frac{1}{2}}(\p\Omega)}\|q\|_{H^{\frac{1}{2}}(\p\Omega)}+\|\pi\|_{L^{2}(\Omega)}\|\nabla q\|_{L^{2}(\Omega)}\\&\lesssim\|g-\pi\cdot N\|_{H^{-\frac{1}{2}}(\p\Omega)}\|q\|_{H^{1}(\Omega)}+\|\pi\|_{L^{2}(\Omega)}\|q\|_{H^{1}(\Omega)}.
\end{aligned}      
\end{align}
Since $q-\bar{q}=q-\frac{1}{|\Omega|}\int_{\Omega}q$ is also a weak solution to \eqref{lemmaforq}, we get
\begin{align}\label{lemmaforqf1}
    \|\nabla(q-\bar{q})\|_{L^{2}(\Omega)}^{2}\lesssim(\|\pi\|_{L^{2}(\Omega)}+\|g-\pi\cdot N\|_{H^{-\frac{1}{2}}(\p\Omega)})\|q-\bar{q}\|_{H^{1}(\Omega)}
\end{align}
Now, invoking Poincar\'{e}-Wirtinger inequality, we have
\begin{align}\label{lemmaforqf2}
    \|q-\bar{q}\|_{L^{2}(\Omega)}\lesssim\|\nabla (q-\bar{q})\|_{L^{2}(\Omega)}.
\end{align}
Combining \eqref{lemmaforqf1} and \eqref{lemmaforqf2}, we obtain
\begin{align}
    \|q-\bar{q}\|_{H^{1}(\Omega)}\lesssim\|\pi\|_{L^{2}(\Omega)}+\|g-\pi\cdot N\|_{H^{-\frac{1}{2}}(\p\Omega)}.
\end{align}

Now, we study our case $\Omega=\T^{2}\times(0,1)$. Notice that
\begin{align}
    \|q\|_{H^{1}(\Omega)}\leq \|q-\bar{q}\|_{H^{1}(\Omega)}+|\bar{q}||\Omega|^{\frac{1}{2}}.
\end{align}
It remains to control $|\bar{q}|$. To this end, let $H$ be a harmonic function satisfying
\begin{align}
    \lap H&=1,\quad \text{in $\Omega$,}\\
    H&=0,\quad\text{on $\Gamma_{1}$,}\\
    \p_{N}H&=1. \quad\text{on $\Gamma_{0}$}
\end{align}
Actually, $H=\frac{1}{2}x_{3}^{2}-1$. A direct computation shows that
\begin{align}
    |\bar{q}|=|\int_{\Omega} q\lap H|=|-\int_{\Omega}\nabla q\nabla H+\int_{\Gamma_{1}}q|\leq\|\nabla q\|_{L^{2}(\Omega)}+\|q\|_{L^{2}(\Gamma_{1})}.
\end{align}
Therefore, we have
\begin{align}
    \|q\|_{H^{1}(\Omega)}\lesssim\|\pi\|_{L^{2}(\Omega)}+\|g-\pi\cdot N\|_{H^{-\frac{1}{2}}(\Omega)}+\|q\|_{L^{2}(\Gamma_{1})}.
\end{align}
\rmk
However, due to the existence of artificial viscosity term, for $q_{ttt}$, we can't use Dirichlet data for $q_{ttt}$ directly. Instead, we use Euler equation to give its Dirichlet boundary data. 
\lem[Pressure estimate]
\begin{align} 
\begin{aligned}
\|q\|_{4.5}+\|q_{t}\|_{3.5}+\|q_{tt}\|_{2.5}
   \lesssim\PP(E^{\kappa})
\end{aligned}   .
\end{align}
\begin{align}\label{estimateforqtttpre}
     \|q_{ttt}\|_{H^{1}(\Omega)}^{2}\lesssim(\PP_{0}+T\PP(E^{\kappa}))\|v_{tttt}\|_{L^{2}(\Omega)}^{2}+\PP(E_{\kappa}^{\prime})+\PP_{0}+T\PP(E^{\kappa})
\end{align}
\pf
Consider
\begin{align}\label{elliforq}
\begin{aligned}
     \p_{i}(g_{\kappa}^{ij}\p_{j}q)&=-\p_{i}(\phi^{\kappa}_{t}g_{\kappa}^{ij}\p_{j}\phi^{\kappa})-\p_{i}((\ak)_{j}^{i}(\nn_{j}\nn_{k}\phi^{\kappa})\nn_{k}\phi^{\kappa})-\p_{i}((\ak)_{k}^{i}v_{t}^{k})\quad \text{in $\Omega$,}\\
    (\ak)_{k}^{3}(\ak)_{k}^{i}\p_{i}q&=-\phi^{\kappa}_{t}(\ak)_{k}^{3}(\ak)_{k}^{j}\p_{j}\phi^{\kappa}-\nn_{k}\phi^{\kappa} (\ak)_{j}^{3}\nn_{j}\nn_{k}\phi^{\kappa}-(\ak)_{k}^{3}v_{t}^{k}\quad\text{on $\Gamma$},\\
    q&=-\sigma (\frac{\sqrt{h}}{\sqrt{h_{\kappa}}}\lap_{h} (\eta)\cdot n_{\kappa})+\kappa\frac{1}{\sqrt{h_{\kappa}}}(1-\bar{\lap})(v\cdot n_{\kappa})\quad \textit{on $\Gamma_{1}$}.
\end{aligned}   
\end{align}
Utilizing Lemma 2.2, we have
\begin{align}
\begin{aligned}
     \|q\|_{4.5}\lesssim\PP(E^{\kappa})(&\|-\p_{i}(\phi^{\kappa}_{t}g_{\kappa}^{ij}\p_{j}\phi^{\kappa})-\p_{i}((\ak)_{j}^{i}(\nn_{j}\nn_{k}\phi^{\kappa})\nn_{k}\phi^{\kappa})-\p_{i}((\ak)_{k}^{i}v_{t}^{k})\|_{2.5}\\&\|-\phi^{\kappa}_{t}(\ak)_{k}^{3}(\ak)_{k}^{j}\p_{j}\phi^{\kappa}-\nn_{k}\phi^{\kappa} (\ak)_{k}^{3}(\ak)_{k}^{j}\p_{j}\nn_{k}\phi^{\kappa}-(\ak)_{k}^{3}v_{t}^{k}\|_{H^{3}(\Gamma)}+\|q\|_{L^{2}(\Omega)}).\end{aligned}  
\end{align}
Using the lower order estimate
\begin{align}
\begin{aligned}
     \|q\|_{H^{1}(\Omega)}&\lesssim\|\phi^{\kappa}_{t}g_{\kappa}^{ij}\p_{j}\phi^{\kappa}\|_{L^{2}(\Omega)}+\|\ak\|_{L^{\infty}}\|\nn\nn\phi^{\kappa}\|_{L^{2}(\Omega
    )}\|\nn\phi^{\kappa}\|_{L^{\infty}}+\|\ak\|_{L^{\infty}(\Omega)}\|v_{t}\|_{L^{2}(\Omega)}\\&+\PP(\|\eta\|_{H^{2}(\p\Omega)}, \kappa\|v\|_{H^{2}(\p\Omega)},\|\p\eta\|_{L^{\infty}(\p\Omega)})\\&\lesssim\PP(E^{\kappa}).
\end{aligned} 
\end{align}
It implies that
\begin{align}
    \|q\|_{4.5}\lesssim\PP(E^{\kappa}).
\end{align}

Differentiating \eqref{elliforq} on time, we find
\begin{align}\label{elliforqt}
\begin{aligned}
   \p_{i}(g_{\kappa}^{ij}\p_{j}q_{t})&=- \p_{i}(\p_{t}g_{\kappa}^{ij}\p_{j}q)+\p_{t}(-\p_{i}(\phi^{\kappa}_{t}g_{\kappa}^{ij}\p_{j}\phi^{\kappa})-\p_{i}((\ak)_{j}^{i}(\nn_{j}\nn_{k}\phi^{\kappa})\nn_{k}\phi^{\kappa})-\p_{i}((\ak)_{k}^{i}v_{t}^{k}))\quad \text{in $\Omega$,}\\
    (\ak)_{k}^{3}(\ak)_{k}^{i}\p_{i}q_{t}&=- \p_{t}((\ak)_{k}^{3}(\ak)_{k}^{i})\p_{i}q-\p_{t}(\phi^{\kappa}_{t}(\ak)_{k}^{3}(\ak)_{k}^{j}\p_{j}\phi^{\kappa}-\nn_{k}\phi^{\kappa} (\ak)_{j}^{3}\nn_{j}\nn_{k}\phi^{\kappa}-(\ak)_{k}^{3}v_{t}^{k})\quad\text{on $\Gamma$},\\
    q_{t}&=-\sigma \p_{t}(\frac{\sqrt{h}}{\sqrt{h_{\kappa}}}\lap_{h} (\eta)\cdot n)+\kappa\frac{1}{\sqrt{h_{\kappa}}}(1-\bar{\lap})(v\cdot n_{\kappa}))\quad \textit{on $\Gamma_{1}$}.
\end{aligned}
\end{align}
Using the estimate of order, it is not difficult to find
\begin{align}
    \|- \p_{i}(\p_{t}g_{\kappa}^{ij}\p_{j}q)\|_{1.5}\lesssim \|\p_{t}g_{\kappa}^{ij}\|_{2.5}\|\p q\|_{2.5}\lesssim\PP(E^{\kappa}),
\end{align}
\begin{align}
    \|\p_{t}(-\p_{i}(\phi^{\kappa}_{t}g_{\kappa}^{ij}\p_{j}\phi^{\kappa})-\p_{i}((\ak)_{j}^{i}(\nn_{j}\nn_{k}\phi^{\kappa})\nn_{k}\phi^{\kappa})-\p_{i}((\ak)_{k}^{i}v_{t}^{k}))\|_{1.5}\lesssim\PP(E^{\kappa}),
\end{align}

\begin{align}
    \begin{aligned}
       & \|- \p_{t}((\ak)_{k}^{3}(\ak)_{k}^{i})\p_{i}q-\p_{t}(\phi^{\kappa}_{t}(\ak)_{k}^{3}(\ak)_{k}^{j}\p_{j}\phi^{\kappa}-\nn_{k}\phi^{\kappa} (\ak)_{j}^{3}\nn_{j}\nn_{k}\phi^{\kappa}-(\ak)_{k}^{3}v_{t}^{k})\|_{H^{2}(\p\Omega)}\\\lesssim &\|\p q\|_{2.5}\times\|\p_{t}(\ak\ak)\|_{2.5}+\PP(E^{\kappa})\\\lesssim &\PP(E^{\kappa}).
    \end{aligned}
\end{align}
\begin{align}
    \|- \p_{t}g_{\kappa}^{ij}\p_{j}q-\p_{t}((\phi^{\kappa}_{t}g_{\kappa}^{ij}\p_{j}\phi^{\kappa})-\p_{i}((\ak)_{j}^{i}(\nn_{j}\nn_{k}\phi^{\kappa})\nn_{k}\phi^{\kappa})-\p_{i}((\ak)_{k}^{i}v_{t}^{k}))\|_{L^{2}(\Omega)}\lesssim\PP(E^{\kappa}),
\end{align}
\begin{align}
\begin{aligned}
      \|q_{t}\|_{L^{2}(\Gamma_{1})}&\lesssim \| -\sigma \p_{t}(\frac{\sqrt{h}}{\sqrt{h_{\kappa}}}\lap_{h} (\eta)\cdot n)+\kappa\frac{1}{\sqrt{h_{\kappa}}}(1-\bar{\lap})(v\cdot  n _{\kappa})\|_{L^{2}(\Gamma_{1})})\|_{L^{2}(\Gamma_{1})}\\&\lesssim\PP(\|v\|_{H^{2}(\Gamma_{1})}\|\eta\|_{H^{2}(\Gamma_{1})},\|\p\eta\|_{L^{\infty}(\Gamma_{1})},\|\p v\|_{L^{\infty}(\Gamma_{1})})+\kappa\|\p_{t}(v\cdot n )\|_{H^{2}(\p\Omega)}\times\|\frac{1}{\sqrt{h_{\kappa}}}\|_{L^{\infty}(\Gamma_{1})}\\&\lesssim\PP(E^{\kappa}).
\end{aligned}
\end{align}
After a similar argument, it gives us
\begin{align}
    \|q_{t}\|_{3.5}\lesssim \PP(E^{\kappa}).
\end{align}
Differentiating \eqref{elliforqt} on time, we can find that
\begin{align}\label{elliforqtt}
\begin{aligned}
   \p_{i}(g_{\kappa}^{ij}\p_{j}q_{tt})&=-\p_{i}(\p_{t}g_{\kappa}^{ij}\p_{j}q_{t})- \p_{i}\p_{t}(\p_{t}g_{\kappa}^{ij}\p_{j}q)-\p_{t}^{2}(\p_{i}(\phi^{\kappa}_{t}g_{\kappa}^{ij}\p_{j}\phi^{\kappa})\\&-\p_{i}((\ak)_{j}^{i}(\nn_{j}\nn_{k}\phi^{\kappa})\nn_{k}\phi^{\kappa})-\p_{i}((\ak)_{k}^{i}v_{t}^{k}))\quad \text{in $\Omega$,}\\
    (\ak)_{k}^{3}(\ak)_{k}^{i}\p_{i}q_{tt}&=- \p_{t}((\ak)_{k}^{3}(\ak)_{k}^{i})\p_{i}q_{t}- \p_{t}(\p_{t}((\ak)_{k}^{3}(\ak)_{k}^{i})\p_{i}q)\\&-\p_{t}^{2}(\phi^{\kappa}_{t}(\ak)_{k}^{3}(\ak)_{k}^{j}\p_{j}\phi^{\kappa}-\nn_{k}\phi^{\kappa} (\ak)_{j}^{3}\nn_{j}\nn_{k}\phi^{\kappa}-(\ak)_{k}^{3}v_{t}^{k})\quad\text{on $\Gamma$},\\
    q_{tt}&=-\sigma \p_{t}^{2}(\frac{\sqrt{h}}{\sqrt{h_{\kappa}}}\lap_{h} (\eta)\cdot  n )+\kappa\frac{1}{\sqrt{h_{\kappa}}}(1-\bar{\lap})(v\cdot  n _{\kappa})\quad \textit{on $\Gamma_{1}$}.
\end{aligned}
\end{align}
Applying a similar argument
\begin{align}
    \|q_{tt}\|_{2.5}\lesssim \PP(E^{\kappa}).
\end{align}

Differentiating \eqref{elliforqtt} on time, we have 
\begin{align}\label{elliforqttt}
\begin{aligned}
   \p_{i}(g_{\kappa}^{ij}\p_{j}q_{ttt})&=- \p_{i}(\p_{t}g_{\kappa}^{ij}\p_{j}q_{tt})-\p_{i}\p_{t}(\p_{t}g_{\kappa}^{ij}\p_{j}q_{t})- \p_{i}\p_{t}^{2}(\p_{t}g_{\kappa}^{ij}\p_{j}q)\\&-\p_{t}^{3}(\p_{i}(\phi^{\kappa}_{t}g_{\kappa}^{ij}\p_{j}\phi^{\kappa})-\p_{i}((\ak)_{j}^{i}(\nn_{j}\nn_{k}\phi^{\kappa})\nn_{k}\phi^{\kappa})-\p_{i}((\ak)_{k}^{i}v_{t}^{k}))\quad \text{in $\Omega$,}\\
    (\ak)_{k}^{3}(\ak)_{k}^{i}\p_{i}q_{ttt}&=- \p_{t}((\ak)_{k}^{3}(\ak)_{k}^{i})\p_{i}q_{tt}- \p_{t}(\p_{t}((\ak)_{k}^{3}(\ak)_{k}^{i})\p_{i}q_{t})- \p_{t}^{2}(\p_{t}((\ak)_{k}^{3}(\ak)_{k}^{i})\p_{i}q)\\&-\p_{t}^{3}(\phi^{\kappa}_{t}(\ak)_{k}^{3}(\ak)_{k}^{j}\p_{j}\phi^{\kappa}-\nn_{k}\phi^{\kappa} (\ak)_{j}^{3}\nn_{j}\nn_{k}\phi^{\kappa}-(\ak)_{k}^{3}v_{t}^{k})\quad\text{on $\Gamma$}.
\end{aligned}
\end{align}
Observe that
\begin{align}
    g-\pi\cdot N=0,
\end{align}
\begin{align}
\begin{aligned}
      \|\pi\|_{L^{2}(\Omega)}&\lesssim  \|\p_{t}g_{\kappa}^{ij}\p_{j}q_{tt}\|_{L^{2}(\Omega)}+\|\p_{t}(\p_{t}g_{\kappa}^{ij}\p_{j}q_{t})\|_{L^{2}(\Omega)}+ \|\p_{t}^{2}(\p_{t}g_{\kappa}^{ij}\p_{j}q)\|_{L^{2}(\Omega)}\\&+\|\p_{t}^{3}(\phi^{\kappa}_{t}g_{\kappa}^{ij}\p_{j}\phi^{\kappa})\|_{L^{2}(\Omega)}+\|\p_{t}^{3}((\ak)_{j}^{i}(\nn_{j}\nn_{k}\phi^{\kappa})\nn_{k}\phi^{\kappa})\|_{L^{2}(\Omega)}+\|\p_{t}^{3}((\ak)_{k}^{i}v_{t}^{k}))\|_{L^{2}(\Omega)}\\&\lesssim \|\p_{t}g_{\kappa}^{ij}\|_{L^{\infty}(\Omega)}\|\p_{j}q_{tt}\|_{L^{2}(\Omega)}+\|(\ak)_{k}^{i}\|_{L^{\infty}(\Omega)}\|v_{tttt}\|_{L^{2}(\Omega)}+\|\p_{t}(\ak)\|_{L^{\infty}(\Omega)}\|\p_{t}^{3}v\|_{L^{2}(\Omega)}\\&+\PP(E_{\kappa}^{\prime})+\PP_{0}+T\PP(E^{\kappa})\\&\lesssim\PP(E^{\kappa}).
\end{aligned} 
\end{align}
Now, we come back to the Euler equation:

Taking $\p_{t}^{3}$
\begin{align}   \p_{t}^{4}v_{i}+\nn_{i}q_{ttt}+\sum_{k=1}^{3}C_{3}^{k}\p_{t}^{k}(\ak)_{i}^{j}\p_{t}^{3-k}\p_{j}q=-\p_{t}^{3}(\phi^{\kappa}_{t}\nn_{i}\phi^{\kappa})-\p_{t}^{3}(\nn_{i}\nn_{k}\phi^{\kappa}\nn_{k}\phi^{\kappa})
\end{align}
Now, testing $x_{3}(\ak)_{i}^{3}\p_{t}^{3}q$, ($x_{3}$ itself is a periodic function)
\begin{align}
    \begin{aligned}
        \int_{\Omega}\frac{1}{2}x_{3}g_{\kappa}^{3j}\p_{j}(|q_{ttt}|^{2})&=\int_{\Omega}-\p_{t}^{4}v_{i}x_{3}(\ak)_{i}^{3}\p_{t}^{3}q-\int_{\Omega}\sum_{k=1}^{3}C_{3}^{k}\p_{t}^{k}(\ak)_{i}^{j}\p_{t}^{3-k}\p_{j}qx_{3}(\ak)_{i}^{3}\p_{t}^{3}q\\&-\int_{\Omega}\p_{t}^{3}(\phi^{\kappa}_{t}\nn_{i}\phi^{\kappa})x_{3}(\ak)_{i}^{3}\p_{t}^{3}q-\int_{\Omega}\p_{t}^{3}(\nn_{i}\nn_{k}\phi^{\kappa}\nn_{k}\phi^{\kappa})x_{3}(\ak)_{i}^{3}\p_{t}^{3}q
    \end{aligned}
\end{align}
Integrating by part, we will have the desired boundary expression:
\begin{align}
    \begin{aligned}
  \int_{\Gamma_{1}}\frac{1}{2}g_{\kappa}^{33}|q_{ttt}|^{2}&=\int_{\Omega}\frac{1}{2}\p_{j}(x_{3}g_{\kappa}^{3j})q_{ttt}+\int_{\Omega}-\p_{t}^{4}v_{i}x_{3}(\ak)_{i}^{3}\p_{t}^{3}q-\int_{\Omega}\sum_{k=1}^{3}C_{3}^{k}\p_{t}^{k}(\ak)_{i}^{j}\p_{t}^{3-k}\p_{j}qx_{3}(\ak)_{i}^{3}\p_{t}^{3}q\\&-\int_{\Omega}\p_{t}^{3}(\phi^{\kappa}_{t}\nn_{i}\phi^{\kappa})x_{3}(\ak)_{i}^{3}\p_{t}^{3}q-\int_{\Omega}\p_{t}^{3}(\nn_{i}\nn_{k}\phi^{\kappa}\nn_{k}\phi^{\kappa})x_{3}(\ak)_{i}^{3}\p_{t}^{3}q
    \end{aligned}
\end{align}
Now, following from our assumption on $g_{\kappa}^{33}$ and Young's inequality, it is not difficult to find
\begin{align}\label{estimateforqtttbdy}
\begin{aligned}
     \|q_{ttt}\|_{L^{2}(\Gamma_{1})}^{2}&\lesssim\|\p_{j}(x_{3}g_{\kappa}^{3j})\|_{L^{2}(\Omega)}\|q_{ttt}\|_{L^{2}(\Omega)}+\|\p_{t}^{4}v_{i}\|_{L^{2}(\Omega)}\|x_{3}(\ak)_{i}^{3}\|_{L^{\infty}(\Omega)}\|\p_{t}^{3}q\|_{L^{2}(\Omega)}\\&+\|\p_{t}(\ak)_{i}^{j}x_{3}(\ak)_{i}^{3}\|_{L^{\infty}(\Omega)}\|\p_{t}^{2}\p_{j}q\|_{L^{2}(\Omega)}\|\p_{t}^{3}q\|_{L^{2}(\Omega)}\\&+\|\p_{t}^{2}(\ak)_{i}^{j}x_{3}(\ak)_{i}^{3}\|_{L^{\infty}(\Omega)}\|\p_{t}\p_{j}q\|_{L^{2}(\Omega)}\|\p_{t}^{3}q\|_{L^{2}(\Omega)}\\&+\|\p_{t}^{3}(\ak)_{i}^{j}\|_{L^{3}(\Omega)}\|x_{3}(\ak)_{i}^{3}\p_{j}q\|_{L^{6}(\Omega)}\|\p_{t}^{3}q\|_{L^{2}(\Omega)}+\|\p_{t}^{3}(\phi^{\kappa}_{t}\nn_{i}\phi^{\kappa})\|_{L^{2}(\Omega)}\|x_{3}(\ak)_{i}^{3}\|_{L^{\infty}(\Omega)}\|\p_{t}^{3}q\|_{L^{2}(\Omega)}\\&+\|\p_{t}^{3}(\nn_{i}\nn_{k}\phi^{\kappa}\nn_{k}\phi^{\kappa})\|_{L^{2}(\Omega)}\|x_{3}(\ak)_{i}^{3}\|_{L^{\infty}(\Omega)}\|\p_{t}^{3}q\|_{L^{2}(\Omega)}\\&\lesssim\delta\|q_{ttt}\|_{L^{2}(\Omega)}^{2} +C(\delta)(\PP(E_{\kappa}^{\prime})+\PP_{0}+T\PP(E^{\kappa}))\\&+C(\delta)(\|\p_{t}^{4}v_{i}\|_{L^{2}(\Omega)}^{2}\|x_{3}(\ak)_{i}^{3}\|_{L^{\infty}(\Omega)}^{2}+\|\p_{t}(\ak)_{i}^{j}x_{3}(\ak)_{i}^{3}\|_{L^{\infty}(\Omega)}^{2}\|\p_{t}^{2}\p_{j}q\|_{L^{2}(\Omega)}^{2})\\&\lesssim\delta\|q_{ttt}\|_{L^{2}(\Omega)}^{2} +C(\delta)\PP(E^{\kappa}).
\end{aligned}  
\end{align}
Thus, 
\begin{align}
\begin{aligned}
     \|q_{ttt}\|_{H^{1}(\Omega)}^{2}&\lesssim\|\pi\|_{L^{2}(\Omega)}^{2}+\|q_{ttt}\|_{L^{2}(\Gamma_{1})}^{2}\\&\lesssim\delta\|q_{ttt}\|_{L^{2}(\Omega)}^{2} +C(\delta)\PP(E^{\kappa})
\end{aligned}   
\end{align}
Choosing appropriate small $\delta$, we have
\begin{align}
     \|q_{ttt}\|_{H^{1}(\Omega)}\lesssim\PP(E^{\kappa}).
\end{align}
Now, thanks to $\|q_{tt}\|_{1}\lesssim\PP_{0}+T\PP(E^{\kappa})$, and $\|\p_{t}^{3}v\|_{L^{2}(\Omega)}\lesssim \PP_{0}+T\PP(E^{\kappa})$, we rewrite
\begin{align}
\begin{aligned}
     \|q_{ttt}\|_{H^{1}(\Omega)}^{2}&\lesssim\|\pi\|_{L^{2}(\Omega)}^{2}+\|q_{ttt}\|_{L^{2}(\Gamma_{1})}^{2}\\&\lesssim \|(\ak)_{k}^{i}\|_{L^{\infty}(\Omega)}^{2}\|v_{tttt}\|_{L^{2}(\Omega)}^{2}+\delta\|q_{ttt}\|_{L^{2}(\Omega)}^{2}+C(\delta)\|\p_{t}^{4}v_{i}\|_{L^{2}(\Omega)}^{2}\|x_{3}(\ak)_{i}^{3}\|_{L^{\infty}(\Omega)}^{2} \\&+C(\delta)(\PP(E_{\kappa}^{\prime})+\PP_{0}+T\PP(E^{\kappa}))
\end{aligned} 
\end{align}
Choosing small $\delta$ again, we have
\begin{align}
     \|q_{ttt}\|_{H^{1}(\Omega)}^{2}\lesssim(\PP_{0}+T\PP(E^{\kappa}))\|v_{tttt}\|_{L^{2}(\Omega)}^{2}+\PP(E_{\kappa}^{\prime})+\PP_{0}+T\PP(E^{\kappa})
\end{align}

We also need the following lemma to control the term $n-n^{\kappa}$.
\lem
\begin{align}
    \sup_{t}(\|n-n^{\kappa}\|_{L^{\infty}(\Gamma_{1})}+\|\pb n-\pb n^{\kappa}\|_{L^{\infty}(\Gamma_{1})}+\|\p_{t}\pb(n-n^{\kappa})\|_{L^{\infty}(\Gamma_{1})})\lesssim \sqrt{\kappa}\PP(E^{\kappa}).
\end{align}
\begin{align}
    \sup_{t}\sum_{k=2}^{3}\|\p_{t}^{k}(n-n^{\kappa})\|_{L^{2}(\Gamma_{1})}\lesssim \kappa\PP(E^{\kappa}),
\end{align}
\begin{align}
    \sup_{t}\|\p_{t}^{4}(n-n^{\kappa})\|_{L^{2}(\Gamma_{1})}\lesssim\kappa\PP(E^{\kappa})+\kappa\PP(E^{\kappa})\|\p_{t}^{3}v\|_{H^{2}(\Gamma_{1})}.
\end{align}
\pf
Recall
\begin{align*}
    J(A)_{i}^{3}&=(\p_{1}\eta_{2}\p_{2}\eta_{3}-\p_{2}\eta_{2}\p_{1}\eta_{3},\p_{2}\eta_{1}\p_{1}\eta_{3}-\p_{1}\eta_{1}\p_{2}\eta_{3},\p_{1}\eta_{1}\p_{2}\eta_{2}-\p_{2}\eta_{1}\p_{1}\eta_{2}),\\ J_{\kappa}(\ak)_{i}^{3}&=(\p_{1}\eta_{2}^{\kappa}\p_{2}\eta_{3}^{\kappa}-\p_{2}\eta_{2}^{\kappa}\p_{1}\eta_{3}^{\kappa},\p_{2}\eta_{1}^{\kappa}\p_{1}\eta_{3}^{\kappa}-\p_{1}\eta_{1}^{\kappa}\p_{2}\eta_{3}^{\kappa},\p_{1}\eta_{1}^{\kappa}\p_{2}\eta_{2}^{\kappa}-\p_{2}\eta_{1}^{\kappa}\p_{1}\eta_{2}^{\kappa}).
\end{align*}
Therefore, we can find a rational function $F$ such that $J(A)_{i}^{3}=F(\p\eta)$ and $J_{\kappa}(\ak)_{i}^{3}=F(\p\eta^{\kappa})$. Define $G(x):\R^{3}-\{0\}\to\R^{3}$ by $G(x)=\frac{x}{|x|}$. Then $n=G(JA_{i}^{3})$ and $n^{\kappa}=G(J_{\kappa}(\ak)_{i}^{3})$.
Therefore,
\begin{align}
\begin{aligned}
    G(JA_{i}^{3})- G(J_{\kappa}(\ak)_{i}^{3})&= \int_{0}^{1}\p_{s}G(J_{\kappa}(\ak)_{i}^{3}+s(J(A)_{i}^{3}-J_{\kappa}(\ak)_{i}^{3}))ds\\&=\int_{0}^{1}DG(J_{\kappa}(\ak)_{i}^{3}+s(J(A)_{i}^{3}-J_{\kappa}(\ak)_{i}^{3}))ds\times (J(A)_{i}^{3}-J_{\kappa}(\ak)_{i}^{3}))\\&=\int_{0}^{1}DG(J(A_{i}^{3})+s(J_{\kappa}(\ak)_{i}^{3}-J(A)_{i}^{3}))ds\times \int_{0}^{1}DF(\pb\eta^{\kappa}+h(\pb\eta-\pb\eta^{\kappa}))dh\times (\pb\eta-\pb\eta^{\kappa}).
\end{aligned}   
\end{align}
Thus, the boundary analysis using the property of the horizontal mollifier yields the desired estimate.

We now present the following list of boundary estimates:
\lem
\begin{align}
    \sup_{t}\sum_{k=0}^{4}\|\p_{t}^{k}n\|_{H^{4-k}(\Gamma_{1})}+\|\p_{t}^{k}n^{\kappa}\|_{H^{4-k}(\Gamma_{1})}+\|\p_{t}^{k}J_{\kappa}\|_{H^{4-k}(\Gamma_{1})}\lesssim\PP(E^{\kappa}),\label{estimatefornongamma1}\\\sup_{t}\|\pb\p_{t}^{4}n^{\kappa}\|_{L^{2}(\Gamma_{1})}+\|\pb\p_{t}^{4}n\|_{L^{2}(\Gamma_{1})}\lesssim\PP(E^{\kappa})+(\PP_{0}+T\PP(E^{\kappa}))\|\p_{t}^{3}v\|_{H^{2}(\Gamma_{1})},\label{estimatefornongamma2}\\\sup_{t}\|\pb^{2}\p_{t}^{3}n^{\kappa}\|_{L^{2}(\Gamma_{1})}+\|\pb^{2}\p_{t}^{3}n\|_{L^{2}(\Gamma_{1})}\lesssim\PP(E^{\kappa})+(\PP_{0}+T\PP(E^{\kappa}))\|\p_{t}^{2}v\|_{H^{3}(\Gamma_{1})},\label{estimatefornongamma3}\\\sup_{t}\| n^{\kappa}\|_{H^{5}(\Gamma_{1})}+\| n\|_{H^{5}(\Gamma_{1})}\lesssim(\PP_{0}+\PP(E^{\kappa}))\|\eta\|_{H^{6}(\Gamma_{1})}+\PP(E^{\kappa})\label{estimateforpb5nkappa}.
\end{align}
\pf Notice that $n\sim \p\eta$. We omit the details.

Now, as a consequence, we rewrite the estimate for $\|q_{ttt}\|_{1}$ as
\cor
\begin{align}
    \sup_{t}\|q_{ttt}\|_{H^{1}(\Omega)}\lesssim\PP(E^{\kappa})+\kappa\|\p_{t}^{3}v\|_{H^{2}(\Gamma_{1})}\times(\PP_{0}+T\PP(E^{\kappa}))+ \kappa\|\p_{t}^{2}v\|_{H^{3}(\Gamma_{1})}\times(\PP_{0}+T\PP(E^{\kappa})).
\end{align}

The following two lemmas provide the error analysis for the boundary coercivity term.
\lem
\begin{align} \label{boundaryapprox1} \sup_{t}\|\sqrt{\kappa}\pb^{5}\eta^{3}-\sqrt{\kappa}\pb^{5}\eta\cdot n_{\kappa})\|_{H^{1}(\Gamma_{1})}^{2}+\sum_{k=0}^{4}\int_{0}^{t}\|\sqrt{\kappa}\p_{t}^{k}\pb^{4-k}v^{3}-\p_{t}^{k}\pb^{4-k}v\cdot n_{\kappa}\|_{H^{1}(\Gamma_{1})}^{2}\lesssim T \PP(E^{\kappa}).
\end{align}
\pf 
The proof follows from the definition of $n_{\kappa}$.
\begin{align}
    \begin{aligned}
      \|\sqrt{\kappa}\pb^{5}\eta^{3}\cdot \int_{0}^{t}\p_{t}n_{\kappa})\|_{H^{1}(\Gamma_{1})} &\lesssim \|\sqrt{\kappa}\pb^{5}\eta\|_{H^{1}(\Gamma_{1})}T\|\p_{t}n_{\kappa}\|_{H^{1.5}(\Gamma_{1})}\\&\lesssim T\PP(E^{\kappa})
    \end{aligned}
\end{align}
\begin{align}
\begin{aligned}
    \sum_{k=0}^{4}\int_{0}^{t}\|\sqrt{\kappa}(\p_{t}^{k}\pb^{4-k}v^{3}-\p_{t}^{k}\pb^{4-k}v\cdot n_{\kappa})\|_{H^{1}(\Gamma_{1})}^{2}&=\sum_{k=0}^{4}\int_{0}^{t}\|\sqrt{\kappa}\p_{t}^{k}\pb^{4-k}v\cdot \int_{0}^{t^{\prime}}\p_{t}n_{\kappa}\|_{H^{1}(\Gamma_{1})}^{2}\\&\lesssim\sum_{k=0}^{4}\int_{0}^{t}(\|\sqrt{\kappa}\p_{t}^{k}\pb^{4-k}v\|_{H^{1}(\Gamma_{1})}^{2} (\int_{0}^{t^{\prime}}\|\p_{t}n_{\kappa}\|_{H^{1.5}(\Gamma_{1})})^{2})\\&\lesssim\sum_{k=0}^{4}\int_{0}^{t}(\|\sqrt{\kappa}\p_{t}^{k}\pb^{4-k}v\|_{H^{1}(\Gamma_{1})}^{2} T\PP(E^{\kappa})\\&\lesssim T\PP(E^{\kappa}).
\end{aligned}    
\end{align}
\lem
\begin{align}\label{boundaryapproxi2}    \sup_{t}\|\pb^{6}\eta^{\lambda}\Pi_{\lambda}^{3}-\pb^{6}\eta^{3}\|_{L^{2}(\Gamma_{1})}^{2}+\sum_{k=0}^{4}\|\p_{t}^{k}\pb^{5-k}\eta^{\lambda}\Pi_{\lambda}^{3}-\p_{t}^{k}\pb^{5-k}\eta^{3}\|_{L^{2}(\Gamma_{1})}^{2}&\lesssim T\PP(E^{\kappa}).
\end{align}

\pf
Notice that when $t=0$, $\pb\eta^{3}=0$ on $\Gamma_{1}$.
\begin{align}
\begin{aligned}
     \sum_{k=0}^{4}\|\p_{t}^{k}\pb^{5-k}\eta^{\lambda}\Pi_{\lambda}^{3}-\p_{t}^{k}\pb^{5-k}\eta^{3}\|_{L^{2}(\Gamma_{1})}&=\sum_{k=0}^{4}\|\p_{t}^{k}\pb^{5-k}\eta^{\lambda}(h^{\alpha\beta}\p_{\alpha}\eta^{3}\p_{\beta}\eta^{l}\delta_{\lambda l})\|_{L^{2}(\Gamma_{1})}\\&=\sum_{k=0}^{4}\|\p_{t}^{k}\pb^{5-k}\eta^{\lambda}(h^{\alpha\beta}(\int_{0}^{t}\p_{\alpha}v^{3})\p_{\beta}\eta^{l}\delta_{\lambda l})\|_{L^{2}(\Gamma_{1})}\\&\lesssim\sum_{k=0}^{4}\|\p_{t}^{k}\eta\|_{H^{5-k}(\Gamma_{1})}\|h^{\alpha\beta}\p_{\beta}\eta^{l}\|_{H^{1.5}(\Gamma_{1})}T\|\pb v\|_{H^{2}(\Gamma_{1})}\\&\lesssim T\PP(E^{\kappa}).
\end{aligned}
\end{align}

Now, we investigate the boundary estimate of $\eta$. The estimate of $I$, $III$, $IV$ is identical to \cite{coutand2007well}. Our analysis aims to explicitly verify the dependence on $E^{\prime}_{\kappa}$. Readers may skip this section and proceed directly to $II$ for the modified part.
\lem[The estimate for the fourth time-differentiated $\kappa$-problem]
\begin{align}\label{fourtimeest}
\begin{aligned}
     &\int_{0}^{t}\int_{\Omega}J_{\kappa}|\p_{t}^{5}\phi|^{2}+\frac{1}{2}\int_{\Omega}J_{\kappa}g_{\kappa}^{ij}\p_{t}^{4}\p_{j}\phi\p_{t}^{4}\p_{i}\phi+\int_{\Omega}J_{k}|\p_{t}^{4}v|^{2}+\sigma\int_{\Gamma_{1}}\sqrt{h}h^{ij}\p_{j}\p_{t}^{3}v^{\lambda}\Pi^{\mu}_{\lambda}\p_{i}\p_{t}^{3}v^{\alpha}\Pi_{\alpha}^{\mu}+\int_{0}^{t}\|\sqrt{\kappa}\p_{t}^{4}v\cdot n^{\kappa}\|_{H^{1}(\Gamma_{1})}^{2}\\&\lesssim\PP(E_{\kappa}^{\prime})+\delta\|\p_{t}^{4}v\|_{L^{2}(\Omega)}^{2}+\delta \|\p_{t}^{3}v\|_{H^{1.5}(\Omega)}^{2}+\delta\int_{0}^{t}\|\sqrt{\kappa}\p_{t}^{4}v\|_{H^{1.5}(\Omega)}^{2}
+\delta\|\p_{t}^{4}\phi\|_{1}^{2}+\delta\int_{0}^{t}\|\p_{t}^{5}\phi\|^{2}+C(\delta)(\PP_{0}+T\PP(E^{\kappa})).
\end{aligned}
\end{align}
\pf
Multiplying $J_{\kappa}$ on \eqref{vkappat}, and taking $\p_{t}^{4}$ then testing with $\p_{t}^{4} v$, we can write
\begin{align}\label{fourthtimeestformu1}
    \underbrace{\int_{0}^{t}\int_{\Omega}\p_{t}^{4}(J_{\kappa} \p_{t}v^{i})\p_{t}^{4}v^{i}}_{I}+\int_{0}^{t}\int_{\Omega}\p_{t}^{4}(J_{\kappa}\nn_{i} q)\p_{t}^{4}v^{i}=\underbrace{\int_{0}^{t}\int_{\Omega}-\p_{t}^{4}(J_{\kappa}\phi^{\kappa}_{t}\nn_{i}\phi^{\kappa})\p_{t}^{4}v^{i}+\int_{0}^{t}\int_{\Omega}-\p_{t}^{4}(J_{\kappa}\nn_{j}\nn_{i}\phi^{\kappa}\nn_{j}\phi^{\kappa})\p_{t}^{4}v^{i}}_{II}.
\end{align}
We calculate the commutator term by term.
\begin{align}
\begin{aligned}
     \int_{0}^{t}\int_{\Omega}\p_{t}^{4}(J_{\kappa} \p_{t}v^{i})\p_{t}^{4}v^{i}&\sim \int_{0}^{t}\int_{\Omega}J_{\kappa}\frac{1}{2}\p_{t}|\p_{t}^{4}v|^{2}+\sum_{\textit{$i_{1}+i_{2}=4$, $i_{1}\geq 1$}}\int_{0}^{t}\int_{\Omega}\p_{t}^{i_{1}}J_{\kappa}\p_{t}^{i_{2}}\p_{t}v^{i}\p_{t}^{4} v^{i}\\&\sim\PP_{0}+\int_{\Omega}J_{k}|\p_{t}^{4}v|^{2}+\int_{0}^{t}\int_{\Omega}\p_{t}J_{\kappa}|\p_{t}^{4}v|^{2}+\sum_{\textit{$i_{1}+i_{2}=4$, $i_{1}\geq 1$}}\int_{0}^{t}\int_{\Omega}\p_{t}^{i_{1}}J_{\kappa}\p_{t}^{i_{2}}\p_{t}v^{i}\p_{t}^{4} v^{i}
\end{aligned}    
\end{align}
From direct calculation
\begin{align}
\begin{aligned}
     \sum_{\textit{$i_{1}+i_{2}=4$, $i_{1}\geq 1$}}\int_{0}^{t}\int_{\Omega}|\p_{t}^{i_{1}}J_{\kappa}\p_{t}^{i_{2}}\p_{t}v^{i}\p_{t}^{4} v^{i}|&\lesssim\int_{0}^{t}\|\p_{t}J_{\kappa}\|_{L^{\infty}(\Omega)}\|\p_{t}^{4}v\|_{L^{2}(\Omega)}^{2}\\&+\sum_{i=2}^{4}\int_{0}^{t}\|\p_{t}^{i}J_{\kappa}\|_{L^{3}(\Omega)}\|\p_{t}^{5-i}v\|_{L^{6}(\Omega)}\|\p_{t}^{4}v\|_{L^{2}(\Omega)}\\&\lesssim\int_{0}^{t}\|\p_{t}J_{\kappa}\|_{2}\|\p_{t}^{4}v\|_{L^{2}(\Omega)}^{2}+\sum_{i=2}^{4}\int_{0}^{t}\|\p_{t}^{i}J_{\kappa}\|_{0.5}\|\p_{t}^{5-i}v\|_{1}\|\p_{t}^{4}v\|_{L^{2}(\Omega)}\\&\lesssim T\PP(E^{\kappa}).
\end{aligned} 
\end{align}
Therefore, we can write
\begin{align}\label{estimateforI}
     I=\int_{\Omega}J_{k}|\p_{t}^{4}v|^{2}+\mathcal{R},
\end{align}
where $|\mathcal{R}|\lesssim \PP_{0}+ T\PP(E^{\kappa})$. Using Piola's formula, and integrating by part, we can show that

\begin{align}
   \begin{aligned}
        \int_{0}^{t}\int_{\Omega}\p_{t}^{4}(J_{\kappa}(\ak)_{i}^{k}\p_{k} q)\p_{t}^{4}v^{i}&= \int_{0}^{t}\int_{\Omega}\p_{k}\p_{t}^{4}(J_{\kappa}(\ak)_{i}^{k} q)\p_{t}^{4}v^{i}\\&=\underbrace{\int_{0}^{t}\int_{\Gamma_{1}} \p_{t}^{4}(J_{\kappa}(\ak)_{i}^{3} q)\p_{t}^{4}v^{i}) }_{III}\underbrace{-\int_{0}^{t}\int_{\Omega}\p_{t}^{4}(J_{\kappa}(\ak)_{i}^{k} q)\p_{t}^{4}\p_{k}v^{i}}_{IV}
   \end{aligned}  
\end{align}
We can show that the estimate in \cite{coutand2007well} is still valid for our model. The third term in \eqref{fourthtimeestformu1} reads
\begin{align}
\begin{aligned}
    |\int_{0}^{t}\int_{\Omega}-\p_{t}^{4}(J_{\kappa}\phi^{\kappa}_{t}\nn_{i}\phi^{\kappa})\p_{t}^{4}v^{i}|&\lesssim \sum_{i_{1}+i_{2}+i_{3}=4}\int_{0}^{t}\int_{\Omega}|\p_{t}^{i_{1}}J_{\kappa}\p_{t}^{i_{2}}\phi^{\kappa}_{t}\p_{t}^{i_{3}}\nn_{i}\phi^{\kappa}\p_{t}^{4}v^{i}|\\&\lesssim \int_{0}^{t}\|\p_{t}^{4}J_{\kappa}\|_{L^{2}(\Omega)}\|\phi^{\kappa}_{t}\|_{L^{\infty}(\Omega)}\|\nn\phi^{\kappa}\|_{L^{\infty}(\Omega)}\|\p_{t}^{4}v\|_{L^{2}(\Omega)}\\&+\int_{0}^{t}\|J_{\kappa}\|_{L^{\infty}(\Omega)}\|\p_{t}^{4}\phi^{\kappa}_{t}\|_{L^{2}(\Omega)}\|\nn\phi^{\kappa}\|_{L^{\infty}(\Omega)}\|\p_{t}^{4}v\|_{L^{2}(\Omega)}\\&+\int_{0}^{t}\|J_{\kappa}\|_{L^{\infty}(\Omega)}\|\p_{t}\phi^{\kappa}_{t}\|_{L^{\infty}(\Omega)}\|\p_{t}^{4}\nn\phi^{\kappa}\|_{L^{2}(\Omega)}\|\p_{t}^{4}v\|_{L^{2}(\Omega)}\\&+\sum_{i_{1}+i_{2}=1}\int_{0}^{t}\int_{\Omega}\|\p_{t}^{3}J_{\kappa}\|_{L^{2}(\Omega)}\|\p_{t}^{i_{1}}\phi^{\kappa}_{t}\|_{L^{\infty}(\Omega)}\|\p_{t}^{i_{2}}\nn\phi^{\kappa}\|_{L^{\infty}(\Omega)}\|\p_{t}^{4}v\|_{L^{2}(\Omega)}\\&+\sum_{i_{1}+i_{2}=2}\int_{0}^{t}\int_{\Omega}\|\p_{t}^{2}J_{\kappa}\|_{L^{\infty}(\Omega)}\|\p_{t}^{i_{1}}\phi^{\kappa}_{t}\|_{L^{4}(\Omega)}\|\p_{t}^{i_{2}}\nn\phi^{\kappa}\|_{L^{4}(\Omega)}\|\p_{t}^{4}v\|_{L^{2}(\Omega)}\\&+\sum_{i_{1}+i_{2}=3}\int_{0}^{t}\int_{\Omega}\|\p_{t}J_{\kappa}\|_{L^{\infty}(\Omega)}\|\p_{t}^{i_{1}}\phi^{\kappa}_{t}\|_{L^{4}(\Omega)}\|\p_{t}^{i_{2}}\nn\phi^{\kappa}\|_{L^{4}(\Omega)}\|\p_{t}^{4}v\|_{L^{2}(\Omega)}\\&\lesssim \delta \int_{0}^{t}|\p_{t}^{5}\phi^{\kappa}|^{2}+C(\delta)T\PP(E^{\kappa})
\end{aligned}  
\end{align}
The fourth term in \eqref{fourthtimeestformu1} leads to a cancellation, which we will see later. Now, we adapt the same strategy to deal with the second term in \eqref{fourthtimeestformu1}.
\begin{align}    III=\underbrace{\int_{0}^{t}\int_{\Gamma_{1}}\p_{t}^{4}(-\sigma\sqrt{h}\lap_{h}(\eta)\cdot n_{\kappa} n_{\kappa}^{i})\p_{t}^{4}v^{i}}_{III_{1}}+\underbrace{\int_{0}^{t}\int_{\Omega}\p_{t}^{4}(\kappa\lap_{0}(v\cdot n_{\kappa})n_{\kappa}^{i})\p_{t}^{4}v^{i}}_{III_{2}}.
\end{align}
We rewrite $III_{1}$ as
\begin{align}
\begin{aligned}
     III_{1}&=\underbrace{\int_{0}^{t}\int_{\Gamma_{1}}\p_{t}^{4}(-\sigma\sqrt{h}\lap_{h}(\eta^{i}))\p_{t}^{4}v^{i}}_{III_{11}}+\underbrace{\int_{0}^{t}\int_{\Gamma_{1}}\p_{t}^{4}(-\sigma\sqrt{h}\lap_{h}(\eta)\cdot  (n_{\kappa}-n)n^{i})\p_{t}^{4}v^{i}}_{III_{12}}\\&\underbrace{+\int_{0}^{t}\int_{\Gamma_{1}}\p_{t}^{4}(-\sigma\sqrt{h}\lap_{h}(\eta) \cdot n_{\kappa}(n_{\kappa}^{i}-n^{i}))\p_{t}^{4}v^{i}}_{III_{13}}
\end{aligned}   
\end{align}
Using \eqref{derivativeofmeancur} and integrating by part, we obtain
\begin{align}
\begin{aligned}
      III_{11}&=\int_{0}^{t}\int_{\Gamma_{1}}-\sigma\p_{t}^{3}(\p_{i}(\sqrt{h}h^{ij}(\delta^{\alpha\lambda}-h^{lk}\p_{l}\eta^{\lambda}\p_{k}\eta^{\alpha})\p_{j}\p_{t}\eta^{\lambda}+\sqrt{h}(h^{ij}h^{kl}-h^{ik}h^{jl})\p_{k}\eta^{\lambda}\p_{t}\p_{l}\eta^{\lambda}\p_{j}\eta^{\alpha}))\p_{t}^{4}v^{\alpha}\\&=\sigma\int_{0}^{t}\int_{\Gamma_{1}}\p_{t}^{3}(\sqrt{h}h^{ij}(\delta^{\alpha\lambda}-h^{lk}\p_{l}\eta^{\lambda}\p_{k}\eta^{\alpha})\p_{j}\p_{t}\eta^{\lambda}+\sqrt{h}(h^{ij}h^{kl}-h^{ik}h^{jl})\p_{k}\eta^{\lambda}\p_{t}\p_{l}\eta^{\lambda}\p_{j}\eta^{\alpha})\p_{i}\p_{t}^{4}v^{\alpha}\\&=\underbrace{\sigma\int_{0}^{t}\int_{\Gamma_{1}}\sqrt{h}h^{ij}(\delta^{\alpha\lambda}-h^{lk}\p_{l}\eta^{\lambda}\p_{k}\eta^{\alpha})\p_{j}\p_{t}^{4}\eta^{\lambda}\p_{i}\p_{t}^{4}v^{\alpha}}_{III_{11,1}}\underbrace{+\sigma\int_{0}^{t}\int_{\Gamma_{1}}\sqrt{h}(h^{ij}h^{kl}-h^{ik}h^{jl})\p_{k}\eta^{\lambda}\p_{t}^{4}\p_{l}\eta^{\lambda}\p_{j}\eta^{\alpha}\p_{i}\p_{t}^{4}v^{\alpha}}_{III_{11,2}}\\&\underbrace{+3\sigma\int_{0}^{t}\int_{\Gamma_{1}}\p_{t}(\sqrt{h}h^{ij}(\delta^{\alpha\lambda}-h^{lk}\p_{l}\eta^{\lambda}\p_{k}\eta^{\alpha}))\p_{j}\p_{t}^{3}\eta^{\lambda}\p_{i}\p_{t}^{4}v^{\alpha}+\p_{t}(\sqrt{h}(h^{ij}h^{kl}-h^{ik}h^{jl})\p_{k}\eta^{\lambda}\p_{j}\eta^{\alpha})\p_{t}^{3}\p_{l}\eta^{\lambda}\p_{i}\p_{t}^{4}v^{\alpha}}_{III_{11,3}}\\&\underbrace{+3\sigma\int_{0}^{t}\int_{\Gamma_{1}}\p_{t}^{2}(\sqrt{h}h^{ij}(\delta^{\alpha\lambda}-h^{lk}\p_{l}\eta^{\lambda}\p_{k}\eta^{\alpha}))\p_{j}\p_{t}^{2}\eta^{\lambda}\p_{i}\p_{t}^{4}v^{\alpha}+\p_{t}^{2}(\sqrt{h}(h^{ij}h^{kl}-h^{ik}h^{jl})\p_{k}\eta^{\lambda}\p_{j}\eta^{\alpha})\p_{t}^{2}\p_{l}\eta^{\lambda}\p_{i}\p_{t}^{4}v^{\alpha}}_{III_{11,4}}\\&\underbrace{+\sigma\int_{0}^{t}\int_{\Gamma_{1}}\p_{t}^{3}(\sqrt{h}h^{ij}(\delta^{\alpha\lambda}-h^{lk}\p_{l}\eta^{\lambda}\p_{k}\eta^{\alpha}))\p_{j}\p_{t}\eta^{\lambda}\p_{i}\p_{t}^{4}v^{\alpha}+\p_{t}^{3}(\sqrt{h}(h^{ij}h^{kl}-h^{ik}h^{jl})\p_{k}\eta^{\lambda}\p_{j}\eta^{\alpha})\p_{t}\p_{l}\eta^{\lambda}\p_{i}\p_{t}^{4}v^{\alpha}}_{III_{11,5}}
\end{aligned}  
\end{align}
We estimate the lower order terms $III_{11,3}$ ,$III_{11,4}$, $III_{11,5}$as
\begin{align}
\begin{aligned}
     III_{11,3}+III_{11,4}+III_{11,5}&\sim\sum_{i_{1}+i_{2}=2}\sigma\int_{\Gamma_{1}}\p_{t}^{i_{1}}\p_{t}(\sqrt{h}h^{ij}(\delta^{\alpha\lambda}-h^{lk}\p_{l}\eta^{\lambda}\p_{k}\eta^{\alpha}))\p_{t}^{i_{2}}\p_{j}\p_{t}\eta^{\lambda}\p_{i}\p_{t}^{3}v^{\alpha}|_{0}^{t}\\&+\sum_{i_{1}+i_{2}=2}\sigma\int_{\Gamma_{1}}\p_{t}^{i_{1}}(\sqrt{h}(h^{ij}h^{kl}-h^{ik}h^{jl})\p_{k}\eta^{\lambda}\p_{j}\eta^{\alpha})\p_{t}^{i_{2}}\p_{l}\p_{t}\eta^{\lambda}\p_{i}\p_{t}^{3}v^{\alpha}|_{0}^{t}\\&+\sigma\sum_{i_{1}+i_{2}=3}\int_{0}^{t}\int_{\Gamma_{1}}\p_{t}^{i_{1}}\p_{t}(\sqrt{h}h^{ij}(\delta^{\alpha\lambda}-h^{lk}\p_{l}\eta^{\lambda}\p_{k}\eta^{\alpha}))\p_{t}^{i_{2}}\p_{j}\p_{t}\eta^{\lambda}\p_{i}\p_{t}^{3}v^{\alpha}\\&+\sigma\sum_{i_{1}+i_{2}=3}\int_{0}^{t}\int_{\Gamma_{1}}\p_{t}^{i_{1}}\p_{t}(\sqrt{h}(h^{ij}h^{kl}-h^{ik}h^{jl})\p_{k}\eta^{\lambda}\p_{j}\eta^{\alpha})\p_{t}^{i_{2}}\p_{l}\p_{t}\eta^{\lambda}\p_{i}\p_{t}^{3}v^{\alpha}
\end{aligned}  
\end{align}
Recall that $\sup_{t}(\sum_{k=0}^{4}\|\p_{t}^{k}\pb\eta\|_{H^{4-k}(\p\Omega)}+\sum_{k=0}^{4}\|\p_{t}^{k}h^{ij}\|_{H^{4-k}(\p\Omega)})\lesssim\PP(E^{\kappa})$, and \\$\sup_{t}(\sum_{k=0}^{3}\|\p_{t}^{k}\p\eta\|_{H^{3-k}(\p\Omega)}+\sum_{k=0}^{3}\|\p_{t}^{k}h^{ij}\|_{H^{3-k}(\p\Omega)})\lesssim\PP(E^{\kappa})+T\PP(E^{\kappa})$. Then we can show that
\begin{align}\label{estiforIII11345}
    \begin{aligned}        |III_{11,3}|+|III_{11,4}|+|III_{11,5}|&\lesssim\delta \|\pb\p_{t}^{3}v\|_{L^{2}(\Gamma_{1})}^{2}+ C(\delta)(\PP_{0}+T\PP(E^{\kappa})).
    \end{aligned}
\end{align}
Now, we consider $III_{11,2}$. There exists an antisymmetric structure on $h^{ij}h^{kl}-h^{ik}h^{jl}$. (In \cite{coutand2007well}, they first observe this structure). i.e., altering the sequence of $il$ or $jk$, the sign will change. We can show that \\
1) Rename $i$ and $l$
\begin{align}
    \sigma\int_{0}^{t}\int_{\Gamma_{1}}\sqrt{h}(h^{ij}h^{kl}-h^{ik}h^{jl})\p_{k}\eta^{\lambda}\p_{l}\p_{t}^{3}v^{\lambda}\p_{j}\eta^{\alpha}\p_{i}\p_{t}^{4}v^{\alpha}=-\sigma\int_{0}^{t}\int_{\Gamma_{1}}\sqrt{h}(h^{ij}h^{kl}-h^{ik}h^{jl})\p_{k}\eta^{\lambda}\p_{i}\p_{t}^{3}v^{\lambda}\p_{j}\eta^{\alpha}\p_{l}\p_{t}^{4}v^{\alpha}
\end{align}
2) Rename $\alpha$ and $\lambda$
\begin{align}
   -\sigma\int_{0}^{t}\int_{\Gamma_{1}}\sqrt{h}(h^{ij}h^{kl}-h^{ik}h^{jl})\p_{k}\eta^{\lambda}\p_{i}\p_{t}^{3}v^{\lambda}\p_{j}\eta^{\alpha}\p_{l}\p_{t}^{4}v^{\alpha}=-\sigma\int_{0}^{t}\int_{\Gamma_{1}}\sqrt{h}(h^{ij}h^{kl}-h^{ik}h^{jl})\p_{k}\eta^{\alpha}\p_{i}\p_{t}^{3}v^{\alpha}\p_{j}\eta^{\lambda}\p_{l}\p_{t}^{4}v^{\lambda}
\end{align}
3) Rename $j$ and $k$
\begin{align}
    -\sigma\int_{0}^{t}\int_{\Gamma_{1}}\sqrt{h}(h^{ij}h^{kl}-h^{ik}h^{jl})\p_{k}\eta^{\alpha}\p_{i}\p_{t}^{3}v^{\alpha}\p_{j}\eta^{\lambda}\p_{l}\p_{t}^{4}v^{\lambda}=\sigma\int_{0}^{t}\int_{\Gamma_{1}}\sqrt{h}(h^{ij}h^{kl}-h^{ik}h^{jl})\p_{j}\eta^{\alpha}\p_{i}\p_{t}^{3}v^{\alpha}\p_{k}\eta^{\lambda}\p_{l}\p_{t}^{4}v^{\lambda}
\end{align}
Consequently,
\begin{align}
      \sigma\int_{0}^{t}\int_{\Gamma_{1}}\sqrt{h}(h^{ij}h^{kl}-h^{ik}h^{jl})\p_{k}\eta^{\lambda}\p_{l}\p_{t}^{3}v^{\lambda}\p_{j}\eta^{\alpha}\p_{i}\p_{t}^{4}v^{\alpha}=\sigma\int_{0}^{t}\int_{\Gamma_{1}}\sqrt{h}(h^{ij}h^{kl}-h^{ik}h^{jl})\p_{j}\eta^{\alpha}\p_{i}\p_{t}^{3}v^{\alpha}\p_{k}\eta^{\lambda}\p_{l}\p_{t}^{4}v^{\lambda}.
\end{align}
However, integrating by part allows us to write
\begin{align}
\begin{aligned}
    \sigma\int_{0}^{t}\int_{\Gamma_{1}}\sqrt{h}(h^{ij}h^{kl}-h^{ik}h^{jl})\p_{j}\eta^{\alpha}\p_{i}\p_{t}^{3}v^{\alpha}\p_{k}\eta^{\lambda}\p_{l}\p_{t}^{4}v^{\lambda}&=\sigma \int_{\Gamma_{1}}\sqrt{h}(h^{ij}h^{kl}-h^{ik}h^{jl})\p_{j}\eta^{\alpha}\p_{i}\p_{t}^{3}v^{\alpha}\p_{k}\eta^{\lambda}\p_{l}\p_{t}^{3}v^{\lambda}|_{0}^{t}\\&-\sigma\int_{0}^{t}\int_{\Gamma_{1}}\p_{t}(\sqrt{h}(h^{ij}h^{kl}-h^{ik}h^{jl})\p_{j}\eta^{\alpha}\p_{k}\eta^{\lambda})\p_{i}\p_{t}^{3}v^{\alpha}\p_{l}\p_{t}^{3}v^{\lambda}\\&-\sigma\int_{0}^{t}\int_{\Gamma_{1}}\sqrt{h}(h^{ij}h^{kl}-h^{ik}h^{jl})\p_{j}\eta^{\alpha}\p_{i}\p_{t}^{4}v^{\alpha}\p_{k}\eta^{\lambda}\p_{l}\p_{t}^{3}v^{\lambda}
\end{aligned}   
\end{align}
Now, we have
\begin{align}\label{iiiestidet1}
\begin{aligned}
     III_{11,2}&=\frac{1}{2}\sigma \int_{\Gamma_{1}}\sqrt{h}(h^{ij}h^{kl}-h^{ik}h^{jl})\p_{j}\eta^{\alpha}\p_{i}\p_{t}^{3}v^{\alpha}\p_{k}\eta^{\lambda}\p_{l}\p_{t}^{3}v^{\lambda}|_{0}^{t}\\&-\frac{1}{2}\sigma\int_{0}^{t}\int_{\Gamma_{1}}\p_{t}(\sqrt{h}(h^{ij}h^{kl}-h^{ik}h^{jl})\p_{j}\eta^{\alpha}\p_{k}\eta^{\lambda})\p_{i}\p_{t}^{3}v^{\alpha}\p_{l}\p_{t}^{3}v^{\lambda}
\end{aligned}   
\end{align}
For the first term of \label{iiiestidet1}, integrating by part on $i$ and $l$, it follows that
\begin{align}
\begin{aligned}
     \frac{1}{2}\sigma \int_{\Gamma_{1}}\sqrt{h}(h^{ij}h^{kl}-h^{ik}h^{jl})\p_{j}\eta^{\alpha}\p_{i}\p_{t}^{3}v^{\alpha}\p_{k}\eta^{\lambda}\p_{l}\p_{t}^{3}v^{\lambda}|_{0}^{t}&=-\frac{1}{2}\sigma \int_{\Gamma_{1}}\p_{i}(\sqrt{h}(h^{ij}h^{kl}-h^{ik}h^{jl})\p_{j}\eta^{\alpha}\p_{k}\eta^{\lambda})\p_{t}^{3}v^{\alpha}\p_{l}\p_{t}^{3}v^{\lambda}|_{0}^{t}\\&-\frac{1}{2}\sigma \int_{\Gamma_{1}}\sqrt{h}(h^{ij}h^{kl}-h^{ik}h^{jl})\p_{j}\eta^{\alpha}\p_{k}\eta^{\lambda}\p_{t}^{3}v^{\alpha}\p_{i}\p_{l}\p_{t}^{3}v^{\lambda}|_{0}^{t}\\&=-\frac{1}{2}\sigma \int_{\Gamma_{1}}\p_{i}(\sqrt{h}(h^{ij}h^{kl}-h^{ik}h^{jl})\p_{j}\eta^{\alpha}\p_{k}\eta^{\lambda})\p_{t}^{3}v^{\alpha}\p_{l}\p_{t}^{3}v^{\lambda}|_{0}^{t}\\&+\frac{1}{2}\sigma \int_{\Gamma_{1}}\p_{l}(\sqrt{h}(h^{ij}h^{kl}-h^{ik}h^{jl})\p_{j}\eta^{\alpha}\p_{k}\eta^{\lambda})\p_{t}^{3}v^{\alpha}\p_{i}\p_{t}^{3}v^{\lambda}|_{0}^{t}\\&+\frac{1}{2}\sigma \int_{\Gamma_{1}}\sqrt{h}(h^{ij}h^{kl}-h^{ik}h^{jl})\p_{j}\eta^{\alpha}\p_{k}\eta^{\lambda}\p_{l}\p_{t}^{3}v^{\alpha}\p_{i}\p_{t}^{3}v^{\lambda}|_{0}^{t}
\end{aligned}   
\end{align}
Now, we rename $i$ and $l$ again,
\begin{align}
\begin{aligned}
    \sigma \int_{\Gamma_{1}}\sqrt{h}(h^{ij}h^{kl}-h^{ik}h^{jl})\p_{j}\eta^{\alpha}\p_{i}\p_{t}^{3}v^{\alpha}\p_{k}\eta^{\lambda}\p_{l}\p_{t}^{3}v^{\lambda}|_{0}^{t}&=-\frac{1}{2}\sigma \int_{\Gamma_{1}}\p_{i}(\sqrt{h}(h^{ij}h^{kl}-h^{ik}h^{jl})\p_{j}\eta^{\alpha}\p_{k}\eta^{\lambda})\p_{t}^{3}v^{\alpha}\p_{l}\p_{t}^{3}v^{\lambda}|_{0}^{t}\\&+\frac{1}{2}\sigma \int_{\Gamma_{1}}\p_{l}(\sqrt{h}(h^{ij}h^{kl}-h^{ik}h^{jl})\p_{j}\eta^{\alpha}\p_{k}\eta^{\lambda})\p_{t}^{3}v^{\alpha}\p_{i}\p_{t}^{3}v^{\lambda}|_{0}^{t}
\end{aligned}    
\end{align}
Hence, we can conclude that
\begin{align}\label{iiiestidet1}
\begin{aligned}
     III_{11,2}&=-\frac{1}{4}\sigma \int_{\Gamma_{1}}\p_{i}(\sqrt{h}(h^{ij}h^{kl}-h^{ik}h^{jl})\p_{j}\eta^{\alpha}\p_{k}\eta^{\lambda})\p_{t}^{3}v^{\alpha}\p_{l}\p_{t}^{3}v^{\lambda}|_{0}^{t}\\&+\frac{1}{4}\sigma \int_{\Gamma_{1}}\p_{l}(\sqrt{h}(h^{ij}h^{kl}-h^{ik}h^{jl})\p_{j}\eta^{\alpha}\p_{k}\eta^{\lambda})\p_{t}^{3}v^{\alpha}\p_{i}\p_{t}^{3}v^{\lambda}|_{0}^{t}\\&-\frac{1}{2}\sigma\int_{0}^{t}\int_{\Gamma_{1}}\p_{t}(\sqrt{h}(h^{ij}h^{kl}-h^{ik}h^{jl})\p_{j}\eta^{\alpha}\p_{k}\eta^{\lambda})\p_{i}\p_{t}^{3}v^{\alpha}\p_{l}\p_{t}^{3}v^{\lambda}
\end{aligned}   
\end{align}
We can also have the following estimate
\begin{align}
    \|\pb(\sqrt{h}(h^{ij}h^{kl}-h^{ik}h^{jl})\p_{j}\eta^{\alpha}\p_{k}\eta^{\lambda})\|_{L^{\infty}(\Gamma_{1})}\lesssim \PP_{0}+T\PP(E^{\kappa}).
\end{align}
Now, we use the interpolation argument again
\begin{align}\label{estiforIII112}
\begin{aligned}
     |III_{11,2}|&\lesssim \delta\|\p_{t}^{3}\pb v\|_{L^{2}(\Gamma_{1})}^{2}+C(\delta)\|\p_{t}^{3}v\|_{L^{2}(\Gamma_{1})}^{2}(\PP_{0}+T\PP(E^{\kappa}))\\&\lesssim \delta\|\p_{t}^{3}\pb v\|_{L^{2}(\Gamma_{1})}^{2}+\|\p_{t}^{3}v\|_{H^{1}(\Omega)}^{2}C(\delta)(\PP_{0}+T\PP(E^{\kappa}))\\&\lesssim \delta\|\p_{t}^{3}\pb v\|_{L^{2}(\Gamma_{1})}^{2}+\|\p_{t}^{3}v\|_{L^{2}(\Omega)}^{\frac{2}{3}}\|\p_{t}^{3}v\|_{H^{1.5}(\Omega)}^{\frac{4}{3}}C(\delta)(\PP_{0}+T\PP(E^{\kappa}))\\&\lesssim\delta\|\p_{t}^{3} v\|_{H^{1.5}(\Omega)}^{2}+C(\delta)(\PP_{0}+T\PP(E^{\kappa})) ,
\end{aligned} 
\end{align}
after choosing an appropriate coefficient when we use Young's inequality. The term $III_{11,1}$ will provide us with the desired improved regularity on the boundary.
\begin{align}
    \begin{aligned}       III_{11,1}&=\sigma\int_{0}^{t}\int_{\Gamma_{1}}\sqrt{h}h^{ij}\Pi_{\alpha}^{\lambda}\p_{j}\p_{t}^{4}\eta^{\lambda}\p_{i}\p_{t}^{4}v^{\alpha}\\&=\sigma\int_{0}^{t}\int_{\Gamma_{1}}\sqrt{h}h^{ij}\p_{j}\p_{t}^{4}\eta^{\lambda}\Pi^{\mu}_{\lambda}\p_{i}\p_{t}^{4}v^{\alpha}\Pi_{\alpha}^{\mu}\\&=\sigma\int_{0}^{t}\int_{\Gamma_{1}}\sqrt{h}h^{ij}\p_{t}(\p_{j}\p_{t}^{3}v^{\lambda}\Pi^{\mu}_{\lambda}\p_{i}\p_{t}^{3}v^{\alpha}\Pi_{\alpha}^{\mu})\\&-\sigma\int_{0}^{t}\int_{\Gamma_{1}}\sqrt{h}h^{ij}\p_{j}\p_{t}^{3}v^{\lambda}\p_{i}\p_{t}^{3}v^{\alpha}\p_{t}(\delta^{\alpha\lambda}-h^{kl}\p_{k}\eta^{\alpha}\p_{l}\eta^{\beta})\\&=\sigma\int_{\Gamma_{1}}\sqrt{h}h^{ij}\p_{j}\p_{t}^{3}v^{\lambda}\Pi^{\mu}_{\lambda}\p_{i}\p_{t}^{3}v^{\alpha}\Pi_{\alpha}^{\mu}|_{0}^{t}\\&-\sigma\int_{0}^{t}\int_{\Gamma_{1}}\p_{t}(\sqrt{h}h^{ij})\p_{j}\p_{t}^{3}v^{\lambda}\Pi^{\mu}_{\lambda}\p_{i}\p_{t}^{3}v^{\alpha}\Pi_{\alpha}^{\mu}\\&-\sigma\int_{0}^{t}\int_{\Gamma_{1}}\sqrt{h}h^{ij}\p_{j}\p_{t}^{3}v^{\lambda}\p_{i}\p_{t}^{3}v^{\alpha}\p_{t}(\delta^{\alpha\lambda}-h^{kl}\p_{k}\eta^{\alpha}\p_{l}\eta^{\beta}).
    \end{aligned}
\end{align}
Since
\begin{align*}
\|\p_{t}(\sqrt{h}h^{ij})\|_{L^{\infty}(\Gamma_{1})}\lesssim\PP(E^{\kappa}),\\
   \| \sqrt{h}h^{ij}\p_{t}(\delta^{\alpha\lambda}-h^{kl}\p_{k}\eta^{\alpha}\p_{l}\eta^{\beta})\|_{L^{\infty}(\Gamma_{1})}\lesssim\PP(E^{\kappa}),
\end{align*}
we can write that
\begin{align}
|-\sigma\int_{0}^{t}\int_{\Gamma_{1}}\p_{t}(\sqrt{h}h^{ij})\p_{j}\p_{t}^{3}v^{\lambda}\Pi^{\mu}_{\lambda}\p_{i}\p_{t}^{3}v^{\alpha}\Pi_{\alpha}^{\mu}|&\lesssim \PP_{0}+T\PP(E^{\kappa})\\
   | -\sigma\int_{0}^{t}\int_{\Gamma_{1}}\sqrt{h}h^{ij}\p_{j}\p_{t}^{3}v^{\lambda}\p_{i}\p_{t}^{3}v^{\alpha}\p_{t}(\delta^{\alpha\lambda}-h^{kl}\p_{k}\eta^{\alpha}\p_{l}\eta^{\beta})|&\lesssim\PP_{0}+T\PP(E^{\kappa}).
\end{align}
Now, we have 
\begin{align}\label{estiforIII111}    III_{11,1}=\sigma\int_{\Gamma_{1}}\sqrt{h}h^{ij}\p_{j}\p_{t}^{3}v^{\lambda}\Pi^{\mu}_{\lambda}\p_{i}\p_{t}^{3}v^{\alpha}\Pi_{\alpha}^{\mu}+\mathcal{R},
\end{align}
where $\mathcal{R}\lesssim \PP_{0}+T\PP(E^{\kappa})$.
Therefore, combining \eqref{estiforIII11345}, \eqref{estiforIII112}, and \eqref{estiforIII111}, we have
\begin{align}    III_{11}=\sigma\int_{\Gamma_{1}}\sqrt{h}h^{ij}\p_{j}\p_{t}^{3}v^{\lambda}\Pi^{\mu}_{\lambda}\p_{i}\p_{t}^{3}v^{\alpha}\Pi_{\alpha}^{\mu}+\mathcal{R},
\end{align}
where $\mathcal{R}\lesssim \delta \|\p_{t}^{3}v\|_{1.5}^{2}+C(\delta)(\PP_{0}+T\PP(E^{\kappa}))$. Now we deal with $III_{12}$ and $III_{13}$.
\begin{align}  
\begin{aligned}    III_{12}&\sim\underbrace{\sum_{k=1}^{2}\int_{0}^{t}\int_{\Gamma_{1}}\p_{t}^{k}(-\sigma\sqrt{h}\lap_{h}(\eta)n^{i}) \p_{t}^{4-k}(n_{\kappa}-n)\p_{t}^{4}v^{i}}_{III_{12,1}}\\&\underbrace{+\int_{0}^{t}\int_{\Gamma_{1}}\p_{t}^{3}(-\sigma\sqrt{h}\lap_{h}(\eta)n^{i}) \p_{t}(n_{\kappa}-n)\p_{t}^{4}v^{i}}_{III_{12,2}}\\&\underbrace{+\int_{0}^{t}\int_{\Gamma_{1}}(-\sigma\sqrt{h}\lap_{h}(\eta)n^{i}) \p_{t}^{4}(n_{\kappa}-n)\p_{t}^{4}v^{i}}_{III_{12,3}}
\end{aligned}
\end{align}
\begin{align}    III_{13}=\int_{0}^{t}\int_{\Gamma_{1}}\p_{t}^{4}(-\sigma\sqrt{h}\lap_{h}(\eta) \cdot n_{\kappa}(n_{\kappa}^{i}-n^{i}))\p_{t}^{4}v^{i}
\end{align}
Further,
\begin{align}
\begin{aligned}    |III_{12,1}|&\lesssim\sum_{k=1}^{2}\int_{0}^{t}\|\p_{t}^{k}(-\sigma\sqrt{h}\lap_{h}(\eta)n^{i})\|_{L^{\infty}(\Gamma_{1})} \|\p_{t}^{4-k}(n_{\kappa}-n)\|_{L^{2}(\Gamma_{1})}\|\p_{t}^{4}v^{i}\|_{L^{2}(\Gamma_{1})}\\&\lesssim\sqrt{\kappa}T^{\frac{1}{2}}\PP(E^{\kappa})(\int_{0}^{t}\|\sqrt{\kappa}\p_{t}^{4}v\|_{L^{2}(\Gamma_{1})}^{2})^{\frac{1}{2}}\\&\lesssim \PP_{0}+T\PP(E^{\kappa}).
\end{aligned}
\end{align}
Similarly,
\begin{align}
    \begin{aligned}        III_{12,2}=\sigma\int_{0}^{t}\int_{\Gamma_{1}}\p_{t}^{3}(h^{ij}\p_{j}\eta n^{i})\p_{i}\p_{t}(n^{\kappa}-n)\p_{t}^{4}v^{i}+\sigma\int_{0}^{t}\int_{\Gamma_{1}}\p_{t}^{3}(h^{ij}\p_{j}\eta n^{i})\p_{t}(n^{\kappa}-n)\p_{i}\p_{t}^{4}v^{i},
    \end{aligned}
\end{align}
\begin{align}
    \begin{aligned}        |III_{12,2}|&\lesssim\sigma\int_{0}^{t}\|\p_{t}^{3}(h^{ij}\p_{j}\eta n^{i})\|_{L^{2}(\Gamma_{1})}\|\p_{i}\p_{t}(n^{\kappa}-n)\|_{L^{\infty}}\|\p_{t}^{4}v^{i}\|_{L^{2}(\Gamma_{1})}\\&+\sigma\int_{0}^{t}\|\p_{t}^{3}(h^{ij}\p_{j}\eta n^{i})\|_{L^{2}(\Gamma_{1})}\|\p_{t}(n^{\kappa}-n)\|_{L^{\infty}}\|\p_{i}\p_{t}^{4}v^{i}\|_{L^{2}(\Gamma_{1})}\\&\lesssim\int_{0}^{t}\PP(E^{\kappa})\sqrt{\kappa}\times\|\p_{t}^{4}v\|_{H^{1}(\Gamma_{1})}\\&\lesssim\delta\int_{0}^{t}\|\sqrt{\kappa}\p_{t}^{4}v\|_{H^{1}(\Gamma_{1})}^{2}+T\PP(E^{\kappa})
    \end{aligned}
\end{align}
For $III_{12,3}$, see our previous lemma and the discussion therein. It's evident that $III_{13}$ can be controlled in a similar fashion. We omit the details. Up to now, we have controlled
\begin{align}\label{estimateforIII1}
III_{1}=\sigma\int_{\Gamma_{1}}\sqrt{h}h^{ij}\p_{j}\p_{t}^{3}v^{\lambda}\Pi^{\mu}_{\lambda}\p_{i}\p_{t}^{3}v^{\alpha}\Pi_{\alpha}^{\mu}+\mathcal{R},
\end{align}
where $|R|\lesssim C(\delta)(\PP_{0}+T\PP(E^{\kappa}))+\delta(\|v_{ttt}\|_{1.5}^{2}+\int_{0}^{t}\|\sqrt{\kappa}\p_{t}^{4}v\|_{H^{1}(\Gamma_{1})}^{2})$. We estimate the viscous term $III_{2}$.
\begin{align}
\begin{aligned}     III_{2}&=\underbrace{\int_{0}^{t}\int_{\Gamma_{1}}\kappa\lap_{0}\p_{t}^{4}(v\cdot n_{\kappa})n_{\kappa}^{i}\p_{t}^{4}v^{i}}_{III_{21}}\\&\underbrace{+\sum_{k=0}^{3}C_{4}^{k}\int_{0}^{t}\int_{\Gamma_{1}}\kappa\lap_{0}\p_{t}^{k}(v\cdot n_{\kappa})\p_{t}^{4-k}n_{\kappa}^{i}\p_{t}^{4}v^{i}}_{III_{22}}.
\end{aligned}  
\end{align}
We estimate $III_{22}$ first. When $k=0,1$, employing H\"{o}lder inequality ($L^{\infty}$-$L^{2}$-$L^{2}$)
\begin{align}
\begin{aligned}
    | C_{4}^{k}\int_{0}^{t}\int_{\Gamma_{1}}-\kappa\lap_{0}(v\cdot n_{\kappa})\p_{t}^{4}n_{\kappa}^{i}\p_{t}^{4}v^{i}|&\leq | C_{4}^{k}\int_{0}^{t}\int_{\Gamma_{1}}-\kappa\p_{t}^{k}(v\cdot n_{\kappa})\p_{t}^{4-k}n_{\kappa}^{i}\p_{t}^{4}v^{i}|\\&+| C_{4}^{k}\int_{0}^{t}\int_{\Gamma_{1}}-\kappa\pb_{l}\p_{t}^{k}(v\cdot n_{\kappa})\pb_{l}\p_{t}^{4-k}n_{\kappa}^{i}\p_{t}^{4}v^{i}|\\&+| C_{4}^{k}\int_{0}^{t}\int_{\Gamma_{1}}-\kappa\pb_{l}\p_{t}^{k}(v\cdot n_{\kappa})\p_{t}^{4-k}n_{\kappa}^{i}\pb_{l}\p_{t}^{4}v^{i}|\\&\lesssim| C_{4}^{k}|\int_{0}^{t}\kappa\|\p_{t}^{k}(v\cdot n_{\kappa})\|_{H^{2.5}(\Gamma_{1})}\|\p_{t}^{4-k}n_{\kappa}^{i}\|_{L^{2}(\Gamma_{1})}\|\p_{t}^{4}v^{i}\|_{H^{1}(\Gamma_{1})}\\&\lesssim\int_{0}^{t}\kappa\PP(E^{\kappa})\|\p_{t}^{4}v\|_{H^{1}(\Gamma_{1})}\\&\lesssim C(\delta)(\PP_{0}+T\PP(E^{\kappa}))+\delta\int_{0}^{t}\|\sqrt{\kappa}\p_{t}^{4}v\|_{H^{1}(\Gamma_{1})}.
\end{aligned} 
\end{align}
When $k=2$, we use H\"{o}lder inequality ($L^{4}$-$L^{4}$-$L^{2}$)
\begin{align}
\begin{aligned}
     | C_{4}^{k}\int_{0}^{t}\int_{\Gamma_{1}}\kappa\lap_{0}\p_{t}^{k}(v\cdot n_{\kappa})\p_{t}^{4-k}n_{\kappa}^{i}\p_{t}^{4}v^{i}|&\leq | C_{4}^{k}\int_{0}^{t}\int_{\Gamma_{1}}\kappa\p_{t}^{k}(v\cdot n_{\kappa})\p_{t}^{4-k}n_{\kappa}^{i}\p_{t}^{4}v^{i}|\\&+| C_{4}^{k}\int_{0}^{t}\int_{\Gamma_{1}}\kappa\pb_{l}\p_{t}^{k}(v\cdot n_{\kappa})\pb_{l}\p_{t}^{4-k}n_{\kappa}^{i}\p_{t}^{4}v^{i}|\\&+| C_{4}^{k}\int_{0}^{t}\int_{\Gamma_{1}}\kappa\pb_{l}\p_{t}^{k}(v\cdot n_{\kappa})\p_{t}^{4-k}n_{\kappa}^{i}\pb_{l}\p_{t}^{4}v^{i}|\\&\lesssim| C_{4}^{k}|\int_{0}^{t}\kappa\|\p_{t}^{k}(v\cdot n_{\kappa})\|_{H^{1.5}(\Gamma_{1})}\|\p_{t}^{4-k}n_{\kappa}^{i}\|_{H^{1.5}(\Gamma_{1})}\|\p_{t}^{4}v^{i}\|_{H^{1}(\Gamma_{1})}\\&\lesssim\int_{0}^{t}\kappa\PP(E^{\kappa})\|\p_{t}^{4}v\|_{H^{1}(\Gamma_{1})}\\&\lesssim C(\delta)(\PP_{0}+T\PP(E^{\kappa}))+\delta\int_{0}^{t}\|\p_{t}^{4}v^{\kappa}\|_{H^{1}(\Gamma_{1})}.
\end{aligned}    
\end{align}
When $k=3$, applying H\"{o}lder inequality ($L^{2}$-$L^{\infty}$-$L^{2}$)
\begin{align}
\begin{aligned}
     |4\int_{0}^{t}\int_{\Gamma_{1}}\kappa\lap_{0}\p_{t}^{3}(v\cdot n_{\kappa})\p_{t}n_{\kappa}^{i}\p_{t}^{4}v^{i}|&\leq | 4\int_{0}^{t}\int_{\Gamma_{1}}\kappa\p_{t}^{3}(v\cdot n_{\kappa})\p_{t}n_{\kappa}^{i}\p_{t}^{4}v^{i}|\\&+| 4\int_{0}^{t}\int_{\Gamma_{1}}\kappa\pb_{l}\p_{t}^{3}(v\cdot n_{\kappa})\pb_{l}\p_{t}n_{\kappa}^{i}\p_{t}^{4}v^{i}|\\&+| 4\int_{0}^{t}\int_{\Gamma_{1}}\kappa\pb_{l}\p_{t}^{3}(v\cdot n_{\kappa})\p_{t}n_{\kappa}^{i}\pb_{l}\p_{t}^{4}v^{i}|\\&\lesssim|\int_{0}^{t}\kappa\|\p_{t}^{3}(v\cdot n_{\kappa})\|_{H^{1}(\Gamma_{1})}\|\p_{t}n_{\kappa}^{i}\|_{L^{\infty}(\Gamma_{1})}\|\p_{t}^{4}v^{i}\|_{H^{1}(\Gamma_{1})}\\&\lesssim\int_{0}^{t}\kappa\PP(E^{\kappa})\|\p_{t}^{4}v\|_{H^{1}(\Gamma_{1})}\\&\lesssim C(\delta)(\PP_{0}+T\PP(E^{\kappa}))+\delta\int_{0}^{t}\|\p_{t}^{4}v^{\kappa}\|_{H^{1}(\Gamma_{1})}.
\end{aligned}    
\end{align}
Consequently, we have
\begin{align}\label{eatimateforIII22}
    |III_{22}|\lesssim C(\delta)(\PP_{0}+T\PP(E^{\kappa}))+\delta\int_{0}^{t}\|\p_{t}^{4}v^{\kappa}\|_{H^{1}(\Gamma_{1})}.
\end{align}
The term $III_{21}$ will contribute another boundary regular term we need.
\begin{align}
\begin{aligned}
     III_{21}&=\int_{0}^{t}\int_{\Gamma_{1}}\kappa\p_{t}^{4}(v\cdot n_{\kappa})n_{\kappa}^{i}\p_{t}^{4}v^{i}+\int_{0}^{t}\int_{\Gamma_{1}}\kappa\p_{t}^{4}\pb_{l}(v\cdot n_{\kappa})\pb_{l}(\p_{t}^{4}v\cdot n_{\kappa})\\&=\int_{0}^{t}\int_{\Gamma_{1}}-\kappa\p_{t}^{4}v\cdot n_{\kappa}\p_{t}^{4}v\cdot n_{\kappa}-\int_{0}^{t}\int_{\Gamma_{1}}\kappa\pb_{l}(\p_{t}^{4}v\cdot n_{\kappa})\pb_{l}(\p_{t}^{4}v\cdot n_{\kappa})\\&+\sum_{k=0}^{3}C_{4}^{k}\int_{0}^{t}\int_{\Gamma_{1}}\kappa\p_{t}^{k}v^{j}\p_{t}^{4-k}n_{\kappa}^{i}\p_{t}^{4}v\cdot n_{\kappa}+\sum_{k=0}^{3}C_{4}^{k}\int_{0}^{t}\int_{\Gamma_{1}}\kappa\pb_{l}(\p_{t}^{k}v^{j}\p_{t}^{4-k} n_{\kappa})\pb_{l} (\p_{t}^{4}v^{i} n_{\kappa}^{i})\\&=\underbrace{\int_{0}^{t}\int_{\Gamma_{1}}\kappa\p_{t}^{4}v\cdot n_{\kappa}\p_{t}^{4}v\cdot n_{\kappa}+\int_{0}^{t}\int_{\Gamma_{1}}\kappa\pb_{l}(\p_{t}^{4}v\cdot n_{\kappa})\pb_{l}(\p_{t}^{4}v\cdot n_{\kappa})}_{III_{21,1}}\\&\underbrace{+\sum_{k=0}^{3}C_{4}^{k}\int_{0}^{t}\int_{\Gamma_{1}}\kappa\p_{t}^{k}v^{j}\p_{t}^{4-k}n_{\kappa}^{i}\p_{t}^{4}v\cdot n_{\kappa}}_{III_{21,2}}\underbrace{+\sum_{k=0}^{3}C_{4}^{k}\int_{0}^{t}\int_{\Gamma_{1}}\kappa\pb_{l}\p_{t}^{k}v^{j}\p_{t}^{4-k} n_{\kappa}\pb_{l} \p_{t}^{4}v^{i} n_{\kappa}^{i}}_{III_{21,3}}\\&\underbrace{+\sum_{k=0}^{3}C_{4}^{k}\int_{0}^{t}\int_{\Gamma_{1}}\kappa\pb_{l}\p_{t}^{k}v^{j}\p_{t}^{4-k} n_{\kappa} \p_{t}^{4}v^{i} \pb_{l}n_{\kappa}^{i}}_{III_{21,4}}\underbrace{+\sum_{k=0}^{3}C_{4}^{k}\int_{0}^{t}\int_{\Gamma_{1}}\kappa\p_{t}^{k}v^{j}\pb_{l}\p_{t}^{4-k} n_{\kappa}\pb_{l} \p_{t}^{4}v^{i} n_{\kappa}^{i}}_{III_{21,5}}\\&\underbrace{+\sum_{k=0}^{3}C_{4}^{k}\int_{0}^{t}\int_{\Gamma_{1}}\kappa\p_{t}^{k}v^{j}\pb_{l}\p_{t}^{4-k} n_{\kappa} \p_{t}^{4}v^{i} \pb_{l}n_{\kappa}^{i}}_{III21,6}.
\end{aligned}   
\end{align}
From direct computation, one can see that
\begin{align}
\begin{aligned}
    |III_{21,2}|&\lesssim\sum_{k=0}^{3}\int_{0}^{t}\kappa\|\p_{t}^{k}v^{j}\p_{t}^{3-k}\p_{t}n_{\kappa}^{i}\|_{L^{2}(\Gamma_{1})}\|\p_{t}^{4}v\|_{L^{2}(\Gamma_{1})}\|n_{\kappa}\|_{L^{\infty}(\Gamma_{1})}\\&\lesssim\int_{0}^{t}\PP(E^{\kappa})\kappa\|\p_{t}^{4}v\|_{L^{2}(\Gamma_{1})}\\&\lesssim\PP_{0}+T\PP(E^{\kappa}),
\end{aligned}    
\end{align}
\begin{align}
    \begin{aligned}
       |III_{21,3}|&\lesssim \sum_{k=0}^{3}|C_{4}^{k}|\int_{0}^{t}\int_{\Gamma_{1}}\kappa\|\p_{t}^{k}\pb_{l}v^{j}\p_{t}^{3-k} \p_{t}n_{\kappa}\|_{L^{2}(\Gamma_{1})}\|\pb_{l} \p_{t}^{4}v^{i}\|_{L^{2}(\Gamma_{1})}\| n_{\kappa}^{i}\|_{L^{\infty}(\Gamma_{1})}\\&\lesssim\int_{0}^{t}\PP(E^{\kappa})\kappa\|\p_{t}^{4}v\|_{H^{1}(\Gamma_{1})}\\&\lesssim\PP_{0}+T\PP(E^{\kappa}),
    \end{aligned}
\end{align}
\begin{align}
    \begin{aligned}
        |III_{21,4}|&\lesssim\sum_{k=0}^{3}|C_{4}^{k}|\int_{0}^{t}\int_{\Gamma_{1}}\kappa\|\p_{t}^{k}\pb_{l}v^{j}\p_{t}^{4-k} n_{\kappa}\|_{L^{2}(\Gamma_{1})} \|\p_{t}^{4}v^{i}\|_{L^{2}(\Gamma_{1})} \|\pb_{l}n_{\kappa}^{i}\|_{L^{\infty}(\Gamma_{1})}\\&\lesssim\int_{0}^{t}\PP(E^{\kappa})\kappa\|\p_{t}^{4}v\|_{L^{2}(\Gamma_{1})}\\&\lesssim\PP_{0}+T\PP(E^{\kappa}).
    \end{aligned}
\end{align}
Now, invoking \eqref{estimatefornongamma2}, we have
\begin{align}    |III_{21,5}|&\lesssim\int_{0}^{t}\kappa \|v^{j}\|_{L^{\infty}(\Gamma_{1})}\|\p_{t}^{4}\pb_{l} n_{\kappa}\|_{L^{2}(\Gamma_{1})}\|\pb_{l} \p_{t}^{4}v^{i}\|_{L^{2}(\Gamma_{1})}\| n_{\kappa}^{i}\|_{L^{\infty}(\Gamma_{1})}\\&+\sum_{k=0}^{2}|C_{4}^{k+1}|\int_{0}^{t}\kappa\|\p_{t}^{k}\p_{t}v^{j}\p_{t}^{3-k}\pb_{l} n_{\kappa}\|_{L^{2}(\Gamma_{1})}\|\pb_{l} \p_{t}^{4}v^{i}\|_{L^{2}(\Gamma_{1})} \|n_{\kappa}^{i}\|_{L^{\infty}(\Gamma_{1})}\\&\lesssim\int_{0}^{t}(\PP(E^{\kappa})+(\PP_{0}+T\PP(E^{\kappa}))\|\p_{t}^{3}v\|_{H^{2}(\Gamma_{1})})\kappa\|\p_{t}^{4}v\|_{H^{1}(\Gamma_{1})}\times(\PP_{0}+T\PP(E^{\kappa}))\\&\lesssim \PP_{0}+T\PP(E^{\kappa})+\int_{0}^{t}(\PP(E^{\kappa})+(\PP_{0}+T\PP(E^{\kappa}))\|\sqrt{\kappa}\p_{t}^{3}v\|_{H^{2}(\Gamma_{1})}\|\sqrt{\kappa}
\p_{t}^{4}v\|_{H^{1}(\Gamma_{1})}\\&\lesssim C(\delta)(\PP_{0}+T\PP(E^{\kappa})+(\int_{0}^{t}\|\sqrt{\kappa}\p_{t}^{3}v\|_{H^{2}(\Gamma_{1})}^{2})^{2}+\delta\int_{0}^{t}\|\sqrt{\kappa}\p_{t}^{4}v\|_{H^{1}(\Gamma_{1})}^{2},
\end{align}
$III_{21,6}$ can be estimated in a similar way,
\begin{align}
    \begin{aligned}       |III_{21,6}|&\lesssim\int_{0}^{t}\kappa \|v^{j}\|_{L^{\infty}(\Gamma_{1})}\|\p_{t}^{4}\pb_{l} n_{\kappa}\|_{L^{2}(\Gamma_{1})}\| \p_{t}^{4}v^{i}\|_{L^{2}(\Gamma_{1})}\| \pb_{l}n_{\kappa}^{i}\|_{L^{\infty}(\Gamma_{1})}\\&+\sum_{k=0}^{2}|C_{4}^{k+1}|\int_{0}^{t}\int_{\Gamma_{1}}\kappa\|\p_{t}^{k}\p_{t}v^{j}\p_{t}^{3-k} \pb_{l}n_{\kappa}\|_{L^{2}(\Gamma_{1})} \|\p_{t}^{4}v^{i}\|_{L^{2}(\Gamma_{1})} \|\pb_{l}n_{\kappa}^{i}\|_{L^{\infty}(\Gamma_{1})}\\&\lesssim C(\delta)(\PP_{0}+T\PP(E^{\kappa})+(\int_{0}^{t}\|\sqrt{\kappa}\p_{t}^{3}v\|_{H^{2}(\Gamma_{1})}^{2})^{2}+\delta\int_{0}^{t}\|\sqrt{\kappa}\p_{t}^{4}v\|_{H^{1}(\Gamma_{1})}^{2}.
    \end{aligned}
\end{align}
Consequently, we have
\begin{align}\label{estimateforIII21}
    III_{21}=\int_{0}^{t}\|\sqrt{\kappa}\p_{t}^{4}v\cdot n^{\kappa}\|_{H^{1}(\Gamma_{1})}+\mathcal{R},
\end{align}
where $\mathcal{R}\lesssim C(\delta)(\PP_{0}+T\PP(E^{\kappa})+(\int_{0}^{t}\|\sqrt{\kappa}\p_{t}^{3}v\|_{H^{2}(\Gamma_{1})}^{2})^{2}+\delta\int_{0}^{t}\|\sqrt{\kappa}\p_{t}^{4}v\|_{H^{1}(\Gamma_{1})}^{2}$.
Thus, combining \eqref{eatimateforIII22} and \eqref{estimateforIII21},
we have
\begin{align}\label{estimateforIII2}
    III_{2}=\int_{0}^{t}\|\sqrt{\kappa}\p_{t}^{4}v\cdot n^{\kappa}\|_{H^{1}(\Gamma_{1})}+\mathcal{R},
\end{align}
where $\mathcal{R}\lesssim C(\delta)(\PP_{0}+T\PP(E^{\kappa})+(\int_{0}^{t}\|\sqrt{\kappa}\p_{t}^{3}v\|_{H^{2}(\Gamma_{1})}^{2})^{2}+\delta\int_{0}^{t}\|\sqrt{\kappa}\p_{t}^{4}v\|_{H^{1}(\Gamma_{1})}^{2}$.
Combining \eqref{estimateforIII1} and \eqref{estimateforIII2}, one can derive that
\begin{align} \label{estimateforIII}  III=\sigma\int_{\Gamma_{1}}\sqrt{h}h^{ij}\p_{j}\p_{t}^{3}v^{\lambda}\Pi^{\mu}_{\lambda}\p_{i}\p_{t}^{3}v^{\alpha}\Pi_{\alpha}^{\mu}+\int_{0}^{t}\|\sqrt{\kappa}\p_{t}^{4}v\cdot n^{\kappa}\|_{H^{1}(\Gamma_{1})}^{2}+\mathcal{R},
\end{align}
where $\mathcal{R}\lesssim C(\delta)(\PP_{0}+T\PP(E^{\kappa})+(\int_{0}^{t}\|\sqrt{\kappa}\p_{t}^{3}v\|_{H^{2}(\Gamma_{1})}^{2})^{2}+\delta\int_{0}^{t}\|\sqrt{\kappa}\p_{t}^{4}v\|_{H^{1}(\Gamma_{1})}^{2}$.

Now, we estimate the interior term $IV$.
\begin{align}\label{expansionforIV}
\begin{aligned}
      IV=\underbrace{-\int_{0}^{t}\int_{\Omega}J_{\kappa} \p_{t}^{4}q(\ak)_{i}^{k}\p_{k}\p_{t}^{4}v^{i}}_{IV_{1}} \underbrace{-\sum_{l=1}^{4}C_{4}^{l}\int_{0}^{t}\int_{\Omega}\p_{t}^{l}(J_{\kappa}(\ak)_{i}^{k})\p_{t}^{4-l}q\p_{t}^{4}\p_{k}v^{i}}_{IV_{{2}}}
\end{aligned}  
\end{align}
We estimate $IV_{2}$ first.
\begin{align}
    \begin{aligned}
        IV_{2}=\underbrace{-\int_{0}^{t}\int_{\Omega}\p_{t}^{4}(J_{\kappa}(\ak)_{i}^{k})q\p_{t}^{4}\p_{k}v^{i}}_{IV_{21}}\underbrace{-\sum_{l=1}^{3}C_{4}^{l}\int_{0}^{t}\int_{\Omega}\p_{t}^{l}(J_{\kappa}(\ak)_{i}^{k})\p_{t}^{4-l}q\p_{t}^{4}\p_{k}v^{i}}_{IV_{22}}
    \end{aligned}
\end{align}
Integrating by part in $IV_{22}$, we have
\begin{align}
    \begin{aligned}
        IV_{22}&=-\sum_{l=l_{1}+l_{2}=1}^{3}C_{4}^{l}C_{l_{1}+l_{2}}^{l_{1}}\int_{0}^{t}\int_{\Omega}\p_{t}^{l_{1}}J_{\kappa}\p_{t}^{l_{2}}(\ak)_{i}^{k}\p_{t}^{4-l_{1}-l_{2}}q\p_{t}^{4}\p_{k}v^{i}\\&=\underbrace{-\sum_{l=l_{1}+l_{2}=1}^{3}C_{4}^{l}C_{l_{1}+l_{2}}^{l_{1}}\int_{0}^{t}\int_{\Gamma_{1}}\p_{t}^{l_{1}}J_{\kappa}\p_{t}^{l_{2}}(\ak)_{i}^{3}\p_{t}^{4-l_{1}-l_{2}}q\p_{t}^{4}v^{i}}_{IV_{22,1}}\\&\underbrace{+\sum_{l_{1}+l_{2}=1}^{3}C_{l_{1}+l_{2}}^{l_{1}}\int_{0}^{t}\int_{\Omega}\p_{k}(\p_{t}^{l_{1}}J_{\kappa}\p_{t}^{l_{2}}(\ak)_{i}^{k})\p_{t}^{4-l_{1}-l_{2}}q\p_{t}^{4}v^{i}}_{IV_{22,2}}\\&\underbrace{+\sum_{l=l_{1}+l_{2}=1}^{3}C_{4}^{l}C_{l_{1}+l_{2}}^{l_{1}}\int_{0}^{t}\int_{\Omega}\p_{t}^{l_{1}}J_{\kappa}\p_{t}^{l_{2}}(\ak)_{i}^{k}\p_{k}\p_{t}^{4-l_{1}-l_{2}}q\p_{t}^{4}v^{i}}_{IV_{22.3}}.
    \end{aligned}
\end{align}
Recall that $\sup_{t}\|q\|_{4.5}+\|q_{t}\|_{3.5}+\|q_{tt}\|_{2.5}\lesssim\PP(E^{\kappa})$, $\|q_{ttt}\|_{1}\lesssim\PP(E^{\kappa})$. We can write
\begin{align}\label{estimateforIV222}
\begin{aligned}
     |IV_{22,2}|&\lesssim\sum_{l=l_{1}+l_{2}=1}^{1}|C_{4}^{l}||C_{l_{1}+l_{2}}^{l_{1}}|\int_{0}^{t}\|\p_{k}(\p_{t}^{l_{1}}J_{\kappa}\p_{t}^{l_{2}}(\ak)_{i}^{k})\|_{L^{\infty}(\Omega)}\|\p_{t}^{3}q\|_{L^{2}(\Omega)}\|\p_{t}^{4}v^{i}\|_{L^{2}(\Omega)}\\&+\sum_{l=l_{1}+l_{2}=2}^{3}|C_{4}^{l}||C_{l_{1}+l_{2}}^{l_{1}}|\int_{0}^{t}\|\p_{k}(\p_{t}^{l_{1}}J_{\kappa}\p_{t}^{l_{2}}(\ak)_{i}^{k})\|_{L^{3}(\Omega)}\|\p_{t}^{4-l_{1}-l_{2}}q\|_{L^{6}(\Omega)}\|\p_{t}^{4}v^{i}\|_{L^{2}(\Omega)}\\&\lesssim\PP_{0}+T\PP(E^{\kappa}),
\end{aligned} 
\end{align}
\begin{align}\label{estimateforIV223}
    \begin{aligned}
       |IV_{22,3}| &\lesssim\sum_{l=l_{1}+l_{2}=1}^{1}|C_{4}^{l}||C_{l_{1}+l_{2}}^{l_{1}}|\int_{0}^{t}\|\p_{t}^{l_{1}}J_{\kappa}\p_{t}^{l_{2}}(\ak)_{i}^{k}\|_{L^{\infty}(\Omega)}\|\p_{k}\p_{t}^{3}q\|_{L^{2}(\Omega)}\|\p_{t}^{4}v^{i}\|_{L^{2}(\Omega)}\\&+\sum_{l=l_{1}+l_{2}=2}^{3}|C_{4}^{l}||C_{l_{1}+l_{2}}^{l_{1}}|\int_{0}^{t}\|\p_{t}^{l_{1}}J_{\kappa}\p_{t}^{l_{2}}(\ak)_{i}^{k}\|_{L^{3}(\Omega)}\|\p_{k}\p_{t}^{4-l_{1}-l_{2}}q\|_{L^{6}(\Omega)}\|\p_{t}^{4}v^{i}\|_{L^{2}(\Omega)}\\&\lesssim\PP_{0}+T\PP(E^{\kappa})
    \end{aligned}
\end{align}
Plunging in the Dirichlet boundary condition of q, we can express $IV_{22,1}$ as
\begin{align}
\begin{aligned}
     IV_{22,1}&=\underbrace{4\int_{0}^{t}\int_{\Gamma_{1}}\p_{t}(J_{\kappa}(\ak)_{i}^{3})\p_{t}^{3}(\sigma(\frac{\sqrt{h}}{\sqrt{h_{\kappa}}} \lap_{h}(\eta) \cdot n_{\kappa}))\p_{t}^{4}v^{i}}_{IV_{22,11}}\\&\underbrace{-4\int_{0}^{t}\int_{\Gamma_{1}}\p_{t}(J_{\kappa}(\ak)_{i}^{3})\p_{t}^{3}(\kappa\frac{1}{\sqrt{h_{\kappa}}}((1-\bar{\lap})(v\cdot n ^{\kappa})))\p_{t}^{4}v^{i}}_{IV_{22,12}}\\&\underbrace{-\sum_{l=l_{1}+l_{2}=2}^{3}C_{4}^{l}C_{l_{1}+l_{2}}^{l_{1}}\int_{0}^{t}\int_{\Gamma_{1}}\p_{t}^{l_{1}}J_{\kappa}\p_{t}^{l_{2}}(\ak)_{i}^{3}\p_{t}^{4-l_{1}-l_{2}}q\p_{t}^{4}v^{i}}_{IV_{22,13}},
\end{aligned}  
\end{align}
Integrating by part in time, we can write $II_{22,11}$ as
\begin{align}
    \begin{aligned}
        IV_{22,11}&=\underbrace{-4\int_{\Gamma_{1}}\p_{t}(J_{\kappa}(\ak)_{i}^{3})\p_{t}^{3}(-\sigma(\frac{\sqrt{h}}{\sqrt{h_{\kappa}}} \lap_{h}(\eta) \cdot n_{\kappa}))\p_{t}^{3}v^{i}|_{0}^{t}}_{IV_{22,111}}\\&\underbrace{-4\int_{0}^{t}\int_{\Gamma_{1}}\p_{t}(\p_{t}(J_{\kappa}(\ak)_{i}^{3})\p_{t}^{3}(-\sigma(\frac{\sqrt{h}}{\sqrt{h_{\kappa}}} \lap_{h}(\eta) \cdot n_{\kappa}))\p_{t}^{3}v^{i}}_{IV_{22,112}}
    \end{aligned}
\end{align}
We expand $IV_{22,111}$ as
\begin{align}
    \begin{aligned}
        IV_{22,111}&=\underbrace{4\sigma\int_{\Gamma_{1}}\p_{t}(J_{\kappa}(\ak)_{i}^{3})\frac{\sqrt{h}}{\sqrt{h_{\kappa}}} n_{\kappa}^{j} \p_{t}^{3}\lap_{h}(\eta^{j})\p_{t}^{3}v^{i}|_{0}^{t}}_{IV_{22,111,1}}\\&\underbrace{+4\sigma C_{3}^{k}\sum_{k=1}^{3}\int_{\Gamma_{1}}\p_{t}(J_{\kappa}(\ak)_{i}^{3})\p_{t}^{k}(\frac{\sqrt{h}}{\sqrt{h_{\kappa}}}n_{\kappa}^{j}) \p_{t}^{3-k}\lap_{h}(\eta^{j})  \p_{t}^{3}v^{i}|_{0}^{t}}_{IV_{22,111,2}}
    \end{aligned}
\end{align}
Integrating by part, one can see that
\begin{align}
\begin{aligned}
     IV_{22,111,1}&=-4\sigma\int_{\Gamma_{1}}\p_{\nu}(\p_{t}(J_{\kappa}(\ak)_{i}^{3})\frac{\sqrt{h}}{\sqrt{h_{\kappa}}}n_{\kappa}^{j})\p_{t}^{3}( h_{\kappa}^{\mu\nu}\p_{\mu}\eta^{j} ) \p_{t}^{3}v^{i}|_{0}^{t}\\&-4\sigma\int_{\Gamma_{1}}\p_{t}(J_{\kappa}(\ak)_{i}^{3})\frac{\sqrt{h}}{\sqrt{h_{\kappa}}}n_{\kappa}^{j}\p_{t}^{3} (h_{\kappa}^{\mu\nu}\p_{\mu}\eta^{j})  \p_{\nu}\p_{t}^{3}v^{i}|_{0}^{t},
\end{aligned} 
\end{align}
Therefore.
\begin{align}\label{estimateforIV221111}
    \begin{aligned}
        |IV_{22,111,1}|&\lesssim \PP_{0}+ \sigma\|\p_{\nu}(\p_{t}(J_{\kappa}(\ak)_{i}^{3})\frac{\sqrt{h}}{\sqrt{h_{\kappa}}}n_{\kappa}^{j})\|_{L^{4}(\Gamma_{1})}\|\p_{t}^{3}( h_{\kappa}^{\mu\nu}\p_{\mu}\eta^{j} )\|_{L^{2}(\Gamma_{1})} \|\p_{t}^{3}v^{i}\|_{L^{4}(\Gamma_{1})}\\&+\sigma\|\p_{t}(J_{\kappa}(\ak)_{i}^{3})\frac{\sqrt{h}}{\sqrt{h_{\kappa}}}n_{\kappa}^{j}\|_{L^{\infty}(\Gamma_{1})}\|\p_{t}^{3} (h_{\kappa}^{\mu\nu}\p_{\mu}\eta^{j})\|_{L^{2}(\Gamma_{1})}\|  \p_{\nu}\p_{t}^{3}v^{i}\|_{L^{2}(\Gamma_{1})}\\&\lesssim\delta\|\p_{t}^{3}v\|_{H^{1}(\Gamma_{1})}^{2}+C(\delta)(\PP_{0}+T\PP(E^{\kappa})).
    \end{aligned}
\end{align}
From a direct calculation,
\begin{align}\label{estimateforIV221112}
    \begin{aligned}
        |IV_{22,111,2}|&\lesssim4\sigma C_{3}^{k+1}\sum_{k=0}^{2}\|\p_{t}(J_{\kappa}(\ak)_{i}^{3})\|_{L^{\infty}(\Gamma_{1})}\|\p_{t}^{k}\p_{t}(\frac{\sqrt{h}}{\sqrt{h_{\kappa}}}n_{\kappa}^{j}) \p_{t}^{2-k}\lap_{h}(\eta^{j}) \|_{L^{2}(\Gamma_{1})} \|\p_{t}^{3}v^{i}\|_{L^{2}(\Gamma_{1})}\\&\lesssim\delta\|\p_{t}^{3}v\|_{H^{1}(\Gamma_{1})}^{2}+C(\delta)(\PP_{0}+T\PP(E^{\kappa})).
    \end{aligned}
\end{align}
Combining \eqref{estimateforIV221111} and \eqref{estimateforIV221112}, one can arrive at
\begin{align}\label{estimateforIV22111}
    |IV_{22,111}|\lesssim\delta\|\p_{t}^{3}v\|_{H^{1}(\Gamma_{1})}^{2}+C(\delta)(\PP_{0}+T\PP(E^{\kappa})).
\end{align}
For $IV_{22,112}$, a direct calculation shows that
\begin{align}
    \begin{aligned}
        IV_{22,112}&=\underbrace{4\sigma\int_{0}^{t}\int_{\Gamma_{1}}\p_{t}^{2}(J_{\kappa}(\ak)_{i}^{3})\p_{t}^{3}(\frac{\sqrt{h}}{\sqrt{h_{\kappa}}} \lap_{h}(\eta) \cdot n_{\kappa})\p_{t}^{3}v^{i}}_{IV_{22,112,1}}\\&\underbrace{+4\sigma\int_{0}^{t}\int_{\Gamma_{1}}\p_{t}(J_{\kappa}(\ak)_{i}^{3})\p_{t}^{4}(\frac{\sqrt{h}}{\sqrt{h_{\kappa}}} \lap_{h}(\eta) \cdot n_{\kappa})\p_{t}^{3}v^{i}}_{IV_{22,112,2}}
    \end{aligned}
\end{align}
\begin{align}\label{estimateforIV221121}
\begin{aligned}
    |IV_{22,112,1}|&\lesssim\sigma\int_{0}^{t}\|\p_{t}^{2}(J_{\kappa}(\ak)_{i}^{3})\|_{L^{\infty}(\Gamma_{1})}\|\p_{t}^{3}\sigma(\frac{\sqrt{h}}{\sqrt{h_{\kappa}}} \lap_{h}(\eta) \cdot n_{\kappa}))\|_{L^{2}(\Gamma_{1})}\|\p_{t}^{3}v^{i}\|_{L^{2}(\Gamma_{1})}\\&\lesssim\PP_{0}+T\PP(E^{\kappa}).
\end{aligned}    
\end{align}
Expanding $IV_{22,112,2}$, and integrating by part, one can write as
\begin{align}
    \begin{aligned}
        IV_{22,112,2}&=4\sigma\int_{0}^{t}\int_{\Gamma_{1}}\p_{t}(J_{\kappa}(\ak)_{i}^{3})\frac{1}{\sqrt{h_{\kappa}}}n_{\kappa}^{j}\p_{t}^{4} \p_{\mu}(h^{\mu\nu}\p_{\nu}\eta^{j}) \p_{t}^{3}v^{i}\\&+4\sigma\sum_{k=1}^{4}C_{4}^{k}\int_{0}^{t}\int_{\Gamma_{1}}\p_{t}(J_{\kappa}(\ak)_{i}^{3})\p_{t}^{k}(\frac{1}{\sqrt{h_{\kappa}}} n_{\kappa}^{j})\p_{t}^{4-k}\p_{\mu}(h^{\mu\nu}\p_{\nu}\eta^{j}) \p_{t}^{3}v^{i}\\&=-4\sigma\int_{0}^{t}\int_{\Gamma_{1}}\p_{\mu}(\p_{t}(J_{\kappa}(\ak)_{i}^{3})\frac{1}{\sqrt{h_{\kappa}}}n_{\kappa}^{j})\p_{t}^{4} (h^{\mu\nu}\p_{\nu}\eta^{j}) \p_{t}^{3}v^{i}\\&-4\sigma\int_{0}^{t}\int_{\Gamma_{1}}\p_{t}(J_{\kappa}(\ak)_{i}^{3})\frac{1}{\sqrt{h_{\kappa}}}n_{\kappa}^{j}\p_{t}^{4} (h^{\mu\nu}\p_{\nu}\eta^{j})\p_{\mu} \p_{t}^{3}v^{i}\\&+4\sigma\sum_{k=1}^{4}C_{4}^{k}\int_{0}^{t}\int_{\Gamma_{1}}\p_{t}(J_{\kappa}(\ak)_{i}^{3})\p_{t}^{k}(\frac{1}{\sqrt{h_{\kappa}}} n_{\kappa}^{j})\p_{t}^{4-k}\p_{\mu}(h^{\mu\nu}\p_{\nu}\eta^{j}) \p_{t}^{3}v^{i}.
    \end{aligned}
\end{align}
Hence,
\begin{align}\label{estimateforIV221122}
    \begin{aligned}
        |IV_{22,112,2}|&\lesssim\sigma\int_{0}^{t}\|\p_{\mu}(\p_{t}(J_{\kappa}(\ak)_{i}^{3})\frac{1}{\sqrt{h_{\kappa}}}n_{\kappa}^{j})\|_{L^{\infty}(\Gamma_{1})}\|\p_{t}^{4} (h^{\mu\nu}\p_{\nu}\eta^{j}) \|_{L^{2}(\Gamma_{1})}\|\p_{t}^{3}v^{i}\|_{L^{2}(\Gamma_{1})}\\&+\sigma\int_{0}^{t}\|\p_{t}(J_{\kappa}(\ak)_{i}^{3})\frac{1}{\sqrt{h_{\kappa}}}n_{\kappa}^{j}\|_{L^{\infty}(\Gamma_{1})}\|\p_{t}^{4} (h^{\mu\nu}\p_{\nu}\eta^{j})\|_{L^{2}(\Gamma_{1})}\|\p_{\mu} \p_{t}^{3}v^{i}\|_{L^{2}(\Gamma_{1})}\\&+\sigma\sum_{k=0}^{3}C_{4}^{k+1}\int_{0}^{t}\|\p_{t}(J_{\kappa}(\ak)_{i}^{3})\|_{L^{\infty}(\Gamma_{1})}\|\p_{t}^{k}\p_{t}(\frac{1}{\sqrt{h_{\kappa}}} n_{\kappa}^{j})\p_{t}^{3-k}\p_{\mu}(h^{\mu\nu}\p_{\nu}\eta^{j})\|_{L^{2}(\Gamma_{1})}\| \p_{t}^{3}v^{i}\|_{L^{2}(\Gamma_{1})}\\&\lesssim\PP_{0}+T\PP(E^{\kappa}).
    \end{aligned}
\end{align}

Hence, combining \eqref{estimateforIV221121} and \eqref{estimateforIV221122}, one can deduce that
\begin{align}\label{estimateforIV22112}
    |IV_{22,112}|&\lesssim\PP_{0}+T\PP(E^{\kappa}).
\end{align}
Combining \eqref{estimateforIV22111} and \eqref{estimateforIV22112}, we have
\begin{align}\label{estimateforIV2211}
    |IV_{22,11}|&\lesssim \delta\|\p_{t}^{3}v\|_{H^{1}(\Gamma_{1})}^{2}+C(\delta)(\PP_{0}+T\PP(E^{\kappa})).
\end{align}
Now, considering $IV_{22,12}$,
\begin{align}
    \begin{aligned}
       IV_{22,12}&= \underbrace{-4\int_{0}^{t}\int_{\Gamma_{1}}\kappa\p_{t}(J_{\kappa}(\ak)_{i}^{3})\frac{1}{\sqrt{h_{\kappa}}}\p_{t}^{3}((1-\bar{\lap})(v\cdot n ^{\kappa}))\p_{t}^{4}v^{i}}_{IV_{22,121}}\\&\underbrace{-4\sum_{k=1}^{3}C_{3}^{k}\int_{0}^{t}\int_{\Gamma_{1}}\kappa\p_{t}(J_{\kappa}(\ak)_{i}^{3})\p_{t}^{k}\frac{1}{\sqrt{h_{\kappa}}}\p_{t}^{3-k}(1-\bar{\lap})(v\cdot n ^{\kappa}))\p_{t}^{4}v^{i}}_{IV_{22,122}}
    \end{aligned}
\end{align}
Integrating by part for $IV_{22,121}$, we have
\begin{align}
    \begin{aligned}        IV_{22,121}&=-4\int_{0}^{t}\int_{\Gamma_{1}}\kappa\p_{t}(J_{\kappa}(\ak)_{i}^{3})\frac{1}{\sqrt{h_{\kappa}}}\p_{t}^{3}(v\cdot n ^{\kappa})\p_{t}^{4}v^{i}\\&-4\int_{0}^{t}\int_{\Gamma_{1}}\kappa\pb_{l}\p_{t}(J_{\kappa}(\ak)_{i}^{3})\frac{1}{\sqrt{h_{\kappa}}}\p_{t}^{3}\pb_{l}(v\cdot n ^{\kappa})\p_{t}^{4}v^{i}\\&-4\int_{0}^{t}\int_{\Gamma_{1}}\kappa\p_{t}(J_{\kappa}(\ak)_{i}^{3})\frac{1}{\sqrt{h_{\kappa}}}\p_{t}^{3}\pb_{l}(v\cdot n ^{\kappa})\pb_{l}\p_{t}^{4}v^{i}.
    \end{aligned}
\end{align}
Thus,
\begin{align}\label{estimateforIV22121}
    \begin{aligned}
       |IV_{22,121}| &\lesssim\int_{0}^{t}\sqrt{\kappa}\|\p_{t}(J_{\kappa}(\ak)_{i}^{3})\frac{1}{\sqrt{h_{\kappa}}}\|_{L^{\infty}(\Gamma_{1})}\|\p_{t}^{3}(v\cdot n ^{\kappa})\|_{L^{2}(\Gamma_{1})}\|\sqrt{\kappa}\p_{t}^{4}v^{i}\|_{L^{2}(\Gamma_{1})}\\&+\int_{0}^{t}\sqrt{\kappa}\|\pb_{l}\p_{t}(J_{\kappa}(\ak)_{i}^{3})\frac{1}{\sqrt{h_{\kappa}}}\|_{L^{\infty}(\Gamma_{1})}\|\p_{t}^{3}\pb_{l}(v\cdot n ^{\kappa})\|_{L^{2}(\Gamma_{1})}\|\sqrt{\kappa}\p_{t}^{4}v^{i}\|_{L^{2}(\Gamma_{1})}\\&+\int_{0}^{t}\sqrt{\kappa}\|\p_{t}(J_{\kappa}(\ak)_{i}^{3})\frac{1}{\sqrt{h_{\kappa}}}\|_{L^{\infty}(\Gamma_{1})}\|\p_{t}^{3}\pb_{l}(v\cdot n ^{\kappa})\|_{L^{2}(\Gamma_{1})}\|\sqrt{\kappa}\pb_{l}\p_{t}^{4}v^{i}\|_{L^{2}(\Gamma_{1})}\\&\lesssim\int_{0}^{t}\sqrt
       {\kappa}\PP(E^{\kappa})\|\sqrt{\kappa}\p_{t}^{4}v\|_{H^{1}(\Gamma_{1})}\\&\lesssim\PP_{0}+T\PP(E^{\kappa}).
    \end{aligned}
\end{align}
The control of $IV_{22,122}$ is straight,
\begin{align}\label{estimateforIV22122}
\begin{aligned}
     |IV_{22,122}|&\lesssim\sum_{k=0}^{2}C_{3}^{k}\int_{0}^{t}\sqrt{\kappa}\|\p_{t}(J_{\kappa}(\ak)_{i}^{3})\|_{L^{\infty}(\Gamma_{1})}\|
     \p_{t}^{k}\p_{t}\frac{1}{\sqrt{h_{\kappa}}}\p_{t}^{2-k}(1-\bar{\lap})(v\cdot n ^{\kappa}))\|_{L^{2}(\Gamma_{1})}\|\sqrt{\kappa}\p_{t}^{4}v^{i}\|_{L^{2}(\Gamma_{1})}\\&\lesssim\int_{0}^{t}\sqrt{\kappa}\PP(E^{\kappa})\|\sqrt{\kappa}\p_{t}^{4}v^{i}\|_{L^{2}(\Gamma_{1})}\\&\lesssim\PP_{0}+T\PP(E^{\kappa}).
\end{aligned}   
\end{align}
Combining \eqref{estimateforIV22121} and \eqref{estimateforIV22122} gives us
\begin{align}\label{estimateforIV2212}    |IV_{22,12}|&\lesssim\PP_{0}+T\PP(E^{\kappa}).
\end{align}
We estimate of $IV_{22,13}$ via integrating by part in time,
\begin{align}
    \begin{aligned}
        IV_{22,13}&=-\underbrace{\sum_{l=l_{1}+l_{2}=2}^{3}C_{4}^{l}C_{l_{1}+l_{2}}^{l_{1}}\int_{\Gamma_{1}}\p_{t}^{l_{1}}J_{\kappa}\p_{t}^{l_{2}}(\ak)_{i}^{3}\p_{t}^{4-l_{1}-l_{2}}q\p_{t}^{3}v^{i}|_{0}^{t}}_{IV_{22,131}}\\&\underbrace{+\sum_{l=l_{1}+l_{2}=2}^{3}C_{4}^{l}C_{l_{1}+l_{2}}^{l_{1}}\int_{0}^{t}\int_{\Gamma_{1}}\p_{t}(\p_{t}^{l_{1}}J_{\kappa}\p_{t}^{l_{2}}(\ak)_{i}^{3})\p_{t}^{4-l_{1}-l_{2}}q\p_{t}^{3}v^{i}}_{IV_{22,132}}\\&\underbrace{+\sum_{l=l_{1}+l_{2}=2}^{3}C_{4}^{l}C_{l_{1}+l_{2}}^{l_{1}}\int_{0}^{t}\int_{\Gamma_{1}}\p_{t}^{l_{1}}J_{\kappa}\p_{t}^{l_{2}}(\ak)_{i}^{3}\p_{t}^{5-l_{1}-l_{2}}q\p_{t}^{3}v^{i}}_{IV_{22,133}}
    \end{aligned}
\end{align}
Similarly,
\begin{align}\label{estimateforIV22131}
    \begin{aligned}        |IV_{22,131}|&\lesssim\PP_{0}+\sum_{l=l_{1}+l_{2}=2}^{2}C_{4}^{l}C_{l_{1}+l_{2}}^{l_{1}}\|\p_{t}^{l_{1}}J_{\kappa}\p_{t}^{l_{2}}(\ak)_{i}^{3}\|_{L^{4}(\Gamma_{1})}\|\p_{t}^{2}(\sigma(\frac{\sqrt{h}}{\sqrt{h_{\kappa}}} \lap_{h}(\eta) \cdot n_{\kappa}))\|_{L^{2}(\Gamma_{1})}\|\p_{t}^{3}v^{i}\|_{L^{4}(\Gamma_{1})}\\&+\sum_{l=l_{1}+l_{2}=2}^{2}C_{4}^{l}C_{l_{1}+l_{2}}^{l_{1}}\|\p_{t}^{l_{1}}J_{\kappa}\p_{t}^{l_{2}}(\ak)_{i}^{3}\|_{L^{4}(\Gamma_{1})}\|\p_{t}^{2}(\kappa\frac{1}{\sqrt{h_{\kappa}}}((1-\bar{\lap})(v\cdot n ^{\kappa}))\|_{L^{2}(\Gamma_{1})}\|\p_{t}^{3}v^{i}\|_{L^{4}(\Gamma_{1})}\\&+\sum_{l=l_{1}+l_{2}=3}^{3}C_{4}^{l}C_{l_{1}+l_{2}}^{l_{1}}\|\p_{t}^{l_{1}}J_{\kappa}\p_{t}^{l_{2}}(\ak)_{i}^{3}\|_{L^{2}(\Gamma_{1})}\|\p_{t}q\|_{H^{1.5}(\Gamma_{1})}\|\p_{t}^{3}v^{i}\|_{L^{2}(\Gamma_{1})}\\&\lesssim\delta\|\p_{t}^{3}v\|_{H^{1}(\Gamma_{1})}+C(\delta)(\PP_{0}+T\PP(E^{\kappa})),
    \end{aligned}
\end{align}
\begin{align}\label{estimateforIV22132}
    \begin{aligned}        |IV_{22,132}|&\lesssim\sum_{l=l_{1}+l_{2}=2}^{3}C_{4}^{l}C_{l_{1}+l_{2}}^{l_{1}}\int_{0}^{t}\|\p_{t}(\p_{t}^{l_{1}}J_{\kappa}\p_{t}^{l_{2}}(\ak)_{i}^{3})\|_{L^{2}(\Gamma_{1})}\|\p_{t}^{4-l_{1}-l_{2}}q\|_{L^{4}(\Gamma_{1})}\|\p_{t}^{3}v^{i}\|_{L^{4}(\Gamma_{1})}\\&\lesssim T\PP(E^{\kappa}),
    \end{aligned}
\end{align}
and
\begin{align}\label{expansionforIV22133}
    \begin{aligned}        IV_{22,133}&=\underbrace{+\sum_{l=l_{1}+l_{2}=3}^{3}C_{4}^{l}C_{l_{1}+l_{2}}^{l_{1}}\int_{0}^{t}\int_{\Gamma_{1}}\p_{t}^{l_{1}}J_{\kappa}\p_{t}^{l_{2}}(\ak)_{i}^{3}\p_{t}^{2}(\sigma\frac{\sqrt{h}}{\sqrt{h_{\kappa}}} \lap_{h}(\eta) \cdot n_{\kappa})\p_{t}^{3}v^{i}}_{IV_{22,133,1}}\\&\underbrace{+\sum_{l=l_{1}+l_{2}=3}^{3}C_{4}^{l}C_{l_{1}+l_{2}}^{l_{1}}\int_{0}^{t}\int_{\Gamma_{1}}\p_{t}^{l_{1}}J_{\kappa}\p_{t}^{l_{2}}(\ak)_{i}^{3}\p_{t}^{5-l_{1}-l_{2}}(\kappa\frac{1}{\sqrt{h_{\kappa}}}(1-\bar{\lap})(v\cdot n ^{\kappa}))\p_{t}^{3}v^{i}}_{IV_{22,133,2}}\\&\underbrace{+\sum_{l=l_{1}+l_{2}=2}^{2}C_{4}^{l}C_{l_{1}+l_{2}}^{l_{1}}\int_{0}^{t}\int_{\Gamma_{1}}\p_{t}^{l_{1}}J_{\kappa}\p_{t}^{l_{2}}(\ak)_{i}^{3}\p_{t}^{3}(\sigma\frac{\sqrt{h}}{\sqrt{h_{\kappa}}} \lap_{h}(\eta) \cdot n_{\kappa})\p_{t}^{3}v^{i}}_{IV_{22,133,3}}\\&\underbrace{+\sum_{l=l_{1}+l_{2}=2}^{2}C_{4}^{l}C_{l_{1}+l_{2}}^{l_{1}}\int_{0}^{t}\int_{\Gamma_{1}}\p_{t}^{l_{1}}J_{\kappa}\p_{t}^{l_{2}}(\ak)_{i}^{3}\p_{t}^{3}(\kappa\frac{1}{\sqrt{h_{\kappa}}}(1-\bar{\lap})(v\cdot n ^{\kappa}))\p_{t}^{3}v^{i}}_{IV_{22,133,4}}.
    \end{aligned}
\end{align}
The first three terms in \eqref{expansionforIV22133} can be estimated directly, which reads
\begin{align}\label{estimateforIV221331}
    \begin{aligned}
       | IV_{22,133,1}|&\lesssim\sum_{l=l_{1}+l_{2}=3}^{3}C_{4}^{l}C_{l_{1}+l_{2}}^{l_{1}}\int_{0}^{t}\|\p_{t}^{l_{1}}J_{\kappa}\p_{t}^{l_{2}}(\ak)_{i}^{3}\|_{L^{4}(\Gamma_{1})}\|\p_{t}^{2}(\sigma\frac{\sqrt{h}}{\sqrt{h_{\kappa}}} \lap_{h}(\eta) \cdot n_{\kappa})\|_{L^{2}(\Gamma_{1})}\|\p_{t}^{3}v^{i}\|_{L^{4}(\Gamma_{1})}\\&\lesssim T\PP(E^{\kappa}),
    \end{aligned}
\end{align}
\begin{align}\label{estimateforIV221332}
    \begin{aligned} |IV_{22,133,2}|&\lesssim\sum_{l=l_{1}+l_{2}=3}^{3}C_{4}^{l}C_{l_{1}+l_{2}}^{l_{1}}\int_{0}^{t}\|\p_{t}^{l_{1}}J_{\kappa}\p_{t}^{l_{2}}(\ak)_{i}^{3}\|_{L^{4}(\Gamma_{1})}\|\p_{t}^{2}(\kappa\frac{1}{\sqrt{h_{\kappa}}}(1-\bar{\lap})(v\cdot n ^{\kappa}))\|_{L^{2}(\Gamma_{1})}\|\p_{t}^{3}v^{i}\|_{L^{2}(\Gamma_{1})}\\&\lesssim T\PP(E^{\kappa}),
    \end{aligned}
\end{align}
\begin{align}\label{estimateforIV221333}
    \begin{aligned}
        |IV_{22,133,3}|&\lesssim\sum_{l=l_{1}+l_{2}=2}^{2}C_{4}^{l}C_{l_{1}+l_{2}}^{l_{1}}\int_{0}^{t}\|\p_{t}^{l_{1}}J_{\kappa}\p_{t}^{l_{2}}(\ak)_{i}^{3}\|_{L^{4}(\Gamma_{1})}\|\p_{t}^{3}(\sigma\frac{\sqrt{h}}{\sqrt{h_{\kappa}}} \lap_{h}(\eta) \cdot n_{\kappa})\|_{L^{2}(\Gamma_{1})}\|\p_{t}^{3}v^{i}\|_{L^{4}(\Gamma_{1})}\\&\lesssim T\PP(E^{\kappa}).
    \end{aligned}
\end{align}
For $IV_{22,133,4}$, we integrate it by part,
\begin{align}
    \begin{aligned}
        IV_{22,133,4}&=\sum_{l=l_{1}+l_{2}=2}^{2}C_{4}^{l}C_{l_{1}+l_{2}}^{l_{1}}\int_{0}^{t}\int_{\Gamma_{1}}\kappa\p_{t}^{l_{1}}J_{\kappa}\p_{t}^{l_{2}}(\ak)_{i}^{3}(\frac{1}{\sqrt{h_{\kappa}}})\p_{t}^{3}(v\cdot n ^{\kappa})\p_{t}^{3}v^{i}\\&+\sum_{l=l_{1}+l_{2}=2}^{2}C_{4}^{l}C_{l_{1}+l_{2}}^{l_{1}}\int_{0}^{t}\int_{\Gamma_{1}}\kappa\pb_{\mu}(\p_{t}^{l_{1}}J_{\kappa}\p_{t}^{l_{2}}(\ak)_{i}^{3}(\frac{1}{\sqrt{h_{\kappa}}}))\p_{t}^{3}\pb_{\mu}(v\cdot n ^{\kappa})\p_{t}^{3}v^{i}\\&+\sum_{l=l_{1}+l_{2}=2}^{2}C_{4}^{l}C_{l_{1}+l_{2}}^{l_{1}}\int_{0}^{t}\int_{\Gamma_{1}}\kappa\p_{t}^{l_{1}}J_{\kappa}\p_{t}^{l_{2}}(\ak)_{i}^{3}(\frac{1}{\sqrt{h_{\kappa}}})\p_{t}^{3}\pb_{\mu}(v\cdot n ^{\kappa})\p_{t}^{3}\pb_{\mu}v^{i}\\&+\sum_{l=l_{1}+l_{2}=2}^{2}C_{4}^{l}C_{l_{1}+l_{2}}^{l_{1}}\sum_{k=1}^{3}C_{3}^{k}\int_{0}^{t}\int_{\Gamma_{1}}\kappa\p_{t}^{l_{1}}J_{\kappa}\p_{t}^{l_{2}}(\ak)_{i}^{3}\p_{t}^{k}(\frac{1}{\sqrt{h_{\kappa}}})\p_{t}^{3-k}(1-\bar{\lap})(v\cdot n ^{\kappa})\p_{t}^{3}v^{i}
    \end{aligned}
\end{align}
Therefore,
\begin{align}\label{estimateforIV221334}
    \begin{aligned}
       |IV_{22,133,4}|&\lesssim \sum_{l=l_{1}+l_{2}=2}^{2}C_{4}^{l}C_{l_{1}+l_{2}}^{l_{1}}\int_{0}^{t}\kappa\|\p_{t}^{l_{1}}J_{\kappa}\p_{t}^{l_{2}}(\ak)_{i}^{3}(\frac{1}{\sqrt{h_{\kappa}}})\|_{L^{\infty}(\Gamma_{1})}\|\p_{t}^{3}(v\cdot n ^{\kappa})\|_{L^{2}(\Gamma_{1})}\|\p_{t}^{3}v^{i}\|_{L^{2}(\Gamma_{1})}\\&+\sum_{l=l_{1}+l_{2}=2}^{2}C_{4}^{l}C_{l_{1}+l_{2}}^{l_{1}}\int_{0}^{t}\kappa\|\pb_{\mu}(\p_{t}^{l_{1}}J_{\kappa}\p_{t}^{l_{2}}(\ak)_{i}^{3}(\frac{1}{\sqrt{h_{\kappa}}}))\|_{L^{4}(\Gamma_{1})}\|\p_{t}^{3}\pb_{\mu}(v\cdot n ^{\kappa})\|_{L^{2}(\Gamma_{1})}\|\p_{t}^{3}v^{i}\|_{L^{4}(\Gamma_{1})}\\&+\sum_{l=l_{1}+l_{2}=2}^{2}C_{4}^{l}C_{l_{1}+l_{2}}^{l_{1}}\int_{0}^{t}\kappa\|\p_{t}^{l_{1}}J_{\kappa}\p_{t}^{l_{2}}(\ak)_{i}^{3}(\frac{1}{\sqrt{h_{\kappa}}})\|_{L^{\infty}(\Gamma_{1})}\|\p_{t}^{3}\pb_{\mu}(v\cdot n ^{\kappa})\|_{L^{2}(\Gamma_{1})}\|\p_{t}^{3}\pb_{\mu}v^{i}\|_{L^{2}(\Gamma_{1})}\\&+\sum_{l=l_{1}+l_{2}=2}^{2}C_{4}^{l}C_{l_{1}+l_{2}}^{l_{1}}\sum_{k=0}^{2}C_{3}^{k+1}\int_{0}^{t}\kappa\|\p_{t}^{l_{1}}J_{\kappa}\p_{t}^{l_{2}}(\ak)_{i}^{3}\|_{L^{\infty}(\Gamma
       _{1})}\|\p_{t}^{k}\p_{t}(\frac{1}{\sqrt{h_{\kappa}}})\p_{t}^{2-k}(1-\bar{\lap})(v\cdot n ^{\kappa})\|_{L^{2}(\Gamma_{1})}\|\p_{t}^{3}v^{i}\|_{L^{2}(\Gamma_{1})}\\&\lesssim T\PP(E^{\kappa}).
    \end{aligned}
\end{align}
Combining \eqref{estimateforIV221331}, \eqref{estimateforIV221332}, \eqref{estimateforIV221333}, and\eqref{estimateforIV221334}, we obtain
\begin{align}\label{estimateforIV22133}
|IV_{22,133}|&\lesssim T\PP(E^{\kappa}).
\end{align}
Combining \eqref{estimateforIV22131}, \eqref{estimateforIV22132}, and \eqref{estimateforIV22133}, we have
\begin{align}\label{estimateforIV2213}
    |IV_{22,13}|&\lesssim\delta\|\p_{t}^{3}v\|_{H^{1}(\Gamma_{1})}+C(\delta)(\PP_{0}+T\PP(E^{\kappa})).
\end{align}
Combining \eqref{estimateforIV2211}, \eqref{estimateforIV2212}, and\eqref{estimateforIV2213}, we have
\begin{align}   \label{estimateforIV221} |IV_{22,1}|&\lesssim\delta\|\p_{t}^{3}v\|_{H^{1}(\Gamma_{1})}+C(\delta)(\PP_{0}+T\PP(E^{\kappa})).
\end{align}
Combining \eqref{estimateforIV221}, \eqref{estimateforIV222}, and \eqref{estimateforIV223}, we have
\begin{align}\label{estimateforIV22}
   |IV_{22}|&\lesssim\delta\|\p_{t}^{3}v\|_{H^{1}(\Gamma_{1})}+C(\delta)(\PP_{0}+T\PP(E^{\kappa})) .
\end{align}
The estimate for $IV_{21}$ is similar to \cite{disconzi2019priori}. The readers may also refer to \cite{gu2024local} for a similar case.
\begin{align}
    \begin{aligned}
        IV_{21}&=\underbrace{-\int_{0}^{t}\int_{\Omega}\p_{t}^{4}(JA_{i}^{k})q\p_{t}^{4}\p_{k}v^{i}}_{IV_{21,1}}\underbrace{-\int_{0}^{t}\int_{\Omega}\p_{t}^{4}(J_{\kappa}(\ak)_{i}^{k}-JA_{i}^{k})q\p_{t}^{4}\p_{k}v^{i}}_{IV_{21,2}},
    \end{aligned}
\end{align}
where $J=\det \nabla\eta$. Thus, we can say $JA_{i}^{k}=Cof(\nabla\eta)$ and $J_{\kappa}(\ak)_{i}^{k}=Cof(\nabla\eta^{\kappa})$. Let $H=Cof$ as the function. Similarly as the way we deal with $n-n^{\kappa}$, one can show that
\begin{align}
    \begin{aligned}
       H(\nabla\eta^{\kappa})-H(\nabla\eta)&=\int_{0}^{1}DH(\nabla\eta+s(\nabla\eta^{\kappa}-\nabla\eta))ds\times(\nabla\eta^{\kappa}-\nabla\eta).
    \end{aligned}
\end{align}
Therefore, we have
\begin{align}
    \begin{aligned}
        \p_{t}^{4}(H(\nabla\eta^{\kappa})-H(\nabla\eta))&=\int_{0}^{1}DH(\nabla\eta+s(\nabla\eta^{\kappa}-\nabla\eta))ds\times\p_{t}^{4}(\nabla\eta^{\kappa}-\nabla\eta)\\&+\sum_{k=1}^{4}C_{4}^{k}\p_{t}^{k}\int_{0}^{1}DH(\nabla\eta+s(\nabla\eta^{\kappa}-\nabla\eta))ds\times\p_{t}^{4-k}(\nabla\eta^{\kappa}-\nabla\eta),
    \end{aligned}
\end{align}
\begin{align}
    \begin{aligned}
        \|\p_{t}^{4}(H(\nabla\eta^{\kappa})-H(\nabla\eta))\|_{L^{2}(\Omega)}&\lesssim\|\int_{0}^{1}DH(\nabla\eta+s(\nabla\eta^{\kappa}-\nabla\eta))ds\|_{L^{\infty}(\Omega)}\times\|\p_{t}^{4}(\nabla\eta^{\kappa}-\nabla\eta)\|_{L^{2}(\Omega)}\\&+4\|\p_{t}\int_{0}^{1}DH(\nabla\eta+s(\nabla\eta^{\kappa}-\nabla\eta))ds\|_{L^{\infty}(\Omega)}\times\|\p_{t}^{3}(\nabla\eta^{\kappa}-\nabla\eta)\|_{L^{2}(\Omega)}\\&+\sum_{k=2}^{3}C_{4}^{k}\|\p_{t}^{k}\int_{0}^{1}DH(\nabla\eta+s(\nabla\eta^{\kappa}-\nabla\eta))ds\|_{L^{2}(\Omega)}\times\|\p_{t}^{4-k}(\nabla\eta^{\kappa}-\nabla\eta)\|_{L^{\infty}(\Omega)}\\&+\|\p_{t}^{4}\int_{0}^{1}DH(\nabla\eta+s(\nabla\eta^{\kappa}-\nabla\eta))ds\|_{L^{2}(\Omega)}\times\|(\nabla\eta^{\kappa}-\nabla\eta)\|_{L^{\infty}(\Omega)}\\&\lesssim\|\int_{0}^{1}DH(\nabla\eta+s(\nabla\eta^{\kappa}-\nabla\eta))ds\|_{H^{2}(\Omega)}\times\kappa\|v_{ttt}\|_{H^{2}(\Omega)}\\&+4\|\p_{t}\int_{0}^{1}DH(\nabla\eta+s(\nabla\eta^{\kappa}-\nabla\eta))ds\|_{H^{2}(\Omega)}\times\kappa\|v_{tt}\|_{H^{2}(\Omega)}\\&+\sum_{k=0}^{1}C_{4}^{k+2}\|\p_{t}^{k}\p_{t}^{2}\int_{0}^{1}DH(\nabla\eta+s(\nabla\eta^{\kappa}-\nabla\eta))ds\|_{L^{2}(\Omega)}\times\kappa\|\p_{t}^{1-k}v\|_{H^{3.5}(\Omega)}\\&+\|\p_{t}^{4}\int_{0}^{1}DH(\nabla\eta+s(\nabla\eta^{\kappa}-\nabla\eta))ds\|_{L^{2}(\Omega)}\times\kappa\|\eta\|_{H^{3.5}(\Omega)}\\&\lesssim \kappa\PP(E^{\kappa})\|v_{ttt}\|_{H^{2}(\Omega)}+\kappa\PP(E^{\kappa}).
    \end{aligned}
\end{align}
Hence,
\begin{align}\label{estimateforIV212}
    \begin{aligned}
        |IV_{21,2}|&\lesssim\int_{0}^{t}\|\p_{t}^{4}(J_{\kappa}(\ak)_{i}^{k}-JA_{i}^{k})\|_{L^{2}(\Omega)}\|q\|_{L^{\infty}(\Omega)}\|\p_{t}^{4}\p_{k}v^{i}\|_{L^{2}(\Omega)}\\&\lesssim\int_{0}^{t}(\kappa\PP(E^{\kappa})\|v_{ttt}\|_{H^{2}(\Omega)}+\kappa\PP(E^{\kappa}))\|\p_{t}^{4}v\|_{H^{1}(\Omega)}\\&\lesssim\delta\int_{0}^{t}\|\sqrt{\kappa}\p_{t}^{4}v\|_{H^{1.5}(\Omega)}^{2}+\int_{0}^{t}\|\p_{t}^{3}v\|_{H^{2}(\Omega)}^{2}+C(\delta)(\PP_{0}+T\PP(E^{\kappa})).
    \end{aligned}
\end{align}
We expand $IV_{21,1}$ first
\begin{align}
\begin{aligned}
     IV_{21,1}&=\underbrace{-\int_{0}^{t}\int_{\Omega}\p_{t}^{4}J A_{i}^{k}q\p_{t}^{4}\p_{k}v^{i}}_{IV_{22,11}}\underbrace{-\int_{0}^{t}\int_{\Omega}J\p_{t}^{4}A_{i}^{k}q\p_{t}^{4}\p_{k}v^{i}}_{IV_{21,12}}\\&\underbrace{-\sum_{k=1}^{3}C_{4}^{k}\int_{0}^{t}\int_{\Omega}\p_{t}^{k}J\p_{t}^{4-k}A_{i}^{k}q\p_{k}\p_{t}^{4}v^{i}}_{IV_{22,13}}.
\end{aligned}
\end{align}
We deal with $IV_{21,13}$ first.
\begin{align}
    \sup_{t}\|q\|_{L^{\infty}(\Omega)}\lesssim\sup_{t}\|q\|_{H^{2}(\Omega)}\lesssim\PP_{0}+\int_{0}^{t}\|q_{t}\|_{H^{2}(\Omega)}\lesssim\PP_{0}+T\PP(E^{\kappa}).
\end{align}
One can compute that
\begin{align}
    \begin{aligned}
        IV_{21,13}=-\sum_{l=0}^{2}C_{4}^{l+1}\int_{\Omega}\p_{t}^{l}\p_{t}J\p_{t}^{2-l}\p_{t}A_{i}^{k}q\p_{k}\p_{t}^{3}v^{i}|_{0}^{t}+\sum_{l=0}^{2}C_{4}^{l+1}\int_{0}^{t}\int_{\Omega}\p_{t}(\p_{t}^{l}\p_{t}J\p_{t}^{2-l}\p_{t}A_{i}^{k}q)\p_{k}\p_{t}^{3}v^{i},
    \end{aligned}
\end{align}
\begin{align}\label{estimateforIV2113}
    \begin{aligned}
         |IV_{21,13}|&\lesssim\PP_{0}+\sum_{l=0}^{2}C_{4}^{l+1}\|\p_{t}^{l}\p_{t}J\p_{t}^{2-l}\p_{t}A_{i}^{k}\|_{L^{2}(\Omega)}\|q\|_{L^{\infty}(\Omega)}\|\p_{k}\p_{t}^{3}v^{i}\|_{L^{2}(\Omega)}\\&+\sum_{l=0}^{2}C_{4}^{l+1}\int_{0}^{t}\|\p_{t}(\p_{t}^{l}\p_{t}J\p_{t}^{2-l}\p_{t}A_{i}^{k})\|_{L^{2}(\Omega)}\|q\|_{L^{\infty}(\Omega)}\|\p_{k}\p_{t}^{3}v^{i}\|_{L^{2}(\Omega)}\\&+\sum_{l=0}^{2}C_{4}^{l+1}\int_{0}^{t}\|\p_{t}^{l}\p_{t}J\p_{t}^{2-l}\p_{t}A_{i}^{k}\|_{L^{2}(\Omega)}\|\p_{t}q\|_{L^{\infty}(\Omega)}\|\p_{k}\p_{t}^{3}v^{i}\|_{L^{2}(\Omega)}\\&\lesssim\delta \|\p_{t}^{3}v\|_{H^{1}(\Omega)}^{2}+C(\delta)(\PP_{0}+T\PP(E^{\kappa})). 
    \end{aligned}
\end{align}
Similarly,
\begin{align}
    \begin{aligned}
       IV_{21,12}&=-\int_{0}^{t}\int_{\Omega}J\p_{t}^{3}(-A_{i}^{l}\p_{l}v^{j}A_{j}^{k})q\p_{t}^{4}\p_{k}v^{i}\\&= \underbrace{\int_{0}^{t}\int_{\Omega}J\p_{t}^{3}\p_{l}v^{j}A_{i}^{l}A_{j}^{k}q\p_{t}^{4}\p_{k}v^{i}}_{IV_{22,121}}\underbrace{+\sum_{k=1}^{3}C_{3}^{k}\int_{0}^{t}\int_{\Omega}Jq\p_{t}^{k}(A_{i}^{l}A_{j}^{k})\p_{t}^{3-k}\p_{l}v^{j}\p_{t}^{4}\p_{k}v^{i}}_{IV_{22,122}}
    \end{aligned}
\end{align}
Applying a similar method, one can find that
\begin{align}
    \begin{aligned}
      IV_{21,122}&=  \sum_{k=0}^{2}C_{3}^{k+1}\int_{\Omega}Jq\p_{t}^{k}\p_{t}(A_{i}^{l}A_{j}^{k})\p_{t}^{2-k}\p_{l}v^{j}\p_{t}^{3}\p_{k}v^{i}\\&- \sum_{k=0}^{2}C_{3}^{k+1}\int_{0}^{t}\int_{\Omega}\p_{t}qJ\p_{t}^{k}\p_{t}(A_{i}^{l}A_{j}^{k})\p_{t}^{2-k}\p_{l}v^{j}\p_{t}^{3}\p_{k}v^{i}\\&- \sum_{k=0}^{2}C_{3}^{k+1}\int_{0}^{t}\int_{\Omega}q\p_{t}J\p_{t}^{k}\p_{t}(A_{i}^{l}A_{j}^{k})\p_{t}^{2-k}\p_{l}v^{j}\p_{t}^{3}\p_{k}v^{i}\\&- \sum_{k=0}^{2}C_{3}^{k+1}\int_{0}^{t}\int_{\Omega}qJ\p_{t}(\p_{t}^{k}\p_{t}(A_{i}^{l}A_{j}^{k})\p_{t}^{2-k}\p_{l}v^{j})\p_{t}^{3}\p_{k}v^{i},
    \end{aligned}
\end{align}
\begin{align}\label{estimateforIV21122}
    \begin{aligned}
        |IV_{21,122}|&\lesssim\PP_{0}+\sum_{k=0}^{2}C_{3}^{k+1}\|J\|_{L^{\infty}(\Omega)}\|q\|_{L^{\infty}(\Omega)}\|\p_{t}^{k}\p_{t}(A_{i}^{l}A_{j}^{k})\p_{t}^{2-k}\p_{l}v^{j}\|_{L^{2}(\Omega)}\|\p_{t}^{3}\p_{k}v^{i}\|_{L^{2}(\Omega)}\\&+ \sum_{k=0}^{2}C_{3}^{k+1}\int_{0}^{t}\|\p_{t}q\|_{L^{\infty}(\Omega)}\|J\|_{L^{\infty}(\Omega)}\|\p_{t}^{k}\p_{t}(A_{i}^{l}A_{j}^{k})\p_{t}^{2-k}\p_{l}v^{j}\|_{L^{2}(\Omega)}\|\p_{t}^{3}\p_{k}v^{i}\|_{L^{2}(\Omega)}\\&+ \sum_{k=0}^{2}C_{3}^{k+1}\int_{0}^{t}\|q\|_{L^{\infty}(\Omega)}\|\p_{t}J\|_{L^{\infty}(\Omega)}\|\p_{t}^{k}\p_{t}(A_{i}^{l}A_{j}^{k})\p_{t}^{2-k}\p_{l}v^{j}\|_{L^{2}(\Omega)}\|\p_{t}^{3}\p_{k}v^{i}\|_{L^{2}(\Omega)}\\&+ \sum_{k=0}^{2}C_{3}^{k+1}\int_{0}^{t}\|q\|_{L^{\infty}}\|J\|_{L^{\infty}(\Omega)}\|\p_{t}(\p_{t}^{k}\p_{t}(A_{i}^{l}A_{j}^{k})\p_{t}^{2-k}\p_{l}v^{j})\|_{L^{2}(\Omega)}\|\p_{t}^{3}\p_{k}v^{i}\|_{L^{2}(\Omega)}\\&\lesssim\delta \|\p_{t}^{3}v\|_{H^{1}(\Omega)}^{2}+C(\delta)(\PP_{0}+T\PP(E^{\kappa})).
    \end{aligned}
\end{align}
Also, we can write
\begin{align}\label{expansionforIV2112}
    \begin{aligned}       IV_{21,121}&=\int_{\Omega}JA_{i}^{l}A_{j}^{k}q\p_{t}^{3}\p_{l}v^{j}\p_{t}^{3}\p_{k}v^{i}-\int_{0}^{t}\int_{\Omega}\p_{t}qJA_{i}^{l}A_{j}^{k}\p_{t}^{3}\p_{l}v^{j}\p_{t}^{3}\p_{k}v^{i}\\&-\int_{0}^{t}\int_{\Omega}q\p_{t}(JA_{i}^{l}A_{j}^{k})\p_{t}^{3}\p_{l}v^{j}\p_{t}^{3}\p_{k}v^{i}-\int_{0}^{t}\int_{\Omega}qJA_{i}^{l}A_{j}^{k}\p_{t}^{4}\p_{l}v^{j}\p_{t}^{3}\p_{k}v^{i}
    \end{aligned}
\end{align}
Notice that the last term in \eqref{expansionforIV2112} coincides with $-IV_{21,121}$, we have
\begin{align}
    \begin{aligned}
        IV_{21,121}&=\frac{1}{2}\int_{\Omega}JA_{i}^{l}A_{j}^{k}q\p_{t}^{3}\p_{l}v^{j}\p_{t}^{3}\p_{k}v^{i}-\frac{1}{2}\int_{0}^{t}\int_{\Omega}\p_{t}qJA_{i}^{l}A_{j}^{k}\p_{t}^{3}\p_{l}v^{j}\p_{t}^{3}\p_{k}v^{i}\\&-\frac{1}{2}\int_{0}^{t}\int_{\Omega}q\p_{t}(JA_{i}^{l}A_{j}^{k})\p_{t}^{3}\p_{l}v^{j}\p_{t}^{3}\p_{k}v^{i}
    \end{aligned}
\end{align}
With the help of Sobolev interpolation theorem, it follows that
\begin{align}\label{estimateforIV21121}
    \begin{aligned}
        |IV_{21,121}|&\lesssim\|JA_{i}^{l}A_{j}^{k}q\|_{L^{\infty}(\Omega)}\|\p_{t}^{3}\p_{l}v^{j}\|_{L^{2}(\Omega)}\|\p_{t}^{3}\p_{k}v^{i}\|_{L^{2}(\Omega)}\\&+\int_{0}^{t}\|\p_{t}q\|_{L^{\infty}(\Omega)}\|JA_{i}^{l}A_{j}^{k}\|_{L^{\infty}(\Omega)}\|\p_{t}^{3}\p_{l}v^{j}\|_{L^{2}(\Omega)}\|\p_{t}^{3}\p_{k}v^{i}\|_{L^{2}(\Omega)}\\&+\int_{0}^{t}\|q\|_{L^{\infty}(\Omega)}\|\p_{t}(JA_{i}^{l}A_{j}^{k})\|_{L^{\infty}(\Omega)}\|\p_{t}^{3}\p_{l}v^{j}\|_{L^{2}(\Omega)}\|\p_{t}^{3}\p_{k}v^{i}\|_{L^{2}(\Omega)}\\&\lesssim(\PP_{0}+T\PP(E^{\kappa}))\|\p_{t}^{3}v\|_{H^{1}(\Omega)}^{2}+\PP_{0}+T\PP(E^{\kappa})\\&\lesssim(\PP_{0}+T\PP(E^{\kappa}))\|\p_{t}^{3}v\|_{L^{2}(\Omega)}^{\frac{2}{3}}\|\p_{t}^{3}v\|_{H^{1.5}(\Omega)}^{\frac{4}{3}}+\PP_{0}+T\PP(E^{\kappa})\\&\lesssim\delta \|\p_{t}^{3}v\|_{H^{1.5}(\Omega)}^{2}+C(\delta)(\PP_{0}+T\PP(E^{\kappa})).
    \end{aligned}
\end{align}
Combining \eqref{estimateforIV21122} and \eqref{estimateforIV21121}, we have
\begin{align}\label{estimateforIV2112}
    |IV_{21,12}|&\lesssim\delta \|\p_{t}^{3}v\|_{H^{1.5}(\Omega)}^{2}+C(\delta)(\PP_{0}+T\PP(E^{\kappa})).
\end{align}
We employ the Jacobi's formula to expand $IV_{21,11}$
\begin{align}
    \begin{aligned}
        IV_{21,11}&=-\int_{0}^{t}\int_{\Omega}\p_{t}^{3}(JA_{l}^{j}\p_{j}v^{l}) A_{i}^{k}q\p_{t}^{4}\p_{k}v^{i}\\&=\underbrace{-\int_{0}^{t}\int_{\Omega}qJA_{l}^{j}A_{i}^{k}\p_{t}^{3}\p_{j}v^{l} \p_{t}^{4}\p_{k}v^{i}}_{IV_{21,111}}\underbrace{-\sum_{s=1}^{3}C_{3}^{s}\int_{0}^{t}\int_{\Omega}A_{i}^{k}q\p_{t}^{s}(JA_{l}^{j})\p_{t}^{3-s}\p_{j}v^{l} \p_{t}^{4}\p_{k}v^{i}}_{IV_{21,112}}
    \end{aligned}
\end{align}
A similar method implies that
\begin{align}
    \begin{aligned}
        IV_{21,112}&=-\sum_{s=0}^{2}C_{3}^{s+1}\int_{\Omega}A_{i}^{k}q\p_{t}^{s}\p_{t}(JA_{l}^{j})\p_{t}^{2-s}\p_{j}v^{l} \p_{t}^{3}\p_{k}v^{i}|_{0}^{t}\\&-\sum_{s=0}^{2}C_{3}^{s+1}\int_{0}^{t}\int_{\Omega}\p_{t}qA_{i}^{k}\p_{t}^{s}(JA_{l}^{j})\p_{t}^{3-s}\p_{j}v^{l} \p_{t}^{3}\p_{k}v^{i}\\&-\sum_{s=0}^{2}C_{3}^{s+1}\int_{0}^{t}\int_{\Omega}q\p_{t}A_{i}^{k}\p_{t}^{s}(JA_{l}^{j})\p_{t}^{3-s}\p_{j}v^{l} \p_{t}^{3}\p_{k}v^{i}\\&-\sum_{s=0}^{2}C_{3}^{s+1}\int_{0}^{t}\int_{\Omega}qA_{i}^{k}\p_{t}(\p_{t}^{s}(JA_{l}^{j})\p_{t}^{3-s}\p_{j}v^{l}) \p_{t}^{3}\p_{k}v^{i},
    \end{aligned}
\end{align}
\begin{align}\label{estimateforIV21112}
    \begin{aligned}
        |IV_{21,112}|&\lesssim\PP_{0}+\sum_{s=0}^{2}C_{3}^{s+1}\|A_{i}^{k}\|_{L^{\infty}(\Omega)}\|q\|_{L^{\infty}(\Omega)}\|\p_{t}^{s}\p_{t}(JA_{l}^{j})\p_{t}^{2-s}\p_{j}v^{l}\|_{L^{2}(\Omega)} \|\p_{t}^{3}\p_{k}v^{i}\|_{L^{2}(\Omega)}\\&+\sum_{s=0}^{2}C_{3}^{s+1}\int_{0}^{t}\|\p_{t}q\|_{L^{\infty}(\Omega)}\|A_{i}^{k}\|_{L^{\infty}(\Omega)}\|\p_{t}^{s}(JA_{l}^{j})\p_{t}^{3-s}\p_{j}v^{l}\|_{L^{2}(\Omega)} \|\p_{t}^{3}\p_{k}v^{i}\|_{L^{2}(\Omega)}\\&+\sum_{s=0}^{2}C_{3}^{s+1}\int_{0}^{t}\|q\|_{L^{\infty}(\Omega)}\|\p_{t}A_{i}^{k}\|_{L^{\infty}(\Omega)}\|\p_{t}^{s}(JA_{l}^{j})\p_{t}^{3-s}\p_{j}v^{l}\|_{L^{2}(\Omega)}\| \p_{t}^{3}\p_{k}v^{i}\|_{L^{2}(\Omega)}\\&+\sum_{s=0}^{2}C_{3}^{s+1}\int_{0}^{t}\|q\|_{L^{\infty}(\Omega)}\|A_{i}^{k}\|_{L^{\infty}(\Omega)}\|\p_{t}(\p_{t}^{s}(JA_{l}^{j})\p_{t}^{3-s}\p_{j}v^{l})\|_{L^{2}(\Omega)}\| \p_{t}^{3}\p_{k}v^{i}\|_{L^{2}(\Omega)}\\&\lesssim\delta \|\p_{t}^{3}v\|_{H^{1}(\Omega)}^{2}+C(\delta)(\PP_{0}+T\PP(E^{\kappa})).
    \end{aligned}
\end{align}
In a similar reason,
\begin{align}\label{expansionofIV21111}
    \begin{aligned}
        IV_{21,111}&=-\int_{\Omega}qJA_{l}^{j}A_{i}^{k}\p_{t}^{3}\p_{j}v^{l} \p_{t}^{3}\p_{k}v^{i}|_{0}^{t}+\int_{0}^{t}\int_{\Omega}\p_{t}qJA_{l}^{j}A_{i}^{k}\p_{t}^{3}\p_{j}v^{l} \p_{t}^{3}\p_{k}v^{i}\\&+\int_{0}^{t}\int_{\Omega}q\p_{t}(JA_{l}^{j}A_{i}^{k})\p_{t}^{3}\p_{j}v^{l} \p_{t}^{3}\p_{k}v^{i}+\int_{0}^{t}\int_{\Omega}qJA_{l}^{j}A_{i}^{k}\p_{t}^{3}\p_{j}v^{l} \p_{t}^{3}\p_{k}v^{i}.
    \end{aligned}
\end{align}
The last term aligns with $-IV_{21,111}$,
\begin{align}
    \begin{aligned}
        IV_{21,111} &=-\frac{1}{2}\int_{\Omega}qJA_{l}^{j}A_{i}^{k}\p_{t}^{3}\p_{j}v^{l} \p_{t}^{3}\p_{k}v^{i}|_{0}^{t}+\frac{1}{2}\int_{0}^{t}\int_{\Omega}\p_{t}qJA_{l}^{j}A_{i}^{k}\p_{t}^{3}\p_{j}v^{l} \p_{t}^{3}\p_{k}v^{i}\\&+\frac{1}{2}\int_{0}^{t}\int_{\Omega}q\p_{t}(JA_{l}^{j}A_{i}^{k})\p_{t}^{3}\p_{j}v^{l} \p_{t}^{3}\p_{k}v^{i}.
    \end{aligned}
\end{align}
Thus,
\begin{align}\label{estimateforIV21111}
    \begin{aligned}
        |IV_{21,111}|&\lesssim\PP_{0}+\|q\|_{L^{\infty}(\Omega)}\|JA_{l}^{j}A_{i}^{k}\|_{L^{\infty}(\Omega)}\|\p_{t}^{3}\p_{j}v^{l}\|_{L^{2}(\Omega)}\| \p_{t}^{3}\p_{k}v^{i}\|_{L^{2}(\Omega)}\\&+\int_{0}^{t}\|\p_{t}q\|_{L^{\infty}(\Omega)}\|JA_{l}^{j}A_{i}^{k}\|_{L^{\infty}(\Omega)}\|\p_{t}^{3}\p_{j}v^{l}\|_{L^{2}(\Omega)} \|\p_{t}^{3}\p_{k}v^{i}\|_{L^{2}(\Omega)}\\&+\int_{0}^{t}\|q\|_{L^{\infty}(\Omega)}\|\p_{t}(JA_{l}^{j}A_{i}^{k})\|_{L^{\infty}(\Omega)}\|\p_{t}^{3}\p_{j}v^{l} \|_{L^{2}(\Omega)}\|\p_{t}^{3}\p_{k}v^{i}\|_{L^{2}(\Omega)}\\&\lesssim\delta \|\p_{t}^{3}v\|_{H^{1.5}(\Omega)}^{2}+C(\delta)(\PP_{0}+T\PP(E^{\kappa})).
    \end{aligned}
\end{align}
Combining \eqref{estimateforIV21112} and \eqref{estimateforIV21111}, we have
\begin{align}\label{estimateforIV2111}
    |IV_{21,11}|&\lesssim\delta \|\p_{t}^{3}v\|_{H^{1.5}(\Omega)}^{2}+C(\delta)(\PP_{0}+T\PP(E^{\kappa})).
\end{align}
Combining \eqref{estimateforIV2113}, \eqref{estimateforIV2112} and \eqref{estimateforIV2111}, we have
\begin{align}\label{estimateforIV211}
    |IV_{21,1}|\lesssim\delta \|\p_{t}^{3}v\|_{H^{1.5}(\Omega)}^{2}+C(\delta)(\PP_{0}+T\PP(E^{\kappa})).
\end{align}
Combining \eqref{estimateforIV212}, \eqref{estimateforIV211}, we have
\begin{align}\label{estimateforIV21}
    |IV_{21}|&\lesssim\delta \|\p_{t}^{3}v\|_{H^{1.5}(\Omega)}^{2}+\delta\int_{0}^{t}\|\sqrt{\kappa}\p_{t}^{4}v\|_{H^{1.5}(\Omega)}^{2}+\int_{0}^{t}\|\p_{t}^{3}v\|_{H^{2}(\Omega)}^{2}+C(\delta)(\PP_{0}+T\PP(E^{\kappa})).
\end{align}
Combining \eqref{estimateforIV22} and \eqref{estimateforIV21}, we have
\begin{align}\label{estimateforIV2}
    |IV_{2}|\lesssim\delta \|\p_{t}^{3}v\|_{H^{1.5}(\Omega)}^{2}+\delta\int_{0}^{t}\|\sqrt{\kappa}\p_{t}^{4}v\|_{H^{1.5}(\Omega)}^{2}+\int_{0}^{t}\|\p_{t}^{3}v\|_{H^{2}(\Omega)}^{2}+C(\delta)(\PP_{0}+T\PP(E^{\kappa})).
\end{align}
Now, we estimate $IV_{1}$, using the incompressible condition,
\begin{align}    IV_{1}=\sum_{l=1}^{4}C_{4}^{l}\int_{0}^{t}\int_{\Omega}J_{\kappa} \p_{t}^{4}q\p_{t}^{l}(\ak)_{i}^{k}\p_{t}^{4-l}\p_{k}v^{i}.
\end{align}
Thus, we can write $IV_{1}$ as
\begin{align}\label{expansionforIV1}
    \begin{aligned}
        IV_{1}&=\underbrace{\int_{0}^{t}\int_{\Omega}J_{\kappa} \p_{t}^{4}q\p_{t}^{4}(\ak)_{i}^{k}\p_{k}v^{i}}_{IV_{11}}\underbrace{+4\int_{0}^{t}\int_{\Omega}J_{\kappa} \p_{t}^{4}q\p_{t}(\ak)_{i}^{k}\p_{t}^{3}\p_{k}v^{i}}_{IV_{12}}\\&\underbrace{+\sum_{l=2}^{3}C_{4}^{l}\int_{0}^{t}\int_{\Omega}J_{\kappa} \p_{t}^{4}q\p_{t}^{l}(\ak)_{i}^{k}\p_{t}^{4-l}\p_{k}v^{i}}_{IV_{13}}
    \end{aligned}
\end{align}
We estimate $IV_{13}$ first. Integrating by part in time,
\begin{align}
   IV_{13} =+\sum_{l=2}^{3}C_{4}^{l}\int_{\Omega}J_{\kappa} \p_{t}^{3}q\p_{t}^{l}(\ak)_{i}^{k}\p_{t}^{4-l}\p_{k}v^{i}|_{0}^{t}-\sum_{l=0}^{1}C_{4}^{l+2}\int_{0}^{t}\int_{\Omega} \p_{t}^{3}q\p_{t}(J_{\kappa}\p_{t}^{l}\p_{t}^{2}(\ak)_{i}^{k}\p_{t}^{1-l}\p_{t}\p_{k}v^{i}).
\end{align}
We can show that
\begin{align}\label{estimateforIV13}
    \begin{aligned}
        |IV_{13}|&\lesssim\PP_{0}+\sum_{l=2}^{3}C_{4}^{l}\|J_{\kappa}\|_{L^{\infty}(\Omega)} \|\p_{t}^{3}q\|_{L^{2}(\Omega)}\|\p_{t}^{l}(\ak)_{i}^{k}\p_{t}^{4-l}\p_{k}v^{i}\|_{L^{2}(\Omega)}\\&+\sum_{l=0}^{1}C_{4}^{l+2}\int_{0}^{t} \|\p_{t}^{3}q\|_{L^{2}(\Omega)}\|\p_{t}(J_{\kappa}\p_{t}^{l}\p_{t}^{2}(\ak)_{i}^{k}\p_{t}^{1-l}\p_{t}\p_{k}v^{i})\|_{L^{2}(\Omega)}\\&\lesssim\delta\|\p_{t}^{3}q\|_{L^{2}(\Omega)}+C(\delta)(\PP_{0}+T\PP(E^{\kappa})).
    \end{aligned}
\end{align}
We deal with the lower order part in $IV_{11}$ in a similar way:
\begin{align}\label{expansionforIV11}
    \begin{aligned}
        IV_{11}&=\underbrace{-\int_{0}^{t}\int_{\Omega}J_{\kappa} \p_{t}^{4}q\p_{t}^{3}\p_{l}v_{\kappa}^{j}(\ak)_{i}^{l}(\ak)_{j}^{k}\p_{k}v^{i}}_{IV_{11,1}}\underbrace{-\sum_{s=1}^{3}C_{3}^{s}\int_{0}^{t}\int_{\Omega}J_{\kappa} \p_{t}^{4}q\p_{t}^{s}((\ak)_{i}^{l}(\ak)_{j}^{k})\p_{t}^{3-s}\p_{l}v_{\kappa}^{j}\p_{k}v^{i}}_{IV_{11,2}}
    \end{aligned}
\end{align}
Integrating by part in time, 
\begin{align}
    \begin{aligned}
     IV_{11,2}   &=\sum_{s=1}^{3}C_{3}^{s}\int_{\Omega}J_{\kappa} \p_{t}^{3}q\p_{t}^{s}((\ak)_{i}^{l}(\ak)_{j}^{k})\p_{t}^{3-s}\p_{l}v_{\kappa}^{j}\p_{k}v^{i}|_{0}^{t}\\&-\sum_{s=1}^{3}C_{3}^{s}\int_{0}^{t}\int_{\Omega} \p_{t}^{3}q\p_{t}(\p_{t}^{s}((\ak)_{i}^{l}(\ak)_{j}^{k})\p_{t}^{3-s}\p_{l}v_{\kappa}^{j}\p_{k}v^{i}J_{\kappa})
    \end{aligned}
\end{align}
using a similar estimate, it follows that
\begin{align}\label{estimateforIV112}
    \begin{aligned}
     |IV_{11,2}|   &\lesssim\PP_{0}+\sum_{s=1}^{3}C_{3}^{s}\|J_{\kappa}\p_{k}v^{i}\|_{L^{\infty}(\Omega)}\| \p_{t}^{3}q\|_{L^{2}(\Omega)}\|\p_{t}^{s}((\ak)_{i}^{l}(\ak)_{j}^{k})\p_{t}^{3-s}\p_{l}v_{\kappa}^{j}\|_{L^{2}(\Omega)}\\&+\sum_{s=1}^{3}C_{3}^{s}\int_{0}^{t} \|\p_{t}^{3}q\|_{L^{2}(\Omega)}\|\p_{t}(\p_{t}^{s}((\ak)_{i}^{l}(\ak)_{j}^{k})\p_{t}^{3-s}\p_{l}v_{\kappa}^{j})\|_{L^{2}(\Omega)}\|\p_{k}v^{i}J_{\kappa}\|_{L^{\infty}(\Omega)}\\&+\sum_{s=1}^{3}C_{3}^{s}\int_{0}^{t} \|\p_{t}^{3}q\|_{L^{2}(\Omega)}\|\p_{t}^{s}((\ak)_{i}^{l}(\ak)_{j}^{k})\p_{t}^{3-s}\p_{l}v_{\kappa}^{j}\|_{L^{2}(\Omega)}\|\p_{t}(\p_{k}v^{i}J_{\kappa})\|_{L^{\infty}(\Omega)}\\&\lesssim\delta\|\p_{t}^{3}q\|_{L^{2}(\Omega)}+C(\delta)(\PP_{0}+T\PP(E^{\kappa})).
    \end{aligned}
\end{align}
Combining \eqref{expansionforIV1}, \eqref{estimateforIV13}, \eqref{expansionforIV11}, and \eqref{estimateforIV112}, we can summarize that
\begin{align}\label{expansionforIV1another}
    IV_{1}=\underbrace{-\int_{0}^{t}\int_{\Omega}J_{\kappa} \p_{t}^{4}q\p_{t}^{3}\p_{l}v_{\kappa}^{j}(\ak)_{i}^{l}(\ak)_{j}^{k}\p_{k}v^{i}}_{IV_{11,1}}\underbrace{-4\int_{0}^{t}\int_{\Omega}J_{\kappa} \p_{t}^{4}q(\ak)_{j}^{k}\p_{k}v_{\kappa}^{i}(\ak)_{i}^{l}\p_{t}^{3}\p_{l}v^{j}}_{IV_{12}}+\mathcal{R},
\end{align}
where $\mathcal{R}\lesssim\delta\|\p_{t}^{3}q\|_{L^{2}(\Omega)}+C(\delta)(\PP_{0}+T\PP(E^{\kappa}))$. We deal with these two terms in a similar way. We only show the $IV_{12}$ as an example
\begin{align}\label{expansionforIV12}
    \begin{aligned}
        IV_{12}=\underbrace{-4\int_{0}^{t}\int_{\Omega} \p_{t}^{4}(J_{\kappa}(\ak)_{i}^{l}q)(\ak)_{j}^{k}\p_{k}v_{\kappa}^{i}\p_{t}^{3}\p_{l}v^{j}}_{IV_{12,1}}\underbrace{+4\sum_{s=1}^{4}C_{4}^{s}\int_{0}^{t}\int_{\Omega} \p_{t}^{s}(J_{\kappa}(\ak)_{i}^{l})\p_{t}^{4-s}q(\ak)_{j}^{k}\p_{k}v_{\kappa}^{i}\p_{t}^{3}\p_{l}v^{j}}_{IV_{12,2}}.
    \end{aligned}
\end{align}
Clearly, we can estimate $IV_{12,2}$ as
\begin{align}\label{estimateforIV122}
    \begin{aligned}
        |IV_{12,2}|&\lesssim\int_{0}^{t} \|\p_{t}(J_{\kappa}(\ak)_{i}^{l})(\ak)_{j}^{k}\p_{k}v_{\kappa}^{i}\|_{L^{\infty}(\Omega)}\|\p_{t}^{3}q\|_{L^{2}(\Omega)}\|\p_{t}^{3}\p_{l}v^{j}\|_{L^{2}(\Omega)}\\&+4\sum_{s=2}^{4}C_{4}^{s}\int_{0}^{t}\| \p_{t}^{s}(J_{\kappa}(\ak)_{i}^{l})\|_{L^{2}(\Omega)}\|\p_{t}^{4-s}q(\ak)_{j}^{k}\p_{k}v_{\kappa}^{i}\|_{L^{\infty}(\Omega)}\|\p_{t}^{3}\p_{l}v^{j}\|_{L^{2}(\Omega)}\\&\lesssim\PP_{0}+T\PP(E^{\kappa}).
    \end{aligned}
\end{align}
For $IV_{12,1}$, integrating by part
\begin{align}\label{expansionforIV121}
    \begin{aligned}
        IV_{12,1}&=\underbrace{4\int_{0}^{t}\int_{\Gamma_{1}} \p_{t}^{4}(\sigma (\sqrt{h}\lap_{h} (\eta)\cdot  n ))(\ak)_{j}^{k}\p_{k}v_{\kappa}^{i}\p_{t}^{3}v^{j}}_{IV_{12,11}}\\&\underbrace{-4\int_{0}^{t}\int_{\Gamma_{1}} \p_{t}^{4}(\kappa(1-\bar{\lap})(v\cdot  n _{\kappa}))(\ak)_{j}^{k}\p_{k}v_{\kappa}^{i}\p_{t}^{3}v^{j}}_{IV_{12,12}}\\&\underbrace{+4\int_{0}^{t}\int_{\Omega} \p_{t}^{4}(J_{\kappa}(\ak)_{i}^{l}\p_{l}q)(\ak)_{j}^{k}\p_{k}v_{\kappa}^{i}\p_{t}^{3}v^{j}}_{IV_{12,13}}\\&\underbrace{+4\int_{0}^{t}\int_{\Omega} \p_{t}^{4}(J_{\kappa}(\ak)_{i}^{l}q)\p_{l}((\ak)_{j}^{k}\p_{k}v_{\kappa}^{i})\p_{t}^{3}v^{j}}_{IV_{12,14}}
    \end{aligned}
\end{align}
We control $IV_{12,11}$ by
\begin{align}
\begin{aligned}
     IV_{12,11}&=-4\sigma\int_{0}^{t}\int_{\Gamma_{1}} \p_{t}^{4}(h^{\mu\nu}\p_{\nu} \eta)\cdot  n (\ak)_{j}^{k}\p_{k}v_{\kappa}^{i}\p_{\mu}\p_{t}^{3}v^{j}\\&-4\sigma\int_{0}^{t}\int_{\Gamma_{1}} \p_{t}^{4}(h^{\mu\nu}\p_{\nu} \eta)\cdot \p_{\mu}( n (\ak)_{j}^{k}\p_{k}v_{\kappa}^{i})\p_{t}^{3}v^{j}\\&+4\sigma\sum_{s=0}^{3}C_{4}^{s}\int_{0}^{t}\int_{\Gamma_{1}} \p_{t}^{s} \p_{\mu}(h^{\mu\nu}\p_{\nu} \eta)\cdot \p_{t}^{4-s} n (\ak)_{j}^{k}\p_{k}v_{\kappa}^{i}\p_{t}^{3}v^{j},
\end{aligned} 
\end{align}
\begin{align}\label{estimateforIV1211}
    \begin{aligned}
        |IV_{12,11}|&\lesssim4\sigma\int_{0}^{t}\| \p_{t}^{4}(h^{\mu\nu}\p_{\nu} \eta)\|_{L^{2}(\Gamma_{1})}\|  n (\ak)_{j}^{k}\p_{k}v_{\kappa}^{i}\|_{L^{\infty}(\Omega)}\|\p_{\mu}\p_{t}^{3}v^{j}\|_{L^{2}(\Gamma_{1}}\\&+4\sigma\int_{0}^{t} \|\p_{t}^{4}(h^{\mu\nu}\p_{\nu} \eta)\|_{L^{2}(\Gamma_{1})}\cdot \|\p_{\mu}( n (\ak)_{j}^{k}\p_{k}v_{\kappa}^{i})\|_{L^{\infty}(\Omega)}\|\p_{t}^{3}v^{j}\|_{L^{2}(\Omega)}\\&+4\sigma\sum_{s=0}^{3}C_{4}^{s}\int_{0}^{t}\|\p_{t}^{s} \p_{\mu}(h^{\mu\nu}\p_{\nu} \eta)\cdot \p_{t}^{4-s} n \|_{L^{2}(\Gamma_{1})}\|(\ak)_{j}^{k}\p_{k}v_{\kappa}^{i}\|_{L^{\infty}(\Omega)}\|\p_{t}^{3}v^{j}\|_{L^{2}(\Omega)}\\&\lesssim\PP_{0}+T\PP(E^{\kappa})
    \end{aligned}
\end{align}
Similarly,
\begin{align}
    \begin{aligned}
        IV_{12,12}&=-4\int_{0}^{t}\int_{\Gamma_{1}} \p_{t}^{4}(\kappa(v\cdot  n _{\kappa}))(\ak)_{j}^{k}\p_{k}v_{\kappa}^{i}\p_{t}^{3}v^{j}\\&-4\int_{0}^{t}\int_{\Gamma_{1}}(\kappa\pb_{l} \p_{t}^{4}(v\cdot  n _{\kappa}))\pb_{l}((\ak)_{j}^{k}\p_{k}v_{\kappa}^{i})\p_{t}^{3}v^{j}\\&-4\int_{0}^{t}\int_{\Gamma_{1}}(\kappa\pb_{l} \p_{t}^{4}(v\cdot  n _{\kappa}))(\ak)_{j}^{k}\p_{k}v_{\kappa}^{i}\pb_{l}\p_{t}^{3}v^{j},
    \end{aligned}
\end{align}
Since 
\begin{align}
    \begin{aligned}
        \|\sqrt{\kappa}\p_{t}^{4}(v\cdot n_{\kappa})\|_{H^{1}(\Gamma_{1})}&\lesssim\|\sqrt{\kappa}\p_{t}^{4}v\|_{H^{1}(\Gamma_{1})}\|n_{\kappa}\|_{H^{2}(\Gamma_{1})}+\|v\|_{H^{2}(\Gamma_{1})}\|\sqrt{\kappa}\p_{t}^{4}n_{\kappa}\|_{H^{1}(\Gamma_{1})}+\sum_{k=0}^{2}\|\sqrt{\kappa}\p_{t}^{k}\p_{t}v\p_{t}^{2-k}\p_{t}n_{\kappa}\|_{H^{1}(\Gamma_{1})}\\&\lesssim(\|\sqrt{\kappa}\p_{t}^{4}v\|_{H^{1}(\Gamma_{1})}+\|\sqrt{\kappa}\p_{t}^{3}v\|_{H^{2}(\Omega)})(\PP_{0}+T\PP(E^{\kappa}))+\PP(E^{\kappa}),
    \end{aligned}
\end{align}
we can calculate that
\begin{align}\label{estimateforIV1212}
    \begin{aligned}
        |IV_{12,12}|&\lesssim4\int_{0}^{t}\sqrt{\kappa} \|\p_{t}^{4}(\sqrt{\kappa}(v\cdot  n _{\kappa}))\|_{L^{2}(\Gamma_{1})}\|(\ak)_{j}^{k}\p_{k}v_{\kappa}^{i}\|_{L^{\infty}(\Gamma_{1})}\|\p_{t}^{3}v^{j}\|_{L^{2}(\Gamma_{1})}\\&+4\int_{0}^{t}\sqrt{\kappa}(\|\sqrt{\kappa}\pb_{l} \p_{t}^{4}(v\cdot  n _{\kappa}))\|_{L^{2}(\Gamma_{1})}\|\pb_{l}((\ak)_{j}^{k}\p_{k}v_{\kappa}^{i})\|_{L^{\infty}(\Gamma_{1})}\|\p_{t}^{3}v^{j}\|_{L^{2}(\Gamma_{1})}\\&+4\int_{0}^{t}\sqrt{\kappa}\|\sqrt{\kappa}\pb_{l} \p_{t}^{4}(v\cdot  n _{\kappa})\|_{L^{2}(\Gamma_{1}}\|(\ak)_{j}^{k}\p_{k}v_{\kappa}^{i}\|_{L^{\infty}(\Gamma_{1})}\|\pb_{l}\p_{t}^{3}v^{j}\|_{L^{2}(\Gamma_{1})}\\&\lesssim\PP_{0}+T\PP(E^{\kappa}).
    \end{aligned}
\end{align}
Now, plugging Euler equation into $IV_{12,13}$, we obtain
\begin{align}\label{expansionforIV1213}
    \begin{aligned}
        IV_{12,13}&=\underbrace{-4\int_{0}^{t}\int_{\Omega} \p_{t}^{4}(J_{\kappa}\p_{t}v_{i})(\ak)_{j}^{k}\p_{k}v_{\kappa}^{i}\p_{t}^{3}v^{j}}_{IV_{12,131}}\underbrace{-4\int_{0}^{t}\int_{\Omega} \p_{t}^{4}(J_{\kappa}\p_{t}\phi^{\kappa}\nn_{i}\phi^{\kappa})(\ak)_{j}^{k}\p_{k}v_{\kappa}^{i}\p_{t}^{3}v^{j}}_{IV_{12,132}}\\&\underbrace{-4\int_{0}^{t}\int_{\Omega} \p_{t}^{4}(J_{\kappa}\nn_{l}\nn_{i}\phi^{\kappa}\nn_{l}\phi^{\kappa})(\ak)_{j}^{k}\p_{k}v_{\kappa}^{i}\p_{t}^{3}v^{j}}_{IV_{12,133}}
    \end{aligned}
\end{align}
Using our estimate for $\phi$,
\begin{align}
    \begin{aligned}
        IV_{12,133}&=-4\int_{\Omega} \p_{t}^{3}\nn_{l}\nn_{i}\phi^{\kappa} J_{\kappa}\nn_{l}\phi^{\kappa}(\ak)_{j}^{k}\p_{k}v_{\kappa}^{i}\p_{t}^{3}v^{j}|_{0}^{t}\\&+4\int_{0}^{t}\int_{\Omega} \p_{t}^{3}\nn_{l}\nn_{i}\phi^{\kappa}\p_{t}( J_{\kappa}\nn_{l}\phi^{\kappa}(\ak)_{j}^{k}\p_{k}v_{\kappa}^{i})\p_{t}^{3}v^{j}\\&+4\int_{0}^{t}\int_{\Omega} \p_{t}^{3}\nn_{l}\nn_{i}\phi^{\kappa}\p_{t}( J_{\kappa}\nn_{l}\phi^{\kappa}(\ak)_{j}^{k}\p_{k}v_{\kappa}^{i})\p_{t}^{3}v^{j}\\&-\sum_{s=0}^{3}4\int_{0}^{t}\int_{\Omega} \p_{t}^{s}(\nn_{l}\nn_{i}\phi^{\kappa})\p_{t}^{4-s}(J_{\kappa}\nn_{l}\phi^{\kappa})(\ak)_{j}^{k}\p_{k}v_{\kappa}^{i}\p_{t}^{3}v^{j},
    \end{aligned}
\end{align}
\begin{align}\label{estimateforIV12133}
    \begin{aligned}
        |IV_{12,133}|&\lesssim\PP_{0}+4\| \p_{t}^{3}\nn_{l}\nn_{i}\phi^{\kappa}\|_{L^{2}(\Omega)}\| J_{\kappa}\nn_{l}\phi^{\kappa}(\ak)_{j}^{k}\p_{k}v_{\kappa}^{i}\|_{L^{\infty}(\Omega)}\|\p_{t}^{3}v^{j}\|_{L^{2}(\Omega)}\\&+4\int_{0}^{t}\|\p_{t}^{3}\nn_{l}\nn_{i}\phi^{\kappa}\|_{L^{2}(\Omega)}\|\p_{t}( J_{\kappa}\nn_{l}\phi^{\kappa}(\ak)_{j}^{k}\p_{k}v_{\kappa}^{i})\|_{L^{\infty}(\Omega)}\|\p_{t}^{3}v^{j}\|_{L^{2}(\Omega)}\\&+4\int_{0}^{t} \|\p_{t}^{3}\nn_{l}\nn_{i}\phi^{\kappa}\|_{L^{2}(\Omega}\| J_{\kappa}\nn_{l}\phi^{\kappa}(\ak)_{j}^{k}\p_{k}v_{\kappa}^{i}\|_{L^{\infty}(\Omega)}\|\p_{t}^{4}v^{j}\|_{L^{2}(\Omega)}\\&+\sum_{s=0}^{3}4\int_{0}^{t}\| \p_{t}^{s}(\nn_{l}\nn_{i}\phi^{\kappa})\p_{t}^{3-s}\p_{t}(J_{\kappa}\nn_{l}\phi^{\kappa})\|_{L^{2}(\Omega)}\|(\ak)_{j}^{k}\p_{k}v_{\kappa}^{i}\|_{L^{\infty}(\Omega)}\|\p_{t}^{3}v^{j}\|_{L^{2}(\Omega}\\&\lesssim\PP(E_{\kappa}^{\prime})+\PP_{0}+T\PP(E^{\kappa}).
    \end{aligned}
\end{align}
Similarly.
\begin{align}
    \begin{aligned}
          IV_{12,132}&=-4\int_{\Omega} \p_{t}^{3}\p_{t}\phi^{\kappa} J_{\kappa}\nn_{i}\phi^{\kappa}(\ak)_{j}^{k}\p_{k}v_{\kappa}^{i}\p_{t}^{3}v^{j}|_{0}^{t}\\&+4\int_{0}^{t}\int_{\Omega} \p_{t}^{4}\phi^{\kappa}\p_{t}( J_{\kappa}\nn_{i}\phi^{\kappa}(\ak)_{j}^{k}\p_{k}v_{\kappa}^{i})\p_{t}^{3}v^{j}\\&+4\int_{0}^{t}\int_{\Omega} \p_{t}^{4}\phi^{\kappa}\p_{t}( J_{\kappa}\nn_{i}\phi^{\kappa}(\ak)_{j}^{k}\p_{k}v_{\kappa}^{i})\p_{t}^{3}v^{j}\\&-\sum_{s=0}^{3}4\int_{0}^{t}\int_{\Omega} \p_{t}^{s}(\p_{t}\phi^{\kappa})\p_{t}^{4-s}(J_{\kappa}\nn_{i}\phi^{\kappa})(\ak)_{j}^{k}\p_{k}v_{\kappa}^{i}\p_{t}^{3}v^{j},
    \end{aligned}
\end{align}
\begin{align}\label{estimateforIV12132}
    \begin{aligned}
       |IV_{12,132}|&\lesssim\PP_{0}+4\| \p_{t}^{4}\phi^{\kappa}\|_{L^{2}(\Omega)}\| J_{\kappa}\nn_{i}\phi^{\kappa}(\ak)_{j}^{k}\p_{k}v_{\kappa}^{i}\|_{L^{\infty}(\Omega)}\|\p_{t}^{3}v^{j}\|_{L^{2}(\Omega)}\\&+4\int_{0}^{t} \|\p_{t}^{4}\phi^{\kappa}\|_{L^{2}(\Omega)}\|\p_{t}( J_{\kappa}\nn_{i}\phi^{\kappa}(\ak)_{j}^{k}\p_{k}v_{\kappa}^{i})\|_{L^{\infty}(\Omega)}\|\p_{t}^{3}v^{j}\|_{L^{2}(\Omega)}\\&+4\int_{0}^{t} \|\p_{t}^{4}\phi^{\kappa}\|_{L^{2}(\Omega}\| J_{\kappa}\nn_{i}\phi^{\kappa}(\ak)_{j}^{k}\p_{k}v_{\kappa}^{i}\|_{L^{\infty}(\Omega)}\|\p_{t}^{4}v^{j}\|_{L^{2}(\Omega)}\\&+\sum_{s=0}^{3}4\int_{0}^{t}\| \p_{t}^{s}(\p_{t}\phi^{\kappa})\p_{t}^{3-s}\p_{t}(J_{\kappa}\nn_{l}\phi^{\kappa})\|_{L^{2}(\Omega)}\|(\ak)_{j}^{k}\p_{k}v_{\kappa}^{i}\|_{L^{\infty}(\Omega)}\|\p_{t}^{3}v^{j}\|_{L^{2}(\Omega}\\&\lesssim\PP_{0}+T\PP(E^{\kappa}). 
    \end{aligned}
\end{align}
After integrating by part in time, a similar argument shows that
\begin{align}
    \begin{aligned}
        IV_{12,131}&=-4\int_{\Omega} J_{\kappa}\p_{t}^{4}v_{i}(\ak)_{j}^{k}\p_{k}v_{\kappa}^{i}\p_{t}^{3}v^{j}|_{0}^{t}+4\int_{0}^{t}\int_{\Omega}\p_{t}^{4}v_{i} \p_{t}(J_{\kappa}(\ak)_{j}^{k}\p_{k}v_{\kappa}^{i})\p_{t}^{3}v^{j}+4\int_{0}^{t}\int_{\Omega}\p_{t}^{4}v_{i} J_{\kappa}(\ak)_{j}^{k}\p_{k}v_{\kappa}^{i}\p_{t}^{4}v^{j}\\&-4\sum_{s=1}^{4}C_{4}^{s}\int_{0}^{t}\int_{\Omega} \p_{t}^{s}J_{\kappa}\p_{t}^{4-s}\p_{t}v_{i}(\ak)_{j}^{k}\p_{k}v_{\kappa}^{i}\p_{t}^{3}v^{j},
    \end{aligned}
\end{align}
\begin{align}\label{estimateforIV12131}
    \begin{aligned}
        |IV_{12,131}|&\lesssim\PP_{0}+4 \|J_{\kappa}(\ak)_{j}^{k}\p_{k}v_{\kappa}^{i}\|_{L^{\infty}(\Omega)}\|\p_{t}^{4}v_{i}\|_{L^{2}(\Omega)}\|\p_{t}^{3}v^{j}\|_{L^{2}(\Omega)}\\&+4\int_{0}^{t}\|\p_{t}^{4}v_{i}\|_{L^{2}(\Omega)}\| \p_{t}(J_{\kappa}(\ak)_{j}^{k}\p_{k}v_{\kappa}^{i})\|_{L^{\infty}(\Omega)}\|\p_{t}^{3}v^{j}\|_{L^{2}(\Omega)}\\&+4\int_{0}^{t}\|\p_{t}^{4}v_{i}\|_{L^{2}(\Omega)}\| J_{\kappa}(\ak)_{j}^{k}\p_{k}v_{\kappa}^{i}\|_{L^{\infty}(\Omega}\|\p_{t}^{4}v^{j}\|_{L6{2}(\Omega)}\\&+4\sum_{s=0}^{3}C_{4}^{s+1}\int_{0}^{t} \|\p_{t}^{s}\p_{t}J_{\kappa}\p_{t}^{3-s}\p_{t}v_{i}\|_{L^{2}(\Omega)}\|(\ak)_{j}^{k}\p_{k}v_{\kappa}^{i}\|_{L^{\infty}(\Omega)}\|\p_{t}^{3}v^{j}\|_{L^{2}(\Omega)}\\&\lesssim\delta\|\p_{t}^{4}v\|_{L^{2}(\Omega)}^{2}+C(\delta)(\PP_{0}+T\PP(E^{\kappa})).
    \end{aligned}
\end{align}
Combining \eqref{expansionforIV1213}, \eqref{estimateforIV12133}, \eqref{estimateforIV12132}, \eqref{estimateforIV12131}, we have
\begin{align}   \label{estimateforIV1213} |IV_{1213}|\lesssim\PP(E_{\kappa}^{\prime})+\delta\|\p_{t}^{4}v\|_{L^{2}(\Omega)}^{2}+C(\delta)(\PP_{0}+T\PP(E^{\kappa})).
\end{align}
The estimate for $IV_{12,14}$ is in a similar fashion. Integrating by part in time,
\begin{align}
\begin{aligned}
     IV_{12,14}&=4\int_{\Omega} J_{\kappa}(\ak)_{i}^{l}\p_{t}^{3}q\p_{l}((\ak)_{j}^{k}\p_{k}v_{\kappa}^{i})\p_{t}^{3}v^{j}|_{0}^{t}\\&-4\int_{0}^{t}\int_{\Omega}\p_{t}^{3}q \p_{t}(J_{\kappa}(\ak)_{i}^{l}\p_{l}((\ak)_{j}^{k}\p_{k}v_{\kappa}^{i}))\p_{t}^{3}v^{j}\\&-4\int_{0}^{t}\int_{\Omega} \p_{t}^{3}q J_{\kappa}(\ak)_{i}^{l}\p_{l}((\ak)_{j}^{k}\p_{k}v_{\kappa}^{i})\p_{t}^{4}v^{j}\\&+4\int_{0}^{t}\int_{\Omega} \p_{t}^{4}(J_{\kappa}(\ak)_{i}^{l})q\p_{l}((\ak)_{j}^{k}\p_{k}v_{\kappa}^{i})\p_{t}^{3}v^{j}\\&+4\sum_{s=1}^{3}C_{4}^{s}\int_{0}^{t}\int_{\Omega} \p_{t}^{s}(J_{\kappa}(\ak)_{i}^{l})\p_{t}^{4-s}q\p_{l}((\ak)_{j}^{k}\p_{k}v_{\kappa}^{i})\p_{t}^{3}v^{j},
\end{aligned}   
\end{align}
\begin{align}\label{estimateforIV1214}
    \begin{aligned}
        |IV_{12,14}|&\lesssim\PP_{0}+4\|\p_{t}^{3}q\|_{L^{6}(\Omega)} \|J_{\kappa}(\ak)_{i}^{l}\p_{l}((\ak)_{j}^{k}\p_{k}v_{\kappa}^{i})\|_{L^{3}(\Omega)}\|\p_{t}^{3}v^{j}\|_{L^{2}(\Omega)}\\&+4\int_{0}^{t}\|\p_{t}^{3}q\|_{L^{6}(\Omega)}\| \p_{t}(J_{\kappa}(\ak)_{i}^{l}\p_{l}((\ak)_{j}^{k}\p_{k}v_{\kappa}^{i}))\|_{L^{3}(\Omega)}\|\p_{t}^{3}v^{j}\|_{L^{2}((\Omega)}\\&+4\int_{0}^{t}\|\p_{t}^{3}q\|_{L^{2}(\Omega)}\| J_{\kappa}(\ak)_{i}^{l}\p_{l}((\ak)_{j}^{k}\p_{k}v_{\kappa}^{i})\|_{L^{\infty}(\Omega)}\|\p_{t}^{4}v^{j}\|_{L^{2}(\Omega)}\\&+4\int_{0}^{t} \|\p_{t}^{4}(J_{\kappa}(\ak)_{i}^{l})\|_{L^{2}(\Omega)}\|q\p_{l}((\ak)_{j}^{k}\p_{k}v_{\kappa}^{i})\|_{L^{\infty}(\Omega)}\|\p_{t}^{3}v^{j}\|_{L^{2}(\Omega)}\\&+4\sum_{s=1}^{3}C_{4}^{s}\int_{0}^{t}\| \p_{t}^{s}(J_{\kappa}(\ak)_{i}^{l})\|_{L^{3}(\Omega)}\|\p_{t}^{4-s}q\|_{L^{2}(\Omega)}\|\p_{l}((\ak)_{j}^{k}\p_{k}v_{\kappa}^{i})\|_{L^{\infty}(\Omega)}\|\p_{t}^{3}v^{j}\|_{L^{2}(\Omega)}\\&\lesssim\delta\|\p_{t}^{3}q\|_{H^{1}(\Omega)}^{2}+C(\delta)(\PP_{0}+T\PP(E^{\kappa})).
    \end{aligned}
\end{align}
Hence, combining \eqref{expansionforIV121}, \eqref{estimateforIV1211}, \eqref{estimateforIV1212}, \eqref{estimateforIV1213} and \eqref{estimateforIV1214}, we have
\begin{align}\label{estimateforIV121}    |IV_{121}|\lesssim\PP(E_{\kappa}^{\prime})+\delta\|\p_{t}^{4}v\|_{L^{2}(\Omega)}^{2}+\delta\|\p_{t}^{3}q\|_{H^{1}(\Omega)}^{2}+C(\delta)(\PP_{0}+T\PP(E^{\kappa})).
\end{align}
Combining \eqref{expansionforIV12}, \eqref{estimateforIV122} and \eqref{estimateforIV121}, we have
\begin{align}\label{estimateforIV12}
    |IV_{12}|\lesssim\PP(E_{\kappa}^{\prime})+\delta\|\p_{t}^{4}v\|_{L^{2}(\Omega)}^{2}+\delta\|\p_{t}^{3}q\|_{H^{1}(\Omega)}^{2}+C(\delta)(\PP_{0}+T\PP(E^{\kappa})).
\end{align}
Since $IV_{11,1}$ has a similar structure as $IV_{12}$, we have
\begin{align}\label{estimateforIV111}
    |IV_{11,1}|\lesssim\PP(E_{\kappa}^{\prime})+\delta\|\p_{t}^{4}v\|_{L^{2}(\Omega)}^{2}+\delta\|\p_{t}^{3}q\|_{H^{1}(\Omega)}^{2}+C(\delta)(\PP_{0}+T\PP(E^{\kappa})).
\end{align}

Combining \eqref{expansionforIV1another}, \eqref{estimateforIV12} and \eqref{estimateforIV111}, we have
\begin{align}\label{estimateforIV1}
    |IV_{1}|\lesssim\PP(E_{\kappa}^{\prime})+\delta\|\p_{t}^{4}v\|_{L^{2}(\Omega)}^{2}+\delta\|\p_{t}^{3}q\|_{H^{1}(\Omega)}^{2}+C(\delta)(\PP_{0}+T\PP(E^{\kappa})).
\end{align}
Combining \eqref{expansionforIV}, \eqref{estimateforIV2} and \eqref{estimateforIV1}, we have
\begin{align}
    |IV|\lesssim\PP(E_{\kappa}^{\prime})+\delta\|\p_{t}^{4}v\|_{L^{2}(\Omega)}^{2}+\delta\|\p_{t}^{3}q\|_{H^{1}(\Omega)}^{2}+\delta \|\p_{t}^{3}v\|_{H^{1.5}(\Omega)}^{2}+\delta\int_{0}^{t}\|\sqrt{\kappa}\p_{t}^{4}v\|_{H^{1.5}(\Omega)}^{2}+C(\delta)(\PP_{0}+T\PP(E^{\kappa})).
\end{align}
Now, invoking \eqref{estimateforqtttpre}, choosing appropriate $\delta$,
we can write
\begin{align}\label{estimateforIV}
\begin{aligned}
     |IV|&\lesssim  \PP(E_{\kappa}^{\prime})+\delta\|\p_{t}^{4}v\|_{L^{2}(\Omega)}^{2}+\delta \|\p_{t}^{3}v\|_{H^{1.5}(\Omega)}^{2}+\delta\int_{0}^{t}\|\sqrt{\kappa}\p_{t}^{4}v\|_{H^{1.5}(\Omega)}^{2}+C(\delta)(\PP_{0}+T\PP(E^{\kappa})).
\end{aligned}   
\end{align}
It's left for us to control $II$.
\begin{align}\label{expansionforII}
    II=\underbrace{\int_{0}^{t}\int_{\Omega}-\p_{t}^{4}(J_{\kappa}\phi^{\kappa}_{t}\nn_{i}\phi^{\kappa})\p_{t}^{4}v^{i}}_{II_{1}}\underbrace{+\int_{0}^{t}\int_{\Omega}-\p_{t}^{4}(J_{\kappa}\nn_{j}\nn_{i}\phi^{\kappa}\nn_{j}\phi^{\kappa})\p_{t}^{4}v^{i}}_{II_{2}}
\end{align}
We expand $II_{1}$ as
\begin{align}
    \begin{aligned}
        II_{1}&=\int_{0}^{t}\int_{\Omega}-J_{\kappa}\p_{t}^{5}\phi^{\kappa}\nn_{i}\phi^{\kappa}\p_{t}^{4}v^{i}-\sum_{s=1}^{4}C_{4}^{s}\int_{0}^{t}\int_{\Omega}\p_{t}^{s}(J_{\kappa}\nn_{i}\phi^{\kappa})\p_{t}^{4-s}\phi^{\kappa}_{t}\p_{t}^{4}v^{i}.
    \end{aligned}
\end{align}
From direct calculation, we have
\begin{align}\label{estimateforII1}
\begin{aligned}
    |II_{1}|&\lesssim\int_{0}^{t}\|J_{\kappa}\nn_{i}\phi^{\kappa}\|_{L^{\infty}(\Omega)}\|\p_{t}^{5}\phi^{\kappa}\|_{L^{2}(\Omega)}\|\p_{t}^{4}v^{i}\|_{L^{2}(\Omega)}\\&+\sum_{s=0}^{3}C_{4}^{s+1}\int_{0}^{t}\|\p_{t}^{s}\p_{t}(J_{\kappa}\nn_{i}\phi^{\kappa})\p_{t}^{3-s}\phi^{\kappa}_{t}\|_{L^{2}(\Omega)}\|\p_{t}^{4}v^{i}\|_{L^{2}(\Omega)}\\&\lesssim\delta\int_{0}^{t}\|\p_{t}^{5}\phi^{\kappa}\|_{L^{2}(\Omega)}^{2}+C(\delta)(\PP_{0}+T\PP(E^{\kappa})).
\end{aligned}
\end{align}
For $II_{2}$, we take its leading order as
\begin{align}\label{expansionforII2}
    \begin{aligned}
        II_{2}&=\int_{0}^{t}\int_{\Omega}-J_{\kappa}\p_{t}^{4}\nn_{j}\nn_{i}\phi^{\kappa}\nn_{j}\phi^{\kappa}\p_{t}^{4}v^{i}-\sum_{s=1}^{4}C_{4}^{s}\int_{0}^{t}\int_{\Omega}\p_{t}^{s}(J_{\kappa}\nn_{j}\phi^{\kappa})\p_{t}^{4-s}\nn_{j}\nn_{i}\phi^{\kappa}\p_{t}^{4}v^{i}\\&=\underbrace{-\int_{0}^{t}\int_{\Omega}(J_{\kappa}\nn_{j}\p_{t}^{4}\nn_{i}\phi^{\kappa})\nn_{j}\phi^{\kappa}\p_{t}^{4}v^{i}}_{II_{21}}\underbrace{-\sum_{s=1}^{4}C_{4}^{s}\int_{0}^{t}\int_{\Omega}J_{\kappa}\p_{t}^{s}(\ak)_{j}^{l}\p_{t}^{4-s}\p_{l}\nn_{i}\phi^{\kappa}\nn_{j}\phi^{\kappa}\p_{t}^{4}v^{i}}_{II_{22}}\\&\underbrace{-\sum_{s=1}^{4}C_{4}^{s}\int_{0}^{t}\int_{\Omega}\p_{t}^{s}(J_{\kappa}\nn_{j}\phi^{\kappa})\p_{t}^{4-s}\nn_{j}\nn_{i}\phi^{\kappa}\p_{t}^{4}v^{i}}_{II_{23}}.
    \end{aligned}
\end{align}
 The control of $II_{22}$ and $II_{23}$ reads
\begin{align}\label{estimateforII22}
    \begin{aligned}
        |II_{22}|&\lesssim\sum_{s=0}^{3}C_{4}^{s+1}\int_{0}^{t}\|J_{\kappa}\nn_{j}\phi^{\kappa}\|_{L^{\infty}(\Omega)}\|\p_{t}^{s}\p_{t}(\ak)_{j}^{l}\p_{t}^{3-s}\p_{l}\nn_{i}\phi^{\kappa}\|_{L^{2}(\Omega)}\|\p_{t}^{4}v^{i}\|_{L^{2}(\Omega)}\\&\lesssim\PP_{0}+T\PP(E^{\kappa}),
    \end{aligned}
\end{align}
and
\begin{align}\label{estimateforII23}
    \begin{aligned}
        |II_{23}|&\lesssim\sum_{s=0}^{3}C_{4}^{s}\int_{0}^{t}\|\p_{t}^{s}\p_{t}(J_{\kappa}\nn_{j}\phi^{\kappa})\p_{t}^{3-s}\nn_{j}\nn_{i}\phi^{\kappa}\|_{L^{2}(\Omega)}\|\p_{t}^{4}v^{i}\|_{L^{2}(\Omega)}\\&\lesssim\PP_{0}+T\PP(E^{\kappa}).
    \end{aligned}
\end{align}

Combining \eqref{expansionforII2}, \eqref{estimateforII22}, and \eqref{estimateforII23}, we have
\begin{align}\label{estimateforII2}
    II_{2}=\underbrace{-\int_{0}^{t}\int_{\Omega}(J_{\kappa}\nn_{j}\p_{t}^{4}\nn_{i}\phi^{\kappa})\nn_{j}\phi^{\kappa}\p_{t}^{4}v^{i}}_{II_{21}}+\mathcal{R},
\end{align}
where $\mathcal{R}\lesssim\PP_{0}+T\PP(E^{\kappa})$. Combining \eqref{expansionforII}, \eqref{estimateforII1} and \eqref{estimateforII2}, we have
\begin{align}\label{estimateforII}    II=-\int_{0}^{t}\int_{\Omega}(J_{\kappa}\nn_{j}\p_{t}^{4}\nn_{i}\phi^{\kappa})\nn_{j}\phi^{\kappa}\p_{t}^{4}v^{i}+\mathcal{R},
\end{align}
where $\mathcal{R}\lesssim\delta\int_{0}^{t}\|\p_{t}^{5}\phi\|_{L^{2}(\Omega)}^{2}+C(\delta)(\PP_{0}+T\PP(E^{\kappa}))$.
Combining \eqref{fourthtimeestformu1}, \eqref{estimateforI}, \eqref{estimateforIII}, \eqref{estimateforIV}, \eqref{estimateforII}, we have
\begin{align} \label{keycancellation2}   \int_{\Omega}J_{k}|\p_{t}^{4}v|^{2}+\sigma\int_{\Gamma_{1}}\sqrt{h}h^{ij}\p_{j}\p_{t}^{3}v^{\lambda}\Pi^{\mu}_{\lambda}\p_{i}\p_{t}^{3}v^{\alpha}\Pi_{\alpha}^{\mu}+\int_{0}^{t}\|\sqrt{\kappa}\p_{t}^{4}v\cdot n^{\kappa}\|_{H^{1}(\Gamma_{1})}^{2}=-\int_{0}^{t}\int_{\Omega}(J_{\kappa}\nn_{j}\p_{t}^{4}\nn_{i}\phi^{\kappa})\nn_{j}\phi^{\kappa}\p_{t}^{4}v^{i}+\mathcal{R},
\end{align}
where $|\mathcal{R}|\lesssim\PP(E_{\kappa}^{\prime})+\delta\|\p_{t}^{4}v\|_{L^{2}(\Omega)}^{2}+\delta \|\p_{t}^{3}v\|_{H^{1.5}(\Omega)}^{2}+\delta\int_{0}^{t}\|\sqrt{\kappa}\p_{t}^{4}v\|_{H^{1.5}(\Omega)}^{2}+\delta\int_{0}^{t}\|\p_{t}^{5}\phi\|^{2}+C(\delta)(\PP_{0}+T\PP(E^{\kappa}))$.
Now, combining \eqref{keycancellation1} and \eqref{keycancellation2}, we have
\begin{align}
\begin{aligned}
    &  \int_{0}^{t}\int_{\Omega}(J_{\kappa}\nn_{j}\p_{t}^{4}\nn_{i}\phi^{\kappa})\nn_{j}\phi^{\kappa}\p_{t}^{4}v^{i}+\int_{0}^{t}\int_{\Omega}J_{\kappa}\nn_{l}\p_{t}^{4}v_{\kappa}^{\mu}\nn_{l}\phi\p_{t}^{4}\nn_{\mu}\phi\\&= -\int_{0}^{t}\int_{\Omega}J_{\kappa}|\p_{t}^{5}\phi|^{2}-\frac{1}{2}\int_{\Omega}J_{\kappa}g_{\kappa}^{ij}\p_{t}^{4}\p_{j}\phi\p_{t}^{4}\p_{i}\phi-\int_{\Omega}J_{k}|\p_{t}^{4}v|^{2}\\&-\sigma\int_{\Gamma_{1}}\sqrt{h}h^{ij}\p_{j}\p_{t}^{3}v^{\lambda}\Pi^{\mu}_{\lambda}\p_{i}\p_{t}^{3}v^{\alpha}\Pi_{\alpha}^{\mu}-\int_{0}^{t}\|\sqrt{\kappa}\p_{t}^{4}v\cdot n^{\kappa}\|_{H^{1}(\Gamma_{1})}^{2}+\mathcal{R},
\end{aligned}
\end{align}
where $|\mathcal{R}|\lesssim\PP(E_{\kappa}^{\prime})+\delta\|\p_{t}^{4}v\|_{L^{2}(\Omega)}^{2}+\delta \|\p_{t}^{3}v\|_{H^{1.5}(\Omega)}^{2}+\delta\int_{0}^{t}\|\sqrt{\kappa}\p_{t}^{4}v\|_{H^{1.5}(\Omega)}^{2}
+\delta\|\p_{t}^{4}\phi\|_{1}^{2}+\delta\int_{0}^{t}\|\p_{t}^{5}\phi\|^{2}+C(\delta)(\PP_{0}+T\PP(E^{\kappa}))$.\\
Notice that
\begin{align}
\begin{aligned}
     & | \int_{0}^{t}\int_{\Omega}(J_{\kappa}\nn_{j}\p_{t}^{4}\nn_{i}\phi^{\kappa})\nn_{j}\phi^{\kappa}\p_{t}^{4}v^{i}+\int_{0}^{t}\int_{\Omega}J_{\kappa}\nn_{l}\p_{t}^{4}v_{\kappa}^{\mu}\nn_{l}\phi\p_{t}^{4}\nn_{\mu}\phi|\\\leq &|\int_{0}^{t}\int_{\Omega}J_{\kappa}\p_{t}^{4}\nn_{i}\phi\nn_{j}\phi\nn_{j}\p_{t}^{4}(v^{i}-v_{\kappa}^{i})|+|\int_{0}^{t}\int_{\Omega}J_{\kappa}\p_{t}^{4}\nn_{i}(\phi-\phi^{\kappa})\nn_{j}\phi\nn_{j}\p_{t}^{4}v^{i}|\\&+|\int_{0}^{t}\int_{\Omega}J_{\kappa}\p_{t}^{4}\nn_{i}\phi^{\kappa}\nn_{j}(\phi-\phi^{\kappa})\nn_{j}\p_{t}^{4}v^{i}|\\+&|\int_{0}^{t}\int_{\p\Omega}(\ak)_{j}^{3}J_{\kappa}\p_{t}^{4}\nn_{i}\phi^{\kappa}\nn_{j}(\phi^{\kappa}-\phi)\p_{t}^{4}v^{i}|+|\int_{0}^{t}\int_{\Omega}(J_{\kappa}\p_{t}^{4}\nn_{i}\phi^{\kappa})\nn_{j}\nn_{j}\phi^{\kappa}\p_{t}^{4}v^{i}|\\\lesssim&\int_{0}^{t}\|\nn_{j}\phi(\ak)_{j}^{l}J_{\kappa}\|_{L^{\infty}(\Omega)}\|\p_{t}^{4}\nn_{i}\phi\|_{L^{2}(\Omega)}\|\p_{l}\p_{t}^{4}(v^{i}-v_{\kappa}^{i})\|_{L^{2}(\Omega)}\\&+\int_{0}^{t}\|J_{\kappa}\nn_{j}\phi\|_{L^{\infty}(\Omega)}\|\p_{t}^{4}\nn_{i}(\phi-\phi^{\kappa})\|_{L^{2}(\Omega)}\|\nn_{j}\p_{t}^{4}v^{i}\|_{L^{2}(\Omega)}\\+&\int_{0}^{t}\|J_{\kappa}\|_{L^{\infty}(\Omega)}\|\p_{t}^{4}\nn_{i}\phi^{\kappa}\|_{L^{2}(\Omega)}\|\nn_{j}(\phi-\phi^{\kappa})\|_{L^{\infty}(\Omega)}\|\nn_{j}\p_{t}^{4}v^{i}\|_{L^{2}(\Omega)}\\+&\int_{0}^{t}\|(\ak)_{j}^{3}J_{\kappa}\|_{L^{\infty}(\p\Omega)}\|\p_{t}^{4}\nn_{i}\phi^{\kappa}\|_{L^{2}(\p\Omega)}\|\nn_{j}(\phi^{\kappa}-\phi)\|_{L^{\infty}(\p\Omega)}\|\p_{t}^{4}v^{i}\|_{L^{2}(\p\Omega)}+\PP_{0}+T\PP(E^{\kappa})\\\lesssim &\int_{0}^{t}\PP(E^{\kappa})\|\sqrt{\kappa}\p_{t}^{4}v\|_{H^{1.5}(\Omega)}+\int_{0}^{t}(\PP_{0}+T\PP(E^{\kappa}))\|\p_{t}^{4}\phi\|_{H^{1.5}(\Omega)}\|\sqrt{\kappa}\p_{t}^{4}v\|_{H^{1}(\Omega)}+\int_{0}^{T}\PP(E^{\kappa})\|\phi\|_{H^{4}(\Omega)}\|\kappa\p_{t}^{4}v\|_{H^{1}(\Omega)}\\+&\int_{0}^{t}\kappa(\PP_{0}+T\PP(E^{\kappa}))\|\p_{t}^{4}\nn\phi^{\kappa}\|_{H^{1}(\Omega)}\|\phi\|_{H^{4}(\Omega)}\|\p_{t}^{4}v\|_{H^{1}(\Omega)}+\PP_{0}+T\PP(E^{\kappa})\\\lesssim&\delta\int_{0}^{t}\|\sqrt{\kappa}\p_{t}^{4}v\|_{H^{1.5}(\Omega)}^{2}+C(\delta)(\PP_{0}+T\PP(E^{\kappa}))+(\PP_{0}+T\PP(E^{\kappa}))\|\sqrt{\kappa}\p_{t}^{3}v\|_{H^{2}(\Omega)}^{2}\\\lesssim &\delta\int_{0}^{t}\|\sqrt{\kappa}\p_{t}^{4}v\|_{H^{1.5}(\Omega)}^{2}+C(\delta)(\PP_{0}+T\PP(E^{\kappa}))+\PP(E_{\kappa}^{\prime}),
\end{aligned} 
\end{align}
where the last two steps we used \eqref{estimateforsqrtkappapt4nablaphih1}. Hence, we have shown
\begin{align}
\begin{aligned}
     &\int_{0}^{t}\int_{\Omega}J_{\kappa}|\p_{t}^{5}\phi|^{2}+\frac{1}{2}\int_{\Omega}J_{\kappa}g_{\kappa}^{ij}\p_{t}^{4}\p_{j}\phi\p_{t}^{4}\p_{i}\phi+\int_{\Omega}J_{k}|\p_{t}^{4}v|^{2}+\sigma\int_{\Gamma_{1}}\sqrt{h}h^{ij}\p_{j}\p_{t}^{3}v^{\lambda}\Pi^{\mu}_{\lambda}\p_{i}\p_{t}^{3}v^{\alpha}\Pi_{\alpha}^{\mu}+\int_{0}^{t}\|\sqrt{\kappa}\p_{t}^{4}v\cdot n^{\kappa}\|_{H^{1}(\Gamma_{1})}^{2}\\&\lesssim\PP(E_{\kappa}^{\prime})+\delta\|\p_{t}^{4}v\|_{L^{2}(\Omega)}^{2}+\delta \|\p_{t}^{3}v\|_{H^{1.5}(\Omega)}^{2}+\delta\int_{0}^{t}\|\sqrt{\kappa}\p_{t}^{4}v\|_{H^{1.5}(\Omega)}^{2}
+\delta\|\p_{t}^{4}\phi\|_{1}^{2}+\delta\int_{0}^{t}\|\p_{t}^{5}\phi\|^{2}+C(\delta)(\PP_{0}+T\PP(E^{\kappa})).
\end{aligned}
\end{align}

Now, we can show the estimate for $E_{\kappa}^{\prime}$ in a similar way, but self-contained.
\lem[The estimate for the three time-differentiated $\kappa$-problem]
\begin{align}\label{threetimeest}
    \begin{aligned}        &\sigma\int_{\Gamma_{1}}\sqrt{h}h^{ij}\p_{j}\pb v_{tt}^{\lambda}\Pi_{\lambda}^{\mu}\p_{i}\pb v_{tt}^{\alpha}\Pi_{\alpha}^{\mu}+\int_{0}^{t}\|\sqrt{\kappa}\pb v_{ttt}\cdot n^{\kappa}\|_{H^{1}(\Gamma_{1})}^{2}\\&\lesssim C(\delta)(\PP_{0}+T\PP(E^{\kappa}))+(\int_{0}^{t}\|\sqrt{\kappa}\pb v_{tt}\|_{H^{2}(\Omega)}^{2})+\delta \int_{0}^{t}\|\sqrt{\kappa}\p_{t}^{3}\pb v\|_{H^{1}(\Gamma_{1})}^{2}.
    \end{aligned}
\end{align}
\lem[The estimate for the two time-differentiated $\kappa$-problem]
\begin{align}\label{twotimeest}
 \begin{aligned}  &\sigma\int_{\Gamma_{1}}\sqrt{h}h^{ij}\p_{j}\pb^{2} v_{t}^{\lambda}\Pi_{\lambda}^{\mu}\p_{i}\pb^{2} v_{t}^{\alpha}\Pi_{\alpha}^{\mu}+\int_{0}^{t}\|\sqrt{\kappa}\pb^{2} v_{tt}\cdot n^{\kappa}\|_{H^{1}(\Gamma_{1})}^{2}\\&\lesssim C(\delta)(\PP_{0}+T\PP(E^{\kappa}))+(\int_{0}^{t}\|\sqrt{\kappa}\pb^{2} v_{t}\|_{H^{2}(\Omega)}^{2})+\delta \int_{0}^{t}\|\sqrt{\kappa}\p_{t}^{2}\pb^{2} v\|_{H^{1}(\Gamma_{1})}^{2}   .
 \end{aligned}   
\end{align}
\lem[The estimate for the time-differentiated $\kappa$-problem]
\begin{align}\label{onetimeest}
\begin{aligned}  &\sigma\int_{\Gamma_{1}}\sqrt{h}h^{ij}\p_{j}\pb^{3} v^{\lambda}\Pi_{\lambda}^{\mu}\p_{i}\pb^{3} v^{\alpha}\Pi_{\alpha}^{\mu}+\int_{0}^{t}\|\sqrt{\kappa}\pb^{3} v_{t}\cdot n^{\kappa}\|_{H^{1}(\Gamma_{1})}^{2}\\&\lesssim C(\delta)(\PP_{0}+T\PP(E^{\kappa}))+(\int_{0}^{t}\|\sqrt{\kappa}\pb^{3} v\|_{H^{2}(\Omega)}^{2})+\delta \int_{0}^{t}\|\sqrt{\kappa}\p_{t}\pb^{3} v\|_{H^{1}(\Gamma_{1})}^{2}.
\end{aligned}   
\end{align}
\lem[The estimate for the $\kappa$-problem]
\begin{align}\label{zerotimeest}
 \begin{aligned}     &\sigma\int_{\Gamma_{1}}\sqrt{h}h^{ij}\p_{j}\pb^{4} \eta^{\lambda}\Pi_{\lambda}^{\mu}\p_{i}\pb^{4} \eta^{\alpha}\Pi_{\alpha}^{\mu}+\int_{0}^{t}\|\sqrt{\kappa}\pb^{4} v\cdot n^{\kappa}\|_{H^{1}(\Gamma_{1})}^{2}\\&\lesssim C(\delta)(\PP_{0}+T\PP(E^{\kappa}))+\delta \int_{0}^{t}\|\sqrt{\kappa}\pb^{4} v\|_{H^{1}(\Gamma_{1})}^{2}.
 \end{aligned}  
\end{align}
\pf
Using the symbol $\pb_{A}=\p_{t}$ or $\pb$ (at least one $\pb$ in $\pb_{A}$ here), multiplying $J_{\kappa}$ on \eqref{vkappat}, and taking $\pb_{A}^{3}\pb$ then testing with $\pb_{A}^{3}\pb v$, we can write
\begin{align}\label{pbaest}
    \begin{aligned}        &\underbrace{\int_{0}^{t}\int_{\Omega}\pb_{A}^{3}\pb(J_{\kappa} \p_{t}v^{i})\pb_{A}^{3}\pb v^{i}}_{I^{\prime}}\underbrace{+\int_{0}^{t}\int_{\Omega}\pb_{A}^{3}\pb(J_{\kappa}\nn_{i} q)\pb_{A}^{3}\pb v^{i}}_{II^{\prime}}\\&=\underbrace{\int_{0}^{t}\int_{\Omega}-\pb_{A}^{3}\pb(J_{\kappa}\phi^{\kappa}_{t}\nn_{i}\phi^{\kappa})\pb_{A}^{3}\pb v^{i}+\int_{0}^{t}\int_{\Omega}-\pb_{A}^{3}\pb(J_{\kappa}\nn_{j}\nn_{i}\phi^{\kappa}\nn_{j}\phi^{\kappa})\pb_{A}^{3}\pb v^{i}}_{III^{\prime}}.
    \end{aligned}
\end{align}
We can show that 
\begin{align}\label{estimateforIIIprime}
    \begin{aligned}
      |III^{\prime}|&\lesssim  \int_{0}^{t}\|\pb_{A}^{3}(\pb J_{\kappa}\phi^{\kappa}_{t}\nn_{i}\phi^{\kappa})\|_{L^{2}(\Omega)}\|\pb_{A}^{3}\pb v^{i}\|_{L^{2}(\Omega)}+\int_{0}^{t}\|\pb_{A}^{3}(J_{\kappa}\pb\phi^{\kappa}_{t}\nn_{i}\phi^{\kappa})\|_{L^{2}(\Omega)}\|\pb_{A}^{3}\pb v^{i}\|_{L^{2}(\Omega)}\\&+\int_{0}^{t}\|\pb_{A}^{3}(J_{\kappa}\phi^{\kappa}_{t}\pb\nn_{i}\phi^{\kappa})\|_{L^{2}(\Omega)}\|\pb_{A}^{3}\pb v^{i}\|_{L^{2}(\Omega)}+\int_{0}^{t}\|\pb_{A}^{3}\pb(J_{\kappa}\nn_{j}\nn_{i}\phi^{\kappa}\nn_{j}\phi^{\kappa})\|_{H^{-\frac{1}{2}}(\Omega)}\|\pb_{A}^{3}\pb v^{i}\|_{H^{\frac{1}{2}}(\Omega)}\\&\lesssim\PP_{0}+T\PP(E^{\kappa})+\int_{0}^{t}\|\pb_{A}^{3}(J_{\kappa}\nn_{j}\nn_{i}\phi^{\kappa}\nn_{j}\phi^{\kappa})\|_{H^{\frac{1}{2}}(\Omega)}\|\pb_{A}^{3} v^{i}\|_{H^{1.5}(\Omega)}\\&\lesssim\PP_{0}+T\PP(E^{\kappa}).
    \end{aligned}
\end{align}
For $I^{\prime}$, it's evident that
\begin{align}
    \begin{aligned}
        I^{\prime}&=\int_{0}^{t}\int_{\Omega}J_{\kappa} \pb_{A}^{3}\pb\p_{t}v^{i}\pb_{A}^{3}\pb v^{i}+\sum_{s=0}^{3}C_{4}^{s+1}\int_{0}^{t}\int_{\Omega}\pb_{A}^{s}\pb_{A}J_{\kappa} \pb_{A}^{3-s}\p_{t}v^{i}\pb_{A}^{3}\pb v^{i}\\&=\underbrace{\frac{1}{2}\int_{\Omega}J_{\kappa} \pb_{A}^{3}\pb v^{i}\pb_{A}^{3}\pb v^{i}|_{0}^{t}}_{I_{1}^{\prime}}\underbrace{+\sum_{s=0}^{3}C_{4}^{s+1}\int_{0}^{t}\int_{\Omega}\pb_{A}^{s}\pb_{A}\p_{t}J_{\kappa} \pb_{A}^{3-s}v^{i}\pb_{A}^{3}\pb v^{i}-\frac{1}{2}\int_{0}^{t}\int_{\Omega}\p_{t}J_{\kappa} \pb_{A}^{3}\pb v^{i}\pb_{A}^{3}\pb  v^{i}}_{I_{2}^{\prime}}.
    \end{aligned}
\end{align}
Since
\begin{align}
    \begin{aligned}
    |I_{2}^{\prime}|&  \lesssim  \sum_{s=0}^{3}C_{4}^{s+1}\int_{0}^{t}\|\pb_{A}^{s}\pb_{A}\p_{t}J_{\kappa} \pb_{A}^{3-s}v^{i}\|_{L^{2}(\Omega)}\|\pb_{A}^{3}\pb v^{i}\|_{L^{2}(\Omega)}+\frac{1}{2}\int_{0}^{t}\int_{\Omega}\|\p_{t}J_{\kappa}\|_{L^{\infty}(\Omega)}\| \pb_{A}^{3}\pb v^{i}\|_{L^{2}(\Omega)}\|\pb_{A}^{3}\pb v^{i}\|_{L^{2}(\Omega)}\\&\lesssim\PP_{0}+T\PP(E^{\kappa}),
    \end{aligned}
\end{align}
we can express $I^{\prime}$ as
\begin{align}\label{estimateforIprime}
    I^{\prime}=\frac{1}{2}\int_{\Omega}J_{\kappa} \pb_{A}^{3}\pb v^{i}\pb_{A}^{3}\pb v^{i}+\mathcal{R},
\end{align}
where $\mathcal{R}\lesssim\PP_{0}+T\PP(E^{\kappa})$. Now, the estimate of $II^{\prime}$ can be written as
\begin{align}\label{expansionforIIprime}
    \begin{aligned}        II^{\prime}&=\underbrace{\int_{0}^{t}\int_{\Gamma_{1}}\pb_{A}^{3}\pb(J_{\kappa}(\ak)_{i}^{3} q)\pb_{A}^{3}\pb v^{i}}_{II_{1}^{\prime}}\underbrace{-\int_{0}^{t}\int_{\Omega}\pb_{A}^{3}\pb(J_{\kappa}(\ak)_{i}^{j} q)\pb_{A}^{3}\pb\p_{j} v^{i}}_{II_{2}^{\prime}}.
    \end{aligned}
\end{align}
For the second term in \eqref{expansionforIIprime} can be controlled as
\begin{align}
    \begin{aligned}       II_{2}^{\prime}&=\int_{0}^{t}\int_{\Omega}J_{\kappa}(\ak)_{i}^{j} \pb_{A}^{3}\pb q\pb_{A}^{3}\pb\p_{j} v^{i}+\sum_{s=1}^{3}C_{4}^{s}\int_{0}^{t}\int_{\Omega}\pb_{A}^{s}(J_{\kappa}(\ak)_{i}^{j} )\pb_{A}^{3-s}\pb q\pb_{A}^{3}\pb\p_{j} v^{i}\\&+\sum_{s=1}^{3}C_{4}^{s}\int_{0}^{t}\int_{\Omega}\pb_{A}^{s}\pb (J_{\kappa}(\ak)_{i}^{j} )\pb_{A}^{3-s} q\pb_{A}^{3}\pb\p_{j} v^{i}\\&=-\sum_{s=1}^{4}C_{4}^{s}\int_{0}^{t}\int_{\Omega}J_{\kappa}\pb_{A}^{3}\pb q \pb_{A}^{s}(\ak)_{i}^{j}\pb_{A}^{4-s}\p_{j} v^{i}+\sum_{s=1}^{3}C_{4}^{s}\int_{0}^{t}\int_{\Omega}\pb_{A}^{s}(J_{\kappa}(\ak)_{i}^{j} )\pb_{A}^{3-s}\pb q\pb_{A}^{3}\pb\p_{j} v^{i}\\&+\sum_{s=1}^{3}C_{4}^{s}\int_{0}^{t}\int_{\Omega}\pb_{A}^{s}\pb (J_{\kappa}(\ak)_{i}^{j} )\pb_{A}^{3-s} q\pb_{A}^{3}\pb\p_{j} v^{i}.
    \end{aligned}
\end{align}
Moreover,
\begin{align}\label{estimateforII2prime}
    \begin{aligned}
       |II_{2}^{\prime}|&\lesssim \sum_{s=0}^{3}C_{4}^{s+1}\int_{0}^{t}\|J_{\kappa}\|_{L^{\infty}(\Omega)}\|\pb_{A}^{3}\pb q\|_{L^{2}(\Omega)}\| \pb_{A}^{s}\pb_{A}(\ak)_{i}^{j}\pb_{A}^{3-s}\p_{j} v^{i}\|_{L^{2}(\Omega)}\\&+\sum_{s=0}^{2}C_{4}^{s+1}\int_{0}^{t}\|\pb_{A}^{s}\p_{t}(J_{\kappa}(\ak)_{i}^{j} )\pb_{A}^{2-s}\pb q\|_{H^{\frac{1}{2}}(\Omega)}\|\pb_{A}^{3}\pb\p_{j} v^{i}\|_{H^{-\frac{1}{2}}(\Omega)}\\&+\sum_{s=0}^{2}C_{4}^{s+1}\int_{0}^{t}\|\pb_{A}^{s}\pb_{A}\pb (J_{\kappa}(\ak)_{i}^{j} )\pb_{A}^{3-s} q\|_{H^{\frac{1}{2}}(\Omega)}\|\pb_{A}^{3}\pb\p_{j} v^{i}\|_{H^{-\frac{1}{2}}(\Omega)}\\&\lesssim \PP_{0}+T\PP(E^{\kappa})+\int_{0}^{t}\PP(E^{\kappa})\|\pb_{A}^{3}\p v\|_{H^{\frac{1}{2}}(\Omega)}\\&\lesssim \PP_{0}+T\PP(E^{\kappa}).
    \end{aligned}
\end{align}
Now, the boundary term reads
\begin{align}
    \begin{aligned}
        II_{1}^{\prime}&=\underbrace{\int_{0}^{t}\int_{\Gamma_{1}}\pb_{A}^{3}\pb(-\sigma \sqrt{h}\lap_{h}(\eta)\cdot n_{\kappa}n_{\kappa}^{i})\pb_{A}^{3}\pb v^{i}}_{II_{11}^{\prime}}\underbrace{+\int_{0}^{t}\int_{\Gamma_{1}}\pb_{A}^{3}\pb(\kappa\lap_{0}(v\cdot n_{\kappa})n_{\kappa}^{i})\pb_{A}^{3}\pb v^{i}}_{II_{12}^{\prime}}
    \end{aligned}
\end{align}
The estimate is similar to $III_{1}$ and $III_{2}$. Here we omit some details. (Replace $\p_{t}^{4}$ by $\pb_{A}^{3}\pb$, every estimate is still valid.) We can write
\begin{align}\label{estimateforII1prime}
    II_{1}^{\prime}=\sigma\int_{\Gamma_{1}}\sqrt{h}h^{ij}\p_{j}\pb_{A}^{3}\pb \eta^{\lambda}\Pi_{\lambda}^{\mu}\p_{i}\pb_{A}^{3}\pb \eta^{\alpha}\Pi_{\alpha}^{\mu}+\int_{0}^{t}\|\sqrt{\kappa}\pb_{A}^{3}\pb v\cdot n^{\kappa}\|_{H^{1}(\Gamma_{1})}^{2}+\mathcal{R},
\end{align}
where $\mathcal{R}\lesssim C(\delta)(\PP_{0}+T\PP(E^{\kappa}))+(\int_{0}^{t}\|\sqrt{\kappa}\pb_{A}^{3}\pb \eta\|_{H^{2}(\Omega)}^{2})+\delta \int_{0}^{t}\|\sqrt{\kappa}\pb_{A}^{3}\pb v\|_{H^{1}(\Gamma_{1})}^{2}$. Thus, combining \eqref{expansionforIIprime}, \eqref{estimateforII2prime} and \eqref{estimateforII1prime}, we have
\begin{align} \label{estimateforIIprime}   II^{\prime}=\sigma\int_{\Gamma_{1}}\sqrt{h}h^{ij}\p_{j}\pb_{A}^{3}\pb \eta^{\lambda}\Pi_{\lambda}^{\mu}\p_{i}\pb_{A}^{3}\pb \eta^{\alpha}\Pi_{\alpha}^{\mu}+\int_{0}^{t}\|\sqrt{\kappa}\pb_{A}^{3}\pb v\cdot n^{\kappa}\|_{H^{1}(\Gamma_{1})}^{2}+\mathcal{R},
\end{align}
where $\mathcal{R}\lesssim C(\delta)(\PP_{0}+T\PP(E^{\kappa}))+(\int_{0}^{t}\|\sqrt{\kappa}\pb_{A}^{3}\pb \eta\|_{H^{2}(\Omega)}^{2})+\delta \int_{0}^{t}\|\sqrt{\kappa}\pb_{A}^{3}\pb v\|_{H^{1}(\Gamma_{1})}^{2}$. Hence, combining \eqref{pbaest}, \eqref{estimateforIIIprime}, \eqref{estimateforIprime} and \eqref{estimateforIIprime}, we have
\begin{align}
 \begin{aligned}    &\sigma\int_{\Gamma_{1}}\sqrt{h}h^{ij}\p_{j}\pb_{A}^{3}\pb \eta^{\lambda}\Pi_{\lambda}^{\mu}\p_{i}\pb_{A}^{3}\pb \eta^{\alpha}\Pi_{\alpha}^{\mu}+\int_{0}^{t}\|\sqrt{\kappa}\pb_{A}^{3}\pb v\cdot n^{\kappa}\|_{H^{1}(\Gamma_{1})}^{2}\\&\lesssim C(\delta)(\PP_{0}+T\PP(E^{\kappa}))+(\int_{0}^{t}\|\sqrt{\kappa}\pb_{A}^{3}\pb \eta\|_{H^{2}(\Omega)}^{2})+\delta \int_{0}^{t}\|\sqrt{\kappa}\pb_{A}^{3}\pb v\|_{H^{1}(\Gamma_{1})}^{2} .
 \end{aligned}   
\end{align}
Notice that $\sup\|\sqrt{\kappa}\eta\|_{6.5}$ has been contained in energy. The lemmas follow immediately. 

The proof of the following lemma is identical to that in \cite{coutand2007well}. Here, we simply provide the details regarding the dependence on the energy.
\lem
\begin{align}\label{estimateforkappa6eta}
\begin{aligned}
     &\frac{1}{2}\int_{\Gamma_{1}}\kappa\pb^{4}\pb_{\mu}\eta\cdot n ^{\kappa}h^{\mu\nu}\pb^{4}\pb_{\nu}\eta\cdot n_{\kappa}+\frac{1}{2}\int_{\Gamma_{1}}\kappa\pb^{4}\pb_{\mu}\p_{\gamma}\eta\cdot n ^{\kappa}h^{\mu\nu}\pb^{4}\pb_{\gamma}\pb_{\nu}\eta\cdot n_{\kappa}\\&+\int_{0}^{t}\int_{\Gamma_{1}}\sigma(\sqrt{h}h^{ij}\Pi_{\lambda}^{l}\pb^{4} \p_{\mu}\p_{j}\eta^{\lambda})( h^{\mu\nu}\p_{i}\pb^{4}\pb_{\nu}\eta^{k}\Pi_{k}^{l})\lesssim\delta\int_{0}^{t}\|\Pi_{\lambda}^{l}\pb^{6}\eta^{\lambda}\|_{L^{2}(\Gamma_{1})}^{2}+C(\delta)T\PP(E^{\kappa}).
\end{aligned} 
\end{align}

\pf
We utilize the boundary condition of $q$:
\begin{align}\label{conditionforq}
   \sqrt{h_{\kappa}} q n_{\kappa}^{l}=-\sigma(\sqrt{h} \lap_{h}(\eta) \cdot n_{\kappa})n_{\kappa}^{l}+\kappa((1-\bar{\lap})(v\cdot n_{\kappa}))n_{\kappa}^{l}
\end{align}
Taking $\pb^{4}\pb_{\mu}$ on \eqref{conditionforq} and testing $h^{\mu\nu}\pb^{4}\pb^{\nu}\eta^{k}\Pi_{k}^{l}$, we have
\begin{align}\label{expansionforIVVVIprime}
\begin{aligned}
     \underbrace{\int_{0}^{t}\int_{\Gamma_{1}} \pb^{4}\pb_{\mu}(\sqrt{h_{\kappa}} q n_{\kappa}^{l})h^{\mu\nu}\pb^{4}\pb_{\nu}\eta^{k}\Pi_{k}^{l}}_{IV^{\prime}}&=\underbrace{\int_{0}^{t}\int_{\Gamma_{1}}-\sigma\pb^{4}\pb_{\mu}(\sqrt{h} \lap_{h}(\eta) \cdot n_{\kappa}n_{\kappa}^{l})h^{\mu\nu}\pb^{4}\pb_{\nu}\eta^{k}\Pi_{k}^{l}}_{V^{\prime}}\\&\underbrace{+\int_{0}^{t}\int_{\Gamma_{1}}\kappa\pb^{4}\pb_{\mu}((1-\bar{\lap})(v\cdot n ^{\kappa})n_{\kappa}^{l})h^{\mu\nu}\pb^{4}\pb_{\nu}\eta^{k}\Pi_{k}^{l}}_{VI^{\prime}}.
\end{aligned}   
\end{align}
We analyze $V^{\prime}$ first.
\begin{align}\label{expansionforVprime}
    \begin{aligned}
        V^{\prime}&=\underbrace{\int_{0}^{t}\int_{\Gamma_{1}}-\sigma\pb^{4}\pb_{\mu}( \p_{i}(\sqrt{h}h^{ij}\p_{j}(\eta^{l})) )h^{\mu\nu}\pb^{4}\pb_{\nu}\eta^{k}\Pi_{k}^{l}}_{V_{1}^{\prime}}\\&\underbrace{+\int_{0}^{t}\int_{\Gamma_{1}}-\sigma\pb^{4}\pb_{\mu}( \p_{i}(\sqrt{h}h^{ij}\p_{j}(\eta)) \cdot (n_{\kappa}-n)n_{\kappa}^{l})h^{\mu\nu}\pb^{4}\pb_{\nu}\eta^{k}\Pi_{k}^{l}}_{V_{2}^{\prime}}\\&\underbrace{+\int_{0}^{t}\int_{\Gamma_{1}}-\sigma\pb^{4}\pb_{\mu}( \p_{i}(\sqrt{h}h^{ij}\p_{j}(\eta)) \cdot n_{\kappa}(n_{\kappa}^{l}-n^{l})h^{\mu\nu}\pb^{4}\pb_{\nu}\eta^{k}\Pi_{k}^{l}}_{V_{3}^{\prime}}
    \end{aligned}
\end{align}
Similar to $III$, we have
\begin{align}\label{expansionforV1prime}
    \begin{aligned}
        V_{1}^{\prime}&=\underbrace{\int_{0}^{t}\int_{\Gamma_{1}}-\sigma\pb^{4} \p_{i}(\sqrt{h}h^{ij}(\delta^{l\lambda}-h^{\alpha\beta}\p_{\alpha}\eta^{\lambda}\p_{\beta}\eta^{l})(\p_{\mu}\p_{j}\eta^{\lambda})) h^{\mu\nu}\pb^{4}\pb_{\nu}\eta^{k}\Pi_{k}^{l}}_{V_{11}^{\prime}}\\&\underbrace{+\int_{0}^{t}\int_{\Gamma_{1}}-\sigma\pb^{4} \p_{i}((h^{ij}h^{\alpha\beta}-h^{i\alpha}h^{j\beta})\p_{\alpha}\eta^{\lambda}\p_{\beta\mu}^{2}\eta^{\lambda}\p_{j}\eta^{l}) h^{\mu\nu}\pb^{4}\pb_{\nu}\eta^{k}\Pi_{k}^{l}}_{V_{12}^{\prime}}.
    \end{aligned}
\end{align}
$V_{11}^{\prime}$ will contribute a positive term as
\begin{align}\label{expansionforV11prime}
    \begin{aligned}
       V_{11}^{\prime}&= \int_{0}^{t}\int_{\Gamma_{1}}\sigma\pb^{4} (\sqrt{h}h^{ij}\Pi_{\lambda}^{l}\p_{\mu}\p_{j}\eta^{\lambda})\p_{i}( h^{\mu\nu}\pb^{4}\pb_{\nu}\eta^{k}\Pi_{k}^{l})\\&=\int_{0}^{t}\int_{\Gamma_{1}}\sigma\pb^{4} (\sqrt{h}h^{ij}\Pi_{\lambda}^{l}\p_{\mu}\p_{j}\eta^{\lambda})( h^{\mu\nu}\p_{i}\pb^{4}\pb_{\nu}\eta^{k}\Pi_{k}^{l})\\&+\int_{0}^{t}\int_{\Gamma_{1}}\sigma\pb^{4} (\sqrt{h}h^{ij}\Pi_{\lambda}^{l}\p_{\mu}\p_{j}\eta^{\lambda})\p_{i}( h^{\mu\nu}\Pi_{k}^{l})\pb^{4}\pb_{\nu}\eta^{k}\\&=\underbrace{\int_{0}^{t}\int_{\Gamma_{1}}\sigma(\sqrt{h}h^{ij}\Pi_{\lambda}^{l}\pb^{4} \p_{\mu}\p_{j}\eta^{\lambda})( h^{\mu\nu}\p_{i}\pb^{4}\pb_{\nu}\eta^{k}\Pi_{k}^{l})}_{V_{11,1}^{\prime}}\\&\underbrace{+\sum_{s=1}^{4}C_{4}^{s}\int_{0}^{t}\int_{\Gamma_{1}}\sigma\pb^{s}(\sqrt{h}h^{ij}\Pi_{\lambda}^{l})(\pb^{4-s} \p_{\mu}\p_{j}\eta^{\lambda})( h^{\mu\nu}\p_{i}\pb^{4}\pb_{\nu}\eta^{k}\Pi_{k}^{l})}_{V_{11,2}^{\prime}}\\&\underbrace{+\int_{0}^{t}\int_{\Gamma_{1}}\sigma\pb^{4} (\sqrt{h}h^{ij}\Pi_{\lambda}^{l}\p_{\mu}\p_{j}\eta^{\lambda})\p_{i}( h^{\mu\nu}\Pi_{k}^{l})\pb^{4}\pb_{\nu}\eta^{k}}_{V_{11,3}^{\prime}}.
    \end{aligned}
\end{align}
The lower orders can be estimated as
\begin{align}\label{estimateforV112prime}
    \begin{aligned}
        |V_{11,2}^\prime|&\lesssim\sum_{s=0}^{3}C_{4}^{s+1}\int_{0}^{t}\sigma\|\pb^{s}\pb(\sqrt{h}h^{ij}\Pi_{\lambda}^{l})(\pb^{3-s} \p_{\mu}\p_{j}\eta^{\lambda})\|_{L^{2}(\Gamma_{1})}\| h^{\mu\nu}\|_{L^{\infty}(\Gamma_{1})}\|\p_{i}\pb^{4}\pb_{\nu}\eta^{k}\Pi_{k}^{l}\|_{L^{2}(\Gamma_{1})}\\&\lesssim\int_{0}^{t}\PP(E^{\kappa})\|\p_{i}\pb^{4}\pb_{\nu}\eta^{k}\Pi_{k}^{l}\|_{L^{2}(\Gamma_{1})}\\&\lesssim\delta\int_{0}^{t}\|\p_{i}\pb^{4}\pb_{\nu}\eta^{k}\Pi_{k}^{l}\|_{L^{2}(\Gamma_{1})}^{2}+C(\delta)T\PP(E^{\kappa}),
    \end{aligned}
\end{align}
and
\begin{align}
    \begin{aligned}
        V_{11,3}^{\prime}&=\int_{0}^{t}\int_{\Gamma_{1}}\sigma (\sqrt{h}h^{ij}\Pi_{\lambda}^{l}\pb^{4}\p_{\mu}\p_{j}\eta^{\lambda})\p_{i}( h^{\mu\nu}\Pi_{k}^{l})\pb^{4}\pb_{\nu}\eta^{k}\\&+\sum_{s=1}^{4}C_{4}^{s}\int_{0}^{t}\int_{\Gamma_{1}}\sigma\pb^{s} (\sqrt{h}h^{ij})\pb^{4-s}(\Pi_{\lambda}^{l}\p_{\mu}\p_{j}\eta^{\lambda})\p_{i}( h^{\mu\nu}\Pi_{k}^{l})\pb^{4}\pb_{\nu}\eta^{k},
    \end{aligned}
\end{align}
\begin{align}\label{estimateforV113prime}
    \begin{aligned}
        |V_{11,3}^{\prime}|&\lesssim\int_{0}^{t}\sigma (\|\sqrt{h}h^{ij}\|_{L^{\infty}(\Gamma_{1})}\|\Pi_{\lambda}^{l}\pb^{4}\p_{\mu}\p_{j}\eta^{\lambda}\|_{L^{2}(\Gamma_{1})}\|\p_{i}( h^{\mu\nu}\Pi_{k}^{l})\|_{L^{\infty}(\Gamma_{1})}\|\pb^{4}\pb_{\nu}\eta^{k}\|_{L^{2}(\Gamma_{1})}\\&+\sum_{s=0}^{3}C_{4}^{s+1}\int_{0}^{t}\sigma\|\pb^{s}\pb (\sqrt{h}h^{ij})\pb^{3-s}(\Pi_{\lambda}^{l}\p_{\mu}\p_{j}\eta^{\lambda})\|_{L^{2}(\Gamma_{1})}\|\p_{i}( h^{\mu\nu}\Pi_{k}^{l})\|_{L^{\infty}(\Gamma_{1})}\|\pb^{4}\pb_{\nu}\eta^{k}\|_{L^{2}(\Gamma_{1})}\\&\lesssim\delta\int_{0}^{t}\|\Pi_{\lambda}^{l}\pb^{4}\p_{\mu}\p_{j}\eta^{\lambda}\|_{L^{2}(\Gamma_{1})}^{2}+C(\delta)T\PP(E^{\kappa}),
    \end{aligned}
\end{align}
after application of Young's inequality. Combining \eqref{expansionforV11prime}, \eqref{estimateforV112prime} and \eqref{estimateforV113prime}, we have
\begin{align}\label{estimateforV11prime}
    V_{11}^{\prime}=\underbrace{\int_{0}^{t}\int_{\Gamma_{1}}\sigma(\sqrt{h}h^{ij}\Pi_{\lambda}^{l}\pb^{4} \p_{\mu}\p_{j}\eta^{\lambda})( h^{\mu\nu}\p_{i}\pb^{4}\pb_{\nu}\eta^{k}\Pi_{k}^{l})}_{V_{11,1}^{\prime}}+\mathcal{R},
\end{align}
where $|\mathcal{R}|\lesssim\delta\int_{0}^{t}\|\Pi_{\lambda}^{l}\pb^{6}\eta^{\lambda}\|_{L^{2}(\Gamma_{1})}^{2}+C(\delta)T\PP(E^{\kappa})$. Now, the estimate of $V_{12}^{\prime}$ needs the following geometric relationship on boundary.
\begin{align}
    \pb_{k}\eta^{i}\Pi_{i}^{l}=0,\quad \textit{$k=1$ or $2$.}
\end{align}
Utilizing this property, we can express $V_{12}^{\prime}$ as
\begin{align}\label{expansionforV12prime}
    \begin{aligned}
        V_{12}^{\prime}&=\int_{0}^{t}\int_{\Gamma_{1}}\sigma\pb^{4} ((h^{ij}h^{\alpha\beta}-h^{i\alpha}h^{j\beta})\p_{\alpha}\eta^{\lambda}\p_{\beta\mu}^{2}\eta^{\lambda}\p_{j}\eta^{l}) h^{\mu\nu}\p_{i}\pb^{4}\pb_{\nu}\eta^{k}\Pi_{k}^{l}\\&+\int_{0}^{t}\int_{\Gamma_{1}}\sigma\pb^{4} ((h^{ij}h^{\alpha\beta}-h^{i\alpha}h^{j\beta})\p_{\alpha}\eta^{\lambda}\p_{\beta\mu}^{2}\eta^{\lambda}\p_{j}\eta^{l}) \p_{i}(h^{\mu\nu}\Pi_{k}^{l})\pb^{4}\pb_{\nu}\eta^{k}\\&=\underbrace{\int_{0}^{t}\int_{\Gamma_{1}}\sigma (h^{ij}h^{\alpha\beta}-h^{i\alpha}h^{j\beta})\p_{\alpha}\eta^{\lambda}\pb^{4}\p_{\beta\mu}^{2}\eta^{\lambda}\p_{j}\eta^{l} h^{\mu\nu}\p_{i}\pb^{4}\pb_{\nu}\eta^{k}\Pi_{k}^{l}}_{V_{12,1}^{\prime}}\\&\underbrace{+\sum_{s=1}^{4}C_{4}^{s}\int_{0}^{t}\int_{\Gamma_{1}}\sigma\pb^{s} ((h^{ij}h^{\alpha\beta}-h^{i\alpha}h^{j\beta})\p_{\alpha}\eta^{\lambda}\p_{j}\eta^{l})\pb^{4-s}(\p_{\beta\mu}^{2}\eta^{\lambda}) h^{\mu\nu}\p_{i}\pb^{4}\pb_{\nu}\eta^{k}\Pi_{k}^{l}}_{V_{12,2}^{\prime}}\\&\underbrace{+\int_{0}^{t}\int_{\Gamma_{1}}\sigma\pb^{4} ((h^{ij}h^{\alpha\beta}-h^{i\alpha}h^{j\beta})\p_{\alpha}\eta^{\lambda}\p_{\beta\mu}^{2}\eta^{\lambda}\p_{j}\eta^{l}) \p_{i}(h^{\mu\nu}\Pi_{k}^{l})\pb^{4}\pb_{\nu}\eta^{k}}_{V_{12,3}^{\prime}}.
    \end{aligned}
\end{align}
Notice that $V_{12,1}^{\prime}=0$. The rest of two terms can be controlled as:
\begin{align}\label{estimateforV122prime}
    \begin{aligned}
       | V_{12,2}^{\prime}|&\lesssim\sum_{s=0}^{3}C_{4}^{s+1}\int_{0}^{t}\sigma\|\pb^{s} \pb((h^{ij}h^{\alpha\beta}-h^{i\alpha}h^{j\beta})\p_{\alpha}\eta^{\lambda}\p_{j}\eta^{l})\pb^{3-s}(\p_{\beta\mu}^{2}\eta^{\lambda})\|_{L^{2}(\Gamma_{1})} \|h^{\mu\nu}\|_{L^{\infty}(\Omega)}\|\p_{i}\pb^{4}\pb_{\nu}\eta^{k}\Pi_{k}^{l}\|_{L^{2}(\Gamma_{1})}\\&\lesssim\delta\int_{0}^{t}\|\p_{i}\pb^{4}\pb_{\nu}\eta^{k}\Pi_{k}^{l}\|_{L^{2}(\Gamma_{1})}^{2}+C(\delta)T\PP(E^{\kappa}),
    \end{aligned}
\end{align}
and using $\Pi_{k}^{l}=\Pi_{\gamma}^{l}\Pi_{k}^{\gamma}$, 
\begin{align}
    \begin{aligned}
        V_{12,3}^{\prime}&=\int_{0}^{t}\int_{\Gamma_{1}}\sigma (h^{ij}h^{\alpha\beta}-h^{i\alpha}h^{j\beta})\p_{\alpha}\eta^{\lambda}\pb^{4}\p_{\beta\mu}^{2}\eta^{\lambda}\p_{j}\eta^{l} \p_{i}(h^{\mu\nu}\Pi_{k}^{l})\pb^{4}\pb_{\nu}\eta^{k}\\&+\sum_{s=1}^{4}C_{4}^{s}\int_{0}^{t}\int_{\Gamma_{1}}\sigma\pb^{s} ((h^{ij}h^{\alpha\beta}-h^{i\alpha}h^{j\beta})\p_{\alpha}\eta^{\lambda}\p_{j}\eta^{l})\pb^{4-s}\p_{\beta\mu}^{2}\eta^{\lambda} \p_{i}(h^{\mu\nu}\Pi_{k}^{l})\pb^{4}\pb_{\nu}\eta^{k}\\&=\int_{0}^{t}\int_{\Gamma_{1}}\sigma (h^{ij}h^{\alpha\beta}-h^{i\alpha}h^{j\beta})\p_{\alpha}\eta^{\lambda}\pb^{4}\p_{\beta\mu}^{2}\eta^{\lambda}\p_{j}\eta^{l} h^{\mu\nu}\p_{i}\Pi_{\gamma}^{l}\Pi_{k}^{\gamma}\pb^{4}\pb_{\nu}\eta^{k}\\&+\sum_{s=1}^{4}C_{4}^{s}\int_{0}^{t}\int_{\Gamma_{1}}\sigma\pb^{s} ((h^{ij}h^{\alpha\beta}-h^{i\alpha}h^{j\beta})\p_{\alpha}\eta^{\lambda}\p_{j}\eta^{l})\pb^{4-s}\p_{\beta\mu}^{2}\eta^{\lambda} \p_{i}(h^{\mu\nu}\Pi_{k}^{l})\pb^{4}\pb_{\nu}\eta^{k}\\&=-\int_{0}^{t}\int_{\Gamma_{1}}\sigma (h^{ij}h^{\alpha\beta}-h^{i\alpha}h^{j\beta})\p_{\alpha}\eta^{\lambda}\pb^{4}\p_{\mu}\eta^{\lambda}\p_{j}\eta^{l} h^{\mu\nu}\p_{i}\Pi_{\gamma}^{l}\Pi_{k}^{\gamma}\pb^{4}\pb_{\beta\nu}^{2}\eta^{k}\\&-\int_{0}^{t}\int_{\Gamma_{1}}\sigma\pb^{4}\p_{\mu}\eta^{\lambda} \pb_{\beta}((h^{ij}h^{\alpha\beta}-h^{i\alpha}h^{j\beta})\p_{\alpha}\eta^{\lambda}\p_{j}\eta^{l} h^{\mu\nu}\p_{i}\Pi_{\gamma}^{l}\Pi_{k}^{\gamma})\pb^{4}\pb_{\nu}\eta^{k}\\&+\sum_{s=1}^{4}C_{4}^{s}\int_{0}^{t}\int_{\Gamma_{1}}\sigma\pb^{s} ((h^{ij}h^{\alpha\beta}-h^{i\alpha}h^{j\beta})\p_{\alpha}\eta^{\lambda}\p_{j}\eta^{l})\pb^{4-s}\p_{\beta\mu}^{2}\eta^{\lambda} \p_{i}(h^{\mu\nu}\Pi_{k}^{l})\pb^{4}\pb_{\nu}\eta^{k},
    \end{aligned}
\end{align}
where we used the orthogonal relationship again. Now, we can show that
\begin{align}\label{estimateforV123prime}
    \begin{aligned}
        |V_{12,3}^{\prime}|&\lesssim\int_{0}^{t}\sigma \|(h^{ij}h^{\alpha\beta}-h^{i\alpha}h^{j\beta})\p_{\alpha}\eta^{\lambda}\|_{L^{\infty}(\Gamma_{1})}\|\pb^{4}\p_{\mu}\eta^{\lambda}\|_{L^{2}(\Gamma_{1})}\|\p_{j}\eta^{l} h^{\mu\nu}\p_{i}\Pi_{\gamma}^{l}\|_{L^{\infty}(\Gamma_{1})}\|\Pi_{k}^{\gamma}\pb^{4}\pb_{\beta\nu}^{2}\eta^{k}\|_{L^{2}(\Gamma_{1})}\\&+\int_{0}^{t}\sigma\|\pb^{4}\p_{\mu}\eta^{\lambda}\|_{L^{2}(\Gamma_{1})}\| \pb_{\beta}((h^{ij}h^{\alpha\beta}-h^{i\alpha}h^{j\beta})\p_{\alpha}\eta^{\lambda}\p_{j}\eta^{l} h^{\mu\nu}\p_{i}\Pi_{\gamma}^{l}\Pi_{k}^{\gamma})\|_{L^{\infty}(\Gamma_{1})}\|\pb^{4}\pb_{\nu}\eta^{k}\|_{L^{2}(\Gamma_{1})}\\&+\sum_{s=0}^{3}C_{4}^{s+1}\int_{0}^{t}\sigma\|\pb^{s} \pb((h^{ij}h^{\alpha\beta}-h^{i\alpha}h^{j\beta})\p_{\alpha}\eta^{\lambda}\p_{j}\eta^{l})\pb^{3-s}\p_{\beta\mu}^{2}\eta^{\lambda}\|_{L^{2}(\Gamma_{1})}\| \p_{i}(h^{\mu\nu}\Pi_{k}^{l})\|_{L^{\infty}(\Gamma_{1})}\|\pb^{4}\pb_{\nu}\eta^{k}\|_{L^{2}(\Gamma_{1})}\\&\lesssim\delta\int_{0}^{t}\|\pb^{4}\pb_{\beta\nu}^{2}\eta^{k}\Pi_{k}^{\gamma}\|_{L^{2}(\Gamma_{1})}^{2}+C(\delta)T\PP(E^{\kappa})
    \end{aligned}
\end{align}
Combining \eqref{expansionforV12prime}, \eqref{estimateforV122prime} and \eqref{estimateforV123prime}, we have
\begin{align}\label{estimateforV12prime}
    \begin{aligned}
        |V_{12}^{\prime}|\lesssim\delta\int_{0}^{t}\|\Pi_{\lambda}^{l}\pb^{6}\eta^{\lambda}\|_{L^{2}(\Gamma_{1})}^{2}+C(\delta)T\PP(E^{\kappa}).
    \end{aligned}
\end{align}
Combining \eqref{expansionforV1prime}, \eqref{estimateforV11prime} and \eqref{estimateforV12prime}, we have
\begin{align}\label{estimateforV1prime}
    V_{1}^{\prime}=\int_{0}^{t}\int_{\Gamma_{1}}\sigma(\sqrt{h}h^{ij}\Pi_{\lambda}^{l}\pb^{4} \p_{\mu}\p_{j}\eta^{\lambda})( h^{\mu\nu}\p_{i}\pb^{4}\pb_{\nu}\eta^{k}\Pi_{k}^{l})+\mathcal{R},
\end{align}
where $|\mathcal{R}|\lesssim\delta\int_{0}^{t}\|\Pi_{\lambda}^{l}\pb^{6}\eta^{\lambda}\|_{L^{2}(\Gamma_{1})}^{2}+C(\delta)T\PP(E^{\kappa})$. Now,
$V_{2}^{\prime}$ and $V_{3}^{\prime}$ can be estimated in a similar fashion. We show $V_{2}^{\prime}$ as an example. With a little bit abuse of notation, we will suppress the index if it is not important in our analysis.
\begin{align}
    \begin{aligned}
        V_{2}^{\prime}&\sim+\int_{0}^{t}\int_{\Gamma_{1}}-\sigma\pb^{4}\pb_{\mu} \p_{i}(\sqrt{h}h^{ij}\p_{j}(\eta^{\lambda})) (n_{\kappa}^{\lambda}-n^{\lambda})n_{\kappa}^{l}h^{\mu\nu}\pb^{4}\pb_{\nu}\eta^{k}\Pi_{k}^{l}\\&+\int_{0}^{t}\int_{\Gamma_{1}}-\sigma\pb^{4} \p_{i}(\sqrt{h}h^{ij}\p_{j}(\eta^{\lambda})) \pb( (n_{\kappa}^{\lambda}-n^{\lambda})n_{\kappa}^{l})h^{.\nu}\pb^{4}\pb_{\nu}\eta^{k}\Pi_{k}^{l}\\&\sim+\int_{0}^{t}\int_{\Gamma_{1}}\sigma\pb^{4}\pb_{\mu} (\sqrt{h}h^{ij}\p_{j}(\eta^{\lambda})) (n_{\kappa}^{\lambda}-n^{\lambda})n_{\kappa}^{l}h^{\mu\nu}\pb^{4}\pb_{\nu}\p_{i}\eta^{k}\Pi_{k}^{l}\\&+\int_{0}^{t}\int_{\Gamma_{1}}\sigma\pb^{4}\pb_{\mu} (\sqrt{h}h^{ij}\p_{j}(\eta^{\lambda})) \p_{i}(n_{\kappa}^{\lambda}-n^{\lambda})n_{\kappa}^{l}h^{\mu\nu}\pb^{4}\pb_{\nu}\eta^{k}\Pi_{k}^{l}\\&+\int_{0}^{t}\int_{\Gamma_{1}}\sigma\pb^{4}\pb_{\mu} (\sqrt{h}h^{ij}\p_{j}(\eta^{\lambda})) (n_{\kappa}^{\lambda}-n^{\lambda})\p_{i}(n_{\kappa}^{l}h^{\mu\nu}\Pi_{k}^{l})\pb^{4}\pb_{\nu}\eta^{k}\\&+\int_{0}^{t}\int_{\Gamma_{1}}-\sigma\pb^{4} \p_{i}(\sqrt{h}h^{ij}\p_{j}(\eta^{\lambda})) \pb( (n_{\kappa}^{\lambda}-n^{\lambda})n_{\kappa}^{l})h^{.\nu}\pb^{4}\pb_{\nu}\eta^{k}\Pi_{k}^{l}
    \end{aligned}
\end{align}
Notice that
\begin{align}
    \begin{aligned}
        \|\pb^{5}(\sqrt{h}h^{ij}\p_{j}(\eta^{\lambda}))\|_{L^{2}(\Gamma_{1})}&\leq \|\pb^{4}(\sqrt{h}h^{ij}(\delta^{\alpha\lambda}-h^{kl}\p_{k}\eta^{\lambda}\p_{l}\eta^{\alpha})\pb\p_{j}\eta^{\alpha}
        ))\|_{L^{2}(\Gamma_{1})}\\&+\|\pb^{4}(\sqrt{h}(h^{ij}h^{kl}-h^{jl}h^{ik})\p_{j}\eta^{\lambda}\p_{k}\eta^{\lambda}\p_{l}\pb\eta^{\lambda})\|_{L^{2}(\Gamma_{1})}\\&\lesssim\PP(E^{\kappa})\|\eta\|_{H^{6}(\Gamma_{1})},
    \end{aligned}
\end{align}
\begin{align}
    \|n-n^{\kappa}\|_{L^{\infty}(\Gamma_{1})}+\|\pb(n-n^{\kappa})\|_{L^{\infty}(\Gamma_{1})}\lesssim \sqrt{\kappa}\PP(E^{\kappa}).
\end{align}
$V_{2}^{\prime}$ can clearly be bounded by
\begin{align}\label{estimateforV2prime}
    \begin{aligned}
        |V_{2}^{\prime}|&\lesssim\int_{0}^{t}\sqrt{\kappa}\PP(E^{\kappa})
    \|\eta\|_{H^{6}(\Gamma_{1})}\|\pb^{6}\eta^{k}\Pi_{k}^{l}\|_{L^{2}(\Gamma_{1})}+\PP_{0}+T\PP(E^{\kappa})\\&\lesssim\delta\int_{0}^{t}\|\pb^{6}\eta^{k}\Pi_{k}^{l}\|_{L^{2}(\Gamma_{1})}^{2}+C(\delta)(\PP_{0}+T\PP(E^{\kappa})).
    \end{aligned}
\end{align}
Similarly, 
\begin{align}\label{estimateforV3prime}
    |V_{3}^{\prime}|\lesssim\delta\int_{0}^{t}\|\pb^{6}\eta^{k}\Pi_{k}^{l}\|_{L^{2}(\Gamma_{1})}^{2}+C(\delta)(\PP_{0}+T\PP(E^{\kappa})).
\end{align}
Combining \eqref{expansionforVprime}, \eqref{estimateforV1prime}, \eqref{estimateforV2prime} and \eqref{estimateforV3prime}, we have
\begin{align}\label{estimateforVprime}
    V^{\prime}=\int_{0}^{t}\int_{\Gamma_{1}}\sigma(\sqrt{h}h^{ij}\Pi_{\lambda}^{l}\pb^{4} \p_{\mu}\p_{j}\eta^{\lambda})( h^{\mu\nu}\p_{i}\pb^{4}\pb_{\nu}\eta^{k}\Pi_{k}^{l})+\mathcal{R},
\end{align}
where $|\mathcal{R}|\lesssim\delta\int_{0}^{t}\|\Pi_{\lambda}^{l}\pb^{6}\eta^{\lambda}\|_{L^{2}(\Gamma_{1})}^{2}+C(\delta)T\PP(E^{\kappa})$. The estimate for $IV^{\prime}$ is quite direct.
\begin{align}
\begin{aligned}
       IV^{\prime}&= -\int_{0}^{t}\int_{\Gamma_{1}} \pb^{4}(\sqrt{h_{\kappa}} q n_{\kappa}^{l})h^{\mu\nu}\pb^{4}\pb_{\mu}\pb_{\nu}\eta^{k}\Pi_{k}^{l}\\&-\int_{0}^{t}\int_{\Gamma_{1}} \pb^{4}(\sqrt{h_{\kappa}} q n_{\kappa}^{l})\pb_{\mu}(h^{\mu\nu}\Pi_{k}^{l})\pb^{4}\pb_{\nu}\eta^{k}
\end{aligned}
\end{align}
Since $\|q\|_{4.5}\lesssim\PP(E^{\kappa})$, we have
\begin{align}\label{estimateforIVprime}
\begin{aligned}
     |IV^{\prime}|&\lesssim\int_{0}^{t}\| \pb^{4}(\sqrt{h_{\kappa}} q n_{\kappa}^{l})\|_{L^{2}(\Gamma_{1})}\|h^{\mu\nu}\|_{L^{\infty}(\Gamma_{1})}\|\pb^{4}\pb_{\mu}\pb_{\nu}\eta^{k}\Pi_{k}^{l}\|_{L^{2}(\Gamma_{1})}\\&+\int_{0}^{t}\int_{\Gamma_{1}} \|\pb^{4}(\sqrt{h_{\kappa}} q n_{\kappa}^{l})\|_{L^{2}(\Gamma_{1})}\|\pb_{\mu}(h^{\mu\nu}\Pi_{k}^{l})\|_{L^{\infty}(\Gamma_{1})}\|\pb^{4}\pb_{\nu}\eta^{k}\|_{L^{2}(\Gamma_{1})}\\&\lesssim\delta\int_{0}^{t}\|\Pi_{\lambda}^{l}\pb^{6}\eta^{\lambda}\|_{L^{2}(\Gamma_{1})}^{2}+C(\delta)T\PP(E^{\kappa}).
\end{aligned}   
\end{align}
Now, we analyze $VI^{\prime}$, which will give us the boundary regularity of $\sqrt{\kappa}\eta$.
\begin{align}\label{expansionforVIprime}
    \begin{aligned}
        VI^{\prime}&=\underbrace{\int_{0}^{t}\int_{\Gamma_{1}}\kappa\pb^{4}\pb_{\mu}(1-\bar{\lap})(v\cdot n ^{\kappa})n_{\kappa}^{l}h^{\mu\nu}\pb^{4}\pb_{\nu}\eta^{k}\Pi_{k}^{l}}_{VI_{1}^{\prime}}\\&\underbrace{+\sum_{s=0}^{4}C_{5}^{s}\int_{0}^{t}\int_{\Gamma_{1}}\kappa\pb^{s}(1-\bar{\lap})(v\cdot n ^{\kappa})(\pb^{5-s}n_{\kappa}^{l})h^{.\nu}\pb^{4}\pb_{\nu}\eta^{k}\Pi_{k}^{l}}_{VI_{2}^{\prime}}.
    \end{aligned}
\end{align}
We estimate $VI_{2}^{\prime}$ first.
\begin{align}\label{expansionforVI2prime}
    \begin{aligned}
        VI_{2}^{\prime}&=\underbrace{\int_{0}^{t}\int_{\Gamma_{1}}\kappa(1-\bar{\lap})(v\cdot n ^{\kappa})(\pb^{5}n_{\kappa}^{l})h^{.\nu}\pb^{4}\pb_{\nu}\eta^{k}\Pi_{k}^{l}}_{VI_{21}^{\prime}}\underbrace{+\sum_{s=1}^{4}C_{5}^{s}\int_{0}^{t}\int_{\Gamma_{1}}\kappa\pb^{s}(1-\bar{\lap})(v\cdot n ^{\kappa})(\pb^{5-s}n_{\kappa}^{l})h^{.\nu}\pb^{4}\pb_{\nu}\eta^{k}\Pi_{k}^{l}}_{VI_{22}^{\prime}},
    \end{aligned}
\end{align}
we control $VI_{21}^{\prime}$ as
\begin{align}\label{estimateforVI21prime}
\begin{aligned}
    |VI_{21}^{\prime}|&\lesssim\int_{0}^{t}\kappa\|(1-\bar{\lap})(v\cdot n ^{\kappa})\|_{L^{\infty}(\Gamma_{1})}\|(\pb^{5}n_{\kappa}^{l})\|_{L^{2}(\Gamma_{1})}\|h^{.\nu}\Pi_{k}^{l}\|_{L^{\infty}(\Gamma_{1})}\|\pb^{4}\pb_{\nu}\eta^{k}\|_{L^{2}(\Gamma_{1})}\\&\lesssim\int_{0}^{t}\kappa\|\eta\|_{H^{6}(\Gamma_{1})}\PP(E^{\kappa})\\&\lesssim \PP_{0}+T\PP(E^{\kappa}).
\end{aligned}  
\end{align}
For $VI_{22}^{\prime}$, we can write
\begin{align}\label{expansionforVI22prime}
    \begin{aligned}
        VI_{22}^{\prime}&=\underbrace{\sum_{s=1}^{4}C_{5}^{s}\int_{0}^{t}\int_{\Gamma_{1}}\kappa\pb^{s}(v\cdot n ^{\kappa})(\pb^{5-s}n_{\kappa}^{l})h^{.\nu}\pb^{4}\pb_{\nu}\eta^{k}\Pi_{k}^{l}}_{V_{22,1}^{\prime}}\\&\underbrace{+\sum_{s=1}^{4}C_{5}^{s}\int_{0}^{t}\int_{\Gamma_{1}}\kappa\pb^{s}\pb_{\gamma}(v\cdot n ^{\kappa})\pb_{\gamma}(\pb^{5-s}n_{\kappa}^{l})h^{.\nu}\pb^{4}\pb_{\nu}\eta^{k}\Pi_{k}^{l}}_{V_{22,2}^{\prime}}\\&\underbrace{+\sum_{s=1}^{4}C_{5}^{s}\int_{0}^{t}\int_{\Gamma_{1}}\kappa\pb^{s}\p_{\gamma}(v\cdot n ^{\kappa})(\pb^{5-s}n_{\kappa}^{l})h^{.\nu}\pb^{4}\pb_{\nu}\pb_{\gamma}\eta^{k}\Pi_{k}^{l}}_{VI_{22,3}^{\prime}}\\&\underbrace{+\sum_{s=1}^{4}C_{5}^{s}\int_{0}^{t}\int_{\Gamma_{1}}\kappa\pb^{s}(v\cdot n ^{\kappa})(\pb^{5-s}n_{\kappa}^{l})\pb_{\gamma}(h^{.\nu}\Pi_{k}^{l})\pb^{4}\pb_{\nu}\eta^{k}}_{V_{22,4}^{\prime}}.
    \end{aligned}
\end{align}
A term by term estimate shows that
\begin{align}\label{estimateforVI221prime}
    \begin{aligned}
        |VI_{22,1}^{\prime}|&\lesssim\sum_{s=0}^{3}C_{5}^{s+1}\int_{0}^{t}\kappa\|\pb^{s}\pb(v\cdot n ^{\kappa})(\pb^{4-s}n_{\kappa}^{l})\|_{L^{2}(\Gamma_{1})}\|h^{.\nu}\Pi_{k}^{l}\|_{L^{\infty}(\Gamma_{1})}\|\pb^{4}\pb_{\nu}\eta^{k}\|_{L^{2}(\Gamma_{1})}\\&\lesssim \PP_{0}+T\PP(E^{\kappa}),
    \end{aligned}
\end{align}
\begin{align}\label{estimateforVI222prime}
    \begin{aligned}
        |VI_{22,2}^{\prime}|&\lesssim C_{5}^{1}\int_{0}^{t}\kappa\|\pb\pb_{\gamma}(v\cdot n ^{\kappa})\|_{L^{\infty}(\Gamma_{1})}\|\pb_{\gamma}(\pb^{4}n_{\kappa}^{l})\|_{L^{2}(\Gamma_{1})}\|h^{.\nu}\Pi_{k}^{l}\|_{L^{\infty}(\Gamma_{1})}\|\pb^{4}\pb_{\nu}\eta^{k}\|_{L^{2}(\Gamma_{1})}\\&+\sum_{s=0}^{1}C_{5}^{s+2}\int_{0}^{t}\kappa\|\pb^{s+2}\pb_{\gamma}(v\cdot n ^{\kappa})(\pb^{1-s}\pb_{\gamma}\pb^{2}n_{\kappa}^{l})\|_{L^{2}(\Gamma_{1})}\|h^{.\nu}\Pi_{k}^{l}\|_{L^{\infty}(\Gamma_{1})}\|\pb^{4}\pb_{\nu}\eta^{k}\|_{L^{2}(\Gamma_{1})}\\&+C_{5}^{4}\int_{0}^{t}\kappa\|\pb^{4}\pb_{\gamma}(v\cdot n ^{\kappa})\|_{L^{2}(\Gamma_{1})}\|\pb_{\gamma}(\pb n_{\kappa}^{l})h^{.\nu}\Pi_{k}^{l}\|_{L^{\infty}(\Gamma_{1})}\|\pb^{4}\pb_{\nu}\eta^{k}\|_{L^{2}(\Gamma_{1})}\\&\lesssim\int_{0}^{t}\kappa\PP(E^{\kappa})\|\eta\|_{H^{6}(\Gamma_{1})}+\int_{0}^{t}\PP(E^{\kappa})+\int_{0}^{t}\kappa (\|v\|_{H^{5}(\Gamma_{1})}+\|\eta\|_{H^{6}(\Gamma_{1})})\PP(E^{\kappa})\\&\lesssim\PP_{0}+T\PP(E^{\kappa}),
    \end{aligned}
\end{align}
\begin{align}\label{estimateforVI223prime}
    \begin{aligned}
        |VI_{22,3}^{\prime}|&\lesssim C_{5}^{4}\int_{0}^{t}\kappa\|\pb^{4}\p_{\gamma}(v\cdot n ^{\kappa})\|_{L^{2}(\Gamma_{1})}\|(\pb n_{\kappa}^{l})h^{.\nu}\|_{L^{\infty}(\Gamma_{1})}\|\pb^{4}\pb_{\nu}\pb_{\gamma}\eta^{k}\Pi_{k}^{l}\|_{L^{2}(\Gamma_{1})}\\&+\sum_{s=0}^{2}C_{5}^{s+1}\int_{0}^{t}\kappa\|\pb^{s}\pb\p_{\gamma}(v\cdot n ^{\kappa})(\pb^{2-s}\pb^{2}n_{\kappa}^{l})\|_{L^{2}(\Gamma_{1})}\|h^{.\nu}\Pi_{k}^{l}\|_{L^{\infty}(\Gamma_{1})}\|\pb^{4}\pb_{\nu}\pb_{\gamma}\eta^{k}\|_{L^{2}(\Gamma_{1})}\\&\lesssim \PP_{0}+T\PP(E^{\kappa}),
    \end{aligned}
\end{align}
\begin{align}\label{estimateforVI224prime}
    \begin{aligned}
        |VI_{22,4}^{\prime}|&\lesssim\sum_{s=0}^{3}C_{5}^{s+1}\int_{0}^{t}\kappa\|\pb^{s}\pb(v\cdot n ^{\kappa})(\pb^{4-s}n_{\kappa}^{l})\|_{L^{2}(\Gamma_{1})}\|\pb_{\gamma}(h^{.\nu}\Pi_{k}^{l})\|_{L^{\infty}(\Gamma_{1})}\|\pb^{4}\pb_{\nu}\eta^{k}\|_{L^{2}(\Gamma_{1})}\\&\lesssim \PP_{0}+T\PP(E^{\kappa}).
    \end{aligned}
\end{align}
Combining 
\eqref{expansionforVI22prime}, \eqref{estimateforVI221prime}, \eqref{estimateforVI222prime}, \eqref{estimateforVI223prime} and \eqref{estimateforVI224prime}, we have
\begin{align}\label{estimateforVI22prime}
    |VI_{22}^{\prime}|\lesssim\PP_{0}+T\PP(E^{\kappa}).
\end{align}
Combining \eqref{expansionforVI2prime}, \eqref{estimateforVI21prime} and \eqref{estimateforVI22prime}, we have
\begin{align}\label{estimateforVI2prime}
    |VI_{2}^{\prime}|\lesssim\PP_{0}+T\PP(E^{\kappa}).
\end{align}
Now, we estimate the major part of $VI^{\prime}$ as
\begin{align}\label{expansionforVI1prime}
    \begin{aligned}
        VI_{1}^{\prime}&=\int_{0}^{t}\int_{\Gamma_{1}}\kappa\pb^{4}\pb_{\mu}(v\cdot n ^{\kappa})n_{\kappa}^{l}h^{\mu\nu}\pb^{4}\pb_{\nu}\eta^{k}\Pi_{k}^{l}+\int_{0}^{t}\int_{\Gamma_{1}}\kappa\pb^{4}\pb_{\mu}\p_{\gamma}(v\cdot n ^{\kappa})n_{\kappa}^{l}h^{\mu\nu}\pb^{4}\pb_{\gamma}\pb_{\nu}\eta^{k}\Pi_{k}^{l}\\&+\int_{0}^{t}\int_{\Gamma_{1}}\kappa\pb^{4}\pb_{\mu}\p_{\gamma}(v\cdot n ^{\kappa})\pb_{\gamma}(n_{\kappa}^{l}h^{\mu\nu}\Pi_{k}^{l})\pb^{4}\pb_{\nu}\eta^{k}\\&=\underbrace{\int_{0}^{t}\int_{\Gamma_{1}}\kappa\pb^{4}\pb_{\mu}v\cdot n ^{\kappa}n_{\kappa}^{l}h^{\mu\nu}\pb^{4}\pb_{\nu}\eta^{k}\Pi_{k}^{l}+\int_{0}^{t}\int_{\Gamma_{1}}\kappa\pb^{4}\pb_{\mu}\p_{\gamma}v\cdot n ^{\kappa}n_{\kappa}^{l}h^{\mu\nu}\pb^{4}\pb_{\gamma}\pb_{\nu}\eta^{k}\Pi_{k}^{l}}_{VI_{11}^{\prime}}\\&\underbrace{+\sum_{s=0}^{4}C_{5}^{s}\int_{0}^{t}\int_{\Gamma_{1}}\kappa\pb^{s}v\cdot\pb^{5-s} n ^{\kappa}n_{\kappa}^{l}h^{.\nu}\pb^{4}\pb_{\nu}\eta^{k}\Pi_{k}^{l}+\sum_{s=1}^{4}C_{6}^{s}\int_{0}^{t}\int_{\Gamma_{1}}\kappa\pb^{s}v\cdot\pb^{6-s} n ^{\kappa}n_{\kappa}^{l}h^{\nu}\pb^{5}\pb_{\nu}\eta^{k}\Pi_{k}^{l}}_{VI_{12}^{\prime}}\\&\underbrace{+\int_{0}^{t}\int_{\Gamma_{1}}\kappa\pb^{4}\pb_{\mu}\p_{\gamma}(v\cdot n ^{\kappa})\pb_{\gamma}(n_{\kappa}^{l}h^{\mu\nu}\Pi_{k}^{l})\pb^{4}\pb_{\nu}\eta^{k}}_{VI_{13}^{\prime}}\underbrace{+\int_{0}^{t}\int_{\Gamma_{1}}\kappa v^{i}\pb^{4}\pb_{\mu}\pb_{\gamma}n_{\kappa}^{i}n_{\kappa}^{l}h^{\mu\nu}\pb^{4}\pb_{\gamma}\pb_{\nu}\eta^{k}\Pi_{k}^{l}}_{VI_{14}^{\prime}}.
    \end{aligned}
\end{align}
The estimate of $VI_{12}^{\prime}$ is similar to the estimate of $VI_{2}^{\prime}$,
\begin{align}\label{estimateforVI12prime}
    \begin{aligned}
        |VI_{12}^{\prime}|&\lesssim\int_{0}^{t}\kappa \|v\|_{L^{\infty}(\Gamma_{1})}\|\pb^{5} n ^{\kappa}\|_{L^{2}(\Gamma_{1})}\|n_{\kappa}^{l}h^{.\nu}\Pi_{k}^{l}\|_{L^{\infty}(\Gamma_{1})}\|\pb^{4}\pb_{\nu}\eta^{k}\|_{L^{2}(\Gamma_{1})}\\&+\sum_{s=0}^{3}C_{5}^{s+1}\int_{0}^{t}\kappa\|\pb^{s+1}v\cdot\pb^{4-s} n ^{\kappa}\|_{L^{2}(\Gamma_{1})}\|n_{\kappa}^{l}h^{.\nu}\Pi_{k}^{l}\|_{L^{\infty}(\Gamma_{1})}\|\pb^{4}\pb_{\nu}\eta^{k}\|_{L^{2}(\Gamma_{1})}\\&+6\int_{0}^{t}\kappa\|\pb v\|_{L^{\infty}(\Gamma_{1})}\|\pb^{5} n ^{\kappa}\|_{L^{2}(\Gamma_{1})}\|n_{\kappa}^{l}h^{\nu}\Pi_{k}^{l}\|_{L^{\infty}(\Gamma_{1})}\|\pb^{5}\pb_{\nu}\eta^{k}\|_{L^{2}(\Gamma_{1})}\\&+\sum_{s=0}^{2}C_{6}^{s+2}\int_{0}^{t}\kappa\|\pb^{s+2}v\cdot\pb^{4-s} n ^{\kappa}\|_{L^{2}(\Gamma_{1})}\|n_{\kappa}^{l}h^{\nu}\Pi_{k}^{l}\|_{L^{\infty}(\Gamma_{1})}\|\pb^{5}\pb_{\nu}\eta^{k}\|_{L^{2}(\Gamma_{1})}\\&\lesssim\PP_{0}+T\PP(E^{\kappa})+\int_{0}^{t}\kappa\PP(E^{\kappa})\|\eta\|_{H^{6}(\Gamma_{1})}^{2}\\&\lesssim \PP_{0}+T\PP(E^{\kappa}).
    \end{aligned}
\end{align}
For $VI_{13}^{\prime}$, we can show that
\begin{align}
    \begin{aligned}       VI_{13}^{\prime}&=-\int_{0}^{t}\int_{\Gamma_{1}}\kappa\pb^{4}\p_{\gamma}(v\cdot n ^{\kappa})\pb_{\mu}\pb_{\gamma}(n_{\kappa}^{l}h^{\mu\nu}\Pi_{k}^{l})\pb^{4}\pb_{\nu}\eta^{k}-\int_{0}^{t}\int_{\Gamma_{1}}\kappa\pb^{4}\p_{\gamma}(v\cdot n ^{\kappa})\pb_{\gamma}(n_{\kappa}^{l}h^{\mu\nu}\Pi_{k}^{l})\pb^{4}\pb_{\mu}\pb_{\nu}\eta^{k}.
    \end{aligned}
\end{align}
Thus,
\begin{align}\label{estimateforVI13prime}
    \begin{aligned}
        |VI_{13}^{\prime}|&\lesssim\int_{0}^{t}\kappa\|\pb^{4}\p_{\gamma}(v\cdot n ^{\kappa})\|_{L^{2}(\Gamma_{1})}\|\pb_{\mu}\pb_{\gamma}(n_{\kappa}^{l}h^{\mu\nu}\Pi_{k}^{l})\|_{L^{\infty}(\Gamma_{1})}\|\pb^{4}\pb_{\nu}\eta^{k}\|_{L^{2}(\Gamma_{1})}\\&+\int_{0}^{t}\kappa\|\pb^{4}\p_{\gamma}(v\cdot n ^{\kappa})\|_{L^{2}(\Gamma_{1})}\|\pb_{\gamma}(n_{\kappa}^{l}h^{\mu\nu}\Pi_{k}^{l})\|_{L^{\infty}(\Gamma_{1})}\|\pb^{4}\pb_{\mu}\pb_{\nu}\eta^{k}\|_{L^{2}(\Gamma_{1})}\\&\lesssim\int_{0}^{t}\kappa(\|v\|_{H^{5}}+\|\eta\|_{H^{6}})\PP(E^{\kappa})\|\eta\|_{H^{6}(\Gamma_{1})}\\&\lesssim \PP_{0}+T\PP(E^{\kappa}).
    \end{aligned}
\end{align}
Now, we estimate the hardest term $VI_{14}^{\prime}$. Writing $v^{i}=v\cdot n^{\kappa}n_{\kappa}^{i}+v\cdot\frac{\p_{1}\eta_{\kappa}}{|\p_{1}\eta_{\kappa}|^{2}}\p_{1}\eta_{\kappa}^{i}+v\cdot\frac{\p_{2}\eta_{\kappa}}{|\p_{2}\eta_{\kappa}|^{2}}\p_{2}\eta_{\kappa}^{i}$, utilizing $n_{\kappa}^{i}n_{\kappa}^{i}=1$ and $n_{\kappa}^{i}\pb_{k}\eta^{i}=0$, $k=1,2$ we have
\begin{align}\label{expansionforVI14prime}
    \begin{aligned}
        VI_{14}^{\prime}&=\underbrace{-\frac{1}{2}\sum_{s=1}^{5}C_{6}^{s}\int_{0}^{t}\int_{\Gamma_{1}}\kappa v\cdot n^{\kappa}\pb^{s}n_{\kappa}^{i}\pb^{6-s}n_{\kappa}^{i}n_{\kappa}^{l}h^{.\nu}\pb^{5}\pb_{\nu}\eta^{k}\Pi_{k}^{l}}_{VI_{14,1}^{\prime}}\\&\underbrace{+\int_{0}^{t}\int_{\Gamma_{1}}\kappa v\cdot\frac{\p_{1}\eta_{\kappa}}{|\p_{1}\eta_{\kappa}|^{2}}\p_{1}\eta_{\kappa}^{i}\pb^{4}\pb_{\mu}\pb_{\gamma}n_{\kappa}^{i}n_{\kappa}^{l}h^{\mu\nu}\pb^{4}\pb_{\gamma}\pb_{\nu}\eta^{k}\Pi_{k}^{l}}_{VI_{14,2}^{\prime}}\\&\underbrace{+\int_{0}^{t}\int_{\Gamma_{1}}\kappa v\cdot\frac{\p_{2}\eta_{\kappa}}{|\p_{2}\eta_{\kappa}|^{2}}\p_{2}\eta_{\kappa}^{i}\pb^{4}\pb_{\mu}\pb_{\gamma}n_{\kappa}^{i}n_{\kappa}^{l}h^{\mu\nu}\pb^{4}\pb_{\gamma}\pb_{\nu}\eta^{k}\Pi_{k}^{l}}_{VI_{14,3}^{\prime}}.
    \end{aligned}
\end{align}
From direct calculation,
\begin{align}\label{estimateforVI141prime}
    \begin{aligned}
        |VI_{14,1}^{\prime}|&\lesssim 6\int_{0}^{t}\kappa \|v\cdot n^{\kappa}\pb n_{\kappa}^{i}\|_{L^{\infty}(\Gamma_{1})}\|\pb^{5}n_{\kappa}^{i}\|_{L^{2}(\Gamma_{1})}\|n_{\kappa}^{l}h^{.\nu}\Pi_{k}^{l}\|_{L^{\infty}(\Gamma_{1})}\|\pb^{5}\pb_{\nu}\eta^{k}\|_{L^{2}(\Gamma_{1})}\\&+\frac{1}{2}\sum_{s=0}^{2}C_{6}^{s+2}\int_{0}^{t}\kappa \|v\cdot n^{\kappa}\|_{L^{\infty}(\Gamma_{1})}\|\pb^{s}\pb^{2}n_{\kappa}^{i}\pb^{4-s}n_{\kappa}^{i}\|_{L^{2}(\Gamma_{1})}\|n_{\kappa}^{l}h^{.\nu}\Pi_{k}^{l}\|_{L^{\infty}(\Gamma_{1})}\|\pb^{5}\pb_{\nu}\eta^{k}\|_{L^{2}(\Gamma_{1})}\\&\lesssim\int_{0}^{t}\kappa\PP(E^{\kappa})\|\eta\|_{H^{6}(\Gamma_{1})}^{2}+\int_{0}^{t}\kappa\PP(E^{\kappa})\|\eta\|_{H^{6}(\Gamma_{1})}\\&\lesssim\PP_{0}+T\PP(E^{\kappa}).
    \end{aligned}
\end{align}
While the structure of $VI_{14,2}^{\prime}$ and $VI_{14,3}^{\prime}$ is the same. It suffices to show $VI_{14,2}^{\prime}$ as an example.
\begin{align}\label{expansionforVI142prime}
    \begin{aligned}
        VI_{14,2}^{\prime}&=\underbrace{-\int_{0}^{t}\int_{\Gamma_{1}}\kappa v\cdot\frac{\p_{1}\eta_{\kappa}}{|\p_{1}\eta_{\kappa}|^{2}}\pb^{4}\pb_{\mu}\pb_{\gamma}\p_{1}\eta_{\kappa}^{i}n_{\kappa}^{i}n_{\kappa}^{l}h^{\mu\nu}\pb^{4}\pb_{\gamma}\pb_{\nu}\eta^{k}\Pi_{k}^{l}}_{VI_{14,21}^{\prime}}\\&\underbrace{-\sum_{s=1}^{5}C_{6}^{s}\int_{0}^{t}\int_{\Gamma_{1}}\kappa v\cdot\frac{\p_{1}\eta_{\kappa}}{|\p_{1}\eta_{\kappa}|^{2}}\pb^{s}\p_{1}\eta_{\kappa}^{i}\pb^{6-s}n_{\kappa}^{i}n_{\kappa}^{l}h^{.\nu}\pb^{5}\pb_{\nu}\eta^{k}\Pi_{k}^{l}}_{VI_{14,22}^{\prime}}.
    \end{aligned}
\end{align}
For $VI_{14,22}^{\prime}$, we can estimate it as
\begin{align}\label{estimateforVI1422prime}
    \begin{aligned}
        |VI_{14,22}^{\prime}|&\lesssim 6\int_{0}^{t}\kappa \|v\cdot\frac{\p_{1}\eta_{\kappa}}{|\p_{1}\eta_{\kappa}|^{2}}\pb\p_{1}\eta_{\kappa}^{i}\|_{L^{\infty}(\Gamma_{1})}\|\pb^{5}n_{\kappa}^{i}\|_{L^{2}(\Gamma_{1})}\|n_{\kappa}^{l}h^{.\nu}\Pi_{k}^{l}\|_{L^{\infty}(\Gamma_{1})}\|\pb^{5}\pb_{\nu}\eta^{k}\|_{L^{2}(\Gamma_{1})}\\&+6\int_{0}^{t}\kappa \|v\cdot\frac{\p_{1}\eta_{\kappa}}{|\p_{1}\eta_{\kappa}|^{2}}\|_{L^{\infty}(\Gamma_{1})}\|\pb^{5}\p_{1}\eta_{\kappa}^{i}\|_{L^{2}(\Gamma_{1})}\|\pb n_{\kappa}^{i}n_{\kappa}^{l}h^{.\nu}\Pi_{k}^{l}\|_{L^{\infty}(\Gamma_{1})}\|\pb^{5}\pb_{\nu}\eta^{k}\|_{L^{2}(\Gamma_{1})}\\&+\sum_{s=0}^{4}C_{6}^{s+2}\int_{0}^{t}\kappa \|v\cdot\frac{\p_{1}\eta_{\kappa}}{|\p_{1}\eta_{\kappa}|^{2}}\|_{L^{\infty}(\Gamma_{1})}\|\pb^{s}\pb^{2}\p_{1}\eta_{\kappa}^{i}\pb^{4-s}n_{\kappa}^{i}\|_{L^{2}(\Gamma_{1})}\|n_{\kappa}^{l}h^{.\nu}\Pi_{k}^{l}\|_{L^{\infty}(\Gamma_{1})}\|\pb^{5}\pb_{\nu}\eta^{k}\|_{L^{2}(\Gamma_{1})}\\&\lesssim\PP_{0}+T\PP(E^{\kappa}).
    \end{aligned}
\end{align}
For $VI_{14,21}^{\prime}$, we first write it into a more symmetric form:
\begin{align}\label{expansionforVI1421prime}
    \begin{aligned}
        VI_{14,21}^{\prime}&=\underbrace{-\int_{0}^{t}\int_{\Gamma_{1}}\kappa v\cdot\frac{\p_{1}\eta_{\kappa}}{|\p_{1}\eta_{\kappa}|^{2}}\pb^{4}\pb_{\mu}\pb_{\gamma}\p_{1}\eta_{\kappa}^{i}n_{\kappa}^{i}n_{\kappa}^{l}h^{\mu\nu}\pb^{4}\pb_{\gamma}\pb_{\nu}\eta^{k}n_{\kappa}^{k}n_{\kappa}^{l}}_{VI_{14,211}^{\prime}}\\&\underbrace{-\int_{0}^{t}\int_{\Gamma_{1}}\kappa v\cdot\frac{\p_{1}\eta_{\kappa}}{|\p_{1}\eta_{\kappa}|^{2}}\pb^{4}\pb_{\mu}\pb_{\gamma}\p_{1}\eta_{\kappa}^{i}n_{\kappa}^{i}n_{\kappa}^{l}h^{\mu\nu}\pb^{4}\pb_{\gamma}\pb_{\nu}\eta^{k}(\Pi_{k}^{l}-(\Pi_{\kappa})_{k}^{l})}_{VI_{14,212}^{\prime}}.
    \end{aligned}
\end{align}
The second term $VI_{14,212}$ can be controlled by a standard $\kappa\|f\|_{H^{s+1}(\Gamma_{1})}\lesssim \|f\|_{s}$. Since $\Pi=n\otimes n$, and we have analyzed that $n$ can be viewed as a smooth function of $\pb\eta$. Thus, we have
\begin{align}\label{linftyestiwrtPi}
    \|\Pi_{k}^{l}-(\Pi_{\kappa})_{k}^{l}\|_{L^{\infty}(\Gamma_{1})}\lesssim\kappa\PP(E^{\kappa})\|\p\eta\|_{H^{2.5}(\Gamma_{1})}.
\end{align}
Now, a standard $L^{2}$ estimate gives us
\begin{align}\label{estimateforVI14212prime}
    \begin{aligned}
        |VI_{14,212}^{\prime}|&\lesssim\int_{0}^{t}\kappa \|v\cdot\frac{\p_{1}\eta_{\kappa}}{|\p_{1}\eta_{\kappa}|^{2}}\|_{L^{\infty}(\Gamma_{1})}\|\pb^{4}\pb_{\mu}\pb_{\gamma}\p_{1}\eta_{\kappa}^{i}\|_{L^{2}(\Gamma_{1})}\|n_{\kappa}^{i}n_{\kappa}^{l}h^{\mu\nu}\|_{L^{\infty}(\Gamma_{1})}\|\pb^{4}\pb_{\gamma}\pb_{\nu}\eta^{k}\|_{L^{2}(\Gamma_{1})}\|(\Pi_{k}^{l}-(\Pi_{\kappa})_{k}^{l})\|_{L^{\infty}(\Gamma_{1})}\\&\lesssim\int_{0}^{t}\kappa\PP(E^{\kappa})\|\eta\|_{H^{7}(\Gamma_{1})}\|\eta\|_{H^{6}(\Gamma_{1})}\kappa\|\eta\|_{H^{3.5}(\Gamma_{1})}\\&\lesssim\int_{0}^{t}\kappa\PP(E^{\kappa})\|\eta\|_{H^{6}(\Gamma_{1})}^{2}\\&\lesssim T\PP(E^{\kappa}).    \end{aligned}
\end{align}
For $VI_{14,211}^{\prime}$, we utilize the double horizontal convolution,
\begin{align}
    \begin{aligned}
        VI_{14,211}^{\prime}&=-\int_{0}^{t}\int_{\Gamma_{1}}\kappa \pb^{4}\pb_{\mu}\pb_{\gamma}\p_{1}\Lambda_{\kappa}\eta^{i}\Lambda_{\kappa}(v\cdot\frac{\p_{1}\eta_{\kappa}}{|\p_{1}\eta_{\kappa}|^{2}}n_{\kappa}^{i}n_{\kappa}^{l}h^{\mu\nu}\pb^{4}\pb_{\gamma}\pb_{\nu}\eta^{k}n_{\kappa}^{k}n_{\kappa}^{l})\\&=-\int_{0}^{t}\int_{\Gamma_{1}}\kappa \pb^{4}\pb_{\mu}\pb_{\gamma}\p_{1}\Lambda_{\kappa}\eta^{i}(v\cdot\frac{\p_{1}\eta_{\kappa}}{|\p_{1}\eta_{\kappa}|^{2}}n_{\kappa}^{i}h^{\mu\nu}n_{\kappa}^{k}\pb^{4}\pb_{\gamma}\pb_{\nu}\Lambda_{\kappa}\eta^{k})\\&-\int_{0}^{t}\int_{\Gamma_{1}}\kappa \pb^{4}\pb_{\mu}\pb_{\gamma}\p_{1}\Lambda_{\kappa}\eta^{i}[\Lambda_{\kappa},(v\cdot\frac{\p_{1}\eta_{\kappa}}{|\p_{1}\eta_{\kappa}|^{2}}n_{\kappa}^{i}h^{\mu\nu}n_{\kappa}^{k})]\pb^{4}\pb_{\gamma}\pb_{\nu}\eta^{k}\\&=\frac{1}{2}\int_{0}^{t}\int_{\Gamma_{1}}\kappa \pb^{4}\pb_{\mu}\pb_{\gamma}\Lambda_{\kappa}\eta^{i}\p_{1}(v\cdot\frac{\p_{1}\eta_{\kappa}}{|\p_{1}\eta_{\kappa}|^{2}}n_{\kappa}^{i}h^{\mu\nu}n_{\kappa}^{k})\pb^{4}\pb_{\gamma}\pb_{\nu}\Lambda_{\kappa}\eta^{k}\\&+\int_{0}^{t}\int_{\Gamma_{1}}\kappa \pb^{4}\pb_{\mu}\pb_{\gamma}\Lambda_{\kappa}\eta^{i}[\Lambda_{\kappa},\p_{1}(v\cdot\frac{\p_{1}\eta_{\kappa}}{|\p_{1}\eta_{\kappa}|^{2}}n_{\kappa}^{i}h^{\mu\nu}n_{\kappa}^{k}]\pb^{4}\pb_{\gamma}\pb_{\nu}\eta^{k})\\&+\int_{0}^{t}\int_{\Gamma_{1}}\kappa \pb^{4}\pb_{\mu}\pb_{\gamma}\Lambda_{\kappa}\eta^{i}[\Lambda_{\kappa},(v\cdot\frac{\p_{1}\eta_{\kappa}}{|\p_{1}\eta_{\kappa}|^{2}}n_{\kappa}^{i}h^{\mu\nu}n_{\kappa}^{k}]\p_{1}\pb^{4}\pb_{\gamma}\pb_{\nu}\eta^{k}).
    \end{aligned}
\end{align}
Invoking \eqref{commutationtypelemma}, we have
\begin{align}\label{estimateforVI14211}
    \begin{aligned}
        |VI_{14,211}^{\prime}|&\lesssim\frac{1}{2}\int_{0}^{t}\kappa \|\pb^{4}\pb_{\mu}\pb_{\gamma}\Lambda_{\kappa}\eta^{i}\|_{L^{2}(\Gamma_{1})}\|\p_{1}(v\cdot\frac{\p_{1}\eta_{\kappa}}{|\p_{1}\eta_{\kappa}|^{2}}n_{\kappa}^{i}h^{\mu\nu}n_{\kappa}^{k})\|_{L^{\infty}(\Gamma_{1})}\|\pb^{4}\pb_{\gamma}\pb_{\nu}\Lambda_{\kappa}\eta^{k}\|_{L^{2}(\Gamma_{1})}\\&+\int_{0}^{t}\kappa \|\pb^{4}\pb_{\mu}\pb_{\gamma}\Lambda_{\kappa}\eta^{i}\|_{L^{2}(\Gamma_{1})}\|[\Lambda_{\kappa},\p_{1}(v\cdot\frac{\p_{1}\eta_{\kappa}}{|\p_{1}\eta_{\kappa}|^{2}}n_{\kappa}^{i}h^{\mu\nu}n_{\kappa}^{k})]\pb^{4}\pb_{\gamma}\pb_{\nu}\eta^{k}\|_{L^{2}(\Gamma_{1})}\\&+\int_{0}^{t}\kappa \|\pb^{4}\pb_{\mu}\pb_{\gamma}\Lambda_{\kappa}\eta^{i}\|_{L^{2}(\Gamma_{1})}\|[\Lambda_{\kappa},(v\cdot\frac{\p_{1}\eta_{\kappa}}{|\p_{1}\eta_{\kappa}|^{2}}n_{\kappa}^{i}h^{\mu\nu}n_{\kappa}^{k})]\p_{1}\pb^{4}\pb_{\gamma}\pb_{\nu}\eta^{k}\|_{L^{2}(\Gamma_{1})}\\&\lesssim\int_{0}^{t}\PP(E^{\kappa})+\int_{0}^{t}\kappa\|\eta\|_{H^{6}(\Gamma_{1})}\|(v\cdot\frac{\p_{1}\eta_{\kappa}}{|\p_{1}\eta_{\kappa}|^{2}}n_{\kappa}^{i}h^{\mu\nu}n_{\kappa}^{k})\|_{H^{2.5}(\Gamma_{1})}\|\eta\|_{H^{6}(\Gamma_{1})}\\&\lesssim \PP_{0}+T\PP(E^{\kappa}).
    \end{aligned}
\end{align}
Combining \eqref{expansionforVI1421prime}, \eqref{estimateforVI14212prime} and \eqref{estimateforVI14212prime}, we have
\begin{align}\label{estimateforVI1421prime}
    |VI_{14,21}^{\prime}|\lesssim\PP_{0}+T\PP(E^{\kappa}).
\end{align}
Combining \eqref{expansionforVI142prime}, \eqref{estimateforVI1422prime} and \eqref{estimateforVI1421prime}, we have
\begin{align}\label{estimateforVI142prime}
    |VI_{14,2}^{\prime}|\lesssim\PP_{0}+T\PP(E^{\kappa}).
\end{align}
Similarly, we have
\begin{align}\label{estimateforVI143prime}
    |VI_{14,3}^{\prime}|\lesssim\PP_{0}+T\PP(E^{\kappa}).
\end{align}
Combining \eqref{expansionforVI14prime}, \eqref{estimateforVI141prime}, \eqref{estimateforVI142prime} and \eqref{estimateforVI143prime}, we have
\begin{align}\label{estimateforVI14prime}
    |VI_{14}^{\prime}|\lesssim \PP_{0}+T\PP(E^{\kappa}).
\end{align}
Now, we focus on $VI_{11}^{\prime}$. Using a similar strategy, we write it into a symmetric form:
\begin{align}\label{expansionforVI11prime}
    \begin{aligned}        VI_{11}^{\prime}&=\underbrace{\int_{0}^{t}\int_{\Gamma_{1}}\kappa\pb^{4}\pb_{\mu}v\cdot n ^{\kappa}n_{\kappa}^{l}h^{\mu\nu}\pb^{4}\pb_{\nu}\eta^{k}n_{\kappa}^{k}n_{\kappa}^{l}+\int_{0}^{t}\int_{\Gamma_{1}}\kappa\pb^{4}\pb_{\mu}\p_{\gamma}v\cdot n ^{\kappa}n_{\kappa}^{l}h^{\mu\nu}\pb^{4}\pb_{\gamma}\pb_{\nu}\eta^{k}n_{\kappa}^{k}n_{\kappa}^{l}}_{VI_{11,1}^{\prime}}\\&\underbrace{+\int_{0}^{t}\int_{\Gamma_{1}}\kappa\pb^{4}\pb_{\mu}v\cdot n ^{\kappa}n_{\kappa}^{l}h^{\mu\nu}\pb^{4}\pb_{\nu}\eta^{k}(\Pi_{k}^{l}-(\Pi_{\kappa})_{k}^{l})+\int_{0}^{t}\int_{\Gamma_{1}}\kappa\pb^{4}\pb_{\mu}\p_{\gamma}v\cdot n ^{\kappa}n_{\kappa}^{l}h^{\mu\nu}\pb^{4}\pb_{\gamma}\pb_{\nu}\eta^{k}(\Pi_{k}^{l}-(\Pi_{\kappa})_{k}^{l})}_{VI_{11,2}^{\prime}}.
    \end{aligned}
\end{align}
For $VI_{11,2}^{\prime}$, we have
\begin{align}
    \begin{aligned}
        |VI_{11,2}^{\prime}|&\lesssim\int_{0}^{t}\int_{\Gamma_{1}}\kappa\|\pb^{4}\pb_{\mu}v\cdot n ^{\kappa}\|_{L^{2}(\Gamma_{1})}\|n_{\kappa}^{l}h^{\mu\nu}\|_{L^{\infty}(\Gamma_{1})}\|\pb^{4}\pb_{\nu}\eta^{k}\|_{L^{2}(\Gamma_{1})}\|(\Pi_{k}^{l}-(\Pi_{\kappa})_{k}^{l})\|_{L^{\infty}(\Gamma_{1})}\\&+\int_{0}^{t}\int_{\Gamma_{1}}\kappa\|\pb^{4}\pb_{\mu}\p_{\gamma}v\cdot n ^{\kappa} \|_{L^{2}(\Gamma_{1})}\|n_{\kappa}^{l}h^{\mu\nu}\|_{L^{\infty}(\Gamma_{1})}\|\pb^{4}\pb_{\gamma}\pb_{\nu}\eta^{k}\|_{L^{2}(\Gamma_{1})}\|(\Pi_{k}^{l}-(\Pi_{\kappa})_{k}^{l})\|_{L^{\infty}(\Gamma_{1})}\\&\lesssim\int_{0}^{t}\kappa^{2}\PP(E^{\kappa})\|\pb^{4}\pb_{\mu}\p_{\gamma}v\cdot n ^{\kappa} \|_{L^{2}(\Gamma_{1})}\|\eta\|_{H^{6}(\Gamma_{1})}\|\eta\|_{H^{3.5}(\Gamma_{1})}.
    \end{aligned}
\end{align}
Now, we solve the elliptic problem:
\begin{align}
     \kappa((1-\bar{\lap})(v\cdot n_{\kappa}))n_{\kappa}^{l} =\sigma(\sqrt{h} \lap_{h}(\eta) \cdot n_{\kappa})n_{\kappa}^{l}-\sqrt{h_{\kappa}} q n_{\kappa}^{l},
\end{align}
\begin{align}
    \|\kappa v\cdot n_{\kappa}\|_{H^{6}(\Gamma_{1})}\lesssim\|\sigma(\sqrt{h} \lap_{h}(\eta) \cdot n_{\kappa})n_{\kappa}^{l}\|_{H^{4}(\Gamma_{1})}+\|\sqrt{h_{\kappa}} q n_{\kappa}^{l}\|_{H^{4}(\Gamma_{1})}\lesssim\PP(E^{\kappa})\|\eta\|_{H^{6}(\Gamma_{1})}+\PP(E^{\kappa}).
\end{align}
Therefore, 
\begin{align}
    \begin{aligned}
        \|\kappa\pb^{6}v\cdot n_{\kappa}\|_{L^{2}(\Gamma_{1})}&\lesssim\|\kappa v\cdot n_{\kappa}\|_{H^{6}(\Gamma_{1})}+\sum_{s=0}^{5}\|\kappa \pb^{s}v \cdot \pb^{6-s}n_{\kappa}\|_{L^{2}(\Gamma_{1})}\\&\lesssim\PP(E^{\kappa})\|\eta\|_{H^{6}(\Gamma_{1})}+\PP(E^{\kappa})+\|\kappa v\|_{H^{5}(\Gamma_{1})}\PP(E^{\kappa})+\kappa\|v\|_{L^{\infty}(\Gamma_{1})}\|\eta^{\kappa}\|_{H^{7}(\Gamma_{1})}\\&\lesssim\PP(E^{\kappa})\|\eta\|_{H^{6}(\Gamma_{1})}+\PP(E^{\kappa})+\|\kappa v\|_{H^{5}(\Gamma_{1})}\PP(E^{\kappa}).
    \end{aligned}
\end{align}
Now, we can write
\begin{align}\label{estimateforVI112prime}
    \begin{aligned}
        |VI_{11,2}^{\prime}|&\lesssim\int_{0}^{t}\kappa\PP(E^{\kappa})(\PP(E^{\kappa})\|\eta\|_{H^{6}(\Gamma_{1})}+\PP(E^{\kappa})+\|\kappa v\|_{H^{5}(\Gamma_{1})}\PP(E^{\kappa}))\|\eta\|_{H^{6}(\Gamma_{1})}\|\eta\|_{H^{3.5}(\Gamma_{1})}\\&\lesssim\PP_{0}+T\PP(E^{\kappa}).
    \end{aligned}
\end{align}
The term $VI_{11,1}^{\prime}$ provide us with the desired estimate:
\begin{align}
    \begin{aligned}
        VI_{11,1}^{\prime}&=\frac{1}{2}\int_{\Gamma_{1}}\kappa\pb^{4}\pb_{\mu}\eta\cdot n ^{\kappa}h^{\mu\nu}\pb^{4}\pb_{\nu}\eta\cdot n_{\kappa}|_{0}^{t}+\frac{1}{2}\int_{\Gamma_{1}}\kappa\pb^{4}\pb_{\mu}\p_{\gamma}\eta\cdot n ^{\kappa}h^{\mu\nu}\pb^{4}\pb_{\gamma}\pb_{\nu}\eta\cdot n_{\kappa}|_{0}^{t}\\&-\frac{1}{2}\int_{0}^{t}\int_{\Gamma_{1}}\kappa\pb^{4}\pb_{\mu}v^{i}\pb^{4}\pb_{\nu}\eta^{k}\p_{t}(n_{\kappa}^{i}h^{\mu\nu}n_{\kappa}^{k})-\frac{1}{2}\int_{0}^{t}\int_{\Gamma_{1}}\kappa\pb^{4}\pb_{\mu}\p_{\gamma}\eta^{i} \pb^{4}\pb_{\gamma}\pb_{\nu}\eta^{k}\p_{t}(n_{\kappa}^{i}h^{\mu\nu}n_{\kappa}^{k}).
    \end{aligned}
\end{align}
Since
\begin{align}
    \begin{aligned}
       &| -\frac{1}{2}\int_{0}^{t}\int_{\Gamma_{1}}\kappa\pb^{4}\pb_{\mu}v^{i}\pb^{4}\pb_{\nu}\eta^{k}\p_{t}(n_{\kappa}^{i}h^{\mu\nu}n_{\kappa}^{k})-\frac{1}{2}\int_{0}^{t}\int_{\Gamma_{1}}\kappa\pb^{4}\pb_{\mu}\p_{\gamma}\eta^{i} \pb^{4}\pb_{\gamma}\pb_{\nu}\eta^{k}\p_{t}(n_{\kappa}^{i}h^{\mu\nu}n_{\kappa}^{k})|\\&\lesssim\int_{0}^{t}\kappa\|\pb^{4}\pb_{\mu}\eta^{i}\|_{L^{2}(\Gamma_{1})}\|\pb^{4}\pb_{\nu}\eta^{k}\|_{L^{2}(\Gamma_{1})}\|\p_{t}(n_{\kappa}^{i}h^{\mu\nu}n_{\kappa}^{k})\|_{L^{\infty}(\Gamma_{1})}\\&+\int_{0}^{t}\kappa\|\pb^{4}\pb_{\mu}\p_{\gamma}\eta^{i} \|_{L^{2}(\Gamma_{1})}\|\pb^{4}\pb_{\gamma}\pb_{\nu}\eta^{k}\|_{L^{2}(\Gamma_{1})}\|\p_{t}(n_{\kappa}^{i}h^{\mu\nu}n_{\kappa}^{k})\|_{L^{\infty}(\Gamma_{1})}\\&\lesssim\PP_{0}+T\PP(E^{\kappa}),
    \end{aligned}
\end{align}
we arrive at 
\begin{align}\label{estimateforVI111prime}
     VI_{11,1}^{\prime}=\frac{1}{2}\int_{\Gamma_{1}}\kappa\pb^{4}\pb_{\mu}\eta\cdot n ^{\kappa}h^{\mu\nu}\pb^{4}\pb_{\nu}\eta\cdot n_{\kappa}+\frac{1}{2}\int_{\Gamma_{1}}\kappa\pb^{4}\pb_{\mu}\p_{\gamma}\eta\cdot n ^{\kappa}h^{\mu\nu}\pb^{4}\pb_{\gamma}\pb_{\nu}\eta\cdot n_{\kappa}+\mathcal{R},
\end{align}
where $|\mathcal{R}|\lesssim \PP_{0}+T\PP(E^{\kappa})$. Combining \eqref{expansionforVI11prime}, \eqref{estimateforVI112prime} and \eqref{estimateforVI111prime}, we have
\begin{align}\label{estimateforVI11prime}
    VI_{11}^{\prime}=\frac{1}{2}\int_{\Gamma_{1}}\kappa\pb^{4}\pb_{\mu}\eta\cdot n ^{\kappa}h^{\mu\nu}\pb^{4}\pb_{\nu}\eta\cdot n_{\kappa}+\frac{1}{2}\int_{\Gamma_{1}}\kappa\pb^{4}\pb_{\mu}\p_{\gamma}\eta\cdot n ^{\kappa}h^{\mu\nu}\pb^{4}\pb_{\gamma}\pb_{\nu}\eta\cdot n_{\kappa}+\mathcal{R},
\end{align}
where $|\mathcal{R}|\lesssim \PP_{0}+T\PP(E^{\kappa})$. Combining \eqref{expansionforVI1prime}, \eqref{estimateforVI12prime}, \eqref{estimateforVI13prime}, \eqref{estimateforVI14prime} and \eqref{estimateforVI11prime}, we have
\begin{align}\label{estimateforVI1prime}
    VI_{1}^{\prime}=\frac{1}{2}\int_{\Gamma_{1}}\kappa\pb^{4}\pb_{\mu}\eta\cdot n ^{\kappa}h^{\mu\nu}\pb^{4}\pb_{\nu}\eta\cdot n_{\kappa}+\frac{1}{2}\int_{\Gamma_{1}}\kappa\pb^{4}\pb_{\mu}\p_{\gamma}\eta\cdot n ^{\kappa}h^{\mu\nu}\pb^{4}\pb_{\gamma}\pb_{\nu}\eta\cdot n_{\kappa}+\mathcal{R},
\end{align}
where $|\mathcal{R}|\lesssim \PP_{0}+T\PP(E^{\kappa})$.
Combining \eqref{expansionforVIprime}, \eqref{estimateforVI2prime} and \eqref{estimateforVI1prime}, we have
\begin{align}\label{estimateforVIprime}
    VI^{\prime}=\frac{1}{2}\int_{\Gamma_{1}}\kappa\pb^{4}\pb_{\mu}\eta\cdot n ^{\kappa}h^{\mu\nu}\pb^{4}\pb_{\nu}\eta\cdot n_{\kappa}+\frac{1}{2}\int_{\Gamma_{1}}\kappa\pb^{4}\pb_{\mu}\p_{\gamma}\eta\cdot n ^{\kappa}h^{\mu\nu}\pb^{4}\pb_{\gamma}\pb_{\nu}\eta\cdot n_{\kappa}+\mathcal{R},
\end{align}
where $|\mathcal{R}|\lesssim \PP_{0}+T\PP(E^{\kappa})$. Combining \eqref{expansionforIVVVIprime}, \eqref{estimateforVprime},\eqref{estimateforIVprime}  and \eqref{estimateforVIprime}, we have
\begin{align}
\begin{aligned}
     &\frac{1}{2}\int_{\Gamma_{1}}\kappa\pb^{4}\pb_{\mu}\eta\cdot n ^{\kappa}h^{\mu\nu}\pb^{4}\pb_{\nu}\eta\cdot n_{\kappa}+\frac{1}{2}\int_{\Gamma_{1}}\kappa\pb^{4}\pb_{\mu}\p_{\gamma}\eta\cdot n ^{\kappa}h^{\mu\nu}\pb^{4}\pb_{\gamma}\pb_{\nu}\eta\cdot n_{\kappa}\\&+\int_{0}^{t}\int_{\Gamma_{1}}\sigma(\sqrt{h}h^{ij}\Pi_{\lambda}^{l}\pb^{4} \p_{\mu}\p_{j}\eta^{\lambda})( h^{\mu\nu}\p_{i}\pb^{4}\pb_{\nu}\eta^{k}\Pi_{k}^{l})=\mathcal{R},
\end{aligned} 
\end{align}
where $|\mathcal{R}|\lesssim\delta\int_{0}^{t}\|\Pi_{\lambda}^{l}\pb^{6}\eta^{\lambda}\|_{L^{2}(\Gamma_{1})}^{2}+C(\delta)T\PP(E^{\kappa})$.

Now, we check that we have finished all energy estimate. Combining
\eqref{Hodgeelliptic},\eqref{apprioriassumption},\eqref{boundaryapprox1},\eqref{boundaryapproxi2}, \eqref{curlest}, \eqref{divest},\eqref{zerotimeest}, \eqref{estimateforkappa6eta}, choosing small $\delta$, we have
\begin{align}
    \sup_{t}\|\eta\|_{5.5}^{2}+\|\sqrt{\kappa}\eta\|_{6.5}^{2}+\int_{0}^{t}\|\sqrt{\kappa}v\|_{5.5}^{2}\lesssim\PP_{0}+T\PP(E^{\kappa}).
\end{align}
Inductively, combining \eqref{Hodgeelliptic},\eqref{apprioriassumption},\eqref{boundaryapprox1},\eqref{boundaryapproxi2}, \eqref{curlest}, \eqref{divest}, \eqref{zerotimeest}\eqref{onetimeest},\eqref{twotimeest},\eqref{threetimeest},and choosing appropriate small $\delta$, we have
\begin{align}\label{estimateforEkappaprime}    E_{\kappa}^{\prime}\lesssim\PP_{0}+T\PP(E^{\kappa}).
\end{align}
Again, combining \eqref{Hodgeelliptic},\eqref{apprioriassumption},\eqref{boundaryapprox1},\eqref{boundaryapproxi2}, \eqref{curlest}, \eqref{divest},\label{estimateforEkappaprime}, \eqref{fourtimeest} and after selecting a enough small $\delta$, we have
\begin{align}\label{estimateforallorderkey}
    \sup_{t}\|\p_{t}^{4}v\|_{L^{2}(\Omega)}^{2}+\|\p_{t}^{3}v\|_{1.5}^{2}+\|\p_{t}^{4}\phi\|_{1}^{2}+\int_{0}^{t}\|\sqrt{\kappa}\p_{t}^{4}v\|_{1.5}^{2}+\int_{0}^{t}\|\p_{t}^{5}\phi\|_{L^{2}(\Omega)}^{2}\lesssim\PP_{0}+T\PP(E^{\kappa}).
\end{align}
Combining
\eqref{phiest1},\eqref{l2estimateforphi}, \eqref{linfinityestimateforphi},\eqref{estimateforEkappaprime} and \eqref{estimateforallorderkey}, we can conclude that
\begin{align}
    E^{\kappa}\lesssim \PP_{0}+T\PP(E^{\kappa}).
\end{align}
Invoking the Gr\"{o}nwall inequality, we prove the uniform in $\kappa$ estimate.
\section{A posteriori estimate}
First, as indicated by the estimate \eqref{l2estimateforphi}, we have:

\begin{align}
\p_{t}^{k} \phi \in L_t^\infty H^{5.5-k}, \quad k = 0, 1, 2, 3.
\end{align}
Applying a similar scheme as above, we can mollify the initial data and establish the solution in higher-order terms. Thus, it suffices to provide the a priori estimates below to improve the regularity, after which a standard approximation scheme follows. 


With this, in \cite{coutand2007well}, the following improvement scheme is established:
We differentiate the boundary equation $qn_{i}=-\sigma\lap_{h}(\eta^{i})$ with $\pb_{A}$, then it becomes
\begin{align}
    \p_{t}^{3}q n_{i}=\lap_{h}(v_{tt})+l.o.t.
\end{align}
The elliptic estimate with curl and divergence analysis will give us 
\begin{align}
    \|v_{tt}\|_{H^{3}(\Omega)}\leq M_{0} .
\end{align}
Similarly, using $\|q_{t}\|_{3.5}\leq  M_{0}$ and $\|q_{tt}\|_{2.5}\leq M_{0}$ we will obtain
\begin{align}
    \|v\|_{H^{5.5}(\Omega)}+\|v_{t}\|_{H^{4.5}(\Omega)}\leq M_{0}.
\end{align}
This will leads to the improvement of regularity of $q$ via elliptic estimate with Neumann boundary:
\begin{align}
    \|q\|_{H^{5}(\Omega)}\leq M_{0}+\|v_{t}\|_{H^{3.5}(\Gamma)}\leq M_{0}.
\end{align}
Now, we rewrite moving boundary as a graph of function $S:\T^{2}\to \R$. Utilizing an identical method in \cite{coutand2007well}, we can obtain 
In our case, the same energy estimate structure gives us
\begin{align}
    \|S\|_{H^{6.5}} \leq M_{0}.
\end{align}

\section{Uniqueness of the solution}
Suppose that $(\eta^{1},v^{1},q^{1},\phi^{1})$ and $(\eta^{2},v^{2},q^{2},\phi^{2})$ are both solutions of \eqref{la2}. We define
\begin{align*}
    w=v^{1}-v^{2},\quad r=q^{1}-q^{2},\quad \xi=\eta^{1}-\eta^{2},\quad \textit{and }z=\phi^{1}-\phi^{2}.  
\end{align*}
Then it follows that
\begin{align}
    \begin{cases}
        \p_{t}r&=w\quad\textit{in $\Omega\times[0,T]$},\\
        \p_{t}w_{i}+(A^{1})_{i}^{j}\p_{j}r&=(A^{2}-A^{1})_{i}^{j}\p_{j}q^{2}-(A^{1})_{k}^{j}\p_{j}((A^{1})_{k}^{l}\p_{l}z\odot(A^{1})_{i}^{m}\p_{m}d^{1})\\&-(A^{1})_{k}^{j}\p_{j}((A^{1})_{k}^{l}\p_{l}d^{2}\odot(A^{1})_{i}^{m}z)-(A^{1}-A^{2})_{k}^{j}\p_{j}((A^{1})_{k}^{l}\p_{l}d^{2}\odot(A^{1})_{i}^{m}\p_{m}d^{2})\quad\textit{in $\Omega\times[0,T]$,}\\&-(A^{2})_{k}^{j}\p_{j}((A^{1}-A^{2})_{k}^{l}\p_{l}d^{2}\odot(A^{1})_{i}^{m}\p_{m}d^{2})-(A^{2})_{k}^{j}\p_{j}((A^{2})_{k}^{l}\p_{l}d^{2}\odot(A^{1}-A^{2})_{i}^{m}\p_{m}d^{2})\\(A^{1})_{i}^{j}\p_{j}w^{i}&=(A^{2}-A^{1})_{i}^{j}v_{j}^{i}\quad \textit{in $\Omega\times[0,T]$,}\\\p_{t}z-\lap_{1}z&=(A^{2}-A^{1})_{i}^{j}\p_{j}((A^{1})_{i}^{k}\p_{k}d^{2})+(A^{2})_{i}^{j}\p_{j}((A^{2}-A^{1})_{i}^{k}\p_{k}d^{2})\\&+z|\nabla_{A^{1}}d^{1}|^{2}+d^{2}((A^{1})_{i}^{k}\p_{k}z\cdot(A^{1})_{i}^{l}\p_{l}d^{1})\quad\textit{in $\Omega\times[0,T]$},\\
        &+d^{2}((A^{1})_{i}^{k}\p_{k}d^{2}\cdot(A^{1})_{i}^{l}\p_{l}z)+d^{2}((A^{1}-A^{2})_{i}^{k}\p_{k}d^{2}\cdot(A^{1})_{i}^{l}\p_{l}d^{2})\\&+d^{2}((A^{2})_{i}^{k}\p_{k}d^{2}\cdot(A^{1}-A^{2})_{i}^{l}\p_{l}d^{2})\\
        r n_{1}&=-\sigma\Pi^{1}h^{1\alpha\beta}\p_{\alpha\beta}\xi-\sigma\sqrt{h^{1}}\lap_{h^{1}-h^{2}}(\eta^{2}),\quad \textit{on $\Gamma\times[0,T]$,}
        \\g_{1}^{3i}\p_{i}z&=(g_{2}^{3i}-g_{1}^{3i})\p_{i}\phi^{2}\quad\textit{on $\Gamma\times[0,T]$,}
        (\xi,w,z)=(0,0,0)\quad\textit{on $\Omega\times\{t=0\}$}.
    \end{cases}
\end{align}
We define
\begin{align}
    H(t)=\sum_{k=0}^{3}(\|\p_{t}^{k}w\|_{4.5-k}^{2}+\int_{0}^{t}\|\p_{t}^{k}z\|_{5.5-k}^{2})
\end{align}
It suffices to prove $H(t)\lesssim T\PP(H)$. From our previous argument, we can set $E^{1}+E^{2}\leq M_{0}$. ($E^{1}$ and $E^{2}$ are the energy defined in Section 4) The basic idea is still to give a $\curl$, $\di$, and boundary estimate. The form is more complicated than original equation, but the proof is similar. To show this, when we give $\curl$, $\di$ estimate, we first write it in the form of $L_{j}\eta^{1}=...$ and $L_{j}\eta^{2}=...$. where $L_{j}$ is $\curl$ or $\di$. Then the argument is quite similar, if we regard $L_{j}$ as a linear operator. Taking the difference, or after integrating by part in time and imposing some derivatives on it. The form of expression and their dependence on initial value has to be similar. We omit the details. For the $\pb_{A}$ estimate, in the absence of artificial viscosity on boundary, the only part needs to be mentioned is the cancellation scheme. That's where used an extra regularity of initial data ($\|u_{0}\|\in H^{5.5}$). Then our solution can be improved in a corresponding way. Now, we don't need to worry about the regularity of coefficient and the pressure term. For $\lap_{g_{1}-g_{2}}$ and the term involving $\phi$, the leading order term still maintain the cancellation. Hence, a similar argument in Section 4 will give us the desired result.

\appendix
\section{Appendix}

In \cite{lin2023nematic}, they proposed a physical boundary condition for nematic liquid crystal in a fixed domain with slippery condition.(See the discussion therein). Simply speaking, their condition can be written as
\begin{align}\label{linbdycondition}
    \phi^{3}=0,\quad
    \nabla_{N}\phi^{1}=\nabla_{N}\phi^{2}=0. 
\end{align}
We don't know whether it is still a physical condition in free boundary problem, but we can give a similar result here. Consider the stress tensor 
\begin{align}
   S_{ij}= p(Id_{3\times 3})_{ij}+\nabla_{i} d \odot\nabla_{j} d-\frac{1}{2}|\nabla d|^{2}(Id_{3\times 3})_{ij}
\end{align}
On boundary, their requirement is roughly to say "free slippery":
\begin{align}
    pN\cdot \tau+\nabla_{N} d \odot\nabla_{\tau} d-\frac{1}{2}|\nabla d|^{2}N\cdot \tau=0.
\end{align}
In our condition, since the boundary is not fixed, we need a new constrain as
\begin{align}
     pN\cdot N+\nabla_{N} d \odot\nabla_{N} d-\frac{1}{2}|\nabla d|^{2}N\cdot N=-\sigma H(\eta).
\end{align}
So, if we still use their condition \eqref{linbdycondition}, our boundary should be written as 
\begin{align}
     p+\nabla_{N} d^{3}\nabla_{N} d^{3}-\frac{1}{2}|\nabla d|^{2}=-\sigma H(\eta).
\end{align}
Their will be an extra term $\nabla_{N} d^{3}\nabla_{N} d^{3}$ compared to our case. However, our work could still handle this situation with a slight modification. Since the boundary condition for $\phi$ is still elliptic, we can obtain a similar estimate for $\phi$. We now provide a sketchy proof. Our energy estimate can start from $\|\p_{t}^{3}v\|_{L_{t}^{\infty}H^{1}(\Omega)}$ and $\|\p_{t}^{4}\phi\|_{L_{t}^{2}H^{1}_{x}(\Omega)}$ instead of $\|\p_{t}^{4}v\|_{L_{t}^{\infty}L^{2}(\Omega)}$ and $\|\p_{t}^{5}\phi\|_{L_{t}^{2}L^{2}_{x}(\Omega)}$. Roughly speaking, if we disregard the mollifier and the associated error analysis, the a priori energy estimate can be shown as follows:
\begin{align}
    \p_{t}^{5}\phi-\p_{t}^{4}\lap_{g}\phi=\p_{t}^{4}(\phi|\nabla\phi|^{2})
\end{align}
Testing the $\p_{t}^{4}\phi$, we have the standard energy estimate. For $\phi^{1}$ and $\phi^{2}$, we consider
\begin{align}
    \p_{t}^{4}\pb\phi-\p_{t}^{3}\pb\lap\phi=\p_{t}^{3}\pb(\phi|\nabla\phi|^{2}).
\end{align}
Testing $\p_{t}^{3}\pb v\cdot\nabla\phi$, we will have a similar cancellation structure for $\phi^{1}$ and $\phi^{2}$.
On the other hand, the Euler equation can be written as
\begin{align}
    \p_{t}v^{j}+\nabla_{j} p=-\p_{t}\phi\nabla_{j}\phi-\nabla_{i}\phi\nabla_{i}\nabla_{j}\phi+\frac{1}{2}\nabla_{j}|\nabla_{N}\phi^{3}|^{2},
\end{align}
where $\nabla_{N}\phi^{3}=\frac{A_{i}^{3}}{|A_{.}^{3}|}A_{i}^{j}\p_{j}\phi^{3}$. For $\phi^{1}$ and $\phi^{2}$, we perform a similar cancellation as in the previous section. For $\phi^{3}$, the leading term reads
\begin{align}
\begin{aligned}
    & \underbrace{\int_{0}^{t}\int_{\Omega}\nabla_{i}\phi^{3}\p_{t}^{3}\pb(\nabla_{i}\phi^{3})\nabla_{j}\p_{t}^{3}\pb v^{j}}_{\dive v=0,l.o.t.}-\int_{0}^{t}\int_{\p\Omega}A_{i}^{3}\p_{3}\phi^{3}\p_{t}^{3}\pb(A_{i}^{3}\p_{3}\phi^{3})A_{j}^{3}\p_{t}^{3}\pb v^{j}\\-&\underbrace{\int_{0}^{t}\int_{\Omega}\nabla_{N}\phi^{3}\p_{t}^{3}\pb\nabla_{N}\phi^{3}\nabla_{j}\p_{t}^{3}\pb v^{j}}_{\dive v=0, l.o.t.}+\int_{0}^{t}\int_{\p\Omega}A_{i}^{3}\p_{3}\phi^{3}\p_{t}^{3}\pb(A_{i}^{3}\p_{3}\phi^{3})A_{j}^{3}\p_{t}^{3}\pb v^{j}\\=&l.o.t.
\end{aligned}   
\end{align}
Thus, we can establish the uniform-$\kappa$ energy estimate in a similar way, with the help of error analysis. Actually, if we don't have the term $\frac{1}{2}\nabla|\nabla_{N}\phi^{3}|^{2}$, we can still cancel the boundary terms via harmonic heat flow. We omit the details. 
\bibliographystyle{plain}
\bibliography{ref}
\end{document}